# Families of Regular Polygons and their Mutations

G.H. Hughes

**Abstract**

Every regular N-gon determines a canonical 'family' of regular polygons which are conforming to the bounds of the 'star polygons' determined by N. This is illustrated below for the regular tetradecagon known as N = 14. These classical star polygons are formed from truncated extended edges of the N-gon and the intersection points ('star points') define the parameters of the cyan family of 'conforming' polygons. Here there are 6 star points on either side of N and in the First Family Theorem of Section 2 each star[k] defines an S[k] 'tile' so left of N there will be 6 family members S[1] to S[6]. When N is even the largest family tile S[N/2-1] is known as D. Since it is congruent to N it can support its own First Family as shown here. The only tile in both families is the central S[N/2-2] tile known as M. Here M is S[5] acting as a surrogate N = 7. In Section 4 we introduce the outer-billiards map $\tau$ and the matching singularity set W.

**Figure 1** The evolution of W and First Families for N = 14 Detail of W

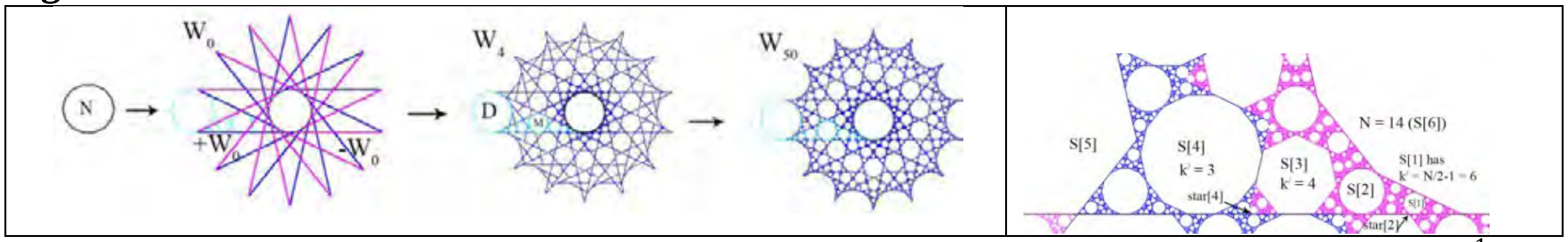

By definition the singularity set W of a convex polygon are the points where neither $\tau$ or $\tau^{-1}$ are defined so W will have clockwise and ccw versions which we call +W and – W. Our convention for N is clockwise so the blue region here is determined by iterating the 14 blue $+W_0$ extended edges to obtain +W and it would be necessary to iterate the magenta 'trailing' edges to obtain -W. By definition $\pm W_0$ are relate by a horizontal reflection T of N, but a regular N-gon is 'chiral' so its 'mirror-image' T(N) cannot be obtained by rotations of $2\pi/N$. Since the evolution of W requires both reflections and rotations the blue +W web will be forever distinct from the magenta –W web as shown in the side-bar – where with apologies +W is now colored magenta.

Here S[4] and S[3] share an edge but they are in distinct invariant regions . Since S[4] is in -W points will rotate locally ccw under $\tau$ while points in S[3] will rotate cw. They share an edge because these rotations are synchronized with 3 steps of S[4] matching 2 of S[3]. In general +W and –W will be zero measure manifolds with disjoint invariant regions which are symmetric left and right of M so for N = 14 the two regions have symmetric counterparts left of S[5].

The Invariance Conjecture of Appendix II says that every boundary between regions will be mediated by a 'dual' pair of edge-sharing S[k]-S[k′] with k′ = N/2-k for N even. The most important example is S[1] and N. They are dual since N is S[N/2-1] relative to the S[k]. As long as one of k or k′ is 'primitive' with gcd(k,N) = 1 these pairs will survive in W. For N even there will be exactly $C_N = \Phi(N)/2$ primitive star[k] on each side of the central M tile and this implies that there will be $C_N$-1 primitive regions right (and left) of M. When N is twice odd, N/2 and N will share the same regions so it is sufficient to assume N even. For N = 14 and N = 7, $C_N = 3$.

One way to distinguish the +Wand –W regions is to use the Jacobian J of N which will have Det(J) ±1 for ccw and cw rotations. Therefore the boundaries of the regions will have Det(J) = 0, but this implies that any map of +W and –W to their union cannot be 'smooth' and it is our contention that the non-smooth ‘solitons’ of this mapping are the S[k] dual tiles. Here $\tau$ commutes with J so it is invariant under ±W and serves as a ‘polarizing’ map. In outer space polarized molecules could mimic our chiral life forms where DNA is -W and RNA is +W.

As noted above, the First Family Theorem (FFT) of Section 2 shows how N constructs each S[k] tile based on star[k] of N. Since the star points progress outwards with star[1] as a vertex of N, this implies an evolution that proceeds right-to-left in a primarily +W clockwise fashion.

**Figure 2** The early cw evolution of the First Family of N = 42

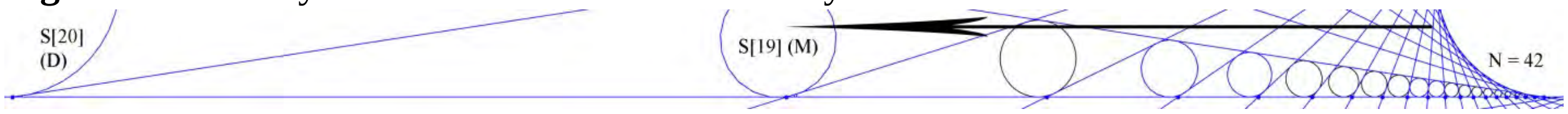

But the –W ccw web is equally important and we will assume a Reflection Principle which implies that the combined web must evolve in both directions. So it is safe to assume that D can generate a symmetric First Family and the only tile in both families will be the S[N/2-2] ‘M’ tile.

**Figure 3** The ccw evolution of the D tile of N = 42

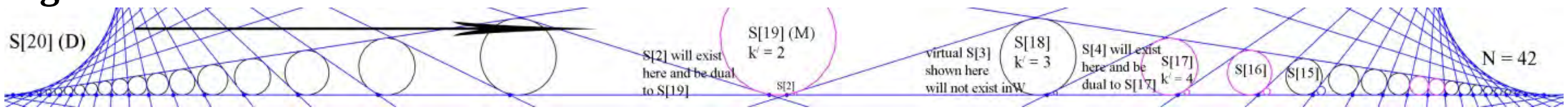

We show here how the interaction of these two webs can create the dual S[k] tiles and the matching invariant regions. Appendix II gives more details about the evolution of N = 42 and many other examples. The twice-even case is similar but the central ‘M’ tile is not an N/2 clone.

A very important feature of the First Family Theorem is that the interaction of N and S[k] has a ‘dual’ version where S[k] can construct N, and the GeneralizedFFT shows how an S[k] can construct an S[k+1] on the left and an S[k′] on the right (if space permits).

Here is an example with S[17] where its left-side star[18] is star[18] of N, so this can be used to construct S[18] just like N. On the right-side every S[k] evolves in the FFT with star[k′] equal to star[k] of N so here star[4] of S[17] is star[17] of N. This star[4] can be used to construct an S[4] clone in the same fashion as it could construct N. In the FFT the midpoint of the new tile is offset by sN/2 from star[k] of N, so here Midpoint S[4] will be offset by sN/2 from star[4] of N so cS[4 ] = Midpoint + {0,hS[4]} and hS[4] is known from the FFT. This new S[4] will share an edge with S[17] because their edges are synchronized with 4 ccw steps of S[17] matching 17 cw steps of S[4]. So their shared edge will align with star[17] of N to define the $I_4$ separatrix.

S[18] has a dual S[3] but it will not survive in W because N/gcd(3,N) = 14 and S[18] will be mutated into the weave of two regular 14-gons which will foster an S[2] with no shared edge.

| **Figure 4** S[17] will foster a dual S[4] clone | **Figure 5** S[18] will be mutated |
|---|---|
| 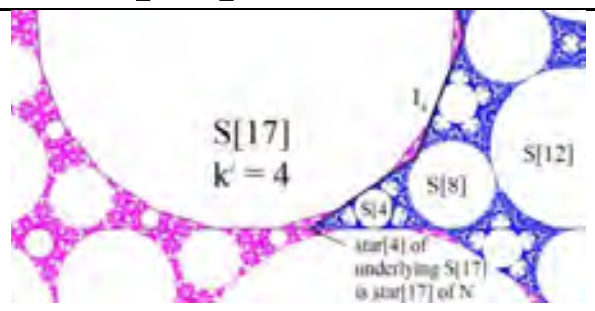 | 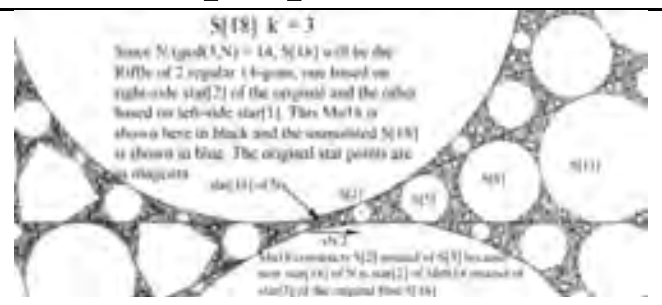 |

**Figure 6** The Invariance Conjecture predicts that N = 42 will have $C_N$-1 = 5 regions right of M

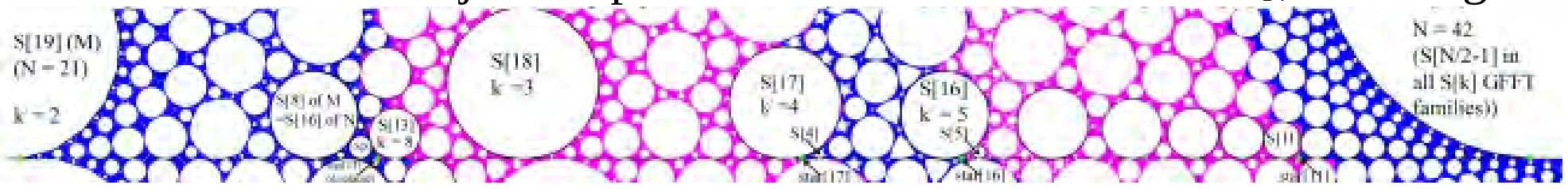

**Organization of the five sections of this paper**
Section 1 (Star Polygons and Star Points) introduces star polygons and their scaling based on ratios of star points. The main results are the Scaling Lemma relating the scaling of N/k to the scaling of N and the Two-Star Lemma for construction of regular polygons.

Section 2 (Conforming Regular Polygons) shows how the S[k] tiles which form the First Families can be constructed from the star[k] points of N. The main results are the First Family Theorem (FFT) and a generalized version (GFFT) for families of the S[k] based on the premise that the S[k] evolve in a multi-step fashion under the outer-billiards map.

Section 3 (Scaling) shows that the S[1] tile of an N-gon can be used as a scaling reference since hS[1]/hS[k] = scale[k] by the First Family Scaling Lemma. We show that the minimal GenScale is inherent in star polygons so every regular N-gon has the potential to support 'generations' of First Families on the edges of it's maximal D tile. The Scaling Field Lemma shows that the primitive scale[k] (with gcd(k,N) = 1) are unit generators of the scaling field $\mathbb{Q}_N^+$.

Section 4 (The Outer Billiards Map ) relates this geometry and scaling to the outer-billiards map $\tau$. Lemma 4.1 uses a height/radius duality between the star[k] points and S[k] centers to show that each S[k] defines a periodic step-k orbit. Example 4.3 shows how the evolution of the web W can be reduced to a 'shear and rotation' - which implies that the S[k] evolve in a multi-step fashion. This explains mutations and also justifies the need for a Generalized FFT to characterize the next generation. But there is no clear path for recursion except for the 8k+2 family where the tiles local to S[2] have a well-defined 'temporal' scaling. In general the Edge Conjecture makes predictions about the mod-4 or mod-8 tiles of S[2] which survive in the 3rd generation. This conjecture distinguishes just 8 classes of regular N-gons of the form 8k + j.

Section 5 has annotated examples of singularity sets and Appendix I describes the 'digital-filter' map of Chau & Lin and a 'dual-center' map of Erik Goetz. Both of these maps reduce to a 'shear and rotation' and their local webs appear to replicate W. Appendix II discusses the large scale web and the invariant regions which arise from the interaction of +W and –W webs.

**Historical Background** In a 1978 article called "*Is the Solar System Stable*?" Jurgen Moser [Mo2] used the landmark KAM Theorem named after A. Kalmogorov,V. Arnold and Moser to show that there is typically a non-zero measure of initial conditions that would lead to a stable solar system. Since the KAM Theorem was very sensitive to continuity, Moser suggested a 'toy model' based on orbits around a polygon where continuity would fail. This is now called the 'outer-billiards' map and in 2007 Richard Schwartz [S1] showed that if the polygon was in a certain class of 'kites', orbits could diverge and stability failed. The special case of a regular polygon was settled earlier in 1989 when F.Vivaldi and A. Shaidanko [VS] showed that all obits are bounded. The author met with Moser at Stanford University in that same year to discuss the 'canonical' structures that always arise in the regular case. Moser suggested that a study of these structures would be an interesting exercise in 'recreational' mathematics and this is what it became over the years, but even this simple case retained a complex topology. Moser was interested in computer applications such as my article in [H1] and I gave him software from the recent 'Summer School' in chaos theory at Los Alamos. See A3 and Chronology.

## Section 1. Star Polygons and Star Points

'Star' polygons or 'stellated' polygons were first studied by Thomas Bradwardine (1290-1349), and later L.Schlafli (1840-1895) and H.S. Coxeter (1907-2003).

The vertices of a regular N-gon with radius r are {rCos[2πk/N], rSin[2πk/N]} for $1 \le k \le N$. A 'star polygon' {p,q} generalizes this by allowing N above to be rational of the form p/q so the vertices are given by {p,q} = {rCos[2πkq/p], rSin[2πkq/p]} for $1 \le k \le p$.

Using the notation of Coxeter [Co] a regular heptagon can be written as {7,1} (or just {7}) and {7,3} is a 'step-3' heptagon formed by joining every third vertex of {7} so the exterior angles are 2π/(7/3) instead of 2π/7.

By the definition above, {14,6} would be the same as {7,3}, but there are two heptagons embedded in N = 14 and a different starting vertex would yield another copy of {7,3}, so a common convention is to define {14,6} using both copies of {7,3} as shown below. This convention guarantees that all the star polygons for {N} will have N vertices.

| {7,1} (a.k.a. N = 7) | {7,3} | {14,6} | {10,2} |
|---|---|---|---|
| 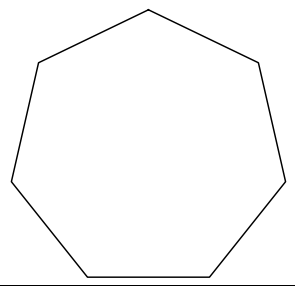 | 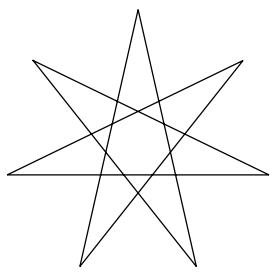 | 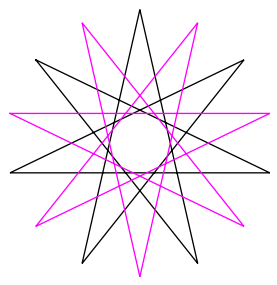 | 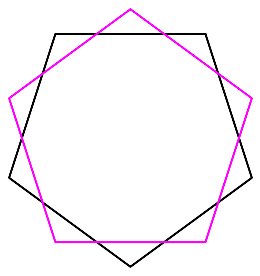 |

The number of 'distinct' star polygons for {N} is the number of integers less than N/2, which we write as $\langle N/2\rangle$. So for a regular N-gon, the 'maximal' star polygon is {N, $\langle N/2\rangle$}.

Our default convention for the 'parent' N-gon will be centered at the origin with 'base' edge horizontal, and the matching {N,1} will be assumed equal to N. In general sN, rN and hN will denote the side, radius and height (apothem). Typically we will use hN as the lone parameter.

**Definition 1.1** The *star points* of a regular N-gon are the intersections of the edges of {N, $\langle N/2\rangle$} with a single extended edge of the N-gon (which will be assumed horizontal). By convention the star points are numbered from star[1] (a vertex of N) outwards to star[$\langle N/2\rangle$] - which is called GenStar[N]. Therefore every star[k] is a vertex of {N,k} embedded in {N,$\langle N/2\rangle$}.

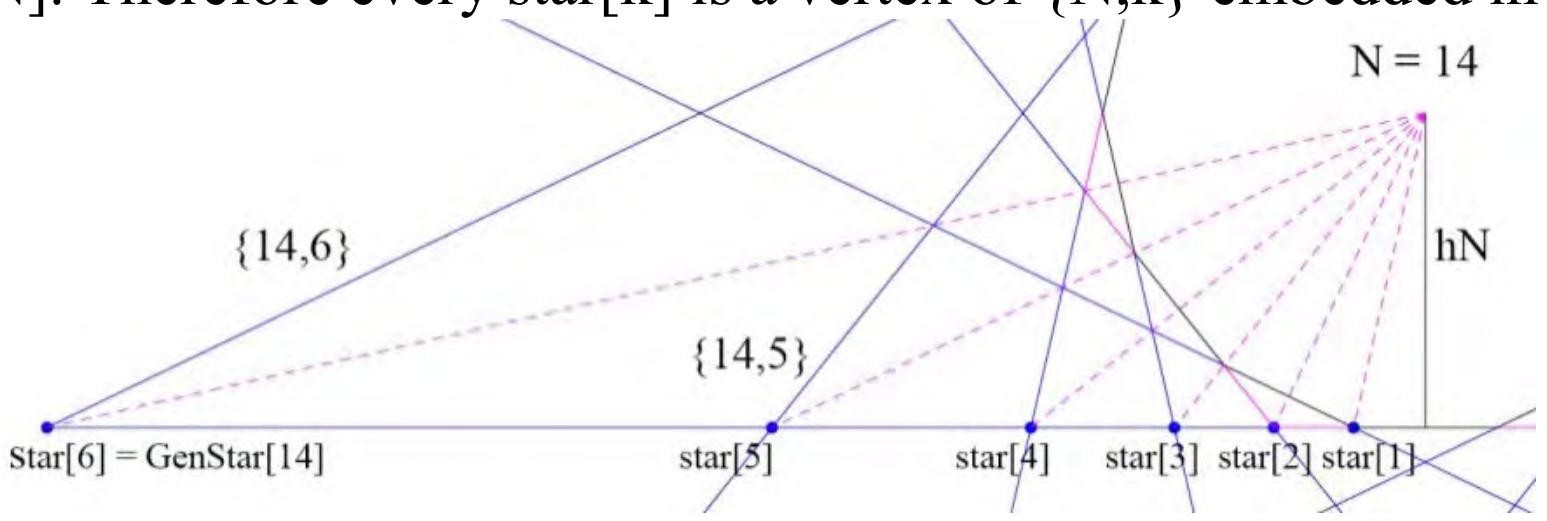

The star points could equally be defined on the positive side of N, but this 'left-side' convention is sometimes convenient. The symmetry of these choices makes it irrelevant which convention is used. In general, the star points of a regular N-gon with apothem hN have coordinates:

$$\text{star}[k] = \{\pm\, hN \cdot s_k, -hN\} \text{ where } s_k = \text{Tan}[k\pi/N] \quad \text{for } 1 \le k < N/2.$$

The *primitive* star[k] and $s_k$ are those with (k,N) = 1, where (k,N) is gcd(k,N). Sometimes we will loosely refer to the $s_k$ as 'star' points.

Note: In the Introduction we noted that there is a long history of interest in trigonometric functions of rational multiples of π. Niven [N] shows that these 'trigonometric numbers' are also algebraic numbers. When C.L. Siegel communicated to S.Chowla [Ch] a proof that the primitive $s_k$ for the cotangent function are linearly independent, he only proved it for prime N. In 1970 Chowla generalized this result using character theory and Dirichlet's $\mathcal{L}$-series, but the general tangent case was only settled recently in [Gi] (1997). Therefore the primitive $s_k$ are linearly independent and it follows that any non-primitive $s_k$ must be a $\mathbb{Q}$-linear combination of the primitive $s_k$ but the coefficients may be very hard to find. See Corollary 3.1 to Lemma 3.3.

Section 3 is devoted to scaling, but here we present the basic definitions and prove the Scaling Lemma. In Section 3 we will use these scales as a basis for the maximal real subfield of $\mathbb{Q}_N$.

**Definition 1.2** The (*canonical) scales* of a regular N-gon are scale[k] = $s_1/s_k$ for $1 \leq k < N/2$ .The *co-scales* are of the form $s_k/s_1$. The *primitive* scales or coscales are those with (k,N) = 1. GenScale[N] is defined to be scale[<N/2>].

By definition scale[1] is always 1 and GenScale is the minimal scale. Since these scales are independent of height, to compare scales for an N-gon and an N/k-gon, the later can be regarded as circumscribed about the N-gon with shared center and height – so the sides can be compared.

**Definition 1.3** If N and M are regular polygons with M = N/k, then ScaleChange(N,M) = ($s_1$ of N)/($s_1$ of M) ≤ 1. This is abbreviated SC(N,M) or just SC(N) < ½ ,when M is N/2.

The case when N and M are related by M = N/2 will play an important role in what follows.

**Definition 1.4** For N even, a regular N/2-gon is the *outer-dual* of N if the N/2-gon can be circumscribed about the N-gon with matching 'base' edge, center and height. In this case sN/s(N/2) will be SC(N) = $s_1/s_2$ = scale[2] of N. When N is twice-odd, this outer-dual will be called the *parity-dual* or *gender-dual* of N. This is the only case where the GenStars match.

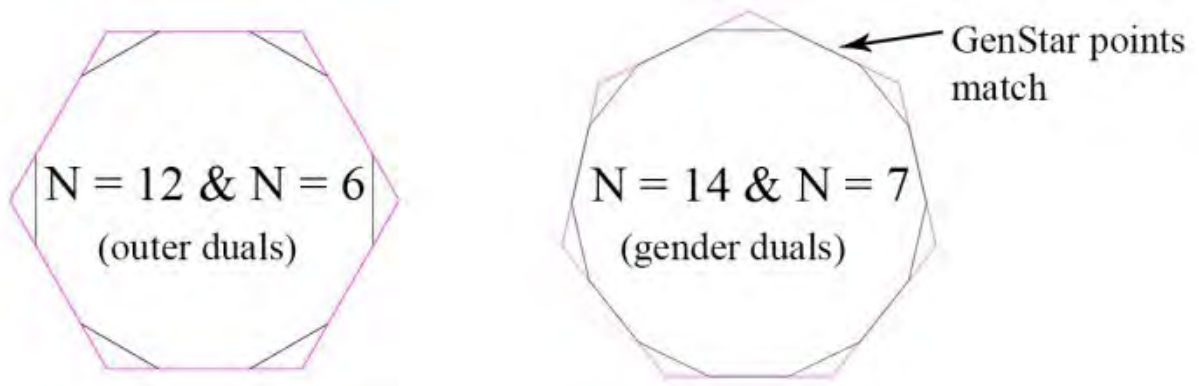

**Lemma 1.1** (Scaling Lemma) Suppose N and M are regular polygons and M = N/k , then

scale[j] of N/k = scale[kj] /scale[k] of N

Proof: There is no loss of generality in assuming that N and N/k are in 'standard position' at the origin with equal heights so N/k will be a 'circumscribed factor polygon' of N. If s is the side of N, then s/SC(N,N/k) will be the side of N/k. The external angle of N/k is 2πk/N so in this position, every kth edge will coincide with an edge of N. Therefore it will share every kth star point with N and by definition the corresponding scales are related by the ratio of the sides of N and N/k, so scale[j] of N/k = (scale[kj] of N)/SC(N,N/k). Since scale[1] of N/k = 1 = scale[k]/SC(N,N/k), it follows that SC(N,N/k) is scale[k] of N. □

Therefore scale[k] of N = 7 is the same as scale[2k]/scale[2] for N = 14 and this twice-odd case is the only nesting where the GenStar points coincide, so GenScale[7] = GenScale[14]/scale[2].

**Lemma 1.2** When N is even, GenScale[N] = $\mathrm{Tan}^2[\pi/N] = s_1^2$ and when N is odd GenScale[N] = $\mathrm{Tan}[\pi/2N]\cdot\mathrm{Tan}[\pi/N] = s_1\cdot s_2$ of 2N.

Proof : For N even, the star points inherit the reflective symmetry of N under $\tan[\pi/2 - \theta] = \cot[\theta]$, so $s_{N/2-k} = 1/s_k$. Setting k = <N/2> = N/2–1, GenScale[N] = $s_1/(1/s_1) = s_1^2$. When N is odd the Scaling Lemma says that GenScale[N] = GenScale[2N]/scale[2] = $s_1^2/(s_1/s_2) = s_1\cdot s_2$ of 2N. □

To define a regular N-gon in a given coordinate system, it is sufficient to know its height (apothem) and its center, but both of these are determined by knowing the co-ordinates of two star points, as any cartographer would know.

**Lemma 1.3 (**Two-Star Lemma) If P is a regular N-gon, any two star points are sufficient to determine the center and height.

Proof. By definition, the star points lie on an extended edge of P. There is no loss of generality in assuming that this extended edge is parallel to the horizontal axis of a known coordinate system with arbitrary center.

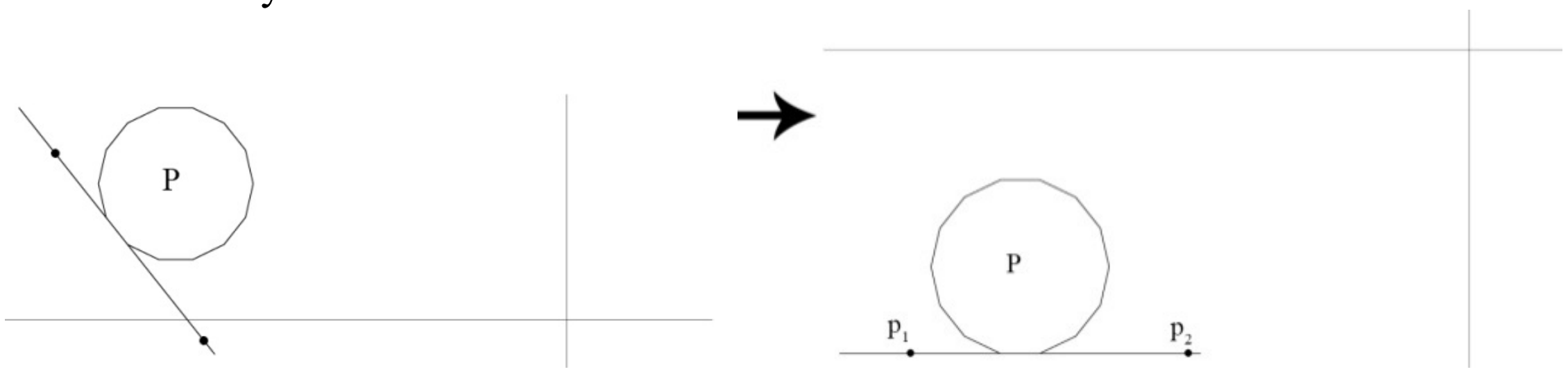

Since all points on this extended edge will have a known second coordinate, we will just need the horizontal coordinates of the star points, namely $p_1$ and $p_2$ with $p_2 > p_1$. Therefore $d = p_2 - p_1$ will be positive. Relative to the polygon P, $p_1 = \pm$ star[j][[1]] = $\pm hP\cdot s_j$ and $p_2 = \pm$ star[k][[1]] = $\pm hP\cdot s_k$ (where x[[1]] is the horizontal coordinate of x). These indices j and k must be known. There are only two cases to consider:

(i) If $p_1$ and $p_2$ are on opposite sides $hP = (p_2 - p_1)/(s_k + s_j)$

(ii) If $p_1$ and $p_2$ are on the same side of P, there is no loss of generality in assuming that it is the right-side of P because star points always exist in their symmetric form with respect to P. In this case we can assume that $1 \leq j < k < N/2$ so $hP = (p_2 - p_1)/(s_k - s_j)$

Now that hP is known, the horizontal displacement of $p_1$ and $p_2$ relative to the center of P are $x = hP\cdot s_j$ and $d + x = hP\cdot s_k$ if both are on the same side or $x = hP\cdot s_j$ and $d - x = hP\cdot s_k$ if they are on opposite sides Of course only one of these displacements is needed to define the center of P. □

**Example 1.1** (The one-elephant case) P below shares two star points with the elephant N = 14 which defines the coordinate system. The shared star points are (right-side) star[3] and star[6] of P which we write as starP[3] and starP[6]. (The points $p_1$ and $p_2$ shown below are the horizontal coodinates of these points relative to MidM.)

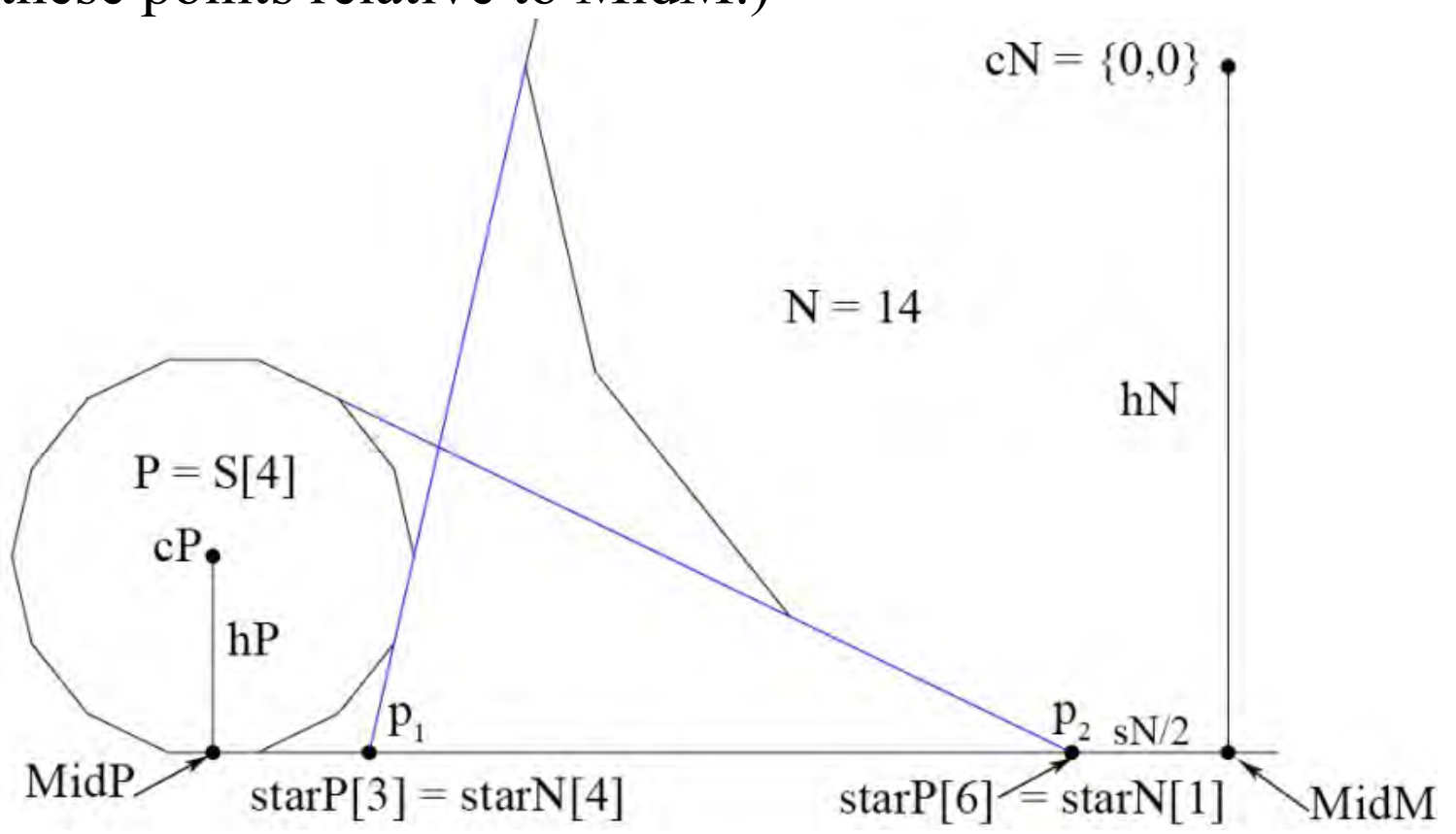

Because of the symmetry between P and N, either one could be used as reference to construct the other, but here we assume that hN is known and it is desired to find hP (and cP) relative to N.

By the Two-Star Lemma hP = $\dfrac{p_2 - p_1}{Tan[6\pi/14] - Tan[3\pi/14]} = hN\dfrac{Tan[\frac{4\pi}{14}] - Tan[\frac{\pi}{14}]}{Tan[\frac{6\pi}{14}] - Tan[\frac{3\pi}{14}]}$

Therefore $\dfrac{hP}{hN} = \dfrac{hS[4]}{hN} = \dfrac{Tan[\frac{4\pi}{14}] - Tan[\frac{\pi}{14}]}{Tan[\frac{6\pi}{14}] - Tan[\frac{3\pi}{14}]} = \dfrac{Tan[\pi/14]}{Tan[3\pi/14]} = \dfrac{s_1}{s_3}$ (which is also $s_1 \cdot s_4$)

Therefore $hP \cdot s_3 = hN \cdot s_1$. By definition, the right side is the horizontal displacement of star[1] from MidN, which is sN/2. This must be equal to the left side which is the horizontal displacement of starP[3] from MidP. Therefore the horizontal displacement of P from starP[3] must be –sN/2. The trigonometric identity above follows from the fact that the displacement from $p_1$ to $p_2$ must be the same when viewed from N or from P.

P here is known as S[4] because it was constructed using star[4] (and star[1]) of N. Clearly this same construction can be carried out for the remaining star[k] points of N = 14 and we will do this below for arbitrary N. Since each S[k] will have a symmetric relationship with N, the displacements will always be –sN/2 and every S[k] will share its star[N/2-1] point with N.

Algebraically the S[k] and N are closely related because hP/hN will always be an element of the 'scaling field' defined by N = 14 (or N = 7). This is the number field generated by 2Cos[2π/7] or GenScale[7] = Tan[π/7]·Tan[π/14]. See Section 3. Using Mathematica:

**AlgebraicNumberPolynomial[ToNumberField[hS[4]]/hN,GenScale[7]],x]** = $\dfrac{x^2 - x + 1}{2}$

So hS[4]/hN = $\dfrac{x^2 - x + 1}{2}$ where x = Tan[π/7]·Tan[π/14] and this scale is an algebraic number.

Since N = 14 and the matching N = 7 have φ(N)/2 = 3 where φ is the Euler totient function, they are classified as 'cubic' polygons. This means that any generator of the scaling field will have minimal degree 3, so the scaling polynomials will be at most quadratic.

**Example 1.2** (The two-elephant case) The tiles below exist in the coordinate space of N = 11. This is a 'quintic' N-gon so the algebra is much more complex than N = 14. Px and DS5 share a star point which is off the page at the right, but they do not share any other star points so it was a challenge to find a second defining star point of Px , even though the parameters and star points of DS5 are known from the First Family Theorem to follow. Since these two elephants are only distantly related, it is unusual for them to share a third tile, which we call Sx. This Sx tile shares extended edges with both Px and DS5 so it could be used to define both of them.

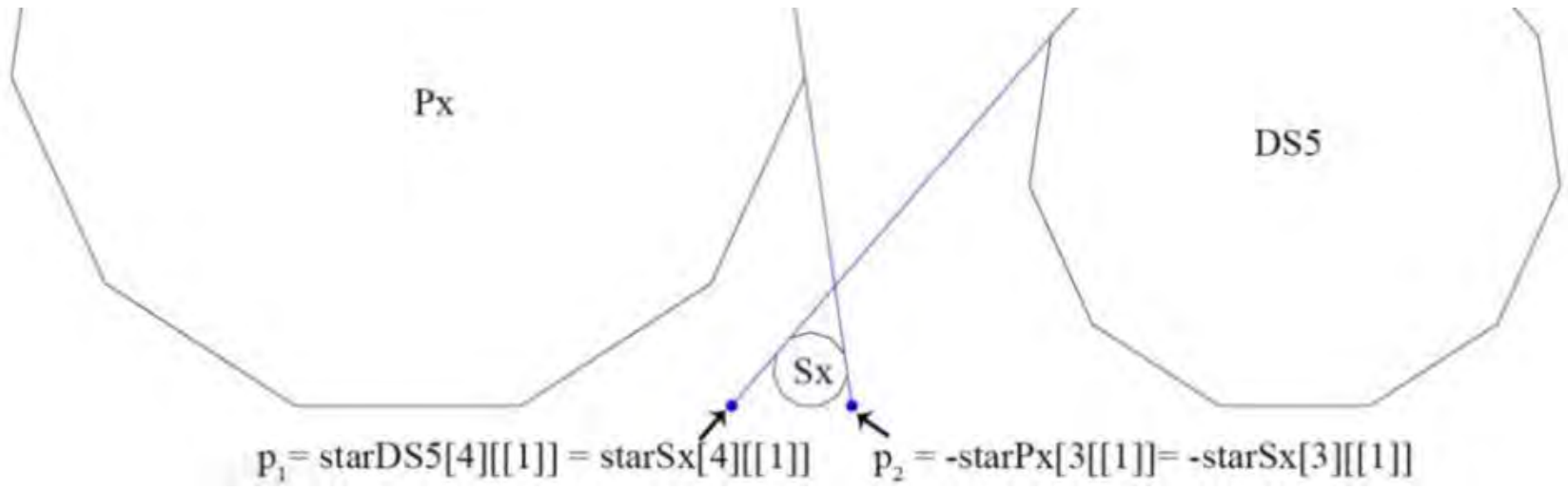

These tiles share the 'base' edge of N = 11. The calculations below assume hN = 1 so the star points have vertical coordinate –1 and we will just need the horizontal coordinates.

$$p_1 = starDS5[4][[1]] = -\tan\left[\frac{\pi}{22}\right]^2\left(2\cot\left[\frac{\pi}{22}\right]-\cot\left[\frac{\pi}{22}\right]^3+\tan\left[\frac{\pi}{11}\right]\right)+2\cot\left[\frac{\pi}{22}\right]\left(-1+2\tan\left[\frac{\pi}{22}\right]^2+\tan\left[\frac{\pi}{22}\right]^3\tan\left[\frac{\pi}{11}\right]\right)-$$
$$\cot\left[\frac{3\pi}{22}\right]\tan\left[\frac{\pi}{22}\right]^2\tan\left[\frac{\pi}{11}\right]\left(2+\tan\left[\frac{\pi}{22}\right]\tan\left[\frac{\pi}{11}\right]\right)\tan\left[\frac{5\pi}{22}\right]+\tan\left[\frac{\pi}{11}\right]\left(1-\tan\left[\frac{\pi}{22}\right]\left(2+\tan\left[\frac{\pi}{22}\right]\tan\left[\frac{\pi}{11}\right]\right)\tan\left[\frac{5\pi}{22}\right]\right)$$

The trigonometric expression for star[3] of Px is more complex because Px shares no canonical scaling with N. Mathematica prefers to do these calculations in 'cyclotomic' form, which will have vanishing complex part here. The cyclotomic form for star[3] of Px is:

$$p_2 = -starPx[3][[1]] = -\frac{\mathrm{i}\left(3-8(-1)^{1/11}+6(-1)^{2/11}+4(-1)^{3/11}+11(-1)^{4/11}+11(-1)^{5/11}+4(-1)^{6/11}+6(-1)^{7/11}-8(-1)^{8/11}+3(-1)^{9/11}\right)}{\left(1+(-1)^{2/11}\right)^2\left(-1+(-1)^{5/11}\right)}$$

By the Two-Star Lemma, hSx = $(p_2-p_1)$/(Tan[4π/11]+Tan[3π/11]) =

$$-(\left(2\left(11\mathrm{i}+2(-1)^{1/22}-5(-1)^{3/22}-5(-1)^{5/22}+2(-1)^{7/22}-3(-1)^{9/22}-2(-1)^{13/22}+2(-1)^{17/22}-11(-1)^{19/22}+3(-1)^{21/22}\right)\right)/$$
$$\left(\left(5-(-1)^{1/11}+6(-1)^{2/11}-(-1)^{3/11}+5(-1)^{4/11}-3(-1)^{5/11}+4(-1)^{6/11}-4(-1)^{7/11}+4(-1)^{8/11}-4(-1)^{9/11}+3(-1)^{10/11}\right)\left(\cot\left[\frac{3\pi}{22}\right]+\cot\left[\frac{5\pi}{22}\right]\right)\right))$$

The ratio hSx/hN will be in the scaling field $S_{11}$ which is generated by x = GenScale[11] = Tan[π/11]·Tan[π/22]. **AlgebraicNumberPolynomial[ToNumberField[hSx/hN, GenScale[11]]** yields p(x) = $-\frac{9}{4}+52x+\frac{63x^2}{2}+x^3-\frac{5x^4}{4}$. This is what we call the 'characteristic' polynomial for hSx. It will be unchanged for any other hN and can also be used for N = 22 since $S_{22}$= $S_{11}$ and the S[9] tile is a surrogate N = 11 with known height. Here inside N = 11, hSx = p(x)·1 ≈ 0.0014424.

To construct an exact Sx using $p_1$ as reference: (i) MidPointSx = {$p_1$ + hSx·Tan[4π/11], –1} ≈ {–6.46283513413515913 44, -1} (ii) cSx = MidPointSx + {0,hSx}
(iii) rSx = RadiusFromHeight[hSx,11] (iv) Sx = RotateVertex[cSx + {0,rSx},11,cSx]

Under the outer- billiards map $\tau$, the edges of tiles such as Sx are part of the singularity set W, which is also known as the 'web'. This set will be defined in Section 4. For a regular N-gon it can be obtained by mapping the extended edges of N under $\tau$ or $\tau^{-1}$. Below is a portion of W in the vicinity of Sx. Note that Sx has a clone obtained by rotation about the center of Px. In the limit this web is probably multi-fractal. See Example 5.4 and [H6]. (Click on the main image or Sx to download larger versions.)

**Example 1.3** The singularity set W local to Sx for N = 11 (two clicks here)

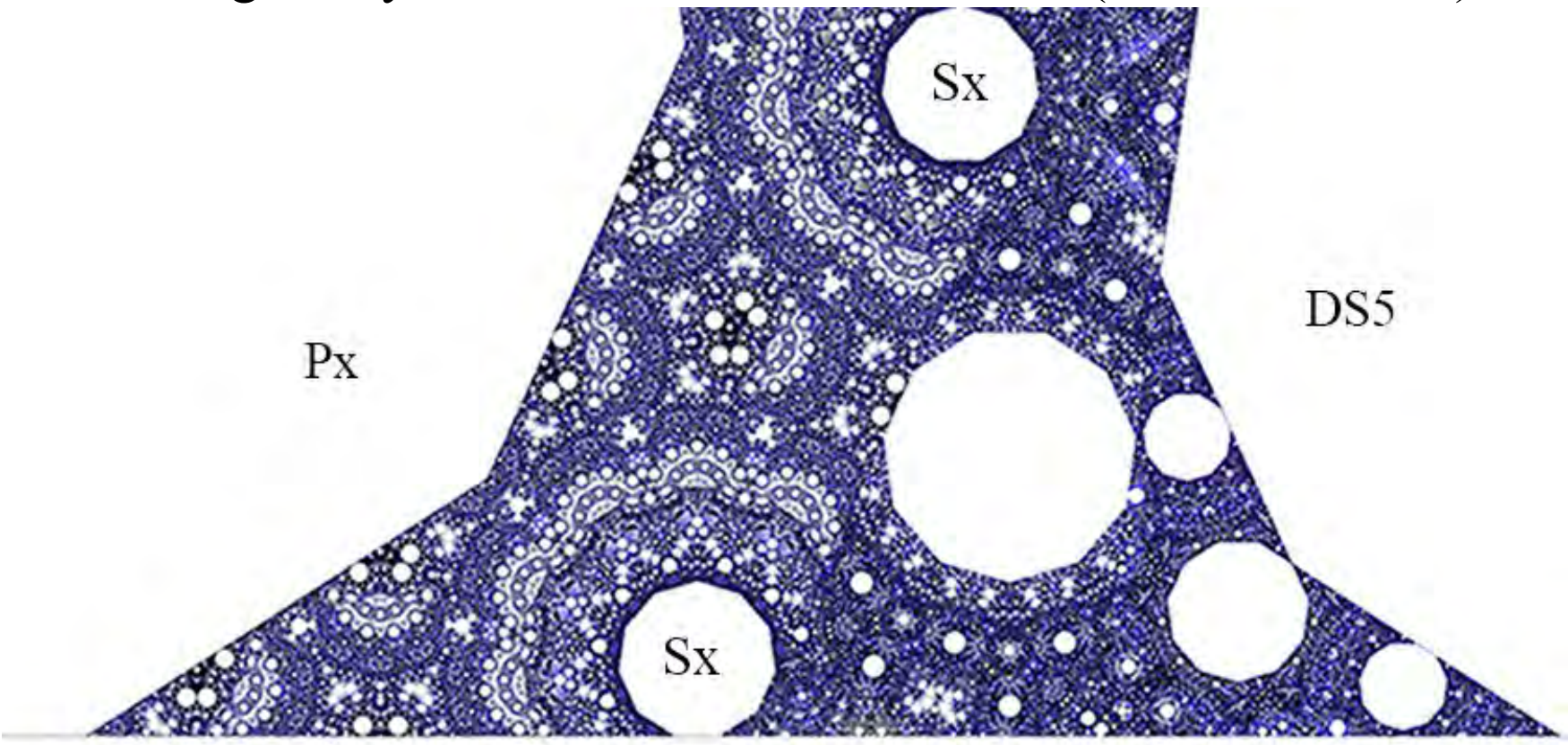

It is our contention that all the polygonal 'tiles' (regular or not) which arise from the outer billiards map of N are defined by scales which lie in the scaling field of N. This is the Scaling Conjecture in Definition 3.2. For N = 11, knowing the exact parameters of '3rd-generation' tiles like Sx allow us to probe deeper into the small-scale structure of N = 11 , which is almost a total mystery. But each new generation may involve a significant increase in computational complexity because the star points of each generation depend on the previous, and these have the same complexity as N. So the true computational complexity of Sx may be on the order of $\varphi(N)^3$ and this implies that it may not be feasible at this time to probe deeper than 10 or 12 generations. Since the generations scale by GenScale[11] = Tan[π/11]·Tan[π/22] ≈ .042217 the 25$^{th}$ generation would be on the order of the Plank scale of $1.6\cdot10^{-35}$ m. In the words of R. Schwartz, "*A case like N = 11 may be beyond the reach of current technology.*"

## Section 2. Conforming Regular Polygons

The S[4] tile from Example 1.1 is 'conforming' to the bounds of the star polygon of N =14 because it shares the 'base' edge of N and the (right-side) GenStar is star[1] of N as shown below. Clearly there are an infinite number of such conforming N-gons for N = 14, but S[4] also shares star[4] (and star[5] ) of N and there are only 6 such 'strongly conforming' tiles for N = 14. These will constitute the nucleus of the First Family of N = 14. It should be clear that only the odd S[k] will have conforming gender duals shown here in magenta.

**Example 2.1** The 6 conforming S[k] tiles of N = 14 in cyan and gender duals in magenta

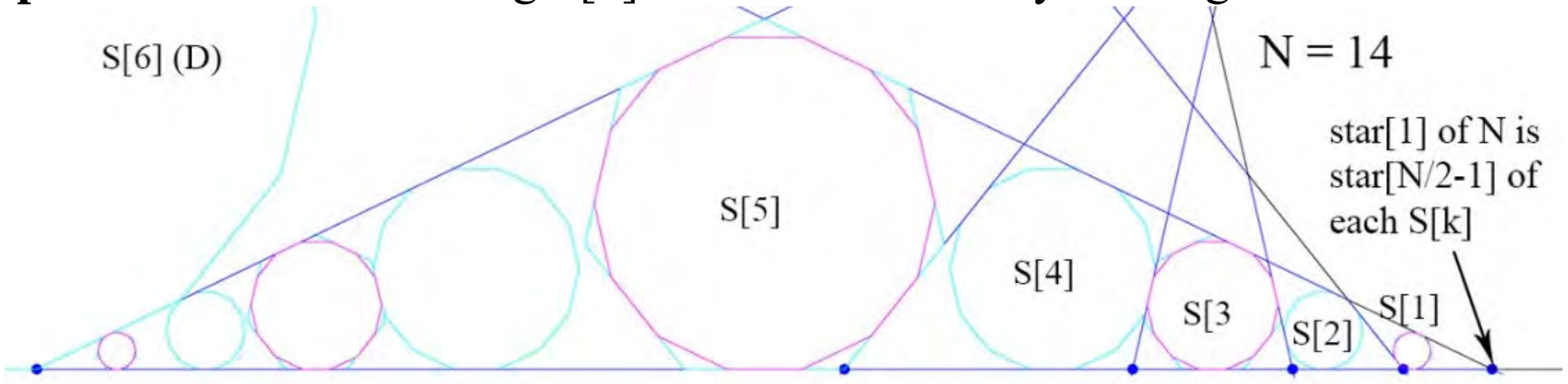

When N is odd, the conforming tiles will again be defined relative to the star[1] point of N as shown below for N = 7. But now these conforming S[k] tiles must be 2N-gons which share their penultimate star point with N, and they cannot be replaced with their duals while still being conforming (Note that S[3] here can serve as N = 14 so the local S[k] are the same as above.) The strongly conforming S[1],S[2] and S[3] tiles of N = 7 will be the nucleus of the First Family.

**Example 2.2** The 3 conforming S[k] tiles of N = 7

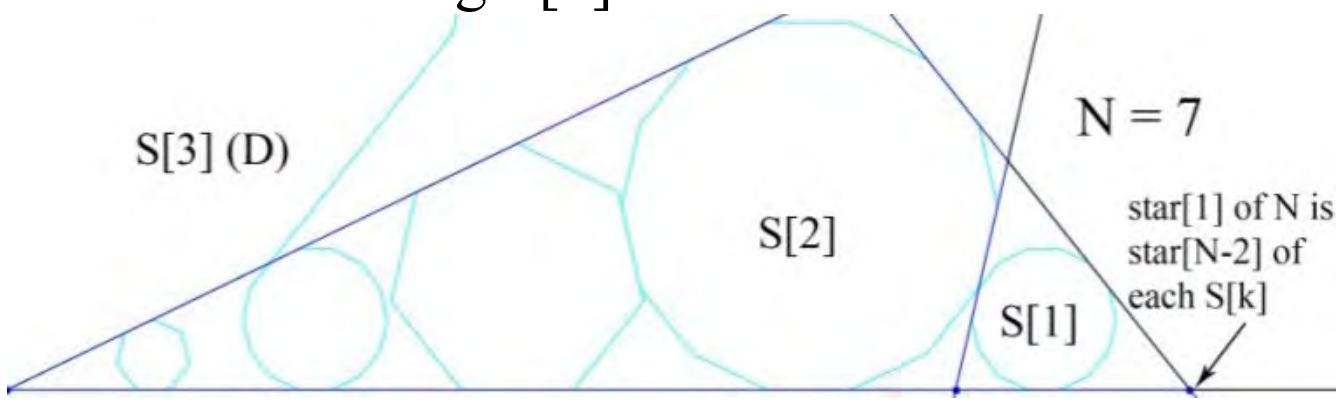

**Definition 2.1** (i) If N is even and P is a regular N-gon or N/2-gon, P is *conforming* relative to N if P shares the same base edge as N and GenStar of P is star[1] of N. (ii) If N is odd and P is a regular 2N-gon then P is *conforming* relative to N if P shares the same base edge as N and star[N–2] of P is star[1] of N. In both cases P is said to be *strongly conforming* if it is conforming and also shares another star point with N.

**Lemma 2.1** (Conformal Replication) Every regular N-gon has a strongly conforming $D_N$ tile which is identical to N for N even and a regular 2N-gon with same side as N for N odd.

Proof: Set star$D_N$[1] = GenStar[N] and center offset sN/2. When N is even this center and side length defines a regular N-gon $D_N$ which is identical to N. By the reflective symmetry of N, $D_N$ must have GenStar[$D_N$] equal to star[1] of N. When N is odd, again use GenStar[N] and offset (sN)/2 but now construct a regular 2N-gon $D_N$. Because the exterior angle of $D_N$ is half of the exterior angle of N, $D_N$ will have star[N-2] equal to star[1] of N so it is strongly conforming. □

**Example 2.3** The left-side and right-side D tiles for N/2 where N = 26

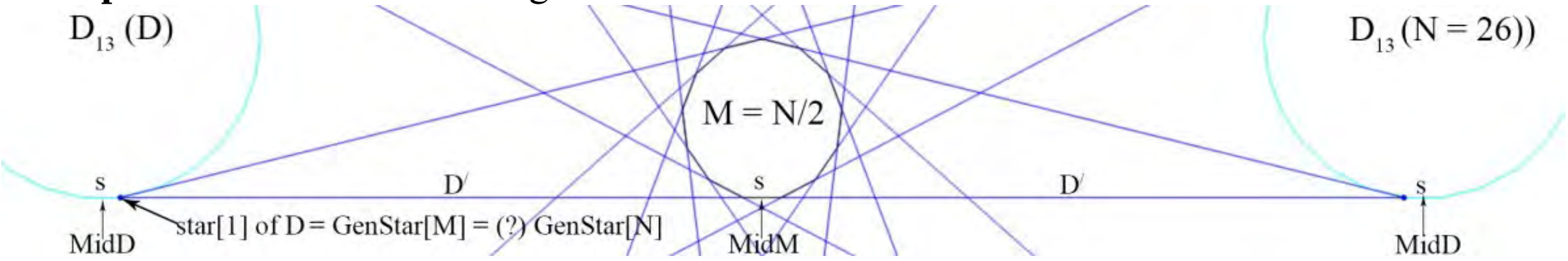

**Lemma 2.2** (Twice Odd Lemma: part 1) Suppose N is twice-odd and M is N/2. For M at the origin define D and $D_R$ to be the left and right-side $D_M$ tiles. Then if the origin is shifted to $D_R$ acting as N, $D_N$ will be D, so N/2 and N share the same D tiles when N is twice-odd.

Proof: It is sufficient to show that GenStar[M] = GenStar[N]. Relative to MidM, GenStar[M] has horizontal displacement $D' = hM\cdot s_{<M/2>}$ as shown above. Relative to N at the origin GenStar[N] has horizontal displacement $hN\cdot s_{N/2-1}$ so it is necessary to show that $hN\cdot s_{N/2-1} = 2D' + s/2 = 2hM\cdot s_{<M/2>} + hN\cdot s_1$. Since N is a D tile relative to M, they have the same side so hM/hN = scale[2] of N, so $hM\cdot s_{<M/2>} = \text{scale}[2]\cdot hN\cdot s_{<M/2>}$ and by definition, $s_k$ of M is $s_{2k}$ of N, so $hM\cdot s_{<M/2>} = \text{scale}[2]\cdot hN\cdot s_{<N/2>} = \text{scale}[2]\cdot hN\cdot s_{N/2-1}$. Therefore we need to show that $s_{N/2-1} = 2\cdot s_{N/2-1}\cdot\text{scale}[2] + s_1 = 2\cdot s_{N/2-1}\cdot s_1/s_2 + s_1$ . For N even $s_{N/2-1} = 1/s_1$ so this reduces to showing that $2/s_2 + s_1 = 1/s_1$ and this is equivalent to the 'half-angle' formula $s_2 = 2s_1/(1- s_1^2)$. □

**Corollary 2.1** (Corollary to Lemma 2.1) For a regular N-gon, there is a conforming regular N-gon or 2N-gon P with $hP \le hD_N$.

Proof. By Lemma 2.1, there is a conformal $D_M$ for any regular N-gon M embedded in N which shares the base edge and star[1] of N. This $D_M$ will also be conformal relative to N because N and M have the same exterior angle so a conforming $D_M$ must exist for any $hD_M \le hD_N$. □

**Theorem 2.1** (First Family Theorem (FFT)) For a regular N-gon every star[k] point defines an S[k] tile which is strongly conforming and has horizontal displacement –sN/2 relative to star[k].

Proof: Suppose that N is even with $p_1$ = starN[k]. Let P be the conforming regular N-gon with center displacement –sN/2 relative to $p_1$ as shown here. Such a P must exist by Corollary 2.1. We will show that P must have starN[k] as a star point and starP[k′] = starN[k] for k′ = N/2–k. (k′ will be called the 'local index'of P.) We will show that when the displacement of MidP from starN[k] is sN/2 then k′ = N/2–k. ( Here k = 4 so k′ =3 so star[4] of N should be star[3] of P.)

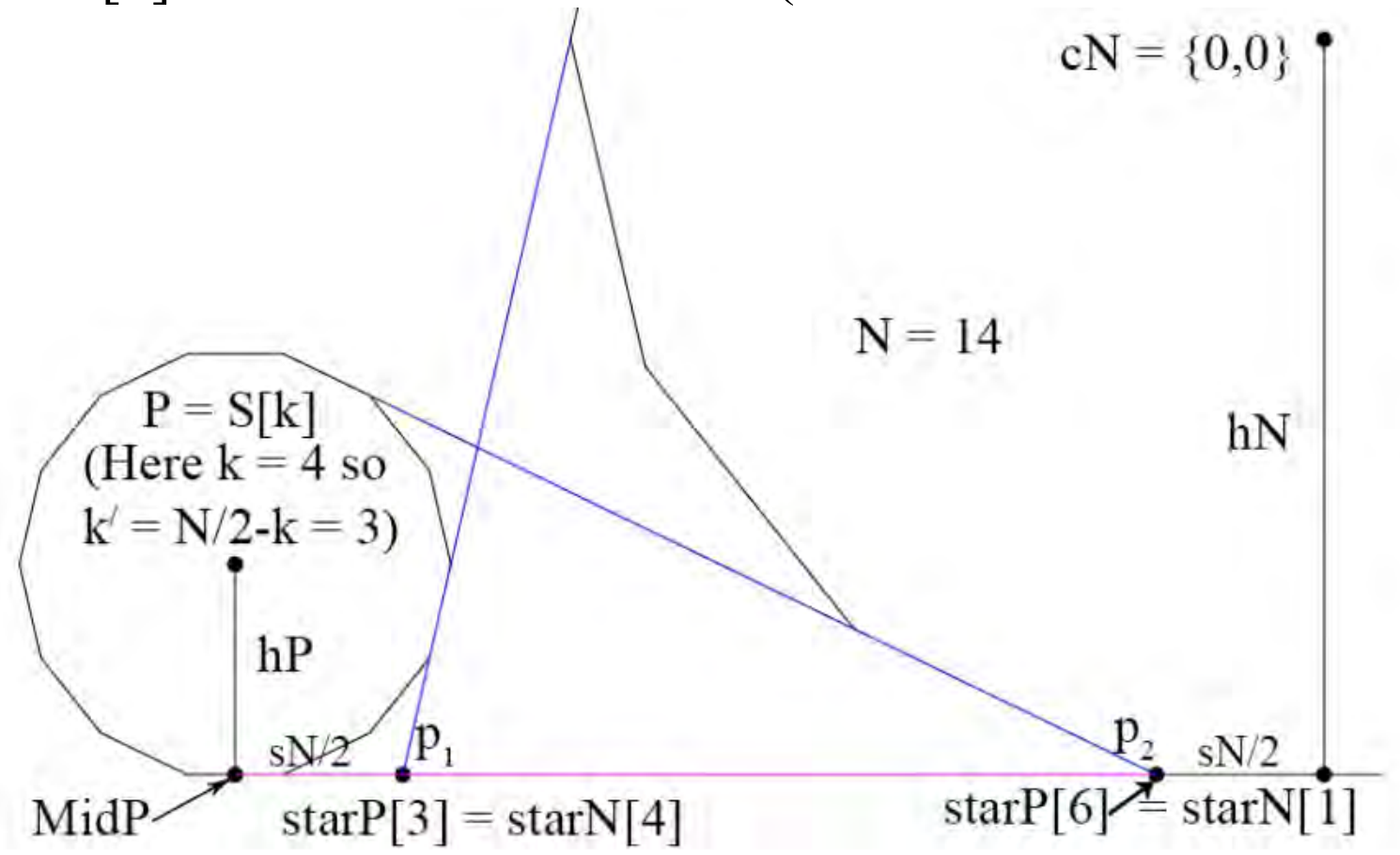

Note that there is symmetry between P and N because these two have equal and opposite displacements of sN/2 from their respective star points. This implies that $hPs_{N/2-1} = hNs_k$ and defines P as S[k] of N (with k = 4 here). If P is known this same formula defines N as S[N/2-1] of P.

Every conforming P has (right-side) starP[N/2–1] = star[1] of N, so they all share a global 'index' of N/2–1, and here we claim that (right-side) starP[k′] is star[k] of N, so the local index of P is k′ = N/2–k. (This 'retrograde' form of the local index is due to the fact that by convention the star points and S[k] are numbered right to left while the star polygon angles increase left to right, so for D at S[N/2–1], k′ = 1 which corresponds to the minimal star angle of $\phi = 2\pi/N$.)

Any horizontal displacement on the base line can be written from the perspective of N or P and these must match. For the magenta interval above from starN[1] to MidP, the perspectives are:
(a) Relative to MidN, it is the displacement of star[k], namely $hN \cdot s_k$
(b) Relative to MidP this displacement is $hP \cdot s_{N/2-1}$ (always the same form for any conforming P)

Therefore $hN \cdot s_k = hP \cdot s_{N/2-1}$ so $\dfrac{hP}{hN} = \dfrac{s_k}{s_{N/2-1}} = \dfrac{s_1}{s_{N/2-k}} = s_1 \cdot s_k$

These last two equalities follow from fact that for N even, $s_{N/2-k} = 1/s_k$ as noted in Lemma 1.2. This implies that $hN \cdot s_1 = hP \cdot s_{N/2-k} = hP \cdot s_{k'}$ and the left side is sN/2, so the displacement of MidP from starP[k′] is sN/2 and this matches the given displacement of starN[k]. Therefore starP[k′] is starN[k] of N as claimed. For k > 1 all the S[k] will be strongly conforming because they share two distinct star points with N - namely star[k] and star[N/2–1] so the theorem is true for N even with the exception of k = 1.

When k = 1 the local index is the same as the global index but the calculations above are valid by setting $p_1 = p_2$, so $\frac{hS[1]}{hN} = \frac{Tan[\pi/N]}{Tan[(N/2-1)\pi/N]} = \frac{s_1}{s_{N/2-1}} = s_1^2$. (Recall that this is GenScale[N])
Therefore S[1] can be constructed without the help of star[2], but S[1] always has star[2] at index 2 because the displacement of star[2] from MidS[1] is $hN\cdot(2s_1 - s_2)$ and this is always $hS[1]\cdot s_2 = hN\cdot s_1^2\cdot s_2$ by the half-angle formula $s_2 = 2s_1/(1-s_1^2)$.

Therefore all of the S[k] tiles are strongly conforming and share the same sN/2 offset so the First Family Theorem is true for N even. The odd case is a simple doubling of the indices from the even case because the conforming P are now 2N-gons with the same sN/2 displacements as the even case. Therefore the parameters of the (unmutated underlying) S[k[ are:

For N even: $\frac{hS[k]}{hN} = \frac{s_1}{s_{N/2-k}} = s_1\cdot s_k$ of N

For N odd: $\frac{hS[k]}{hN} = \frac{s_2}{s_{N-2k}} = s_2\cdot s_{2k}$ of 2N

This completes the First Family Theorem for arbitrary N. □

Definition 2.2 below will define the set of S[k] that make up the First Family for arbitrary N. Then we will prove a Generalized First Family Theorem (GFFT) that describes the local families of the S[k]. This would not be necessary if these local families were just scaled versions of the First Families, but this is only true in the linear or quadratic cases of N = 3,4,5,6,8,10 and 12.

Example 4.3 will show that the S[k] evolve in a multi-step fashion in the web W. The steps will correspond to the local index of S[k] which is k′ = N/2–k for N even or 2k′ for N odd. This occurs because the 'star angles' formed by the extended edges of N are of the form $k\phi'$ with $\phi = 2\pi/N$. These angles determine the recursive evolution of each S[k]. Only the D tiles will have a 'normal' step-1 evolution, as seen below for N = 22 where the local family of D is the same as N.

**Example 2.4** The S[k] tiles of N = 22 in black and (right-side) local families of the S[k] in cyan.

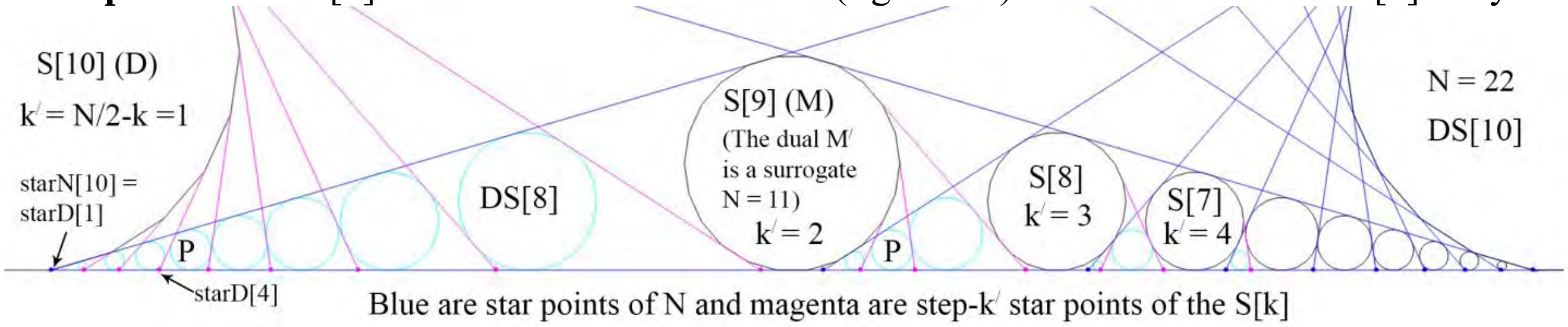

Blue are star points of N and magenta are step-k′ star points of the S[k]

When N is even, $hS[N/2-1]/hN = s_{N/2-1}\cdot s_1 = 1$ so D here is congruent to N, and it is identical to the $D_N$ of the Replication Lemma since it is strongly conforming with displacement sN/2 from GenStar[N]. In this N-even case, the central penultimate S[N/2-2] is called $M_N$ or simply M.

When N is twice-odd M will remain conforming when replaced with the N/2 gender-dual M′. As seen in Example 2.1, this 'mutation' can be applied to any odd S[k] (with even k′) without violating the conformity of the FFT. M at k′ = 2 is special because Lemma 2.2 shows that any D tile of N will have a matching N/2 tile which will share the same D tiles as N. Since M′ and N/2 share their left and right side GenStar points, M′ can be regarded as a surrogate N/2-gon.

**Definition 2.2** (First Family) For any regular N-gon, the First Family Theorem (FFT) defines the strongly conforming S[k] tiles for 1 ≤k ≤ <N/2>. These tiles will be called the First Family Nucleus and as defined earlier the S[<N/2>] tile will be called D and for N even S[N/2-2] will be called M. For N >4, the S[1] & S[2] tiles of D will be called DS[1] & DS[2] or M[1] & D[1].

(i) For N twice-even, the First Family will consist of the S[k] in the First Family Nucleus together with the (right-side) First Family Nucleus of D – called the DS[k].

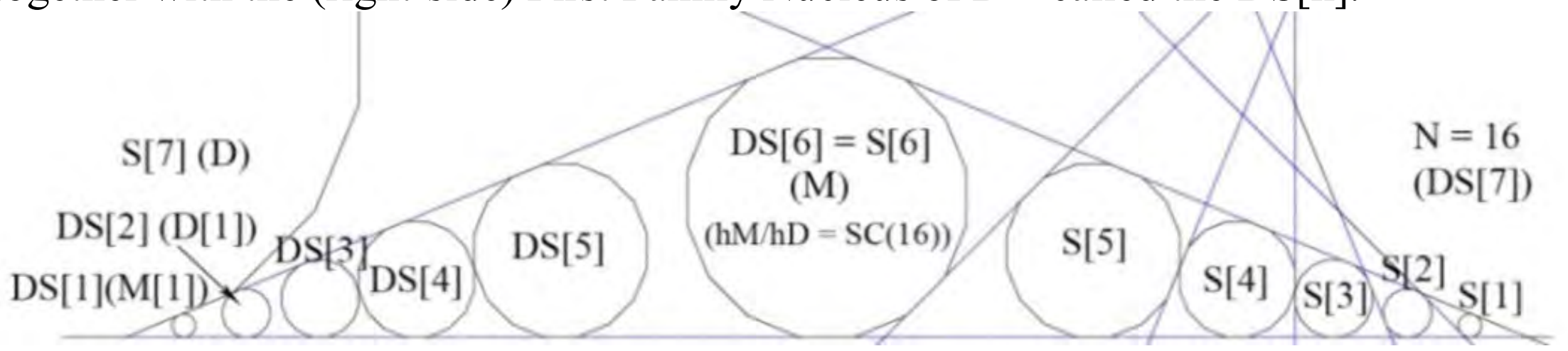

(ii) For N twice-odd, the odd S[k] in the First Family Nucleus will be replaced with their N/2 counterparts, and the First Family will consist of this revised First Family Nucleus together with the (right-side) revised First Family nucleus of D - called the DS[k].

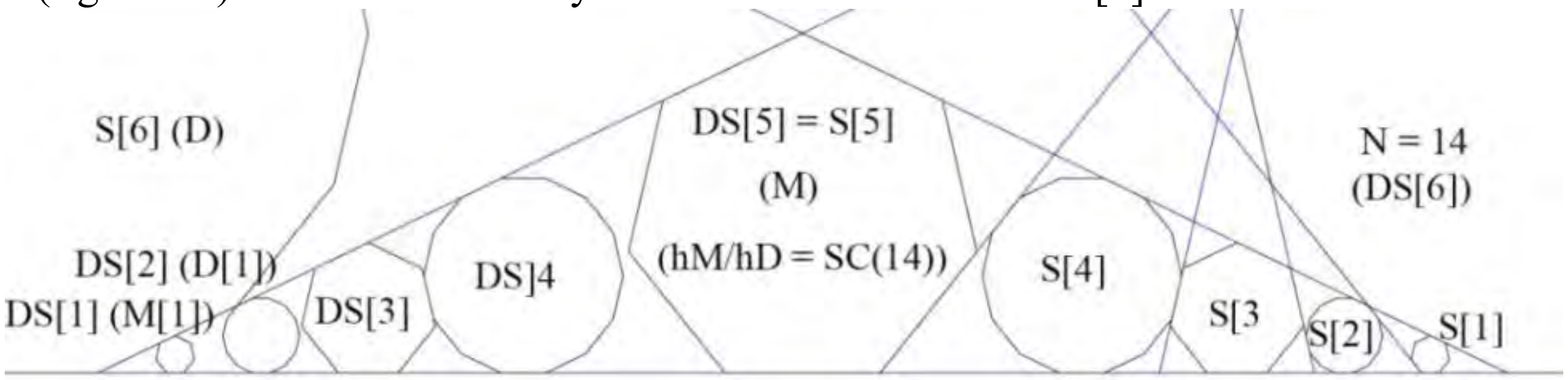

(iii) For N odd, the First Family will consist of the First Family Nucleus together with the (right-side) revised First Family Nucleus of D - called the DS[k]. We also include the reflections of these tiles about cN. With this convention all First Families span from D to matching D

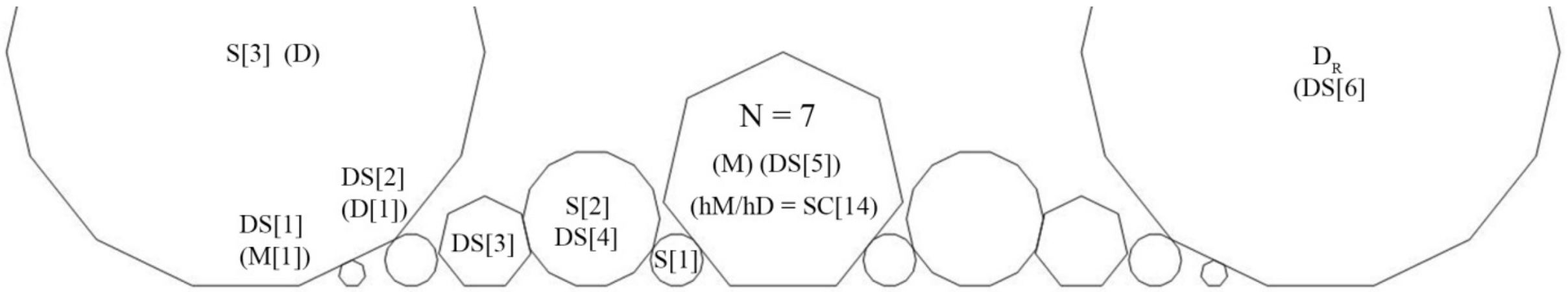

All these families run from D to matching D and involve a fundamental M-D relationship which defines the overall structure of the family. This will define future ‘generations’(if they exist). In general parts (ii) and (iii) will be identical except for scale. The only issue here is the S[1] tile of N = 7 which is technically not part of the First Family of N = 14. This tile is regarded as part of the ‘local family’ of N = 7 inside N = 14. We will define these local families in the GFFT and if desired they can be included as an extended N = 14 family. The mapping that transforms any twice-odd N-gon family to the matching N/2-gon family is a simple gender-change scaling followed by a translation of origin from N to M: **T[x] = TranslationTransform [{0,0}– cS[N/2–2]] [x]/SC[N,N/2]]**. This mapping and its inverse make up the Twice-Odd Lemma 4.2.

To return to the case of N = 22 above, note that the S[2] (P) tile in the local family of M is congruent to the S[4] tile of D or N. It is true in general that the local family of M consists of S[k] which are congruent to S[2k] of N. This is a special case of the GFFT because M has $k' = 2$. For S[8], $k' = 3$, so the local family will be congruent to mod-3 S[k] of N (and by extension this must include the neighbor S[7] as well as N at S[10]). In the web, S[8] will look like a normal N-gon because $\gcd(k',22) = 1$. For N = 18, S[3] at $k' = 6$ will be mutated into two interwoven $N/\gcd(k',N)$-gons at different radii. For N twice-prime the only be mutations are when $k'$ is even.

The GFFT below is a direct consequence of the observation made in the FFT about the symmetry between N and the S[k]. This implies that for N = 22, every S[k] will have N in its local family at S[10] and this step-$k'$ pattern must apply to all local tiles. (We do just the N-even case here.)

**Lemma 2.3** (Generalized FFT (GFFT)) For N-even every S[k] defines a 'family' of secondary tiles conforming to star[$k'$] of S[k] where $k' = N/2–k$. These family tiles will always include N at S[N/2–1] and step down by steps of $k'$. (We show in Section 4 why these are indeed the surviving (right-side) star points of S[k] in the web W.) Therefore the effective star points of S[k] will be subsets of the normal star points of S[k] of the form star[$k''$] where $k'' = (N/2–1) – j\cdot k'$ for $j \geq 0$ and the corresponding tiles will be congruent to S[$k''$] of N with offset sN/2 from star[$k''$] so the $k' = 1$ case will reproduce the (right-side) S[k] of D from the FFT.

Proof: Suppose S[k] is at the origin with star[N/2–1] which is star[1] of some $P_0$ tile and conversely this $P_0$ tile is conforming to star[$k'$] of S[k]. Then if $P_0$ and S[k] have equal (and opposite) displacements from their respective star points (as in the FFT), $hS[k]\cdot s_{N/2-1} = hP_0 \cdot s_k$ and $hS[k] = hN\cdot s_1\cdot s_k$ by the FFT, so $hP_0\cdot s_k = hN\cdot s_1\cdot s_k\cdot s_{N/2-1}$ and $hP_0/hN = s_1\cdot s_{N/2-1} = 1$. Therefore $P_0$ = N and N is congruent to S[N/2–1] of S[k].

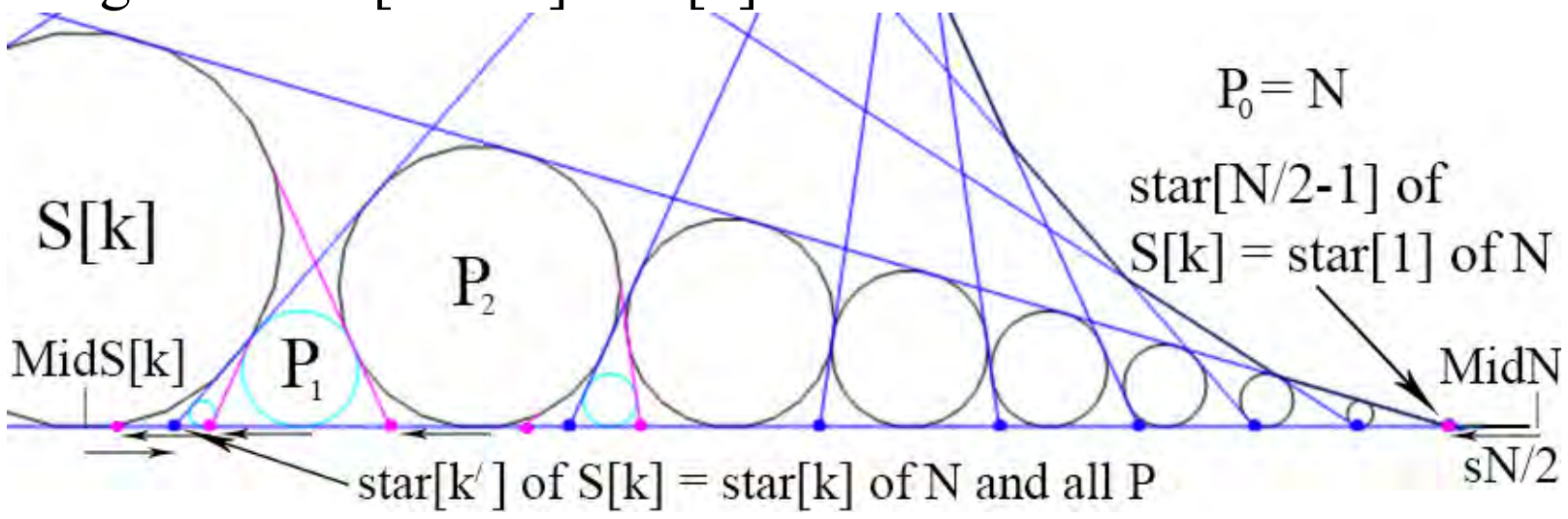

This argument will remain valid for comforming $P_k$ iff $P_k$ and S[k] have equal (and opposite) displacement from their respective starpoints. For S[k] the midpoint displacement from star[$k'$] is defined to be sN/2, so any hypothetical $P_k$ must have the opposite displacement from the defining star[$k''$] point of S[k].Suppose that star[$k''$] is an effective star point of S[k]. The conformal replication lemma says that there is a regular N-gon $P_k$ conforming to star[$k'$] of S[k] with displacement sN/2 from star[$k''$] as shown above.Therefore $hP_k\cdot s_k = hS[k]\cdot s_{k''} = hN\cdot s_1\cdot s_k\cdot s_{k''}$ so $hP_k/hN = s_1\cdot s_{k''}$ and by the FFT, $P_k$ is congruent to S[$k''$] of N. □

Setting j = 1 in the GFFT, $k'' = k–1$ so every S[k] has the right-side neighbor S[k-1] in its GFFT local family. Here we have concentrated on the right-side families of each S[k] with the undestanding that the left-side geometry must match the neighboring S[k+1]. All of these predicted secondary S[k] will be conforming to a star[k] of N, so their centers will be colinear with the origin as illustrated by the magenta lines in the web plots below.

**Example 2.5** In Section 4 we we will merge these First Families and secondary GFFT families with the 'web' of N. This is a portion of the web for N = 22 showing the GFFT families of S[k].

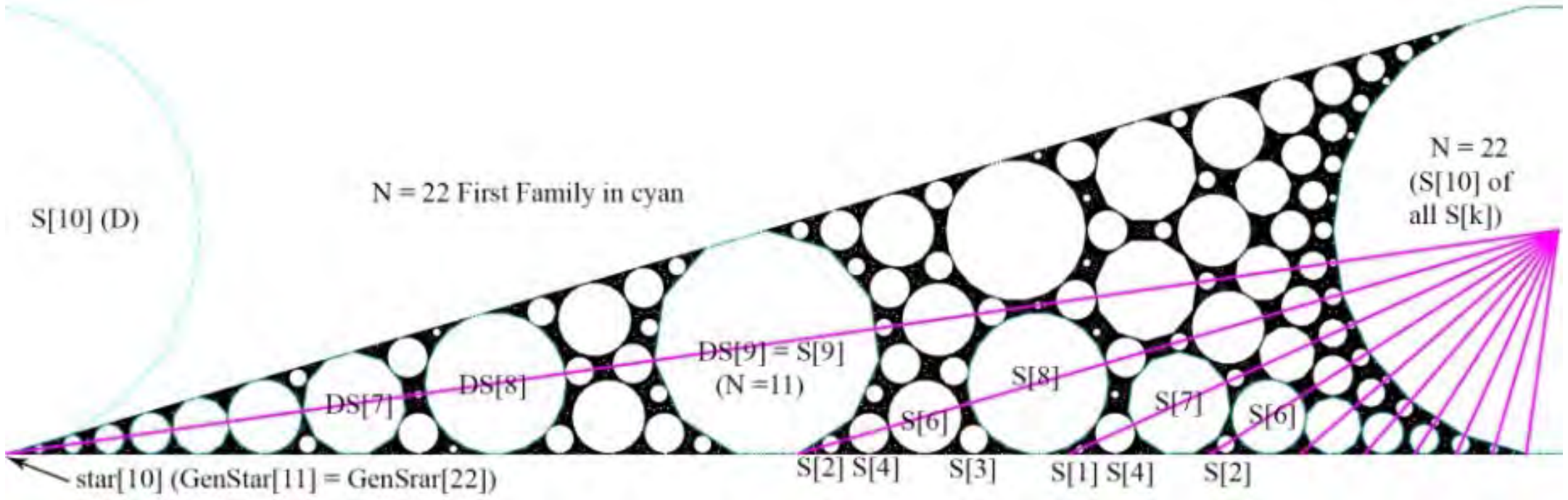

Since N is twice-odd the penultimate S[9] 'M' tile can serve as a scaled copy of N = 11 and we do not need an 'odd' version of the GFFT. Since k′ = N/2–k = 2, S[9] will have a step-2 family counting down from N at S[10]. It is no coincidence that these are also the S[k] predicted by the FFT. In the same fashion, the S[8] tile of N = 22 will have a step-3 right-side family, but the left-side family must be consistent with the right-side of S[9]. Here a copy of S[6] is shared by S[8] and S[9] and this determines the S[3] at the base of S[8] to complete the step-3 left-side family.

**Example 2.6** The local GFFT families of the S[25] tile of N = 60

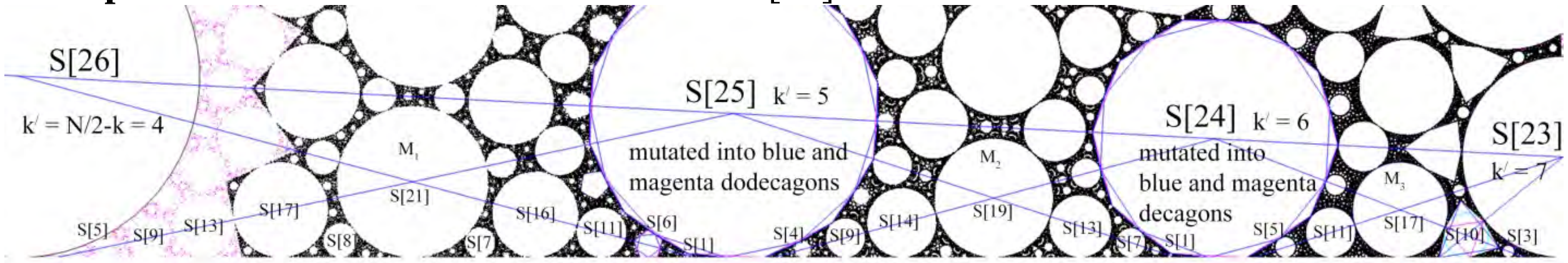

As explained in Appendix II, the web W for every regular N-gon is the disjoint union of invariant regions and the number of such regions is on the order of the algebraic complexity of N, namely $C_N = \varphi(N)/2 - 1$ so for N = 60 there are 7 such regions and the black region shown above is the first such region, with S[25] playing the role of 'shepherd' . See Fig.11 of App. II.

The Invariance Conjecture states that the transition between regions is always mediated by pairs of edge-sharing S[k]-S[k'] tiles with k' = N/2-k. We call them 'dual' pairs and a S[k] is called 'primitive' with gcd[k,N) = 1. For N = 60 the magenta-black dual pair is S[13] and S[17] of S[26] where S[17] has k' = 13 and likewise S[13] has k' = 17 with both indices primitive. The GFFT right-side family of S[26] is step-4 and it must be compatible with S[25] which evolves first, so counting down mod 4 it includes S[17] and S[13] which share an edge. The next transition for N = 60 will be mediated by S[11] and S[19] and $L_{60}$ = {1,13,11,7} describes the first 4 transitions (in order) up to the S[N/2-2] 'M' tile of N. By symmetry the remaining 4 transitions can be determined a simple reflection about cM, so $L_{60}$ above is sufficient to describe all 8 transitions. The most important transition is the first one mediated by S[1] and N. For N even they are always a dual pair since N is DS[N/2-1] and conversely. This is the 'template' for all edge-sharing and the rings of maximal D tiles that surround all regular N-gons are linked by common vertices as well as shared S[1] tiles. It is important to understand that these invariant regions are only those with occur on the scale of the S[k] tiles so they constitute just a first step toward a true theory of invariance - which typically occurs on all scales.

**Section 3. Scaling of a Regular N-gon**

In Section 4 we will introduce the outer-billiards map and prove some basic facts about dynamics, but here we continue to concentrate on the geometry of regular polygons to obtain results that could be applied to any mapping. The most important result of this section will be the Scaling Field Lemma 3.3 which shows that the primitive scales for a regular N-gon will generate the maximal real subfield of the cyclotomic field $\mathbb{Q}_N$ of N. This is what we call the scaling field.

In Section 1 we defined the (canonical) scales of a regular polygon N to be of the form scale[k[ = $s_1/s_k$ for $1 \leq k < N/2$ where $s_k = \text{Tan}[k\pi/N]$. The primitive scales are those with gcd(k,N) = 1

For N = 60 above, there will be 29 scales and the scaling field will be generated by the $\varphi(N)/2$ = 8 primitive scales so N = 60 is said to have 'algebraic complexity' 8. Scale[1] is always 1 and the minimal scale[N/2-1] is called GenerationScale[60] $\approx$ 0.002746. We will see that this is the scale of (ideal) sequences of S[1] and S[2] tiles converging to star[1] of N (or D).

In this twice-even case it is not clear how the geometry of N = 60 is related to N = 30 or any other factor, but for all factors j, $\mathbb{Q}_{N/j} \subseteq \mathbb{Q}_N$. N has a natural embedding inside each factor polygon that preserves center and height, so the star[1] point of N/j aligns with star[N/j] of N and this makes it easy to compare their scales as in Lemma 1.1 earlier.

**Example 3.1** The 9 factor polygons for N = 60

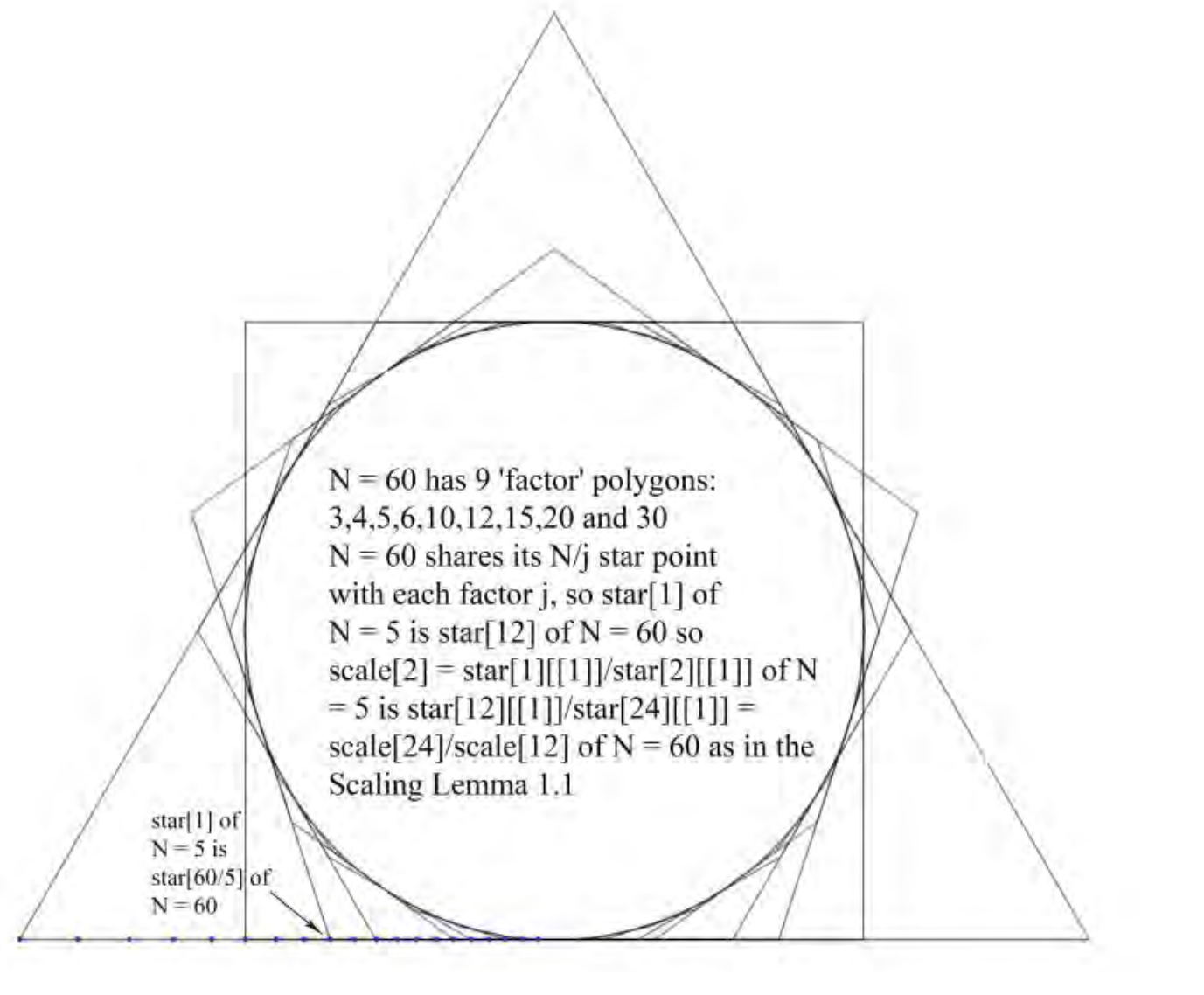

Since the scales are defined as ratios of the $s_k$ 'star' points, any N-gon must share scales with its 'factors' because N can be regarded as being embedded in each factor with shared height and center. For any factor k in this embedding, the star[N/k] point of N must be the star[1] point of N/k and Lemma 1.1 says that scale[j] of N/k = scale[kj]/scale[k] of N.

This is illustrated here with k = 12 and the matching factor polygon N = 5 where scale[2] = $s_1/s_2$ = star[1][[1]]/star[2][[1]] of N = 5. This must also be star[12][[1]]/star[24][[1]] of N = 60 which is scale[24]/scale[12].

For N = 30 with factor 2, scale[j] = scale[2j]/scale[2] of N = 60, so scale[2] of N = 60 converts one scale to the other. This scale[2] relationship between N and N/2 applies to the twice-odd case also and there it guarantees that N/2 is a valid member of the First Family of N.

The First Family Theorem defines the heights of the S[k] tiles in terms of the $s_k$ – which we loosely call 'star' points.

**Lemma 3.1 (First Family Scaling)** For all regular N-gons, hS[1]/hS[k] = scale[k] of N

Proof: For N-even hS[1]/hS[k] = $s_1^2/s_1 \cdot s_k$ = scale[k] of N. For N odd hS[1]/hS[k] = $s_2^2/ s_2 \cdot s_{2k}$ of 2N which is by definition $s_1/s_k$ = scale[k] of N.□

Since D is always the S[<N/2>] tile of N, hS[1]/hD will always be the minimal scale[<N/2>] which we define to be GenScale[N]. When N is even, D is congruent to N, so hS[1]/hN = scale[N/2-1] = GenScale[N] = $s_1^2 = \mathrm{Tan}^2[\pi/N]$. But when N is odd, D is twice-odd so hDS[1]/hD is still $s_1^2$ of D but it does not make sense to compare N-gons with 2N-gons, so our convention is to treat D temporarily as an N/2-gon (with the same side length) by simply dividing by scale[2] of 2N. Then hS[1]/hD = $s_1^2/(s_1/s_2) = s_1 \cdot s_2$ of D = GenScale[N]. This will give the correct scaling for both N and D in the twice-odd case while $s_1^2 = \mathrm{Tan}^2[\pi/N]$ may fail to be a unit.

Based on these results, it may seem that S[1] is the natural choice for a 'next-generaton' tile, but we will discover in Sections 4 and 5 that the most fruitful path to a conforming next generation is by using S[1] and S[2] as a pair and combining their webs together. This will give a combined web which is step-4 relative to S[2] in the even case and step-8 relative to S[2] in the odd case. Both of these webs will have potential to foster next generation tiles on the edges of S[2] playing the role of a second generation D[1]. S[1] will still have a largely autonomous local web.

The following Lemma shows that GenScale scaling is inherent in the star polygons because every N-gon can be regarded as an S[1] of some upscaled D′ tile. We just do the odd case here.

**Lemma 3.2:** Every regular N-gon is the S[1] tile of a D′ tile where hN/hD′ = GenScale[D]

Proof: There is no loss of generality in assuming that D′ shares the horizontal edge of N and has star[1] matching star[1] of the (right-side) D tile of N, and sD′ = $2hN \cdot s_{<N/2>}$.

We show in Lemma 2.2 that MidM is the midpoint of the extended edge of D as shown below, so sN/sD = $2hN \cdot s_1/2hN \cdot s_{<N/2>} = s_1/s_{<N/2>} \equiv$ GenScale[N]. Since this is a gender mismatch, hN/hD′ = GenScale[N]/scale[2] = GenScale[2N] = $s_1^2$ (of 2N) so by the FFT, N has the correct scale to be the S[1] of D′ and by definition it has displacement sD′/2 from star[1], so N is a valid S[1] tile of D′.

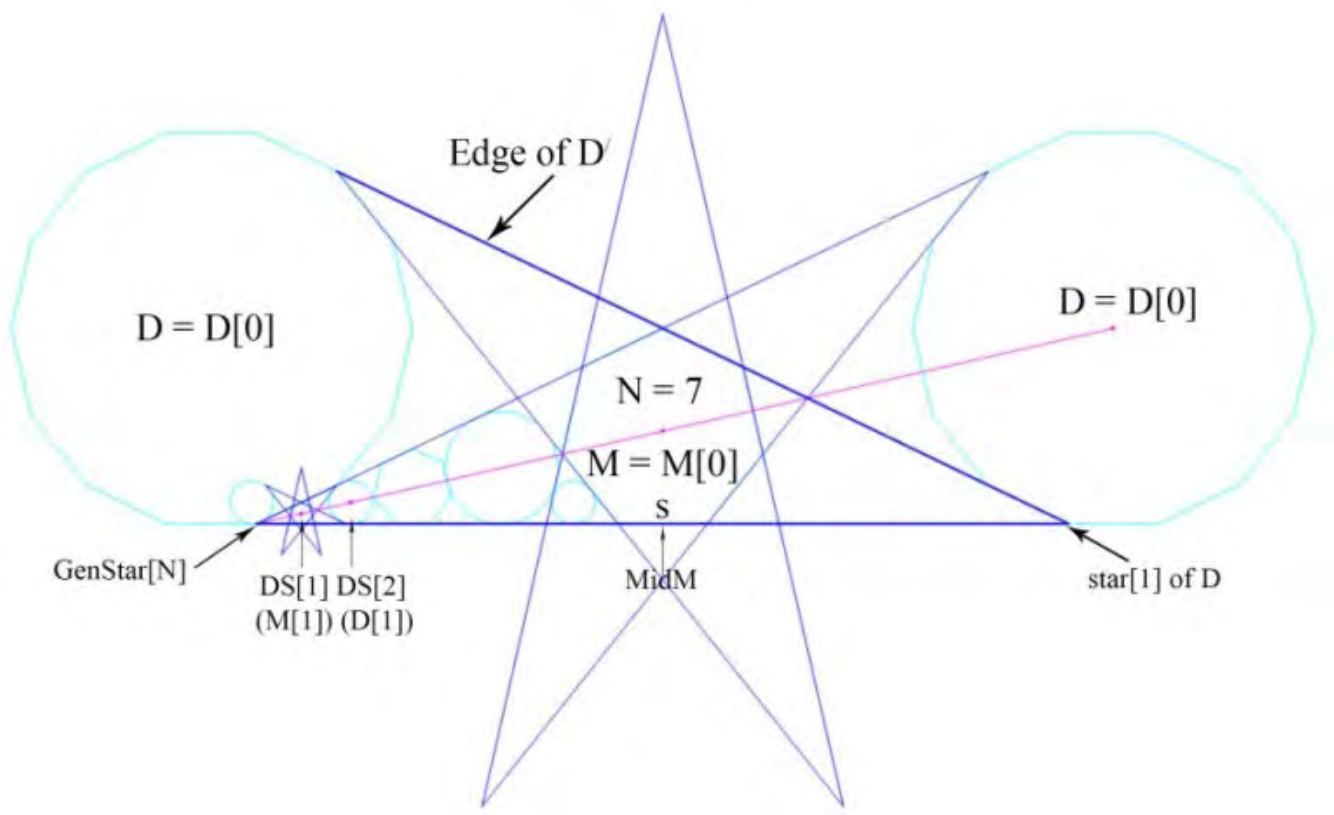

**The cyclotomic field of N**

**Definition 3.1** For a regular N-gon, the cyclotomic field $\mathbb{Q}_N$, is the algebraic number field which can be generated by $\zeta = \exp(2\pi i/N)$ so $\mathbb{Q}_N = \mathbb{Q}(\zeta)$ which can be shown to have order $\varphi(N)$ over $\mathbb{Q}$. Since $\zeta + \zeta^{-1} = 2\cos(2\pi/N)$, and complex conjugation in $\mathbb{Q}_N$ is always an automorphism of order 2, $\mathbb{Q}(\zeta + \zeta^{-1})$ is order $\varphi(N)/2$. This is called the maximal real subfield and denoted $\mathbb{Q}_N^+$. The full field $\mathbb{Q}_N$ is always a quadratic extension of $\mathbb{Q}_N^+$. As a vector space, $\mathbb{Q}_N$ is the direct sum of its real and imaginary parts $\mathbb{Q}_N^+$ and $\mathbb{Q}_N^-$.

Because of the Lemma below $\mathbb{Q}_N^+$ is also called the *'scaling field'* of N – written $S_N$. Since $\mathbb{Q}_N$ is equivalent to $\mathbb{Q}_{N/2}$ for N twice-odd, $S_N = S_{N/2}$.

**Lemma 3.3** (Scaling Field Lemma) For any regular N-gon, the maximal real subfield $\mathbb{Q}_N^+$ has a unit basis consisting of the primitive (canonical) scales.

**Proof**: For a given N we have defined $s_k = \tan(k\pi/N)$ for $1 \le k < N/2$. These 'star points' are classified as 'primitive' if (k,N) = 1 and 'degenerate' otherwise. Here we show that the set of primitive scales, $T = \{t_K = s_1/s_k: (k,N) = 1\}$ is a basis for $\mathbb{Q}_N^+ = \mathbb{Q}_N \cap \mathbb{R}$. $\mathbb{Q}_N$ can be generated by any 'primitive' N-th root of unity of the form $\zeta^k = \exp(2k\pi i/N)$ with (k,N) = 1 so the indices of the primitive scales are also indices of primitive roots of unity.

Since $i\tan(\theta) = (e^{i2\theta}-1)/(e^{i2\theta}+1)$, $is_k = (e^{2k\pi i/N}-1)/(e^{2k\pi i/N}+1) = (\zeta^k-1)/(\zeta^k+1)$ so $is_k$ is in $\mathbb{Q}_N$ and scale[k] = $\tan[\pi/N]/\tan[k\pi/N] = is_1/is_k = \left[\frac{\zeta-1}{\zeta+1}\right]\left[\frac{\zeta^k+1}{\zeta^k-1}\right]$ is in $\mathbb{Q}_N^+$. We will show that when (k,N) = 1, scale[k] is an algebraic integer in $\mathbb{Z}[\zeta]$ and its inverse coscale[k] = $\cot[\pi/N]/\cot[k\pi/N]$ is also an integer. Regrouping terms, scale[k] $=\left[\frac{\zeta-1}{\zeta^k-1}\right]\left[\frac{\zeta^k+1}{\zeta+1}\right]$. The term on the left is called a 'cyclotomic unit' in $\mathbb{Z}[\zeta]$, but we will show from first principles that this product is a unit. Since (k,N) = 1, there is a (rational) integer j such that kj =1(ModN). Therefore

$$\left[\frac{\zeta-1}{\zeta^k-1}\right] = \left[\frac{(\zeta^k)^j-1}{(\zeta^k)-1}\right] = \sum_0^{j-1}(\zeta^k)^i \text{ and when k is odd } \left[\frac{\zeta^k+1}{\zeta+1}\right] = \left[\frac{(-\zeta)^k-1}{(-\zeta)-1}\right] = \sum_0^{k-1}(-\zeta)^i$$

If k is even N must be odd, so repeat the above with $k' = k+N$. This substitution leaves scale[k] and coscale[k] unchanged. Therefore scale[k] is an algebraic integer.

For coscale[k] replace $\zeta$ with $-\zeta$ so when k is odd the two quotients are $\sum_0^{j-1}(-\zeta^k)^i$ and $\sum_0^{k-1}\zeta^i$.

When k is even, again replace k with $k' = k+N$ which will be odd, and the inverse j can be chosen odd so that it will preserve the sign change. Therefore scale[k] is an algebraic unit in $\mathbb{Q}_N^+$ with inverse coscale[k].

To show that the set T of primitive scales forms a basis for $\mathbb{Q}_N^+$, note that $|T| = \varphi(N)/2$ because (k,N) = 1 implies that (N– k, N) = 1. It only remains to show that the primitive scales are independent over $\mathbb{Q}$.

Suppose that $\sum_{1\le k<N/2} a_i t_k = 0$ with (k,N) =1 and $a_i \in \mathbb{Q}$, then $\frac{1}{s_1}\sum_{1\le k<N/2} a_i t_k = \sum_{1\le k<N/2} a_i r_k = 0$ where $r_k = \cot[k\pi/N]$ with (k,N) = 1. This contradicts the Siegel-Chowla result that the primitive $r_k$ are independent over $\mathbb{Q}$. Therefore the primitive scales are a unit basis for $\mathbb{Q}_N^+$. □

**Example 3.2** (N = 11) Set $\zeta = \zeta_{11}$ then scale[4] = $\frac{\mathrm{Tan}[\pi/11]}{\mathrm{Tan}[4\pi/11]} = \left[\frac{\zeta-1}{\zeta+1}\right]\left[\frac{\zeta^4+1}{\zeta^4-1}\right] = \left[\frac{\zeta-1}{\zeta^{15}-1}\right]\left[\frac{\zeta^{15}+1}{\zeta+1}\right]$

$$= \left[\frac{(\zeta^{15})^3-1)}{\zeta^{15}-1}\right]\left[\frac{(-\zeta)^{15}-1}{-\zeta-1}\right] = \sum_{j=0}^{2}(\zeta^{15})^j\cdot\sum_{j=0}^{14}(-\zeta)^j = (1+\zeta^4+\zeta^8)\cdot\sum_{j=0}^{14}(-\zeta)^j$$

coscale[4] = $\frac{Cot[\pi/11]}{Cot[4\pi/11]} = \sum_{j=0}^{2}(-\zeta^{15})^j\cdot\sum_{j=0}^{14}\zeta^j = (1-\zeta^4+\zeta^8)\cdot\sum_{j=0}^{14}\zeta^j$

This shows a relationship between scale[4] and $\zeta^4$ which goes back to star[4] and S[4]. Since $\zeta^4$ is a generator of $\mathbb{Q}_{11}$ this relationship can be made exact:

**AlgebraicNumberPolynomial[ToNumberField[scale[4], Exp[2*Pi*I/11]^4],x]** yields

scale[4] = $-1+2[x^2+x^3-x^5-x^6+x^8+x^9]$ where x = $\zeta^4$ , Since $S_{11} = \mathbb{Q}_{11}^+$ is generated by $\lambda_{11}$ = 2cos(2π/11) the integers in the field $\mathbb{Q}(\lambda_{11})$ must be in the ring $\mathbb{Z}[\lambda_{11}] = \mathbb{Z}[\zeta + 1/\zeta]$ where scale[4] simplifies to $3-4x-2x^2+2x^3$. Among the scales, the best match for $\lambda_{11}$ is GenScale[11] (scale[5]) at $1-6x^2+2x^4$. We will often use GenScale[N] as a 'surrogate' generator for $S_N$ when scaling tiles. This will yield more meaningful polynomials for the S[k].

**Lemma 3.4** (Alternate generators of the scaling field $S_N$)
(i) For all N, $s_1^2 = \tan(\pi/N)^2$ is a generator of the scaling field $S_N$
(ii) When N is twice-even, $s_1^2$ is GenScale[N] which is a unit generator of $S_N$.
(iii)When N is odd , $s_1\cdot s_2$ of 2N is GenScale[N] and it is a unit generator of $S_N$ and $S_{2N}$

Proof (i) Since $\cos(2\theta) = \frac{1-\tan^2(\theta)}{1+\tan^2(\theta)}$ and $\lambda_N$ = 2cos(2π/N) is always a generator of $S_N$, $\tan(\pi/N)^2$ is also a generator. (ii) When N is twice-even $\tan(\pi/N)^2$ is GenScale[N] and it is clearly primitive so it is a unit by Lemma 3.3. (iii) When N is odd $\lambda_N$ = 2/(GenScale[N] +1) so GenScale[N] is a generator of $S_N$ and a primitive scale so it is a unit generator of $S_N$ which is identical to $S_{2N}$. □

Therefore our convention for (unit) generators of $S_N$ will be to use GenScale[N] for N twice-even and when N is odd GenScale[N] will serve as a unit generator for both $S_N$ and $S_{2N}$. (Any primitive scale of N should also be a generator of $S_N$ but the proof is not obvious.)

The classical paper by J.C. Calcutt [Ca] proves the Complexity Theorem for cos, sin and tan from first principles and develops 'Chebyshev-like' polynomials for tan that highlight the close connection between the $s_k$ and $\zeta^k$. For example $\mathbb{Q}(s_1) = \mathbb{Q}(s_k)$ for primitive $s_k$ and the order of $s_1$ over $\mathbb{Q}$ is $\varphi$(N) – except when 4|N it is $\varphi$(N)/2. In [H4] we show that the primitive $s_k$ are the positive Galois conjugates of $s_1$ and $-s_1$. These cases are distinct for N = 24. Click to see them.

John Stillwell [St] used results such as these and the half-angle formulas relating tan to sin and cos to say "*the tan function may be considered more fundamental than either cos or sin*".

The degenerate scales of N are still in $\mathbb{Q}_N^+$ so any such scale will be a linear combination of the primitive scales. This should be true for the degenerate $s_k$ themselves and we prove this fact here. Since the $s_k$ are typically not in $\mathbb{Q}(\zeta)$ we will continue to use the $is_k$ as surrogates.

**Corollary 3.1** For a given N, every $s_j$ for $j < N$, is a linear combination of the primitive $s_k$.

Proof: Suppose $(j,N) = m \geq 1$, then $(j/m,N/m) = 1$ and $s_{j/m}$ of N/m = $s_j$ of N. Since $\mathbb{Q}_{N/m} \subseteq \mathbb{Q}_N$, $is_j \in \mathbb{Q}_N^-$ and the set $iS = \{is_k : (k,N) = 1, k < N/2\}$ forms a (vector space) basis for $\mathbb{Q}_N^-$ because $|iS| = \varphi(N)/2$ and the $is_k$ are independent over $\mathbb{Q}$ (iff the $s_k$ are independent). Therefore $is_j$ is a linear combination of the primitive $is_k$ which implies that $s_j$ is a linear combination of the $s_k$. □

**Example 3.3** For a case like N = 18, known summation formula such as those found in [P] can be used to find the coefficients. For example $\tan(\pi/18) + \tan(\pi/18 + \pi/3) + \tan(\pi/18 + 2\pi/3) = \tan(6\pi/18) = \sqrt{3}$. Note that all the factors on the left side are primitive in N = 18, with k values 1, 7 and 13 (= – 5), so $\tan(\pi/18) + \tan(7\pi/18) - \tan(5\pi/18) = \sqrt{3}$. In general it can be very difficult to find these coefficients and the same applies to coefficients of the scales.

Even though $s_1^2 = \tan^2(\pi/N)$ will always be a generator of $S_N$, it may not be a unit or even an integer when N is twice-odd. Therefore it would not be a good choice of generator. The First Family geometry shows clearly the gender issue discussed earlier in this section.

**Example 3.4** (The First Family for N = 6 ) The First Family Theorem says that hS[1]/hN = scale[2] = GenScale[6] = $\tan^2(\pi/6) = 1/3$, but this is a gender mismatch and not an algebraic integer, so a better choice of scaling is hS[2]/hN = GenScale[3] = $\tan(\pi/6)\cdot\tan(\pi/3) = 1$. This avoids the gender mismatch and assigns the correct scaling to the web W.

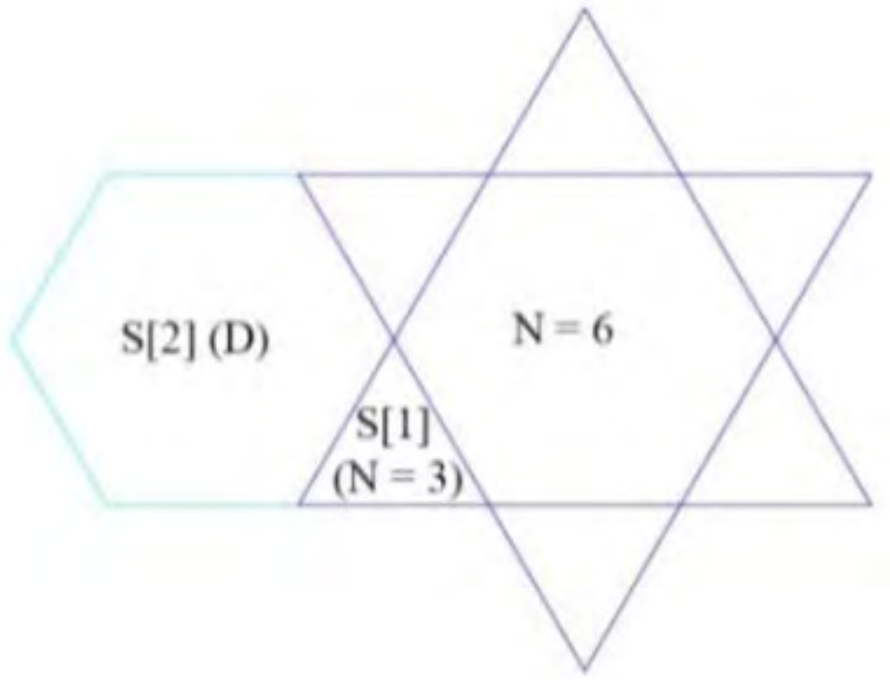

For N = 3, 4 and 6, $\varphi(N)/2 = 1$ so they have linear algebraic complexity. This means that there will be no accumulation points in the outer-billiards web W because the web consists of rays or segments parallel to the sides of N and under the outer-billiards map, these segments are bounded apart by linear combinations of the vertices. When N is regular the vector space determined by the vertices has the same rank as $\mathbb{Q}_N$ - namely $\varphi(N)$. Therefore the coordinate space of N-gons with linear complexity will have rank 2 over $\mathbb{Q}$ and hence affinely rational coordinates. Since the outer-billiards map is itself an affine transformation these 'rational' coordinates are preserved and the web is also affinely rational – with no limit points.

**Definition 3.2** (Canonical polygons) Every regular N-gon defines a coordinate system, and any line segment or polygon P (regular or not) that exists in this co-ordinate system will be called *canonical* relative to N if for every side sP, the ratio sP/sN is in $S_N$. The Scaling Conjecture says that all tiles and line segments which arise in the web W of the outer-billiards map are canonical.

**Example 3.5** Any line segment defined by a linear combination of $s_k$ of N is canonical because if T is such a linear combination $T/s_1$ is a linear combination of dual scales, so it is in $S_N$. Therefore the multi-step S[k] mutations in Section 4 will be canonical because the edges are the difference of two star points. Our convention for scaling N-gons is to scale them relative to another N-gon and in this case side scaling is the same as height scaling.

**Lemma 3.5** For a regular N-gon, the First Family S[k] tiles are canonical.

Proof: For N even: $\dfrac{hS[k]}{hN} = \dfrac{s_1}{s_{N/2-k}} = s_1 \cdot s_k$ of N - which is scale[N/2– k] of N.

For N odd: $\dfrac{hS[k]}{hN} = \dfrac{s_2}{s_{N-2k}} = s_2 \cdot s_{2k}$ of 2N = scale[N–2k]/scale[2] of 2N so it is in $S_{2N} = S_N$. □

Since the S[k] are canonical regular N-gons or N/2 gons, the (ideal) First Families and subsequent families generated by the S[k] will also be canonical regular N-gons or N/2 gons. The linear or quadratic cases like N = 3,4,5,6, 8, 10, and 12 have webs consisting of only scaled First Family tiles with possible multi-step mutations, so all tiles are canonical.

**Definition 3.3 (**Ideal generations at GenStar[N]**)** For any regular N-gon with N > 4 , define D[0] = D and M[0] = M (the penultimate tile of D). Then for any natural number k > 0, define M[k] and D[k] to be DS[1] and DS[2] of D[k–1] so for all k, M[k] will be the penultimate tile of D[k]. The (ideal) *kth generation* of N is defined to be the (ideal) First Family of D[k–1]. Therefore M[k] and D[k] will be 'matriarch' and 'patriarch' of the next generation - which is generation k+1. By Lemma 3.1, hM[k]/hM[k–1] will be GenScale[N], or GenScale[N/2] if N is twice-odd.

**Example 3.6** (Ideal generations at GenStar[N] for N = 9, where by convention M[0] = N)

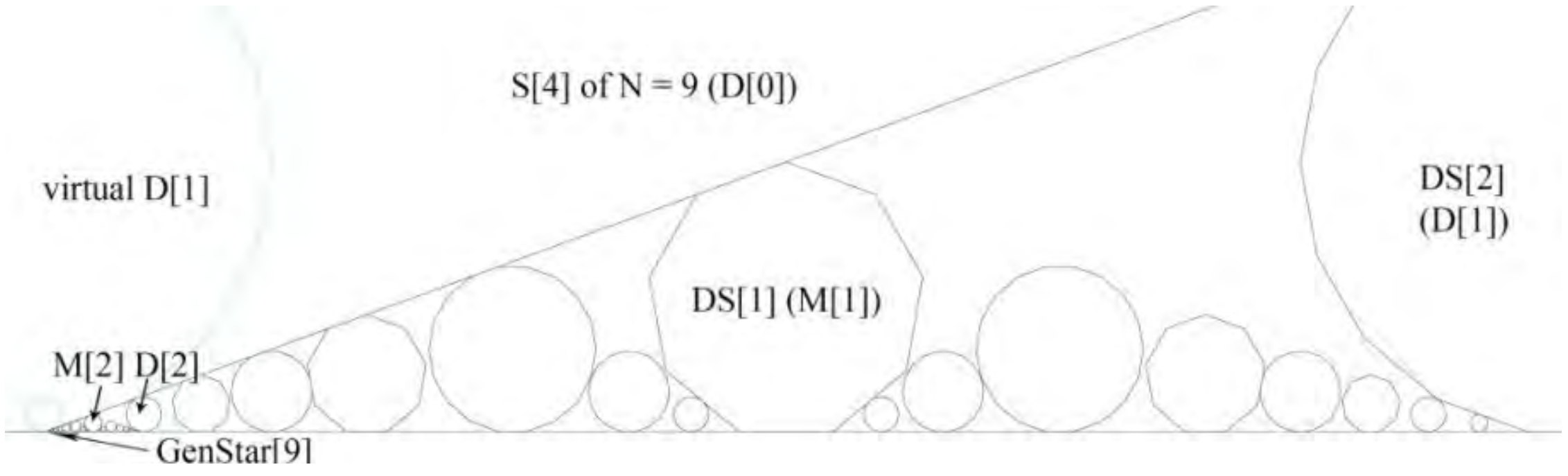

(The Mathematica code for these families is **FFM1 = TranslationTransform[cM[1]][First Family*GenScale]** and **FFM2 = TranslationTransform[cM[2]][FirstFamily*GenScale^2]**.)

## Section 4. The Outer-Billiards map

Except for motivation, the geometry of the first three sections is independent of any mapping, and here we will show connections between this geometry and the outer-billiards map $\tau$. In particular we will prove that for any regular N-gon, the First Family S[k] tiles will exist in the complement of the 'singularity' set W – which is defined in Definition 4.2 below. This implies that the S[k] arise naturally under iterations of $\tau$. This has long been observed in computer simulations, but now the exact parameters of the S[k] are known. See Lemma 4.1.

In truth it is impossible to properly define $\tau$ for 'most' regular N-gon because $\tau$ is a discontinuous mapping and no one knows exactly where it is defined, because we only know the limiting singularity set W up to approximation. However we do know how W evolves and the S[k] are an integral part of the early evolution of the 'web' W.

**Definition 4.1** (The outer-billiards map $\tau$) Suppose that P is a convex polygon in Euclidean space with origin internal to P. For an arbitrary point p external to P there will usually be a unique vertex c of P so that the line from p to c can serve as a (cw. or ccw.) 'tangent' line or 'support' line. Here we show the clockwise case. If this line is extended an equal distance past c, the new point will be a reflection about c, namely $\tau(p) = 2c - p$. If no such c exists, $\tau$ is not defined and p is an element of the singularity set W. All the blue 'trailing edges' of P are in W.

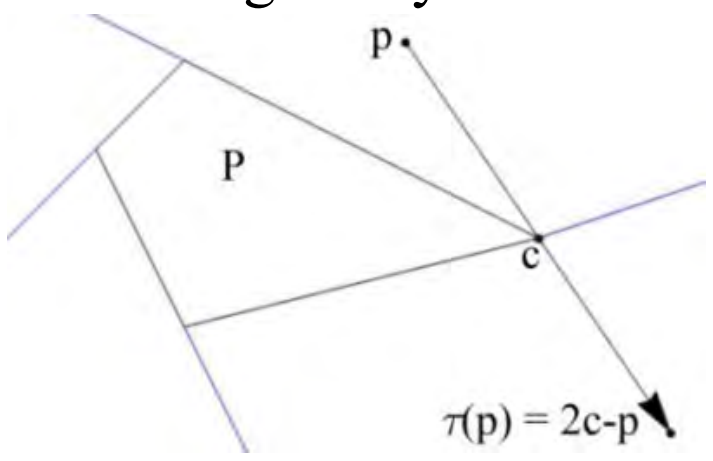

**Example 4.1** (S[2] for N = 7). By Lemma 4.1 to follow, cS[2] will have an orbit that advances two vertices of N on each iteration so it will be period 7. But $\tau$ is an inversion so only the center can have odd period. All other points in S[2] will be period 14 as shown in the right. After 14 iterations the original offset from the center will be negated by the second round of 7 iterations.

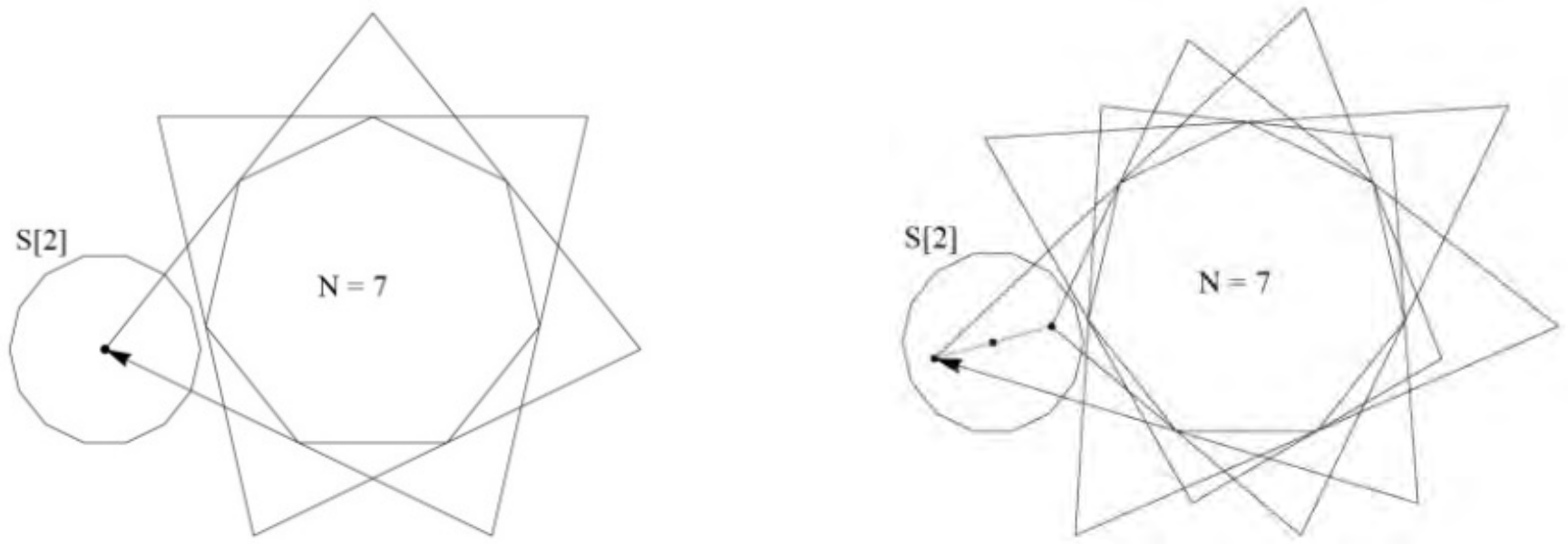

Note that $\tau^k(p)$ will always have the form $2Q + (-1)^k p$ where Q is a sum of k vertices of N with alternating sign. So the center point is period 7 iff $\tau^7(cS[2]) = 2Q - cS[2] = cS[2]$. This has only one solution which is Q = cS[2]. For all points in S[2] (including cS[2]), any 'second round' of 7 iterations will have the same Q (except for sign), so the concatenation yields Q = 0 and $\tau^{14}(p) = p$. This is a 'stable' period 14 orbit for all points except cS[2].

**Lemma 4.1** (Canonical orbits of S[k] centers) For any regular N-gon, cS[k] has an outer-billiards orbit that advances by k vertices on each iteration, so it has period N/(k,N) ≡ N/gcd(k,N).

Proof: We will show that there is a height/radius 'duality' between cS[k] and star[k] defined by **Du[x] = RotationTransform[-Pi/N,{0,0}][x*rN/hN]**. For any N-gon P centered at the origin, Du[P] will map midpoints of edges of P to the adjacent (cw) vertex of P. Therefore Du can be used to map any edge-based (magenta) orbit of a star[k] point to a (green) vertex-based orbit of Du[star[k]] as shown here for N = 10 and N = 7. (The green orbits are shortened.)

Every star[k] point is a vertex of an {N,k} star polygon embedded in {N,<N/2>}. Therefore this star polygon defines an 'edge-orbit' $O_k$ which coincides with the edges of N. This orbit advances k edges on each iteration so it has period N/(k,N). (Since k = 2 here the periods are 5 and 7.) Because {N,k} inherits the rotational symmetry of N, these magenta orbits will extend equal distance on either side of the midpoint of each edge, so $O_k$ also defines an outer-billiards orbit relative to M -which has vertices equal to the midpoints of N. Since Du[M] = N, it follows that Du[$O_k$] will be an outer-billiards orbit around N with period N/(k,N) as shown in green. Every vertex point in this orbit is Du[p] for some midpoint p in $O_k$. The initial point of these green orbits must be Du[star[k]] but it is not obvious that this point is also cS[k]. We will prove this.

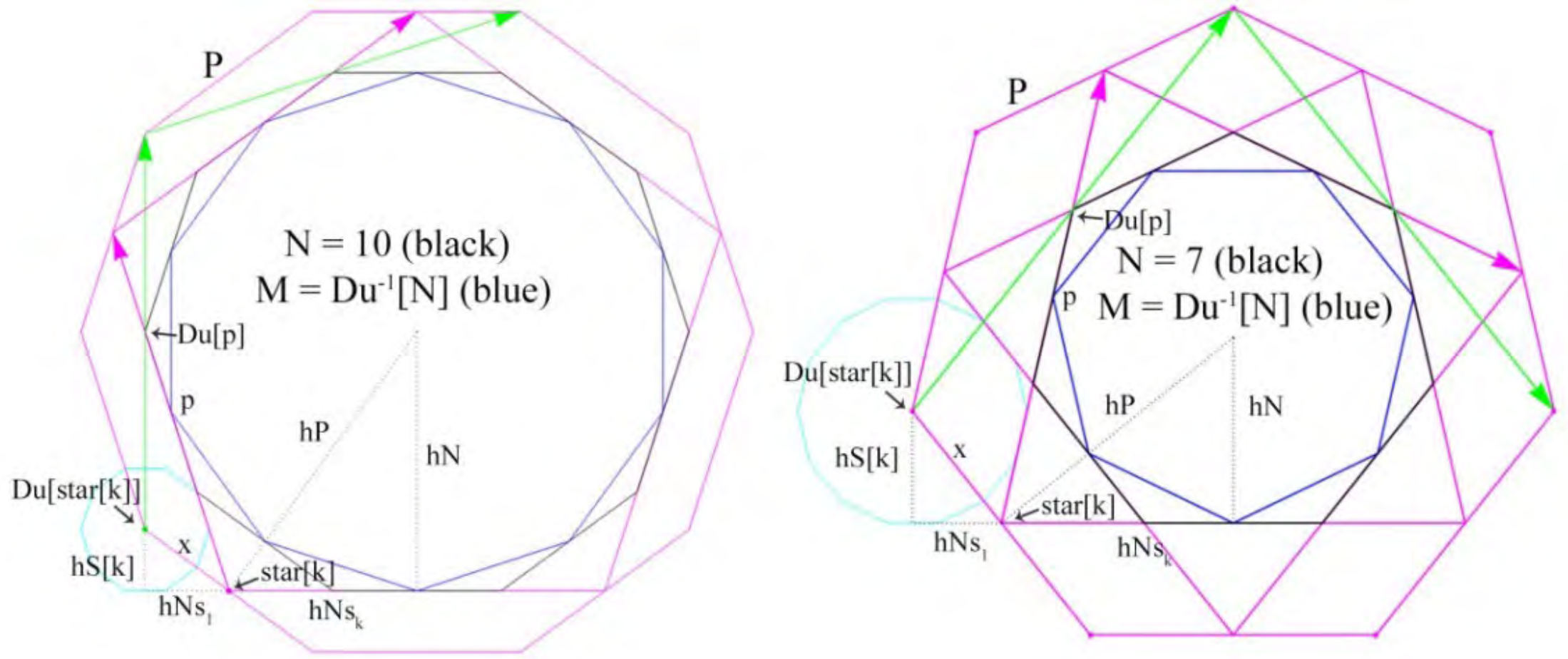

The polygons P were not needed earlier, but they will be useful here. P is defined to be the unique N-gon centered at the origin with star[k] as its midpoint. Therefore Du[star[k]] is a vertex of P and we will show that it must also be cS[k]. The length x shown here is sP/2 which is $hPs_1$ (where $s_k = \tan(k\pi/N) = \tan(k\pi/P)$. From the FFT, the displacement of cS[k] from star[k] is sN/2 = $hN\cdot s_1$. Therefore hS[k] should be $s_1\sqrt{hP^2 - hN^2}$ . But $hP^2 = hN^2(1 + s_k^2)$ so this implies that $hS[k] = s_1\sqrt{hN^2(1+s_k^2) - hN^2} = s_1 \cdot hN \cdot s_k$ .When N is even the FFT says that this is indeed hS[k] but when N is odd S[k] is a 2N-gon and $hS[k] = hN\cdot s_2\cdot s_{2k}$ of 2N. But $s_1$ of N is $s_2$ of 2N and $s_k$ of N is $s_{2k}$ of 2N, so once again $hS[k] = hN\cdot s_2\cdot s_{2k}$ of 2N as in the FFT. Therefore cS[k] is always Du[star[k]] and the period of the orbit of cS[k] will be N/gcd(k,N). □

We say that a regular tile P has period k if its center has period k. Of course all points in P must map together but as noted above, when N is even the off-center points may have 'period doubling' and visit P twice. This can actually be helpful because it can be used to identify the center. Of course any periodic tile will be in $W^C$ (the complement of W) which has full measure.

**Definition 4.2** (The outer-billiards singularity set W of a convex polygon P)
Let $W_0 = \bigcup E_j$ where the $E_j$ are the (open) extended edges of P. The singularity set W is defined to be $\bigcup_{n=-\infty}^{\infty} \tau^n(W_o)$ and for a positive index k this can be partitioned into the level-k (forward) web $W_k = \bigcup_{j=0}^{k} \tau^{-j}(W_0)$ and the level-k (inverse) web $W_k^i = \bigcup_{j=0}^{k} \tau^{j}(W_0)$. We define $-W = \lim_{k\to\infty} W_k$ and $+W = \lim_{k\to\infty} W_k^i$ .When cP is internal to P, a horizontal reflection about cP will convert $\tau$ to $\tau^{-1}$ so in this sense $\tau$ can be regarded as a 'chiral' map where $+W_0$ and $-W_0$ are related by a simple reflection. But under iteration they evolve independently so the geometric relationship between $-W_k$ and $+W_k$ is no longer a simple reflection. Therefor for the Penrose 'kite' of R. Schwartz [S1] he showed that both +W and –W defined fractal paths allowing points to escape to $\pm\infty$. For a regular N-gon the $D_N$ dihedral symmetry group of rotations $r$ and reflections $f$ is non-Abelian since $f \circ r = -r \circ f$ because N has an orientation. So now N is 'self-chiral' independent of $\tau$ and once again $+W_0$ is a simple reflection of $-W_0$ but the long-term relationship between +W and –W is still far from trivial and our Reflection Principle says that the overall web W must evolve both right and left handed.

**Example 4.2** (The +W star polygon web of N = 14 in blue and $-W_0$ in magenta).
For years we were under the delusion that +W would be sufficient and will eventually cover the magenta trailing edges in $-W_0$. This delusion was supported by plots such as these where the (truncated) forward and trailing edges are shown in blue and magenta. Here we generate the level-k (inverse) webs $W_k^i$ by iterating the blue forward edges under $\tau^k$ for k = 0,1,2,3 and 50. At this point the magenta trailing edges are nearly invisible. We did acknowledge that +W and –W were always distinct, but suggested that the union of $W_k$ and$-W_k$ was an optional step to enhance the images and in the limit it would not matter. This is false and these two must evolve together. With apologies here we swapped colors in the sidebar to match the Introduction, so -W is blue.

| The evolution of +W by iterating the 'forward' edges of N | +W and –W with colors swapped |
|---|---|
| 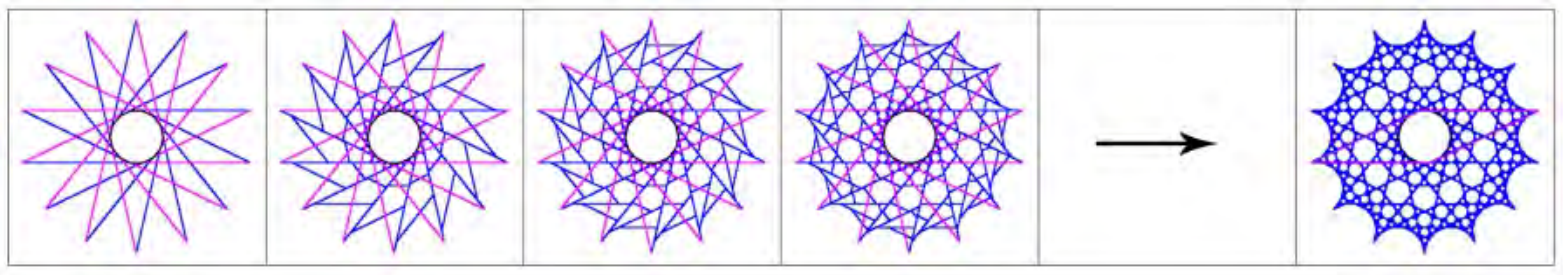 | 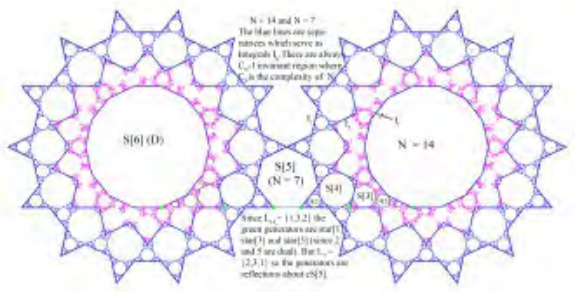 |

Here for N = 14 the S[5] tile can serve as a surrogate N = 7 because N = 14 and N = 7 share the same cyclotomic field . In terms of the invariant regions, the S[k] transition points for N = 7 are just a palindrome of the points for N = 14, so $L_{14}$ = {1,3,2} and $L_7$ = {2,3,1}.

**Lemma 4.2** (Twice-odd Lemma) Lemma 2.1 shows that for N twice-odd the First Families of N and N/2 are related by the scaling and change of origin between M and N given by **T[x] = TranslationTransform[{0,0}– cS[N/2–2]][x]/SC[N,N/2]]**. Here we show that N/2 and N must also have locally equivalent webs.

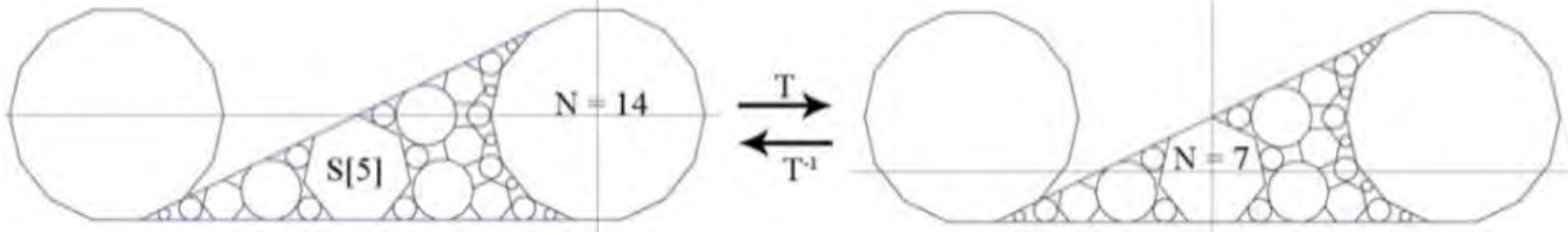

Proof: Because T and $T^{-1}$ are affine transformations they will commute with the affine transformation $\tau$ and hence preserve the web W. □

## Evolution of the Web: Part I

The web W can be formed by iterating the extended 'forward' edges of N under $\tau$ and here we will concentrate on the local web evolution of S[k] tiles with a more detailed analysis to follow in Part II. The main issue is that there are N-domains which map to each other and we want to describe how a single domain evolves. In the graphical iteration of a function of one variable the y = x line is used to swap domain and range. This is called 'cobwebbing'. For rational piecewise isometries such as $\tau$ this cobwebbing typically involves a translation or 'shear' followed by a rotation to swap domain and range and prepare for the next iteration. This is also a recipe for constructing N and the two mappings of the Appendix use a shear and rotation to mimic W.

We will show that W can be generated by an iteration of 'shear and rotate' where the shear is of constant magnitude sN and the rotation angle for each S[k] interval is the 'star[k] angle' $k'\phi$ with $\phi = 2\pi/N$ and $k' = N/2–k$ for N even and $k' = N–2k$ for N odd. The N even case is shown below.

**Example 4.3** N = 22 - showing the $2^{nd}$ iteration of W in a single domain of $\tau^{-1}$

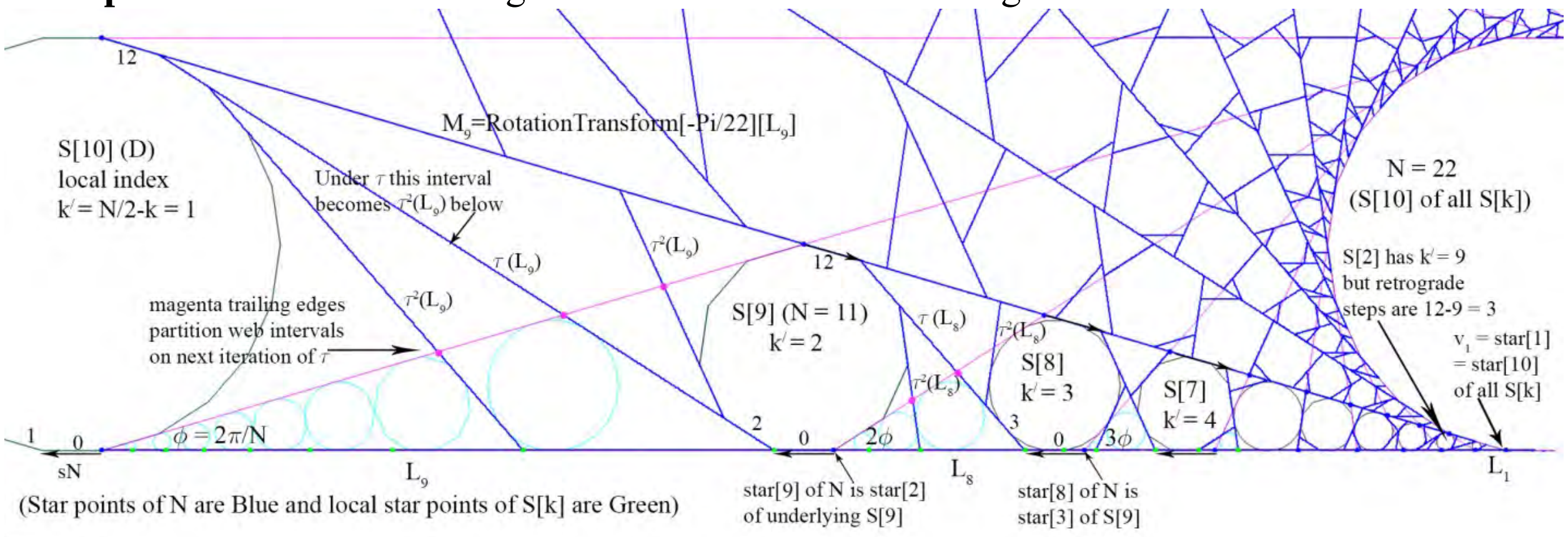

By symmetry it is sufficient to track the evolution of one of the domains of $\tau$ or $\tau^{-1}$. Our canonical choice is the (open) domain of $\tau^{-1}$ shown here between two blue extended edges of N. This is called Domv$_1$ because for all p, $\tau^{-1}(p) = 2v_1–p$ where $v_1$= star[1] of N. This domain of $\tau^{-1}$ will intersect N/2–1 magenta domains of $\tau$ defined by the trailing edges of N and the star[k] points. Our goal is to track the evolution of these star[k] domains under $\tau$. As p approaches the horizontal edge from the top, $\tau^{-1}(p) = 2v_1–p$ becomes a horizontal shear of magnitude $2v_1$ = sN, so under $\tau$, points in Domv$_1$ will feel an outward shear of this magnitude as they approach this edge and by continuity as p approaches the top edge the shears will be reversed.

Our choice of initial 'level-0' web interval will be the 'forward' edges of N such as edge L shown here. This edge runs from star[1] of N to star[1] of D and is partitioned by the star points into open intervals $L_k$ which run from star[k] to star[k+1]. The star points are the only points of discontinuity of $\tau$ on L, so W is by definition the disjoint union of the local webs of the $L_k$.

Even though these blue star points have no image under $\tau$ (or $\tau^{-1}$) by (one-sided) continuity they must feel this same shear so each star[k] and shear can be used to define a 'level-0 base' for the S[k] as subsets of $L_k$. We will see that the S[k] evolve from these level-0 bases so the webs of the $L_k$ will include the S[k] as well as the local webs of the S[k].

Note that the S[k] centers of the First Family Theorem are consistent with these shears. Recall that the local index k′ of S[k] was defined so that star[k] of N is star[k′] of S[k]. This means that the proposed rotation angle $k'\phi$ is the angle between edge L and the next step-k′ edge of S[k].

Therefore rotating the level-0 base edge of S[k] by $k'\phi$ will align a segment of edge L with an edge of S[k]. Under $\tau$, this segment will feel the sN shear and form a portion of the 'next' edge of S[k]. This is repeated recursively to generate the edges of S[k] in a step-k′ fashion.

Any portion of the original level-0 sN segment that extends beyond the bounds of the magenta trailing edges will lie in an adjacent $\tau$-domain and feel a different rotation angle on the next iteration - but the S[k] are conforming to these bounds, so there will always be a portion of the original segment that survives within these bounds and this 'web-cycle' will have period N/gcd(k′,N). By symmetry there must be an $M_k$ cycle on the top edge with this same period.

When N is even k′ = N/2–k so D at GenStar[N] will have local index k′ = 1 and the rotation angle will be $\phi = 2\pi/N$ to match that of N. This is the rotation that will allow the hypothetical edge 1 shown above to be a horizontal edge and hence experience the shear sN. The top and bottom shears are $\phi$ apart because they are on consecutive edges. Therefore for D, the top shear will be relative to edge 12 (N/2+1). The matching S[k] edges will have this same 12-step orientation with respect to the level-0 horizontal base. This value will be even when N is twice-odd so the primary and secondary cycles are synchronized mod 2 and the S[k] are always formed in a redundant step-k′ fashion.

When N is twice-odd, all the odd S[k] will have even k′ and hence web-periods N/2. Therefore the odd S[k] will be N/2-gons in the web W and this explains the 'androgynous' nature of the twice-odd N-gons. This is what we call a 'mutation'. Whenever gcd(k′,N) > 1, the S[k] will be formed from shortened web cycles so S[6] of N = 18 will be mutated in a known fashion.

When N is twice-even, N/2 ±1 is odd and the two web cycles are no longer synchronized mod 2. Therefore mutations will occur iff gcd(k′,N) > 2, because when gcd(k′,N) = 2 the shortened primary cycle will be corrected by the top cycle which is also period N/2. This means that the twice-odd case will have no 'gender' distinction between S[k] for k even or odd. In particular the penultimate S[N/2-2] 'M' tile will be free of mutations. Therefore for N = 12, S[2] and S[3] will be mutated, but S[4] will not. The 8k+4 'family' of N-gons is the only one with mutations in S[2] and this will have a strong influence on future generations. See Example 5.6.

When N is odd the (relative) shears are unchanged from the even case and the star angles are $(\phi/2)(N-2k) = (\phi/2)k'$ where k′ = N–2k is double the index for N even. These will be the local indices of the S[k] (which are now 2N-gons). Once again D will have index 1 and hence rotation angle $\phi/2$. Therefore D will be a 2N-gon with the same side as N. Since k′ = N–2k must be odd, the primary web cycle be odd and it will be shortened iff gcd(N–2k,2N) > 1. The top cycle will also be odd since it is based on edge N+2 of D, so these two cycles will be synchronized mod 2 as in the twice-odd case, but now both cycles are odd relative to the base edge, so there will be no gender-based mutations in the S[k], but of course D will have the full spectrum of gender changes in the DS[k]. In terms of mutations, the S[3], S[5] and S[6] tiles of N = 15 will be mutated. See the mutation conjecture below.

**Mutations of Regular N-gons**

In Example 4.6 we will extend this web analysis to the regions between the S[k] as described in the GFFT but here we summarize what occurs when the step-$k'$ evolution of an S[k] tile has a resonance with N. The simplest example is a resonance of 2 when N is twice-odd. As noted above, the odd S[k] with even $k'$ will have 'unfinished' webs which renders them as N/2-gons.

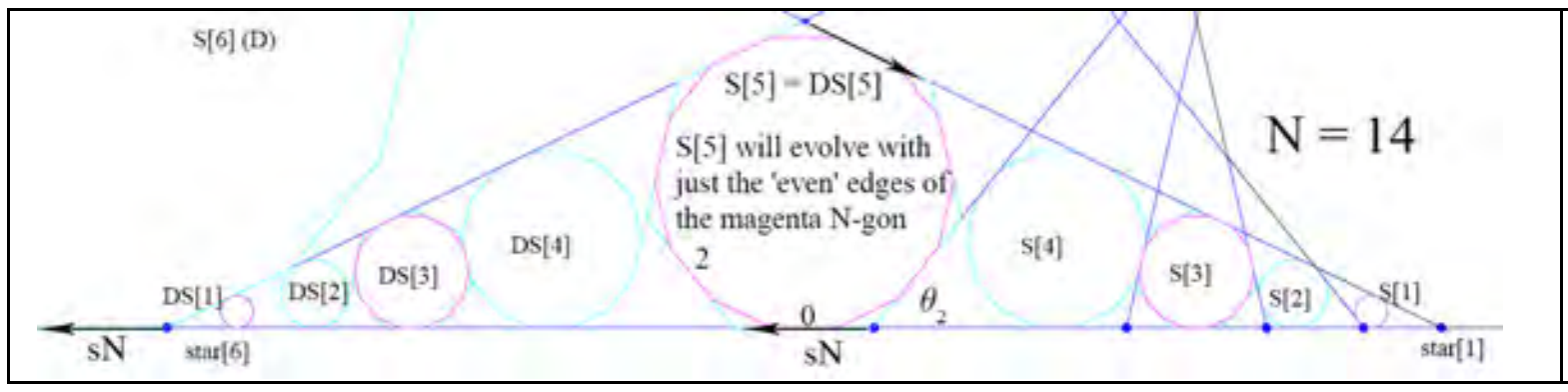

S[1] ,S[3] and S[5] for N = 14 will have gcd($k'$,N) = 2 so the web will only 'see' the even edges of the magenta underlying N-gon.

**Conjecture for mutations of S[k] due to their multi-step evolution**

(i) When N is twice-odd, S[k] will have local index $k'$ = N/2–k so the rotational period of S[k] in the web will be N/gcd($k'$,N) and S[k] will be mutated iff gcd($k'$,N) > 1.

(ii) When N is twice-even the local index of S[k] will be $k'$ = N/2–k as above so the rotational period is identical, but now it takes an odd number of rotational steps to reach the 'retrograde' top edge so the rotational period can be N/2 with no mutation of S[k]. Therefore S[k] will be mutated iff gcd($k'$,N) > 2.

(iii) When N is odd, S[k] will have local index $k'$ = N–2k which is twice the index above because the S[k] are now 2N-gons. The rotational period of the S[k] local web will be 2N/gcd($k'$,2N) and S[k] will be mutated iff gcd($k'$,2N) >1.

(iv) Typically a mutated S[k] will have a base that spans one web cycle, which is gcd(N,$k'$) or gcd(2N,$k'$) star points for N even or N odd. One of these surviving star points should be the minimum of N/2-1-j$k'$ for N even and N-2 -j$k'$ for N odd. The mutation will be an equilateral 2M-gon (or M-gon) where M = N/(gcd($k'$,N) for N even or 2N/(gcd($k'$,2N) for N odd . □

In [H6] we extend this conjecture to the next-generation S[k] tiles of S[2] which are called DS[k]. When N is even, S[2] seems to inherit the $k'$ = N/2–k step structure so the DS[k] and S[k] may share the same mutations. For N = 30 below DS[3], DS[5] and DS[9] have $k'$ = 12 ,10 and 6, so DS[5] is a hexagon while DS[3] and DS[9 ] are decagons, as in the matching S[k]. But when N is odd the S[k] web steps are N-2k so the mutations may not match.

Mutations in the next-generation DS[k] tiles of S[2] for N = 30 in the 8k+6 family

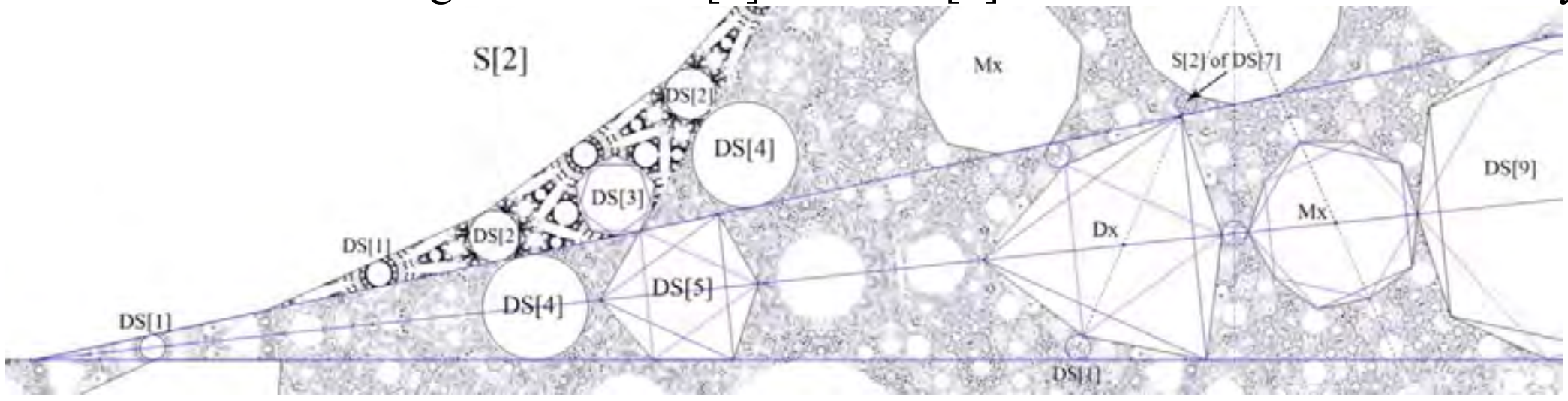

**Example 4.4** (Mutations of S[k] for N = 24) Because N is twice-even the web is more forgiving and S[2] and S[10] escape mutation. For S[8], gcd(N, k′) = 4, so MuS[8] will span 4 star points and the minimal (right-side) star point will be as predicted by the GFFT, namely the minimum of N/2-1-jk′ = 3. This is the smallest star point of the underlying S[8] that survives the local step-4 web and it is an ‘anchor’ for the magenta N/4 –gon shown below. The ‘base length’ includes the edge of S[8], so the matching left-side anchor is star[1] and this will be the reference for the blue hexagon. MuS[8] will be the ‘interleave’ of these two hexagons, shown here in black.

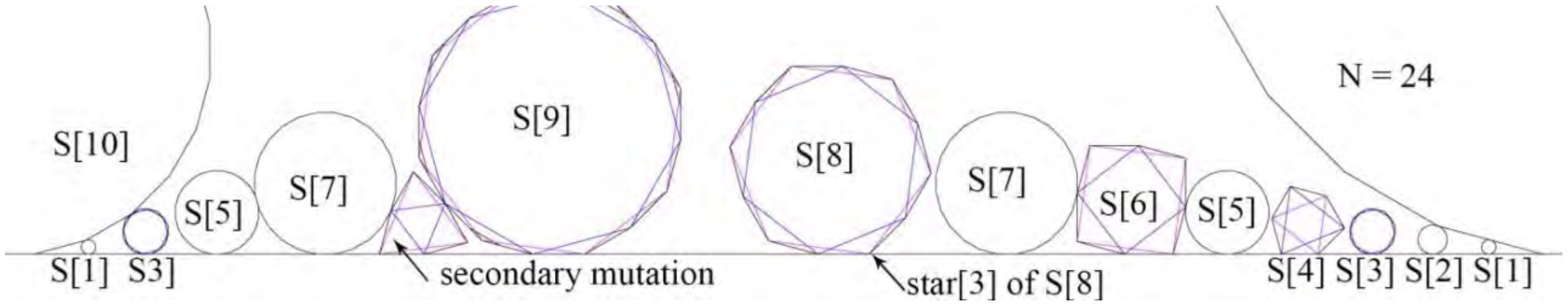

In [H3] we call MuS[8] a ‘semi-regular’ dodecagon because it is equilateral and has dihedral symmetry group $\mathcal{D}_6$ rather than $\mathcal{D}_{12}$. The ‘in-vitro’ version at the origin seems to have bounded dynamics but we do not know if this class of semi-regular polygons is always bounded.

Here are the steps to construct this mutated S[8] based on the known height and center of S[8].
(i) Find the star[3] point: MidS8= {cS[8][[1]],–hN}; star[3] = MidS8+{Tan[3·Pi/24]·hS[8],0};
(ii) M1=RotateVertex[star[3], 6, cS[8]] (magenta); M2 =RotateVertex[star[1],6, cS[8]] (blue)
(iii) MuS8 = Riffle[M1,M2] (black) (This weaves them as in a card shuffle.)

See Example 2.3 for the web. Since these mutations are incomplete local webs, they can be regarded as ‘scaffolding’ for the underlying S[k] tiles. The ideal S[k] will always be in the convex envelope of this scaffolding so in this sense the underlying S[k] are preserved by W.

**The general web**

Because of the multi-step origin of the S[k] the small-scale web may be very different from the orderly First Family structure. The limiting ‘tiles’ could be points - which must have non-periodic orbits, or possibly lines similar to the structures that appear for N = 11 on the right below. No one has ruled out the possibility of limiting regions with non-zero lebesgue measure. This is a long-standing open question in the phase-space geometry of Hamiltonian systems. Any non-limiting tile must be convex with edges parallel to those of N, so it is easy to see that the D tiles are maximal among regular polygons – and in fact rings of these tiles must exist at all radial distances so the dynamics in any finite region must be bounded. See Appendix II.

**Example 4.5** Some non-regular tiles in W (click to enlarge)

| N = 7 | N = 14 | N = 11 | N = 11 | N = 11 | N = 15 |
|---|---|---|---|---|---|
| 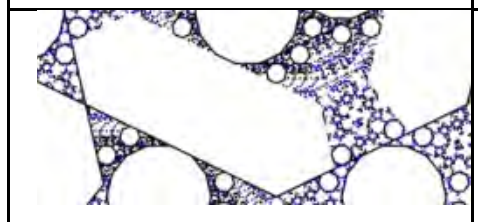 | 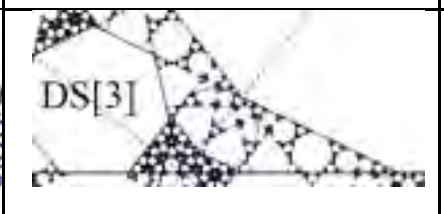 | 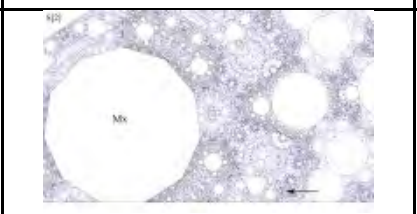 | 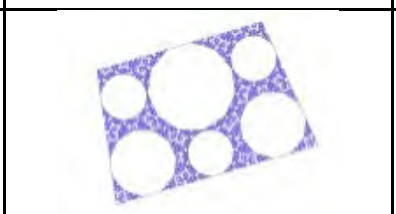 | 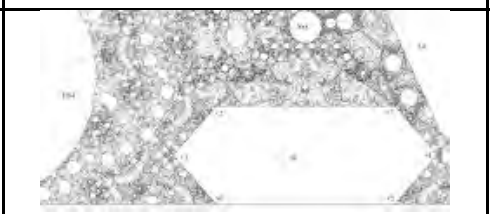 | 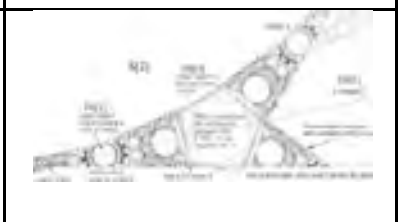 |
| N = 18 | N =36 | N = 41 | N = 48 | N = 60 | N = 200 |
| 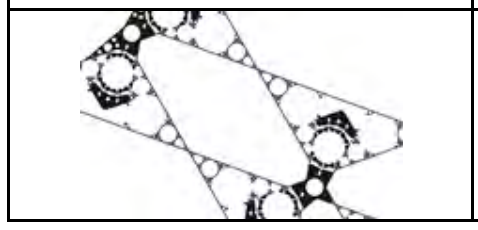 | 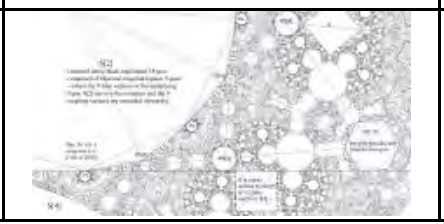 | 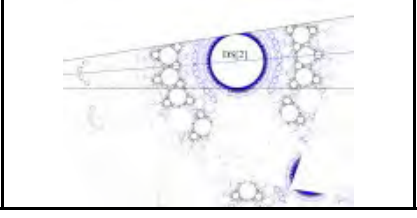 | 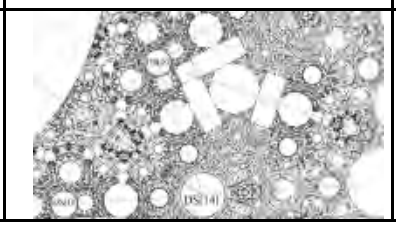 | 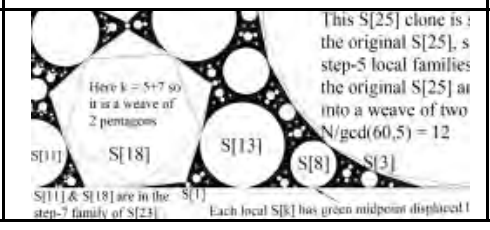 | 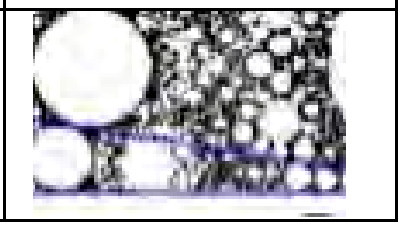 |

**Evolution of the Web : Part II - Local webs of the S[k] and in-situ vs. in-vitro evolution**

Every polygon P can be regarded as an 'in-vitro parent' at the origin, or P could be just a tile in the web or family of another polygon and we want to know how it evolves 'in-situ'. When N is twice-odd, the transformation T of Lemma 4.2 relates the in-situ M tile of N with its in-vitro form as an N/2-gon. When N is even the D tile is a retrograde version of N so the in-situ to in-vitro transformation T is a reflection about cM. For most P the two viewpoints are very different.

The previous web analysis in Part I shows that all S[k] (except D) evolve in a multi-step fashion and this will cause visible mutations when the steps are synchronized with N, but even when the final S[k] look normal, the multi-step web evolution will yield star points which are also multi-step and hence the families will be multi-step as shown by the GFFT and the cyan tiles below.

**Example 4.6** Enlargement of Example 4.3 showing the surviving star points for N = 22

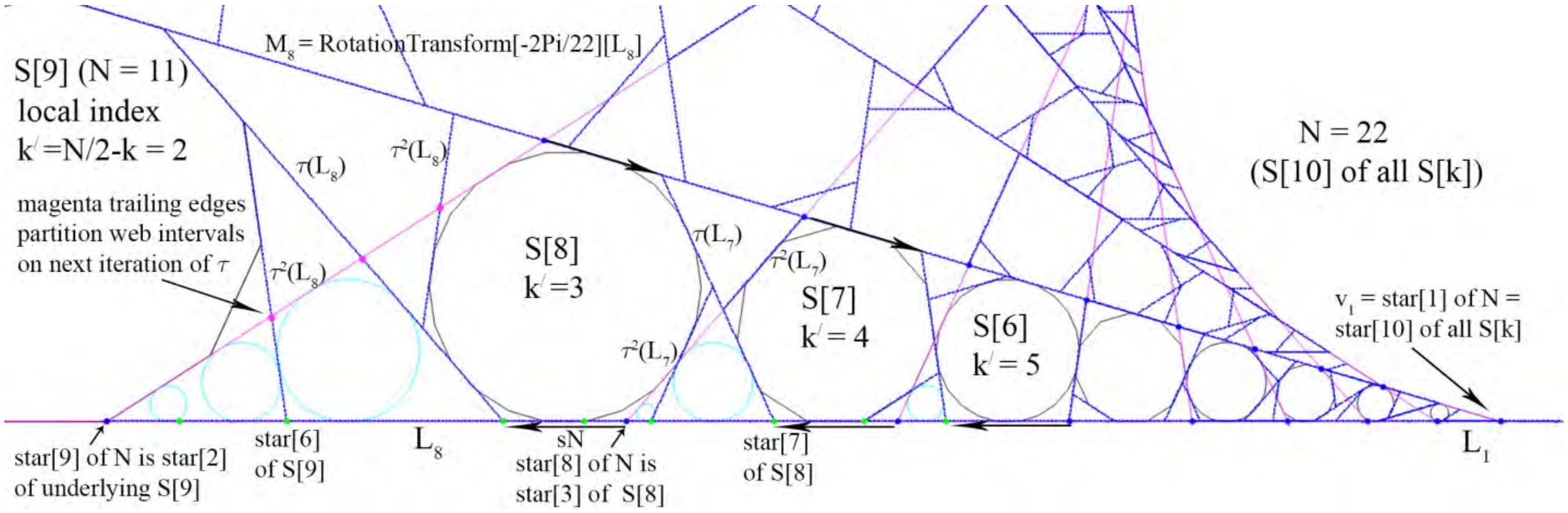

In Example 4.3 the primary interest was the evolution of the S[k] and the mutations that can occur in their web cycles. Here we attempt to describe the local webs of these S[k], because the gobal web W is by definition the disjoint union of these local webs. In this 'in-situ' evolution of the S[k] there is no reason to expect that the right-side geometry will match the left-side, but we can temporarily regard the left-side web of an S[k] to be the right-side web of S[k+1] and concentrate on just the right-side evolution as discussed in Lemma 2.3 (GFFT). The premise of that lemma was that for N even, the 'effective' right-side star points of the S[k] were step $k' =$ N/2– k and Example 4.3 showed that for N even the S[k] evolve in a step $k'$ fashion. This should imply that the surviving star points are also step $k'$ as assumed by the GFFT.

We know that there are two such cycles for each S[k] and each cycle has period N/gcd($k'$,N). In the right-side evolution of an S[k], the primary cycle is the 'top' cycle which is based on edge N/2+1 for N even. It will count down in a skip $k'$ fashion and this explains why the surviving star points are star[$k''$] where $k'' = (N/2-1) - j\cdot k'$ for $j \geq 0$. From this it is easy to construct the local cyan S[$k''$] families as in the GFFT.

Dynamically speaking these families can be regarded as 'satellites'of the S[k] which makes them seondary resonances of N. Since the S[k] have constant step-k orbits around N, the step sequences and periods of the families should reflect this relationship. This is discussed in Appendix C of [H4] and here we give a conjecture about the progression of these sequences.

**Step Sequence Conjecture for S[k] families for N even**

Suppose the local step-k′ 'GFFT family of S[k] consists of $P_j$ for j = 1,..,m where $P_{m+1}$ = S[k-1] and $P_{m+2}$ = N are not part of this immediate familiy. Then the step-sequence of $P_j$ is periodic of the form {(m – j +1)*{k}, k–1} where * is concatenation of copies of {k}.□

For example with N = 22, the local family of S[9] is $P_1,P_2,P_3$ = S[2],S[4],S[6], so m = 3 and the sequence for $P_1$ has the closest match with S[9] at {3*{9},8} = {9,9,9,8}. The sequence for S[4] is {9,9,8}and finally the local S[6] M' tile with {9,8} and inherited step k′ ={11-9 + 11-8} = 5. This can be continued with S[8] which has a local family with {8,8,7} and {8,7}, and S[7] contributes {7}, {7,6}. After S[7] all the local families contain just N.

To analyze orbits like these, we use a tool borrowed from classical analysis. The winding number (rotation number) of an orbit is a measure of the average rotation per iteration

**Definition 4.3** (winding numbers) For any N-gon, suppose a point p has step sequence S = $s_1,s_2$,.. The winding number of S is defined to be $\omega(S) = \frac{1}{N}[\lim_{m\to\infty}\frac{1}{m}\sum_{j=1}^{m} s_j]$.

If S is a periodic step sequence with period k then $\omega(S) = \frac{1}{N}[\frac{1}{k}\sum_{j=1}^{k} s_j]$ so in this case the winding number is 1/N times the mean number of steps in one period of the step sequence. For a regular N-gon Lemma 4.1 says every S[k] will have constant step sequence {k} with winding number k/N and the conjecture above extends this to the local families of the S[k].

In general these winding numbers are irrational and very difficult to compute. A knowledge of these sequences would be a major step toward understanding the dynamics of $\tau$ for any N-gon – regular or not. The only non-trivial regular cases known are quadratic with N = 5 by S. Tabachnikov in [T] (1995) and N = 8 by R. Schwartz in [S2] (2006). See Appendix E of [H3].

**Example 4.7** On the right are winding numbers as relative heights in the vicinity of S[1] for N = 11. It is typical for larger 'resonant' tiles to have relatively small winding numbers. The Sk 'skating rink' has period 338 and ω = 41/169 compared with S[1] at 1/11 (which is minimal).

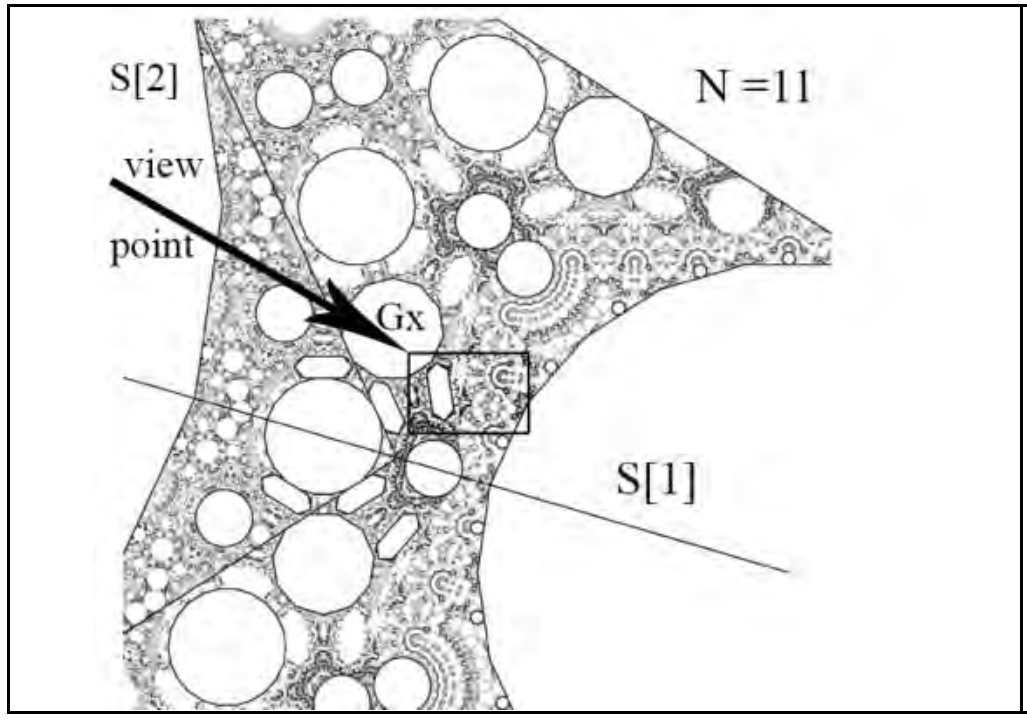

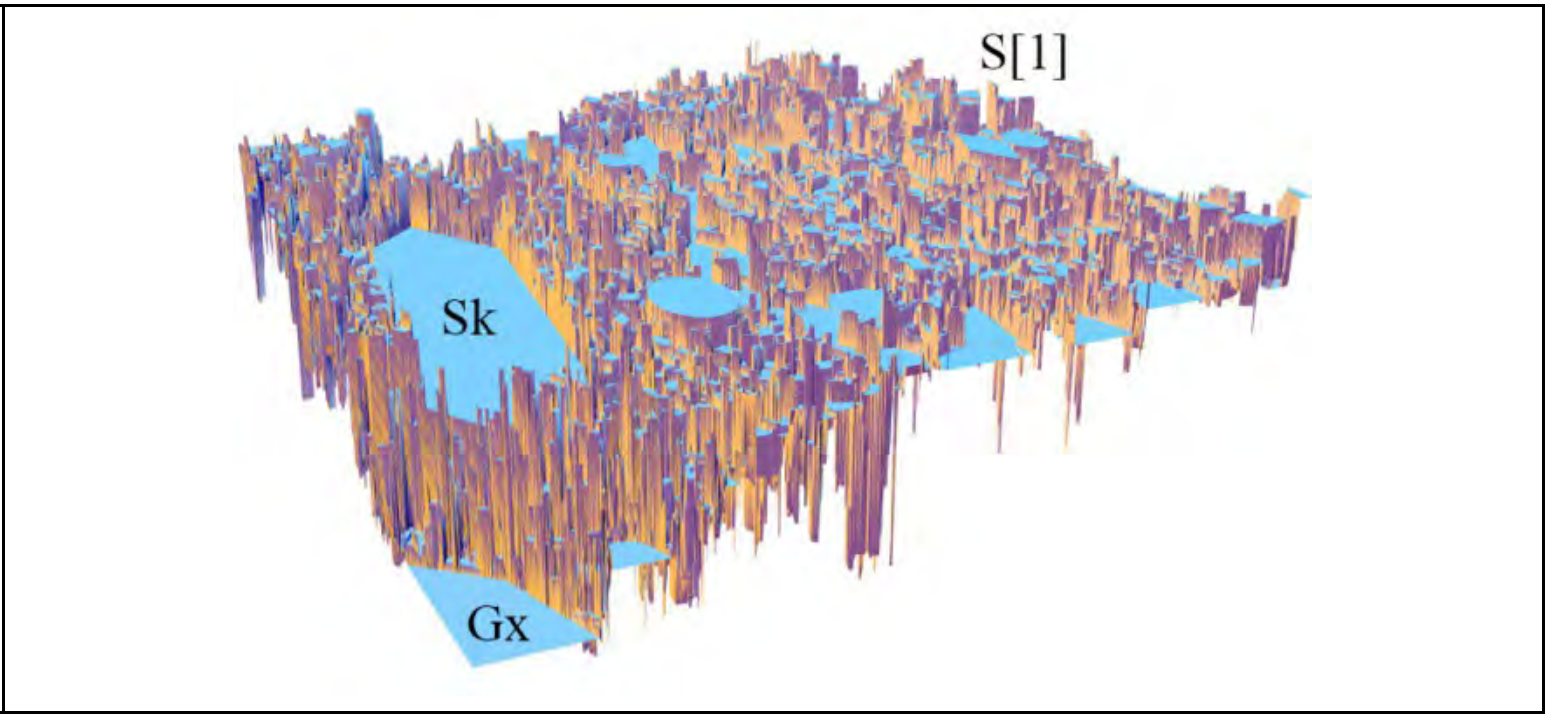

**The Edge Geometry of Regular Polygons**

Our ultimate goal is to understand the limiting web W and the FFT is just a first step. The derivation of these local families in the GFFT was based on the step-k′ web evolution of the S[k] and it was only possible because these S[k] 'parents' are closely related to each other and to N, so every S[k] shares star points with its neighbors. This may no longer true for the 3rd generation and in general there is no clear path toward recursion so the small-scale geometry of the 3rd generation may be very complex. In H[6[ and H[7] we will try to make some progress by focusing on the edge geometry of N. Appendix II shows that for N = 22 shown here there will be 4 invariant region on the right side of S[9] acting as a surrogate N = 11. The first such region is shown here bounded by the edge-sharing pair S[5] and S[6] . Here we can expect the S[k] tiles to share their local webs but we will concentrate on S[1] and S[2].

**Example 4.8** The invariant local geometry of N = 22 extends out to S[5]

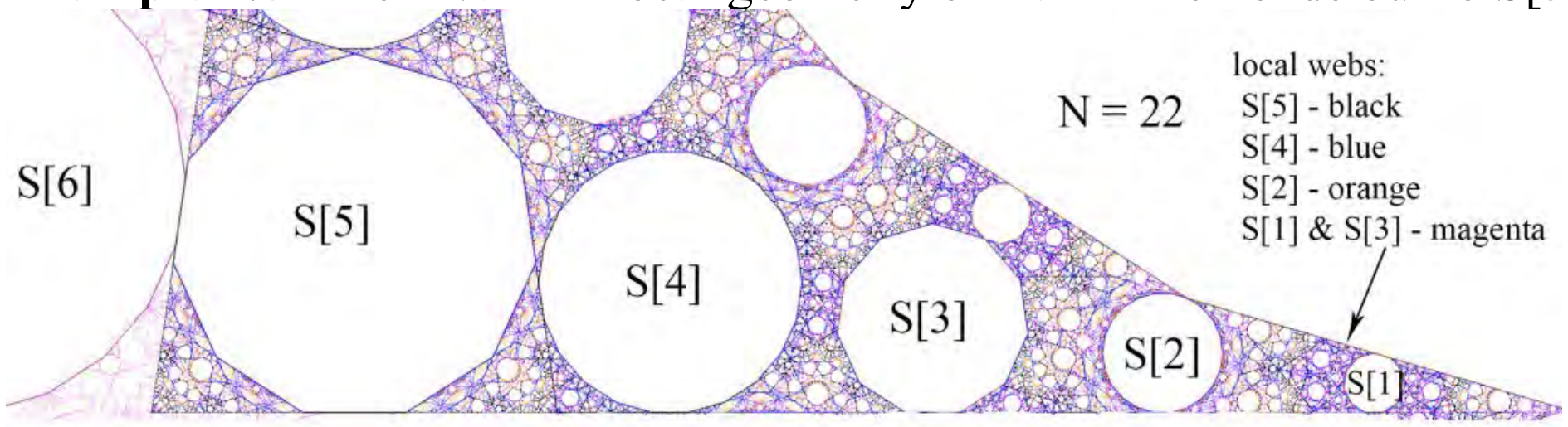

Since S[2] has period N/gcd(2, 22) = 11 it's orbit will decompose into two identical invariant regions and it is sufficient to study just one. There will be further invariant regions on all scales but it will be very difficult to find their boundaries. To probe this region it will fruitful to combine the S[1] and S[2] webs because they complement each other. S[2] here has k′ = N/2–k = 9. That is an awkward step size, but all S[k] are formed from two complementary cycles. For N even the complement cycle has N/2+1– (N/2– k) = k+1 steps. This says that the 'retrograde' S[2] web (which is just a right-side version of W) evolves in a step-3 fashion.

**Fxample 4.9** The combined web evolution of S[1] and S[2] for N = 22 begins with an initial (magenta and blue) interval that spans both tiles. This will be our default initial interval 0.

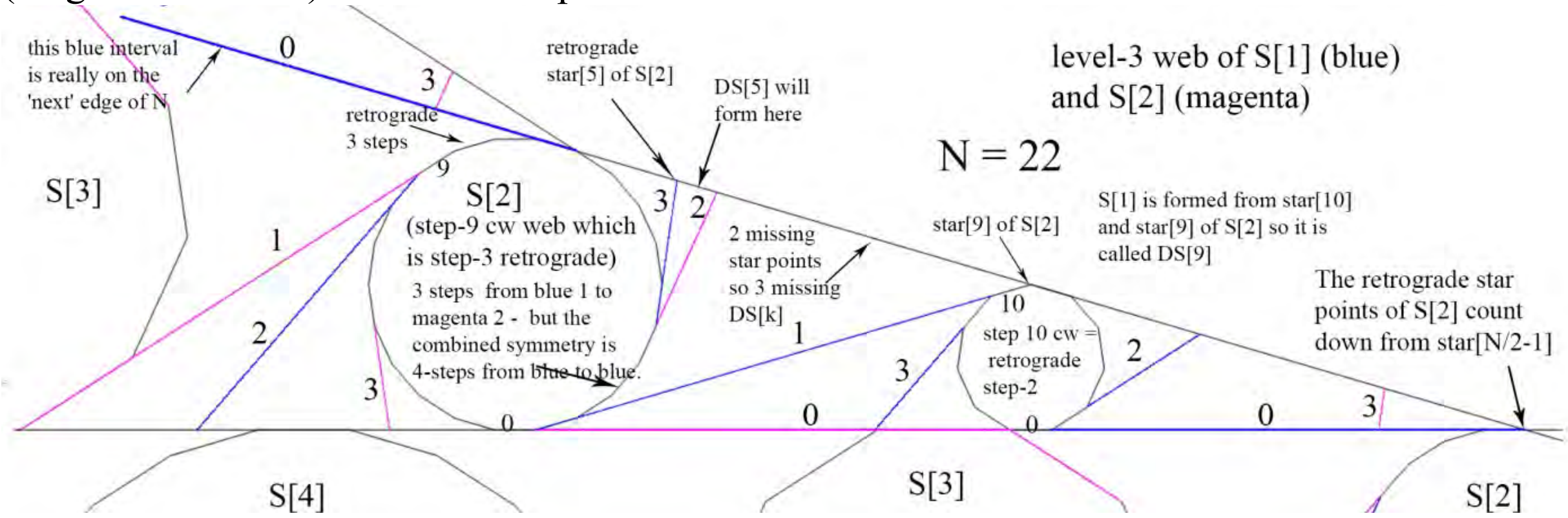

As indicated above for N even, S[1] and S[2] will have step-2 and step-3 ccw webs, and by default we will combine their webs together as shown here. There will be a 1-step offset between these webs so the result is what we call a hybrid step-4 web for S[2]. This extra step is perfect for generating the shared next generation S[k] tile of S[2], which we call the DS[k].

The early web evolution of these DS[k] is predictable because it is based on S[1] which is formed from star[N/2-1] and star[N/2-2] of S[2] – so it is called DS[N/2-2]. Because the combined web of S[1] and S[2] is step-4, the surviving DS[k] will count-down mod-4 to yield DS[5] and DS[1] for N = 22. This is the **Rule of 4 Conjecture** *For N even, counting down from S[1] at DS[N/2-2] the DS[k] will exist at least Mod-4.* Other 'volunteer' tiles may exist, but these DS[k] will evolve early and determine the symmetry of the web local to S[2]. This is illustrated with N = 40 below where the step-4 symmetry of S[2] is partitioned into step-2 increments yielding 4 star points between DS[k] (which here are simply S[k] of S[2]).

**Example 4.10** The web of N = 40, showing the symmetry local to S[2]

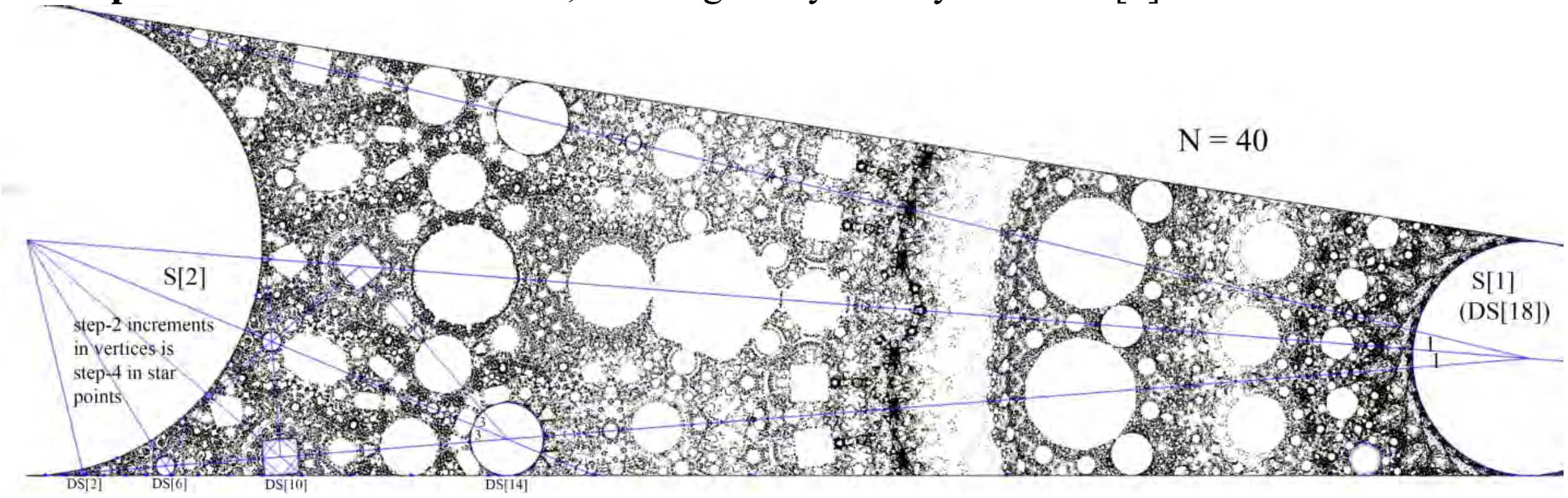

The primary lines of symmetry are determined by S[1] which has k′ = N/2-k = 18 and retrograde step k + 1 = 2. At each effective star point of S[2], the web splits and DS[14] will have k′ = 6 . This determines the local symmetry of DS[14] as shown by the blue lines. DS[10] will have k′ = 10 and be mutated into two squares since N/gcd(10,N) = 4. This implies that for N-even, the critical DS[1] and DS[2] tiles (if they exist) will have webs with k′ = N/2-1 and N/2-2 and Mod[N,N/2-1] is always 2 while Mod[N,N/2-2] is 4, so there two will have step-2 and step-4 limiting webs. This means that DS[2] could play the role of a D[2], but for N even this DS[2] is not predicted to survive except in the 8k family and usually there is no matching M[2. But the 8k+2 family has a predicted DS[3] and this can support D[2]s as well as M[2]s.

**N-odd** case is very similar but now S[1] and S[2] are 2N-gons, so the web steps are doubled from N/2-k to N -2k and the retrograde steps are also doubled at k′= 2k+2. This means that S[1] will be retrograde step 4 and S[2] will be step-6, but the combined web of S[1] and S[2] will be step-8 as shown below. The DS[k] are now based on star[2] of S[2] and they will also be step-8 to match the effective star points. This is what we call the Rule of 8.

**Example 4.11** The level-3 right-side webs of S[1] (blue) and S[2] (magenta) for N = 23

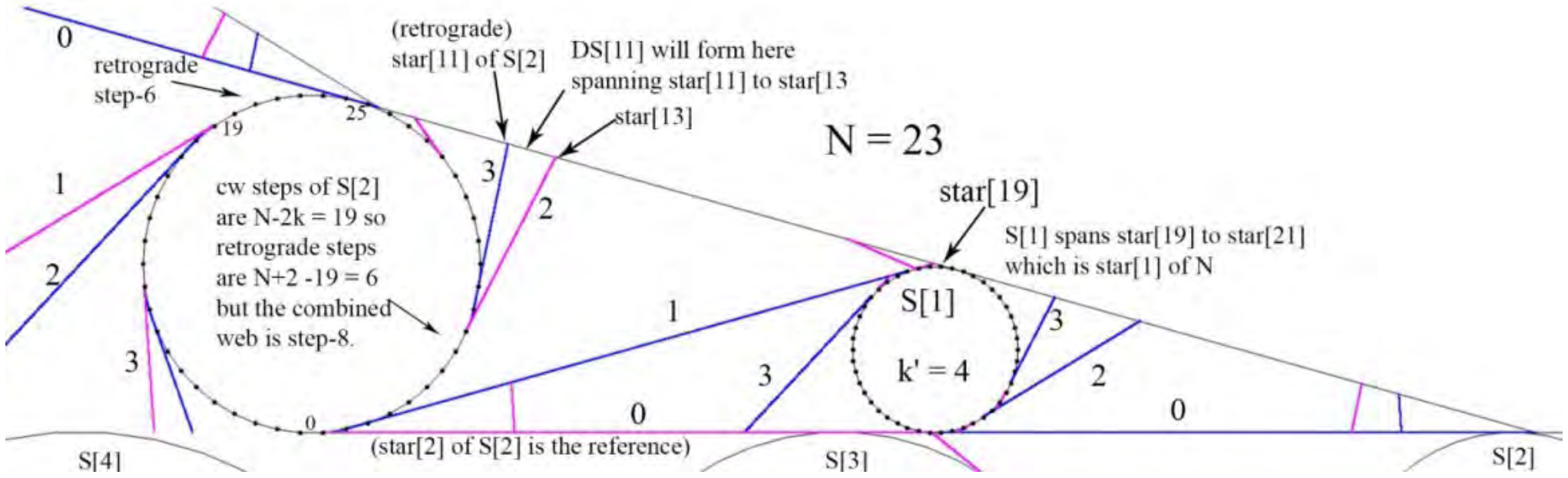

**Rule of 8 Conjecture**: *For N odd, counting down from S[1] at DS[N–4] the DS[k] will exist at least mod 8*. The important caveat is that the DS[k] tiles here are no longer part of the (ideal) First Family of S[2] because they will be based at star[2] of S[2]. They will be called the 'star2 family' of S[2] and they will be D tiles relative to the (virtual) S[k] of N. The S[1] tile will always be DS[N-4] so here it is DS[35] and theis combined S[1]-S[2] web, the effective star points of S[2] will count down mod-8 from 35, to yield DS[27[, DS[19], DS[11] and DS[3].

**Example 4.12** The web symmetry of S[2] local to S[2] for N = 39 in the 8k +7 family

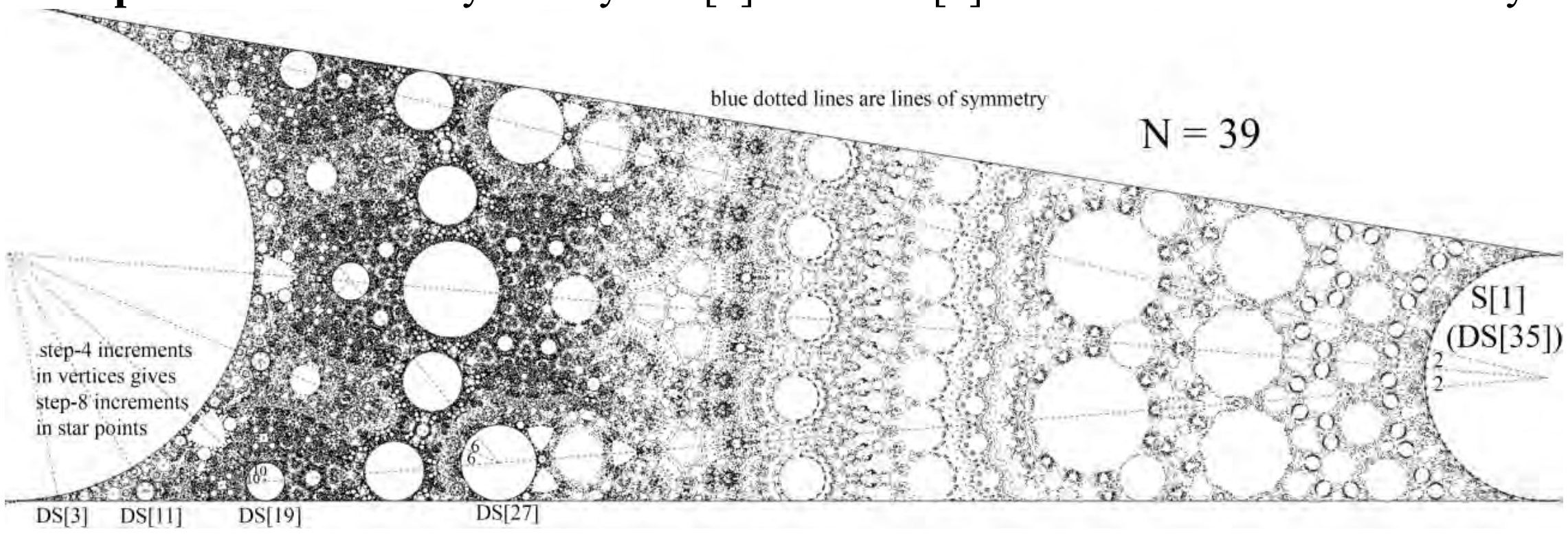

Since S[1] and S[2] are both 2N-gons the web steps of the DS[k] will be k′ = 2(N/2 )- k = N – k just like the even case. So DS[3] and DS[27] will be mutated with k′ = 36 and 12. Both have N/gcd(k′,N) = 13, so they will be the weave of 2 regular 13-gons. But unlike the N even-case these mutations will usually not be shared with S[k] of N because when N is odd k′ =N-2k. For N odd the S[k] tiles are 2N-gons and the GFFT will generate the step-k local families of the S[k]. This is done below for N = 17, but constructing the local DS[k] of S[2] is another issue. Here is a derivation using the Two-Star Lemma that will work for any DS[k]: np = 2N; h = hS[2]*Tan[Pi/np]*Tan[k*Pi/np] /(Tan[Pi/np]/Tan[2*Pi/np]); Mid=StarS2[[k]] + {h*Tan[((np/2)-k)*Pi/np],0}; c =Mid+{0,h};r=RadiusFromHeight[h,np];DS[k]:=RotateCorner[c +{r,0}, np, c];

**Example 4.13** The large-scale geometry for N-odd is a mixture of odd and even geometries as shown here for N = 17 . The black invariant region local to N = 17 extends out only to S[5].

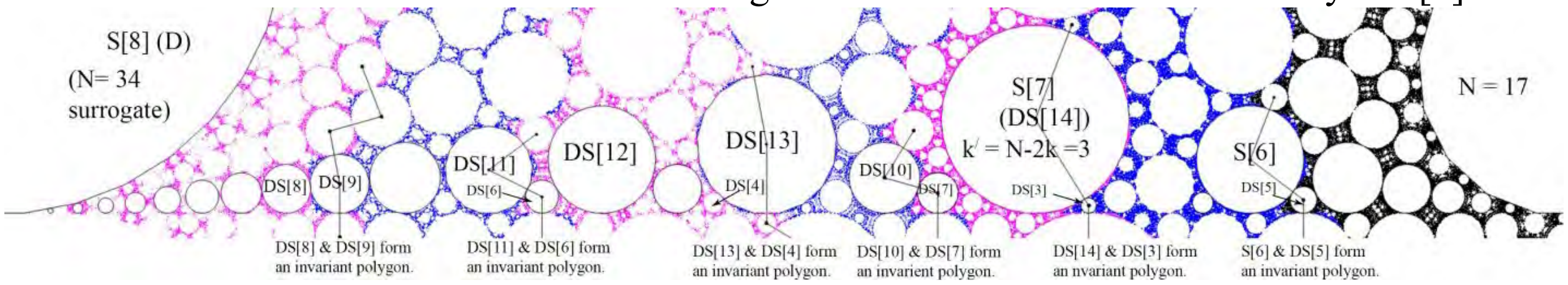

The Twice-Odd Lemma implies that N = 17 and N = 34 share the same web, but that does not imply shared dynamics and we do not know how the dynamics are related aside from the periods of the S[k] of N and DS[k] of D. This plot shows 7 invariant regions that create barriers for the dynamics so there may be very little sharing between N and D. The 8k+2 Conjecture below will show that we have some understanding of the local (magenta) dynamics of D, but all regions are invariant under $\tau$, so no point can have an orbit that visits other regions. This makes it difficult to relate the dynamics of N = 17 and N = 34. (Here we used the GFFT to construct the local step-k families of the S[k] to help understand the origin of the invariant regions. See Appendix II.)

**The 8k+2 Conjecture**

Suppose N is a regular polygon in the 8k+2 family. Since N is even, it is interchangeable with D and Definition 3.3 describes a well-defined (ideal) sequence of M[k] and D[k] tiles converging to star[1] of N which is also known as GenStar[N] or GenStar[N/2]. By reflective symmetry there will be an equivalent sequence converging to star[1] of S[2] which is the DS[2] of N. It is our convention to study the sequence there because S[2] can serve as a 'surrogate' N or D.

In this sequence N itself can be regarded as D[0] and the matching M[0] will be the S[N/2-2] penultimate tile of N. M[1] and D[1] will be the DS[1] and DS[2] of N which are simply the S[1] and S[2] of N since N is even. Then in the ideal sequence of Definition 3.3, for any positive integer k, hM[k+1]/hM[k] = hD[k+1]/hD[k] = GenScale[N/2] since N is twice-odd.

(i) We conjecture that these M[k] and D[k] tiles will always exist and converge to GenStar[N] with geometric scaling GenScale[N/2] and temporal scaling given by (ii) below

(ii) If $D_k$ is the $\tau$-period of D[k], it will satisfy the second-order difference equation

$D_k = nD_{k-1} + (n+1)D_{k-2}$ where n = N/2 and $D_1 = n$ and $D_2 = n^2$

The solution to this equation is $D_k = -\dfrac{n\left((-1)^k - (1+n)^k\right)}{2+n}$.

Therefore the periods of the D[k] for N = 2n will satisfy this closed-form equation.

(iii) The ratio of the periods will be $D_k/D_{k-1} = \dfrac{(-1)^k - (1+n)^k}{(-1)^{1+k} - (1+n)^{-1+k}}$

This sequence of ratios will clearly approach n + 1 = N/2 + 1 in the limit. It will begin with n and exceed n+1 on the second iteration and then alternate low-high relative to the limit.

(iv) Since 'most' M[k] form on the edges of the D[k-1] they will have the same limiting ratios and also satisfy the same basic difference equation as the D[k], but the new initial conditions will be $M_1$ =period of M[1] = n and $M_2$ =period of M[2] = n(3(n-1)/2 + 2 ), (This later period can be explained by noting that as shown below, 'most' D[k] arise in groups of 2 with 3 M[k] in each group plus the 2 M[k] of the single D[k]). The solution is

$$M_k = \frac{n\left((-1)^{1+k} + (-1)^k n + 3(1+n)^k\right)}{2(2+n)}$$ (period of the M[k] for n = 2N)

(v) Each of the predicted DS[k] in the Edge Conjecture will also survive in each generation.

**Note**:Part (v) implies that the entire region local to S[2] should share the geometric scaling of the M[k] and D[k] tiles, but the remaining DS[k] may have very different temporal scaling. Working back from D[2] toward star[1] of S[2] we know both the geometric and temporal scaling so the local fractal dimension should be Log[N/2+1]/Log[1/GenScale[N/2] for N = 8k+2. This begins with 1.2411 for N = 10 and approaches ½ with N. This does not imply that the subsequent generations will be self-similar and evidence says that typically they are not.

**Example 4.14** N = 26 has algebraic complexity 6 and it is the first 8k+2 case which is neither quadratic nor cubic in complexity. We will use it to help explain the difference equations for the evolution of the D[k] and M[k]. The Edge Conjecture predicts that S[1] will be DS[11] and there will also be DS[7]s and DS[3]s as shown below. The solid blue center lines track the local web symmetry and this also applies to the First Family as shown by the S[2] line of symmetry.

The geometry of N = 26 showing clusters of D[2] tiles anchored by S[3]s (not shown). The blue lines are lines of symmetry

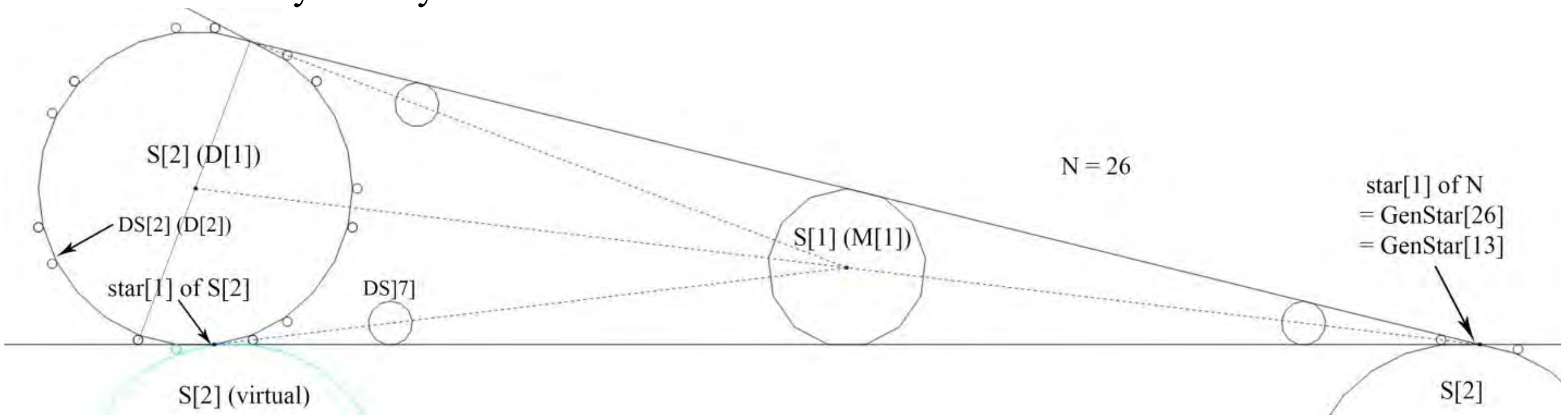

Because the combined web of S[1] and S[2] is step-4 there will be k (from 8k+2) of these S[3] 'clusters' on either side of this line of symmetry. Each cluster has 2 D[2]s and there is a lone D[2] on the line of symmetry for a total of 4k+1 D[2]s. This is the 'n = N/2' in the difference equation. Since the period of any S[k] is N/gcd(k,N), S[2] has period n also. Therefore the period of the D[2]s is $n^2$. Every member of the 8k+2 family will have D[1] at period n and D[2] at period $n^2$. These are the initial conditions of the p[k] equation. (Note that the S[2] local to N is in a different group of S[2]s so its small D[2] is not counted here. However when D[1] is replaced with D[2] to count the next generation D[3]s, these 'outliers' must be counted.)

At D[1] there will be 13 magenta D[3] around each D[2] - each with their outliers.
So there will be 14 outlier D[3]s which are just visible here as magenta dots.

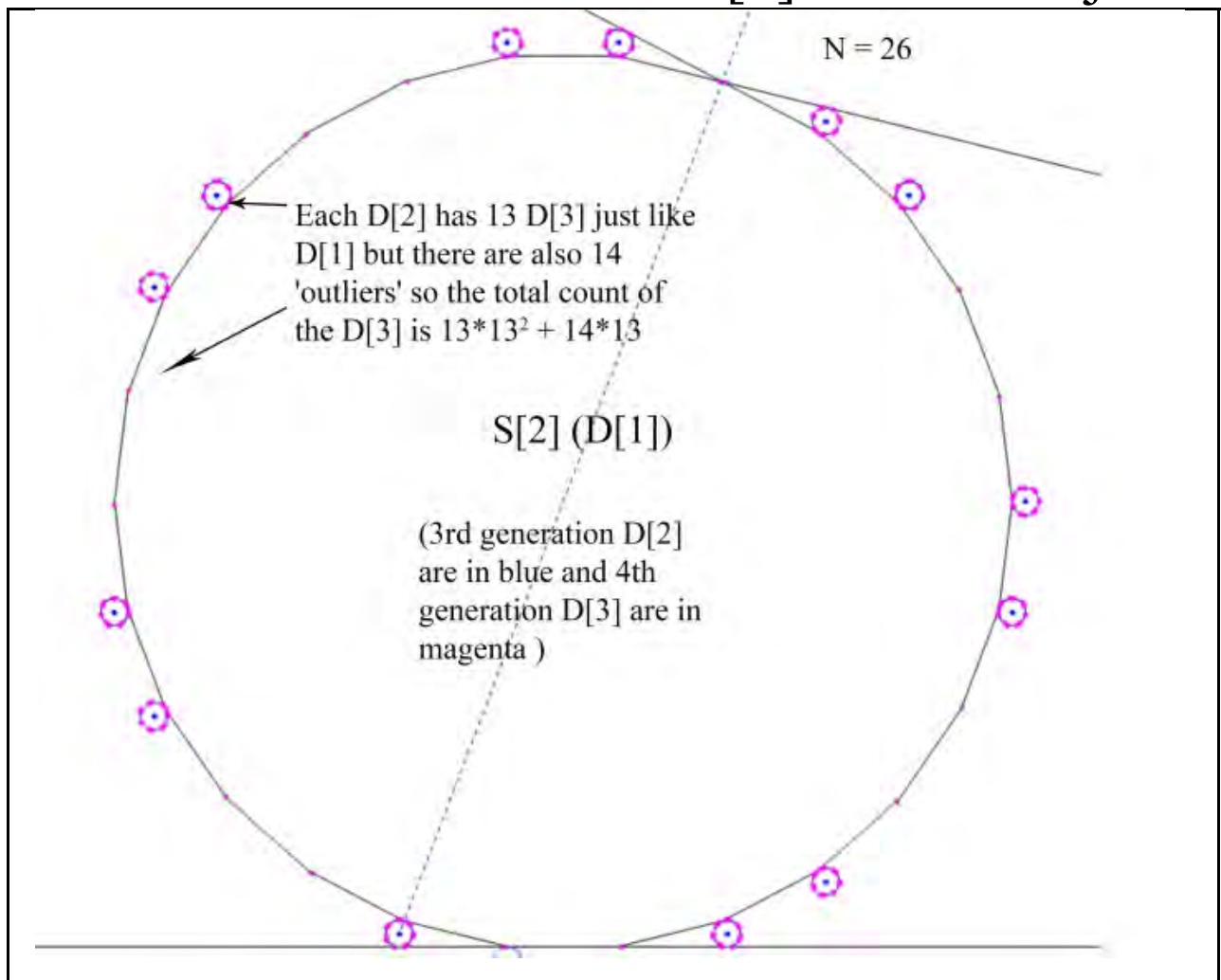

The periods of the first three D[k] are:
D[1] at 13, D[2] at $13^2$ and D[3] at $13^3 + 14\cdot 13$ so the proposed difference equation for the periods $D_k$ of the D[k] is

$$D_k = nD_{k-1} + (n+1)D_{k-2}$$

where n = N/2 and $D_1 = n$ and $D_2 = n^2$
Mathematica's solution is:

**RSolve**[{D[k]= D[k-1]n + D[k-2](n+1), D[1]=n, D[2]=n^2}, D[k], k]

$$D[k] \to -\frac{n\left((-1)^k - (1+n)^k\right)}{2+n}$$

as in the 8k+2 Conjecture

Therefore: $D[n,k] = -\dfrac{n\left((-1)^k - (1+n)^k\right)}{2+n}$ gives the $\tau$-period of each D[k] for N = 2n

**Table**[D[13,k], {k,1,8}] = {13, 169, 2379, 33293, 466115, 6525597, 91358371, 1279017181}

The $\tau$-periods of the M[k] of N = 26 should satisfy $M[n,k]=\dfrac{n\left((-1)^{1+k}+(-1)^{k}n+3(1+n)^{k}\right)}{2(2+n)}$

**Table**[M[13,k], {k,1,9}] ={{13, 260, 3562, 49946, 699166, 9788402, 137037550, 1918525778}

Because 'most' D[k] come in clusters of 2 the D[k]/M[k] ratio will approach 2/3 in the limit. Since the centers of these tiles are known, we have verified the early cases. See N = 34 in [H7].

**Example 4.16** N = 10 is the charter member of the 8k+2 family. Since it has quadratic complexity there is just one non-trivial scale which it shares with N = 5, namely GenScale[5] = $(\sqrt{5}+1)/2$. In Section 5 and [H6] we will look at the shared web of N = 10 and N = 5 and explain why both of these must have temporal scaling 6. Here we use N = 10 to illustrate the fact that even N-gons will have S[2] orbits which 'decompose' into two congruent invariant regions.

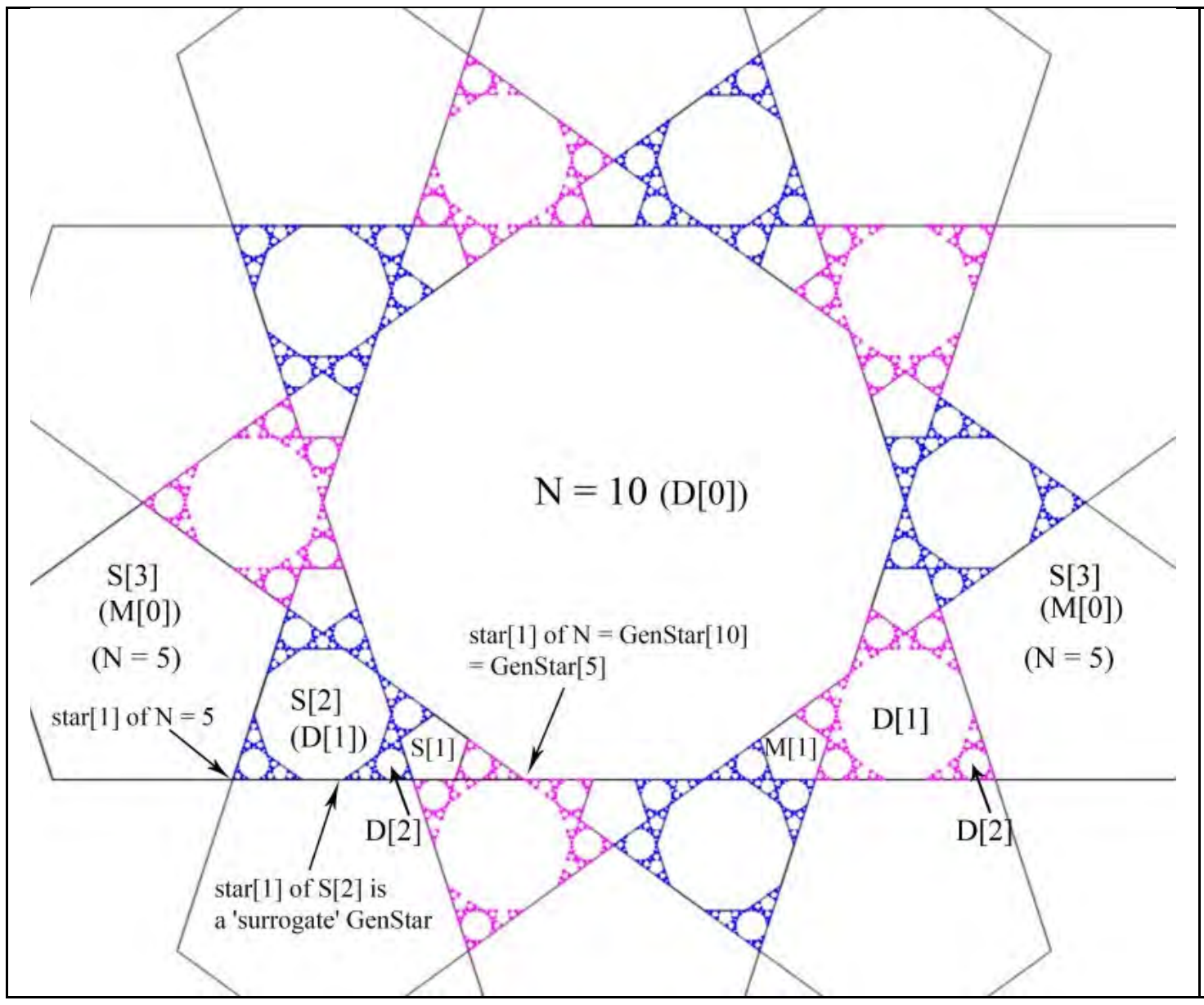

When dealing with dynamics it is necessary to take a global perspective. For all S[k], Lemma 4.1 shows that the period will be N/gcd(k,N), so for N even S[2] will have period N/2 as shown here by the blue decagons. This is why we could ignore the small D[2] 'outlliers' local to star[1] of N. It is sufficient to find the blue periods and this applies to all members of the 8k+2 family. But clearly the embedded N/2 - gon may have very different dynamics and it is only in this special case that we can show that the temporal scaling is unchanged when the center is shifted to N = 5.

This invariance occurs at all scales and these regions are typically 'nested' as shown below.

**Example 4.17** For N = 26, the 13 'even' and 'odd' S[2] still map to each other but only locally. The blue and magenta regions below are the orbits of points 'close' to S[2] and S[1] respectively. On close inspection is it clear that the odd S[2] will have magenta D[k] and M[k].

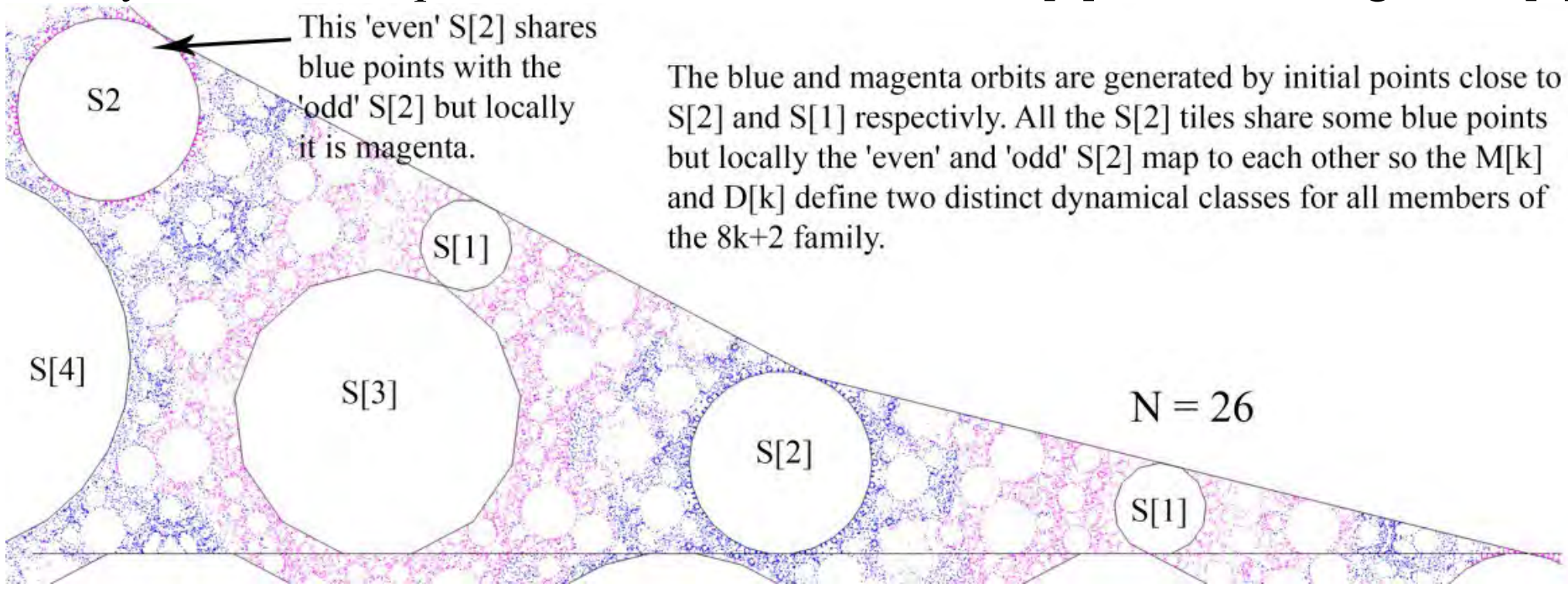

**Table 4.4** A classification of web geometry on the edges of regular N-gons, based on the Rule of 4 for N even (top) and the Rule of 8 for N odd (bottom). (This is the ‘8-fold way’ for N-gons.)

<table>
<tr><td>8k family (8,16,24,..)</td><td>8k+2 family (10,18,26,..)</td><td>8k+4 family (4,12,20,..)</td><td>8k+6 family (6,14,22,..)</td></tr>
<tr><td>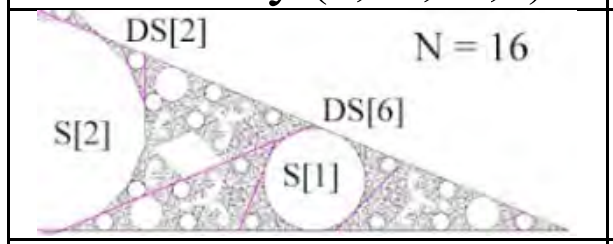
</td><td>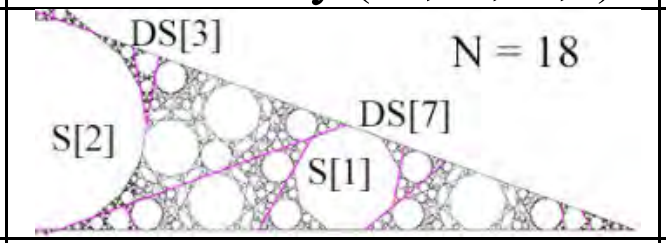
</td><td>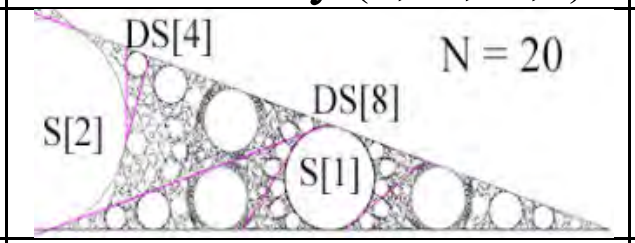
</td><td>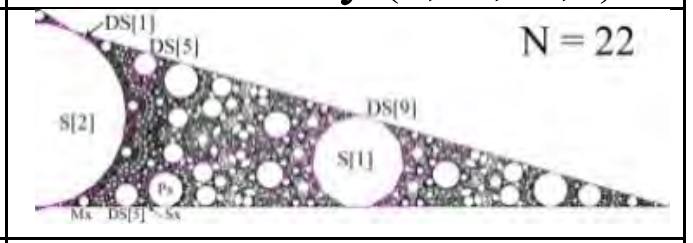
</td></tr>
<tr><td>8k+1 family (9,17,25,..)</td><td>8k+3 family (3,11,19,..)</td><td>8k+5 family (5,13,21,..)</td><td>8k+7 family (7,15,23,..)</td></tr>
<tr><td>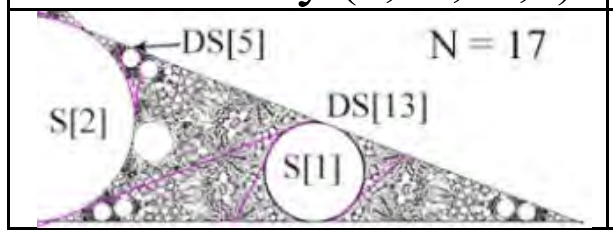
</td><td>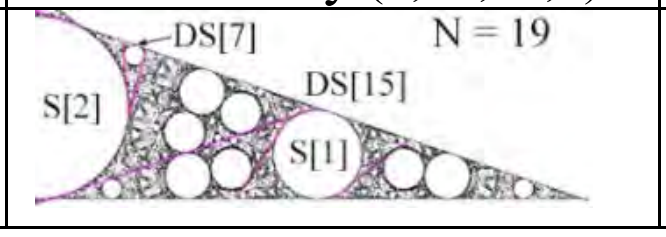
</td><td>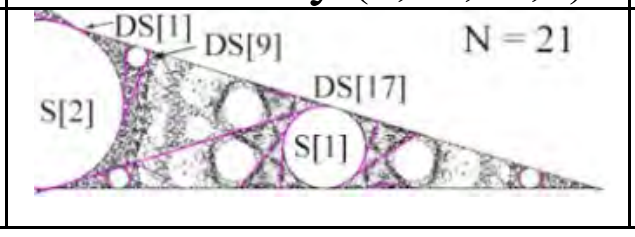
</td><td>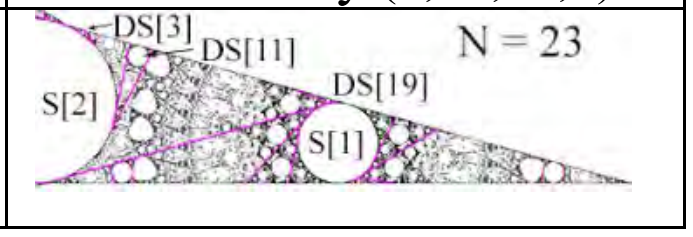
</td></tr>
</table>

These even and odd cases are linked because any odd N-gon can be regarded as embedded in the 2N-case. Therefore the 8k+1 and 8k+5 families can be embedded in 8k+2 families while the 8k +3 and 8k+7 families can be embedded in 8k+6. From a practical standpoint this embedding of an odd N-gon may be of little value because the step-2 web local to N may be very different than the web local to D. See the Twice-odd Lemma 4.2.

The **8k+2 Conjecture** described above was originally called the 4k+1 Conjecture in [H3]. The **8k+4 Conjecture** describes what happens when our 'host' S[2] is mutated. See Example 5.6 . The **N-odd S[2] Conjecture** says that the 8k+1, 8k+5 and 8k+7 families will have a DS[1] or DS[2] with potential for extended structure. The **8k+7 Conjecture** says that the predicted DS[3] will generate dual DS[1]s. The T**wice-even S[1] Conjecture** says that S[1] can support 'step-2' families which are D tiles relative to normal S[k].

Every regular N-gon has a locally invariant web that would be expected to contain about 1/3 or 1/4 of the S[k], so there is a connection between edge geometry and the large scale geometry. Both are driven by the cyclotomic field and the corresponding scaling field $S_N$ with complexity $C_N = \varphi(N)/2$. Hopefully the examples below and the survey in [H6] and [H7] may shed some light on the issue of ‘nature’ (algebraic complexity) vs. ‘nurture’ (web and edge complexity under $\tau$.) Appendix II shows that the locally invariant web is one of $C_N$ -1 invariant regions.

**Table 4.5** Algebraic Complexity of r egular N-gons for N ≤ 50

| $\phi$N/2 | 1 | 2 | 3 | 4 | 5 | 6 | 8 | 9 | 10 | 11 | 12 | 14 | 15 | 18 | 20 | 21 | 23 |
|---|---|---|---|---|---|---|---|---|---|---|---|---|---|---|---|---|---|
| N | 3<br>4<br>6 | 5<br>8<br>10<br>12 | 7<br>9<br>14<br>18 | 15<br>16<br>20<br>24<br>30 | 11<br>22 | 13<br>21<br>26<br>28<br>36<br>42 | 17<br>32<br>40<br>48 | 19<br>27<br>38 | 25<br>44<br>50 | 23<br>46 | 35<br>39<br>45 | 29 | 31 | 37 | 41 | 43<br>49 | 47 |

## Section 5. Examples of Singularity Sets

These examples will include the twice-odd pairs 5 and 10, 7 and 14, 11 and 22 and 13 and 26 as well as the twice-even cases N = 12, 16 and 20.

For N = 10 the 8k+2 Conjecture makes a number of assertions about the web and dynamics local to S[1] and S[2]. The Edge Conjecture makes no DS[k] predictions beyond S[1] at DS[3]. Since N = 10 has quadratic complexity, the web it shares with N = 5 has just one non-trivial scale and should be fractal in nature. As indicated earlier, N = 5 and N = 8 are the only non-trivial regular cases where the dynamics and singularity sets have been studied in detail. In [T] (1995) S. Tabachnikov derived the fractal dimension of W for N = 5 using 'normalization' methods and symbolic dynamics and in [S2] (2006) R. Schwartz used similar methods for N = 8. In [BC] (2011) Bedaride and Cassaigne reproduced Tabachnikov's results in the context of 'language' analysis and showed that N = 5 and N= 10 had equivalent sequences. (The Twice-Odd Lemma implies equivalent webs but not equivalent dynamics.) In [H3] we gave an independent analysis of the temporal scaling based on difference equations and that is reproduced in [H6]. Here we will give the implications of the 8k+2 Conjecture.

**Example 5.1** The temporal scaling of any 8k+2 N-gon (with N = 5 as the lone N/2 case)

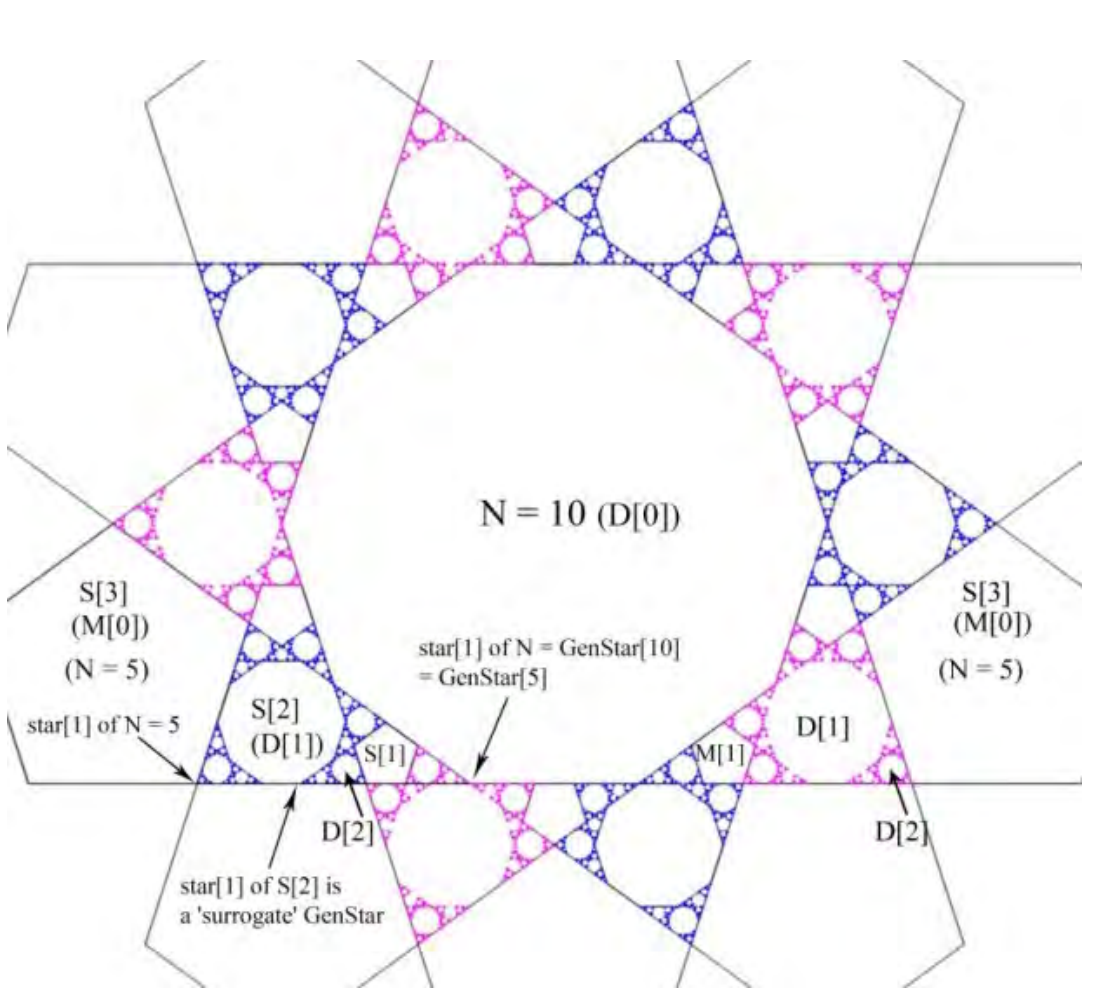

(i) The D[k] of N will have period $D_k$ where $D_k = nD_{k-1} + (n+1)D_{k-2}$ where n = N/2 and $D_1 = n$ and $D_2 = n^2$
(ii) The M[k] of N will have period $M_k$ where $M_k = nM_{k-1} + (n+1)M_{k-2}$ where n = N/2 and $M_1 = n$ and $M_2 = n(3(n-1)/2 + 2)$

(iii) The D[k] of N = 5 will have period $D_k$ where $D_k = nD_{k-1} + (n+1)D_{k-2}$ where n = N and $D_1 = n$ and $D_2 = n(n+2)$
(iv) The M[k] of N = 5 will have period $M_k$ where $M_k = nM_{k-1} + (n+1)M_{k-2}$ where n = N and $M_1 = 2n$ and $M_2 = 2n^2$ (these will be $2D_k$ for N = 10)

The solutions are:

| | |
|---|---|
| $(i)\ n = N/2: D[n,k] = -\dfrac{n\left((-1)^k - (1+n)^k\right)}{2+n}$ | $(ii)\ n = N/2, M[n,k] = \dfrac{n\left((-1)^{1+k} + (-1)^k n + 3(1+n)^k\right)}{2(2+n)}$ |

The first few periods are: $D_k$:{5, 25, 155, 925, 5555, 33325, 199955, 1199725, 7198355} and $M_k$ {5, 40, 230, 1390, 8330, 49990, 299930, 1799590,10797530}

| | |
|---|---|
| $(iii)\ n = 5, D[n,k] = \dfrac{n\left((-1)^k + (-1)^k n + 3(1+n)^k + n(1+n)^k\right)}{(1+n)(2+n)}$ | $(iv)\ n = 5, M[n,k] = -\dfrac{2n\left((-1)^k - (n+1)^k\right)}{2+n}$ |

The first few periods are: $D_k$ {5, 35, 205,1235, 7405, 44435, 266605, 1599635, 9597805} and $M_k$:{10, 50, 310, 1850, 11110, 66650, 399910, 2399450, 14396710} = $2D_k$ for N = 10.

**Example 5.2** Comparison of the edge geometry of the quadratic polygons N = 8, 10 and 12

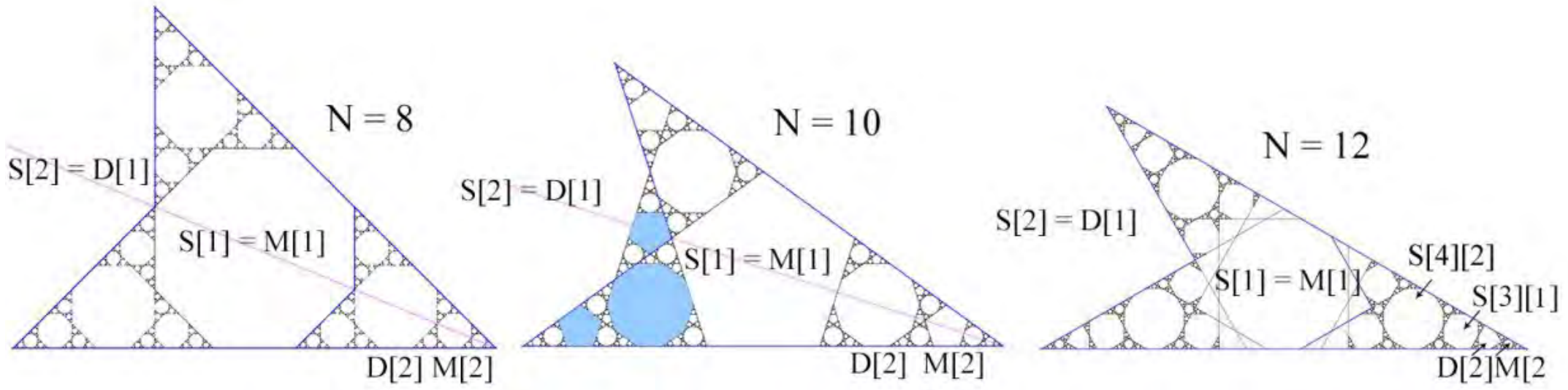

These webs have just one non-trivial primitive geometric scale which is GenScale[N] (or GenScale[N/2] for N = 10). This is hD[k]/hD[k–1] = hM[k]/hM[k–1], so it is also the scale of the darts (or triangles) which are anchored by M[k] or D[k].

In each case the magenta 'renormalization' line shows how the initial dart (or triangle) is mapped to a self-similar version of itself under $\tau^k$ for some k. As expected the N = 8 and N = 12 cases are closely related since their cyclotomic fields are generated by $\{\sqrt{2},i\}$ and $\{\sqrt{3},i\}$.

(i) For N = 8 (&12), the D[k] and M[k] are identical except for size. There are 3 darts in this invariant region and they are anchored by D[k]s and each D[k] is surrounded by 3 M[k]s - so the M[k] have temporal scaling of 9.

(ii) The temporal scaling of N = 10 is 6 as predicted by the 8k+2 conjecture and Example 5.1. Note that the 'next-generation' light-blue region shown here is composed of two overlapping 'towers' containing an M[k] and each M[k[ is surrounded by 3 D[k+1]s for a temporal scaling of 6 for the D[k]. Since these towers form a sequence converging to star[1] of M[1] this helps to explain why the D[k]s overall should have this same scaling. This is a non-trivial fact and as explained earlier, the key issue is the relationship between the decagons and pentagons.

(iii) For N = 12, each dart is anchored by an S[4] and each S[4] is surrounded by 3 S[3]s for a combined scaling of 9, and in the limit each S[3] will account for 3 M[k]s, so the M[k] scale by 27. This case is also not trivial and it is covered in more detail in Example 5.6.

Therefore the similarity (box-counting) dimension of the three webs should be:
(i) N =5 &10: Log[6]/Log[1/GenScale[5]] ≈ 1.2411 where GenScale[5] = Tan[π/5]Tan[π/10]
(ii) N = 8: Log[9]/Log[1/GenScale[8]] ≈ 1.2465 where GenScale[8] = $\text{Tan}[\pi/8]^2$
(iii) N = 12: Log[27]/Log[1/GenScale[12]] ≈ 1.2513 where GenScale[12] = $\text{Tan}[\pi/12]^2$

For compact self-similar sets such as these, the similarity dimension will match the traditional Hausdorff fractal dimension. It is no surprise that these dimensions are increasing, but this applies only to the quadratic family. For the cubic family and beyond, the webs are probably multi-fractal with a spectrum of dimensions. However it is likely that the maximal Hausdorff dimension will increase with the algebraic complexity of N, with limiting value of 2. See[LKV].

**Example 5.3** (N = 7 & N = 14) With N = 7 at the origin the inner invariant region would contain just S[1] and S[2] as show by the dark blue region below. Shifting the origin to D acting as N = 14, the web would be unchanged and the dynamics relative to D would have inner invariant region consisting of DS[1]- DS[3]] which is consistent with N= 7. The major convergent sequences are shown here and it should be clear that the dynamics relative to N = 7 are much more complex than at D. See [H6] for more detail about N = 7 and N = 14.

Below is an overview of the N = 7 and 14 families showing the major convergent sequenves

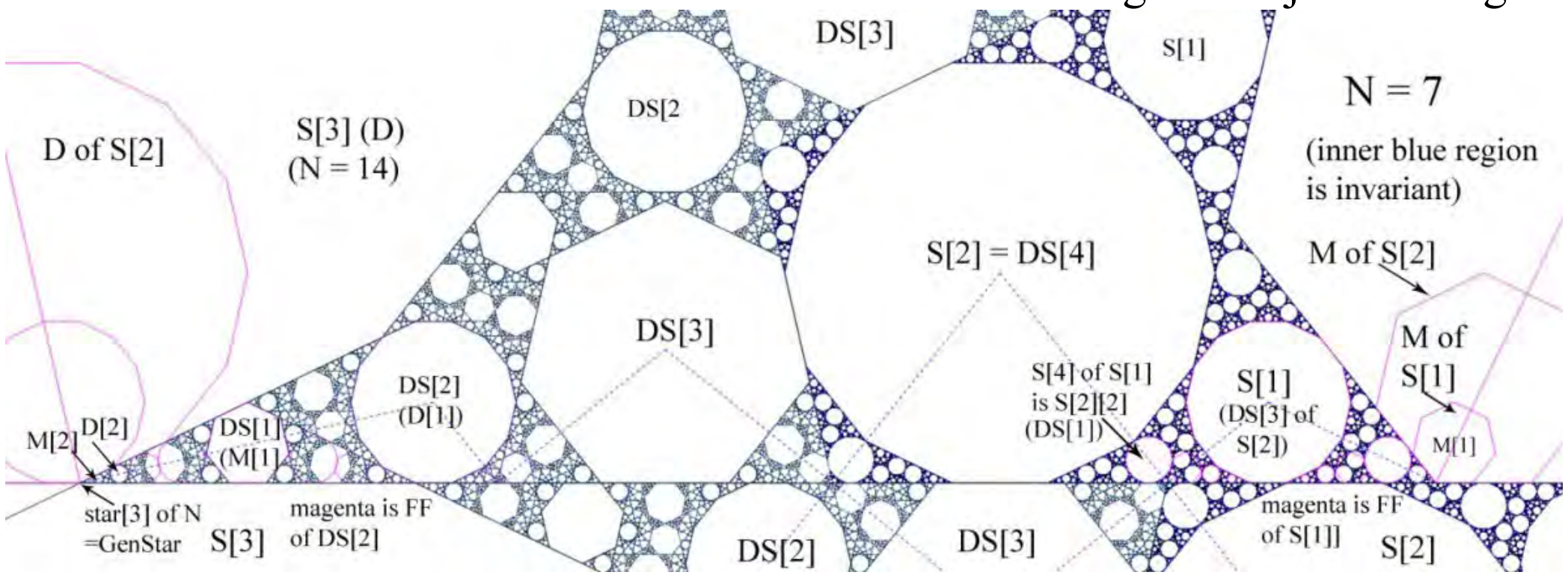

The invariant region around D is dominated by scale[3] = GenScale[7] but conflicts occur with DS[4] = S[2] since hS[1]/hS[2] = scale[2] of N = 7. By the Siegel-Chowla result these scales are non-commensurate. Below is an enlargement of D acting as N = 14 in the 8k+6 family. The Edge Conjecture predicts that DS[1] should exist as an M[2], but as evidenced by N = 22, the matching DS[2] may not exist. Here it appears that the intermittent DS[3] tiles come to the rescue with edges that help to form a DS[2] which can act as a D[2] and foster a 3rd generation at star[1] of D. It appears that this 3rd generation is self-similar to the first so there is potential for sequences of M[k] and D[k] converging to star[1] of D which is also GenStar[14]. The geometric scaling here is simply GenScale[7]$^k$ but the temporal scaling appears to be a multi-fractal as with N = 7. The first few M[k] periods are 28, 98, 2212, 17486, 433468, 3482794, 86639924 with ratios that appear to approach 8 and 25.Perhaps there is a 4th degree difference equation for these.

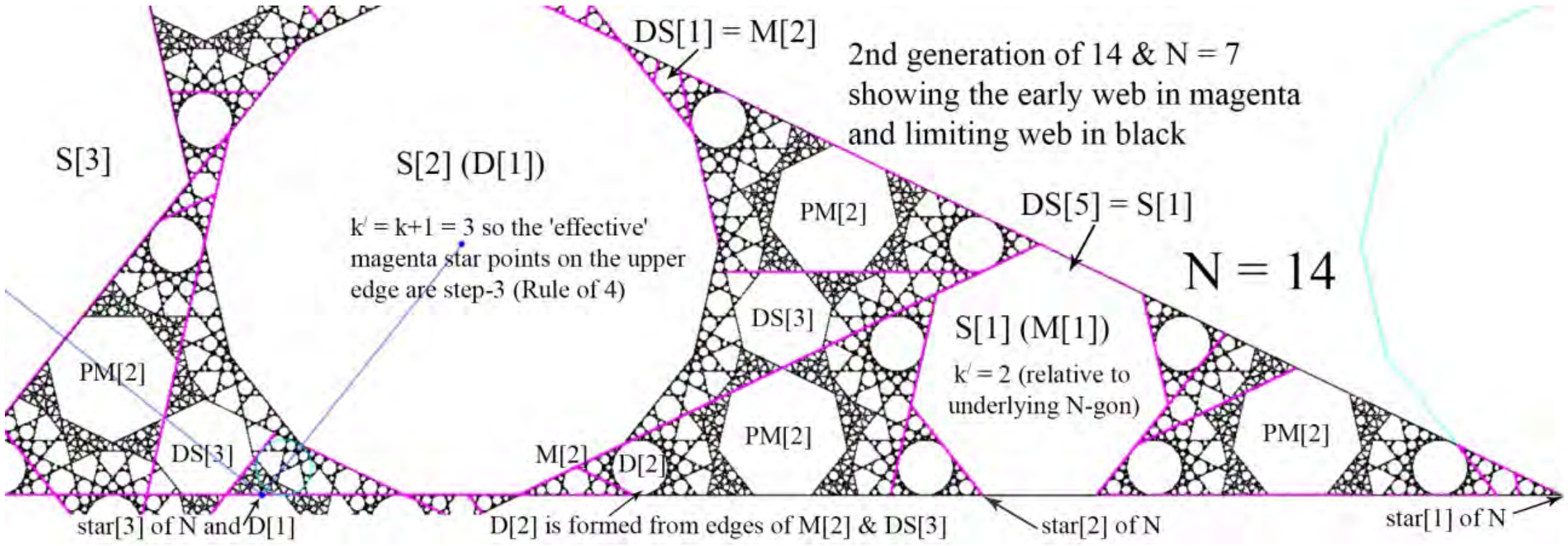

There is a dual orthogonal convergence at star[3] of D[1] with alternating PM and DS[3] tiles on the left and alternating virtual and real D[k] on the right. This yields a local geometry at star[3] of S[2] which is a mixture of single-scale self-similarity and multi-scale dynamics. This seems to be the nature of the limiting geometry for the 'cubic' N-gons N = 7, 9, 14 and 18.

**Example 5.4** (N = 11 & N = 22) These are the only 'quintic' cases. The two-elephant case from Example 1.2 takes place in the 2nd generation so these tiles exist on the edges of D which is a reflection of N = 22. This is an 8k+6 family like N = 14 so there is a DS[1] serving as an M[2], but apparently it is not capable of generating D[2] without the help of a DS[3] and it appears that the weakly conforming Mx tiles that appear in the 2nd generation are remnants of failed DS[3]s.

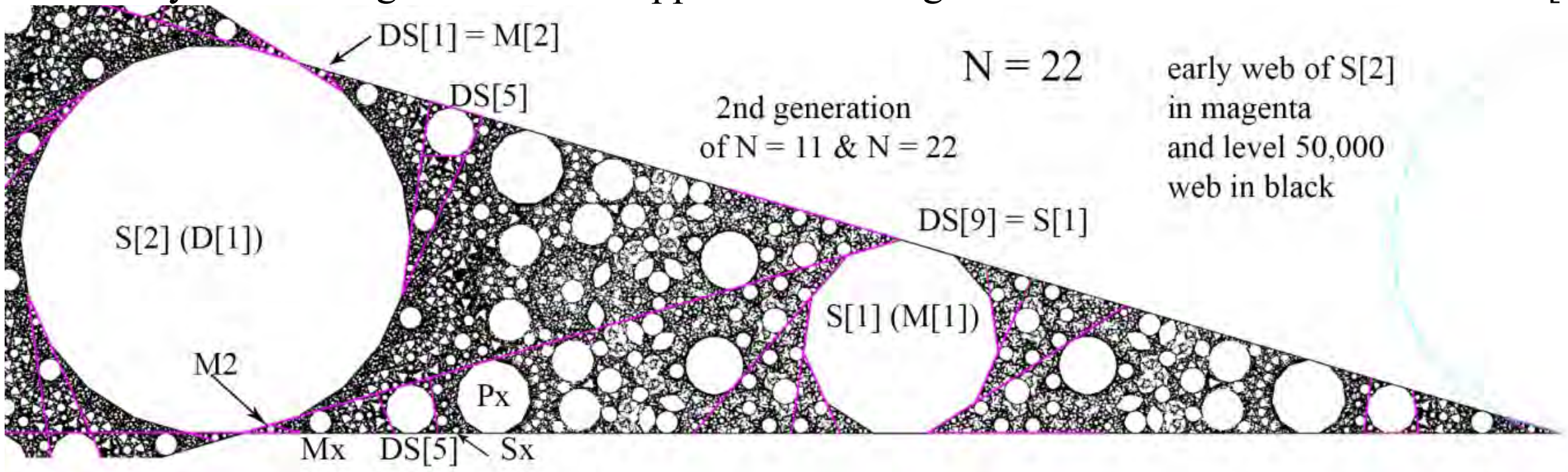

The web plot shows the predicted survival of M[2] and DS[5] along with S[1] at DS[9]. This DS[5] played a role in the earlier construction of Sx and the other 'elephant' was Px which is only weakly conforming to D[1]. In H[6] we derive the parameters of Px along with Mx and Sx. The difficult cases were Mx and Px and it was only possible to derive their parameters by doing exact calculations with the web of S[2]. (These calculations for Mx are repeated in the Appendix to show that the toral digital-filter map and the complex-valued dual-center map give the same parameters. We will sometimes use the dual-center map in these plots because it is very efficient but the Dc map requires a change of origin, so in the plot above, star[1] of N is at the origin.)

**Example 5.5** Detail of the edge geometry of N = 11 local to S[1].

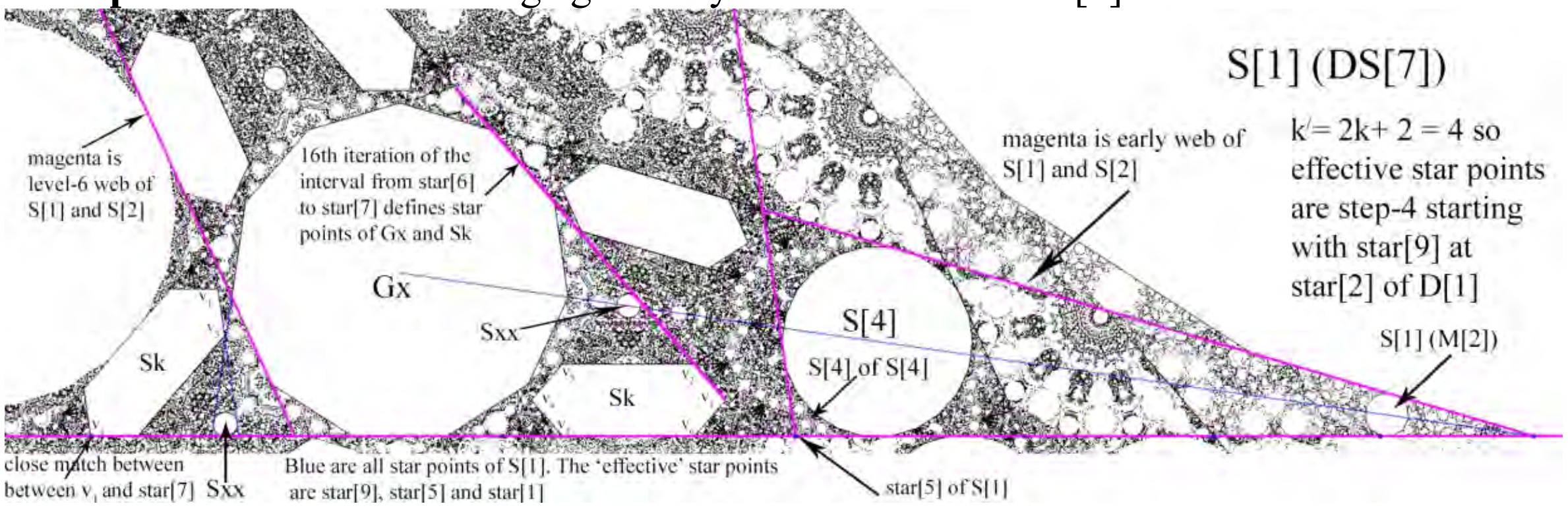

This 8k+3 family has no canonical DS[k] adjacent to S[2], but relative to S[1], the effective star points are mod-4 counting down from star[N-2] at star[2] of S[2]. Therefore the 8k+3 and 8k+7 families will have star[1] effective and for N = 11 this provides support for an S[4] a well as an S[1]. The Gx tile is a weakly conforming volunteer that also forms early in the web and is no doubt responsible for the Sk (See Example 4.7). In [H6] we find the parameters of Gx and the Sk. Using Gx and one of the Sk as 'elephants' it is easy to find the parameters of the small Sxx tile, as shown by the blue extended edges above. This is the first documented 'offspring' of two tiles which are themselves volunteers so the (location) parameters of Sxx are very complex - but hSxx is simple because Sxx is congruent to S[4] of S[4] shown in magenta. Because S[1] itself is congruent to D[1] of N = 22 , there may be subtle connections between Gx, Sk, Sxx and the N = 22 volunteers Mx, Px and Sx.

**Example 5.6** (N =12). N = 12 is a quadratic polygon and the first non-trivial member of the 8k+4 family. This is a very interesting family geometrically and algebraically. Since the mutation criteria for S[k] in the twice-even case is gcd(N/2-k,N) > 2, this family is the only one where S[2] is mutated. This is a simple star[1] to opposite-side star[3] mutation, so S[2] will always be the equilateral Riffle or 'weave' of two N/4-gons. This is compatible with the Edge Conjecture which predicts the survival of a DS[4]. For N = 12 , S[3] is also mutated, but as always the M tile escapes mutation because it is S[N/2-2].

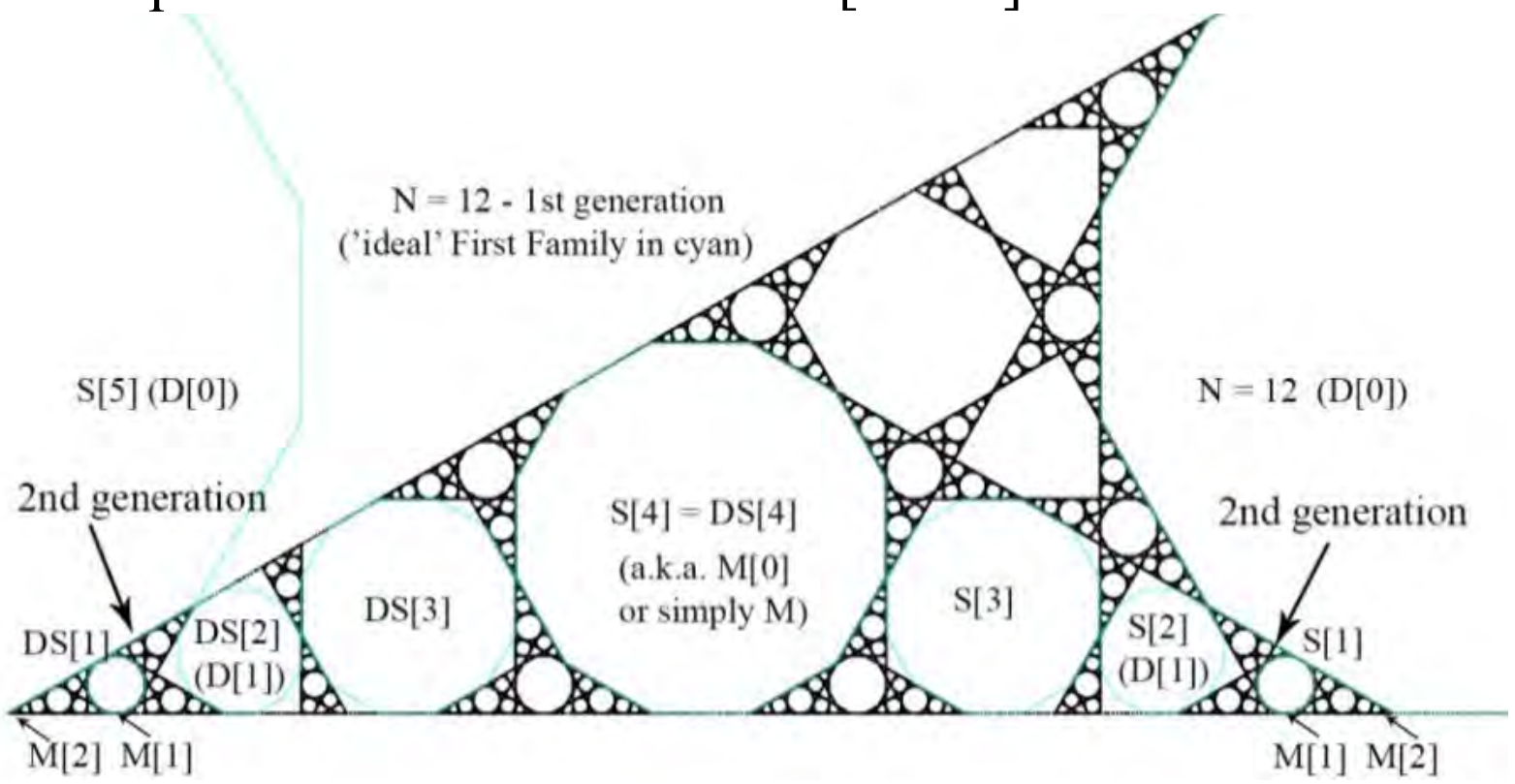

Below is the 2[nd] generation where the predicted DS[4] is S[1], also known as M[1]. The geometry local to M[1] looks like an unfinished web but this is consistent with the fact that S[2] is composed of two regular triangles which have trivial local webs. As N increases in the 8k+ 4 family, this local N/4-web of S[2] will be more complex, but apparently it will never yield 'normal' families for the DS[4]s at the magenta extended vertices. However the blue vertices will be also be vertices of the underlying S[2] and have a more promising geometry. It appears that for all members of this 8k+4 family, the blue star[1] of S[2] will be shared by an alternative 'parent' Px of DS[4]. For N = 12 this Px is the large S[3] tile of N shown here.

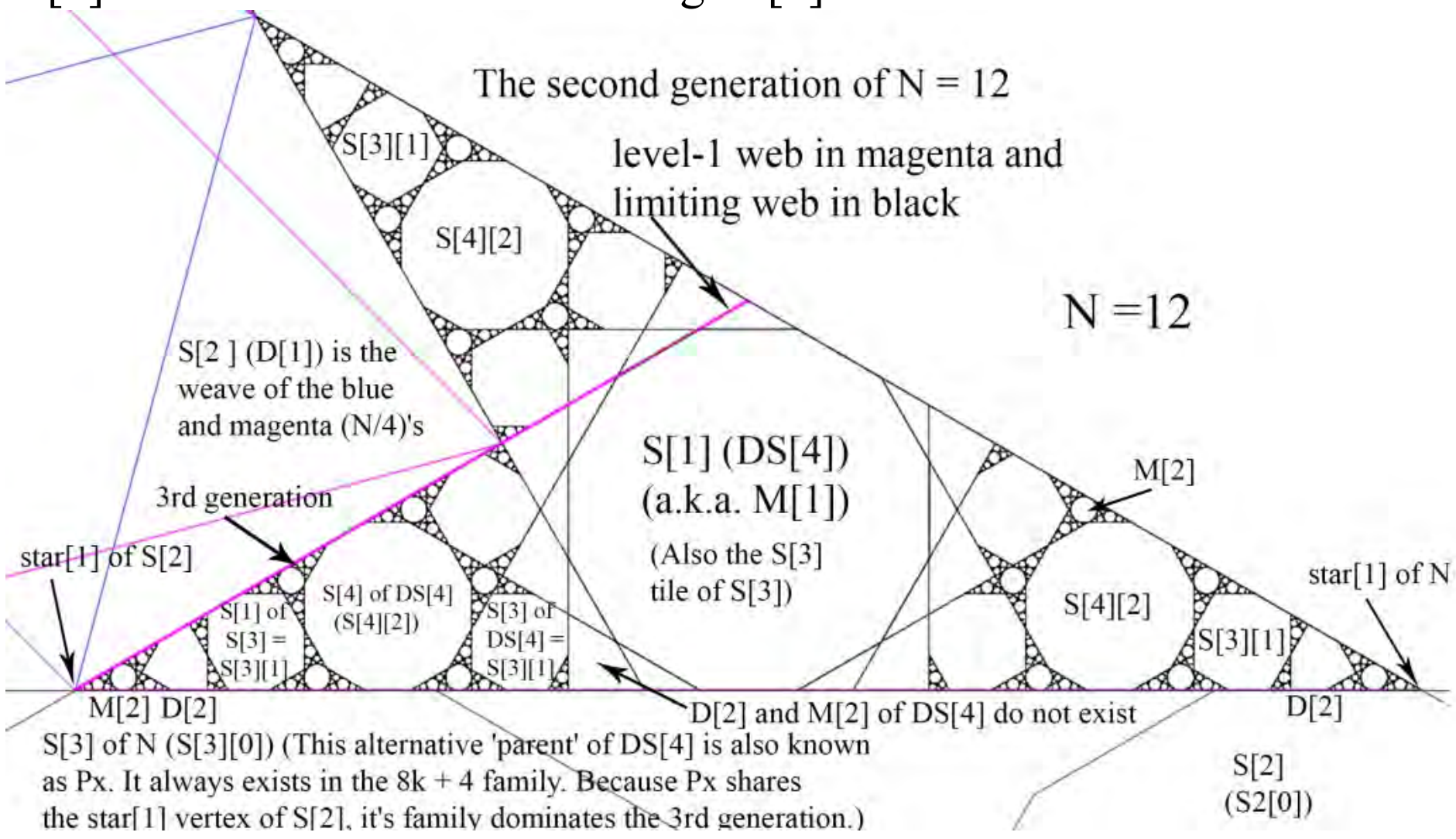

The First Family of S[3] does include DS[4] as well as a next-generation S[3][1] at the S[1] position. This S[3][1] survives the web along with DS[4]. (The S[2] of S[3] does not survive.) This S[3][1] and symmetric copies then generate M[2]'s as their S[3] tiles. We will show that in the limit each S[3][k] will account for 3 M[k+1]s, and since there will 3 of these S[3][1]s in each 'dart' and there are 3 darts in each generation, the overall growth rate should be 27 as expected.

Here the count of M[2]s is a little short with only 24 but that is because the mutation of S[2] destroys potential M[2]s (and D[2]s) which would exist on the edges of DS[4]. We will show why these mutations will play a diminishing role with each new generation.

Below is an enlargement of one dart. The only M[2] that is generated by a (mutated) D[2] is at star[1] of N. The remaining 5 M[2]s are generated as step-3 tiles of the S[3][1] so they have normal webs which will contain M[3]s at the 'next-generation' step-1 positions. Therefore of the 23 'darts' shown here only 3 are affected by the mutations - and this ratio will decrease with each new generation to yield a limiting count of 3 M[k]s for each S[3][k–1] and a limiting temporal scaling of 27 for the M[k]s (and D[k]s). See Example 5.2 for the resulting fractal dimension.

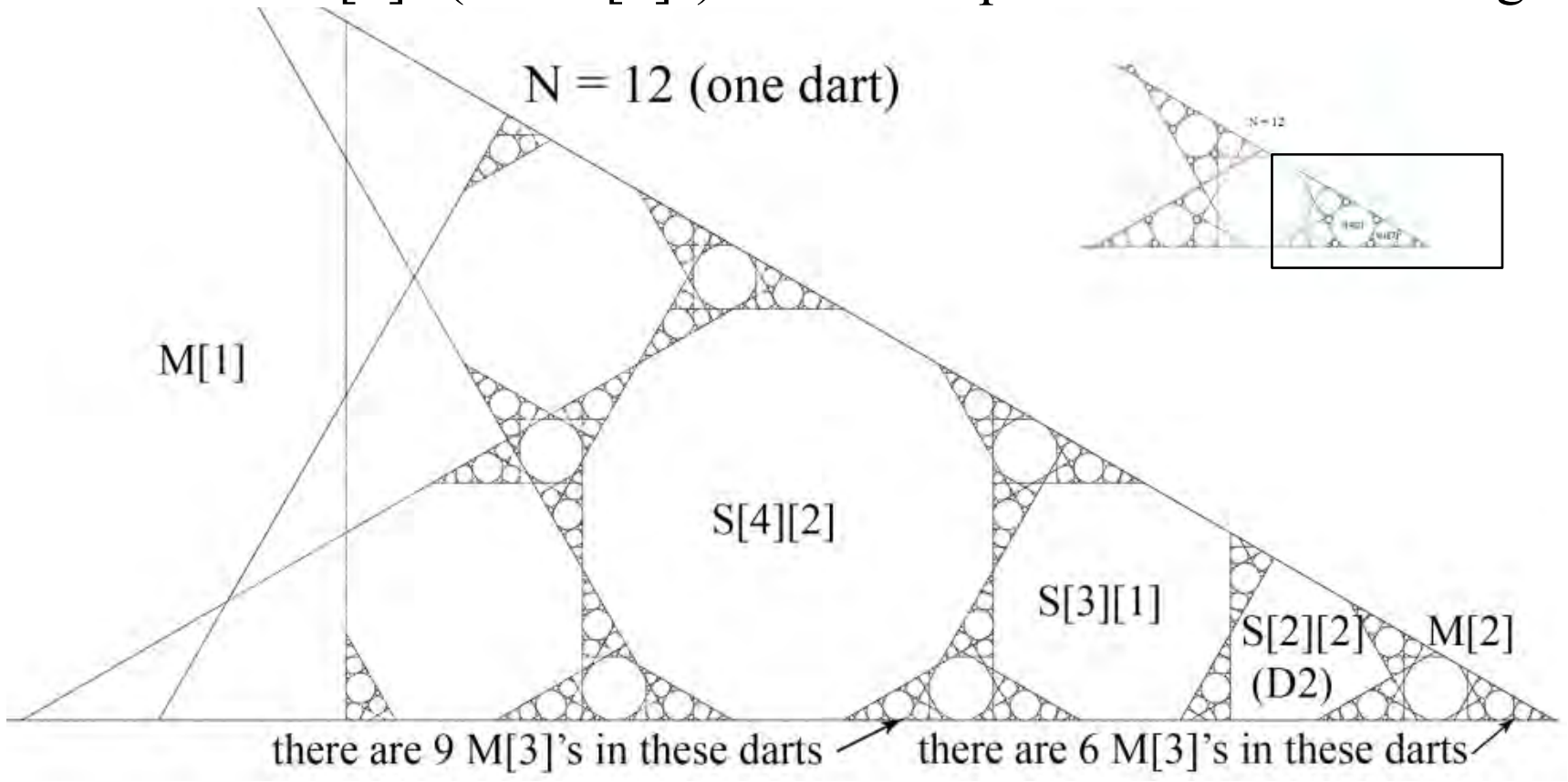

Another way to verify the temporal scaling for quadratic self-similar webs is to simply count the growth of the tiles using the $\tau$-periods. Here the $\tau$-periods of the M[k] in the canonical invariant 'star-region' of Example 4.2 are 60, 942, 28292, 775356, 21055308,.. with ratios of 15, 30, 27.41, 27.15. These are the combined periods of the M[k] at GenStar and their reflections about S[4]. Even though the local web has perfect reflective symmetry with respect to S[4], the dynamics are different and this combined count helps to minimize these differences. The dynamics of any composite N-gon allows for the possible 'decomposition' of expected orbits unto groups of orbits with smaller periods. This makes it difficult to match tile counts with periods, but for self-similar webs, the effect of these exceptions diminishes with each generation and in the limit the $\tau$-ratios will match the geometric ratios.

Returning to the First Family above, note that for N even, there is always an 'extended edge' relationship between N and M, with a variable number of intermediate ties. Here the intermediate tiles are D[1] and M[1], but with the mutation of S[2] in the 2$^{nd}$ generation, these intermediate tiles become virtual tiles in the First Family of DS[4].

The **8k+4 Conjecture** states that these Px 'parents' of DS[4] always exist and come in two versions. N = 12 is the charter member of a mod-16 family 12+16j where Px has DS[4] at S[3 +8j]. The next member is N = 28 where DS[4] is S[11] of Px and this in turn defines Px. The matching mod-16 family is based on N = 20 and here Px is simply a D tile of DS[4] but as evidenced by N = 52 the Px may now be virtual. But these virtual embeddings are guaranteed by the FFT, so this branch may tell us little about the geometry at star[1] of S[2]. Based on N = 44 in the N = 12 branch, even the real Px may have First Families with no actual survivors except DS[4] and the local geometry at DS[4] is a disrupted because of the slight offset from S[2].

**Example 5.7** (N = 26 & N = 13 ) These two have algebraic order 6 along with N = 21,28,36 and 42. Since N = 26 is the first member of the 8k+2 family with complexity that is neither quadratic or cubic, it was used as an example for the 8k+2 Conjecture earlier. Below is what we call a symmetry plot where the blue center lines show local web symmetry. When N is even the combined S[1] and S[2] web will be step-4 and the Edge Conjecture predicts that there will be surviving mod-4 DS[k]. The salient feature of the 8k+2 family is the early occurrence of DS[3]s with their matching DS[1]s and DS[2]s. It appears that both of these are needed so that the D[2]s will inherit the step-4 web structure of S[2] and support sequences of S[k] and M[k] converging to star[1] of S[2]. The 8k+2 Conjecture predicts that tiles like DS[7] will also survive in subsequent generations but this does not imply that volunteers like Px will survive. However for N = 26 the Py appear to survive in the web local to S[1] and their parameters are easy to find.

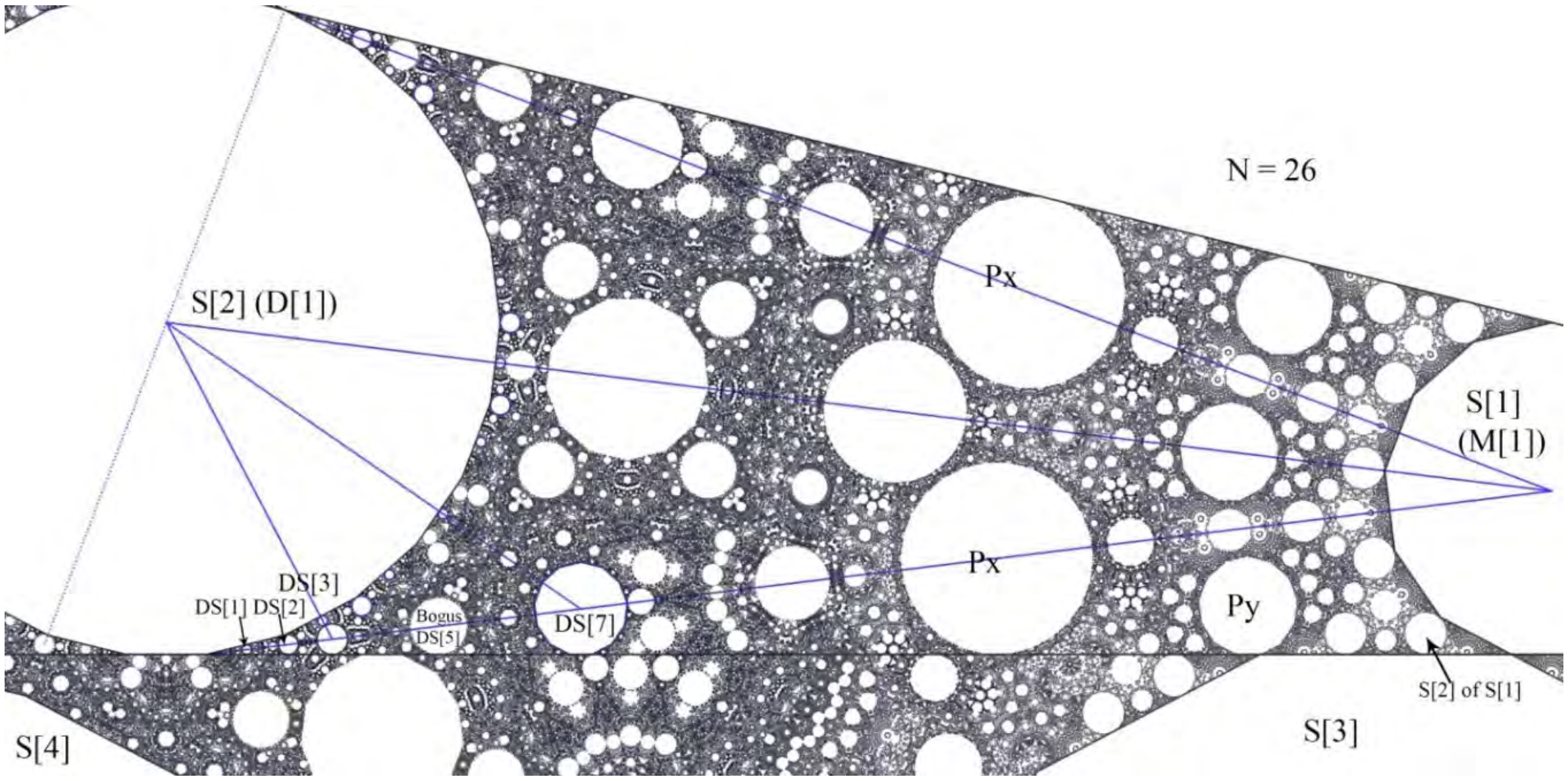

The Twice-Odd Lemma implies that the S[11] (M) tile of N = 26 will be a surrogate N = 13. But this tile is a different invariant region so its dynamics would be expected to be very different and we have no theory that relates this M[0] tile with M[1] here.

Below is a symmetry plot for N = 13 in the 8k+5 family. The Rule of 8 makes no predictions except for S[1] at DS[9] and an isolated DS[1] which could serve as a D[2], but there is no matching M[2] and it seems that the only First Family tiles it supports are S[4]s as vertex tiles.

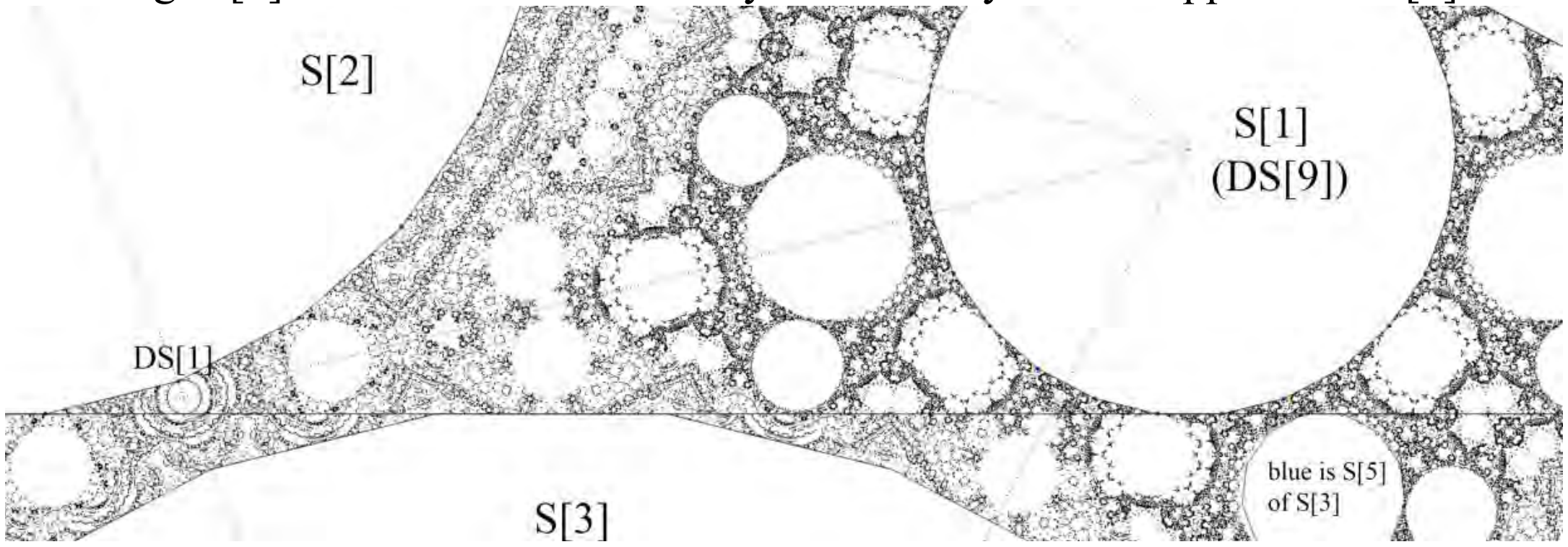

**Example 5.8** - The edge geometry of N = 16.
N = 16 has 'quartic' complexity along with N = 15, 20, 24 and 30. Both N = 16 and N =24 are in the 8k family so the Rule of 4 implies an isolated DS[2] will exist and here it acts as a D[2] and supports a chain of D[k] tiles with an almost-perfect chain of matching M[k].

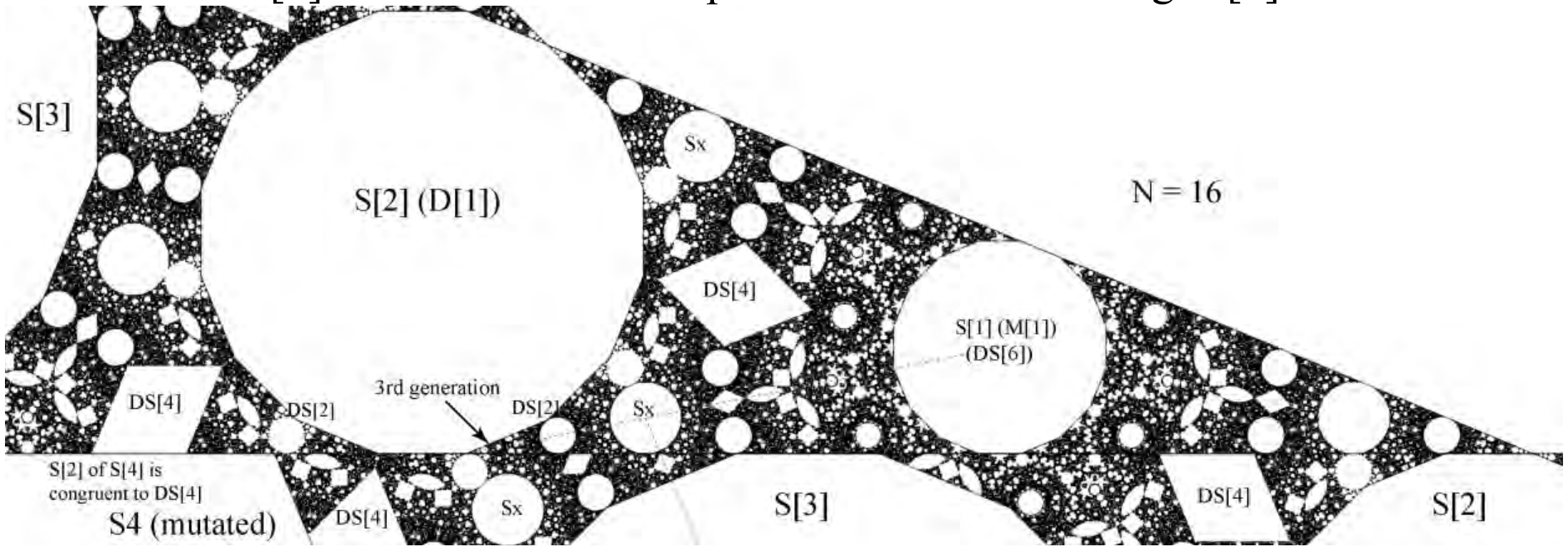

The S[4] of N is mutated since k′ = N/2-k = 4, and gcd(4,N) = 4. Therefore the mutated S[4] will consist of the octagonal weave of two squares. Normally any DS[4] of S[2] would inherit this same mutation, but here it seems the interaction between S[2] and S[4] turn DS[4] into a simple rhombus as shown here. This is an extended edge version of the mutation, which is not unique in the very 'lazy' web. The real surprise here is the chain of D[k] and M[k] which apparently exist even though M[2] is missing. The periods of the first 9 D[k] are: 8, 16·2, 8·57, 16·154, 8·2609, 16·6898, 8·121863, 16·322826, which gives even and odd ratios of about 5.3 and 8.8.

**Example 5.9** - The edge geometry of N = 20 showing a mutated S[2]

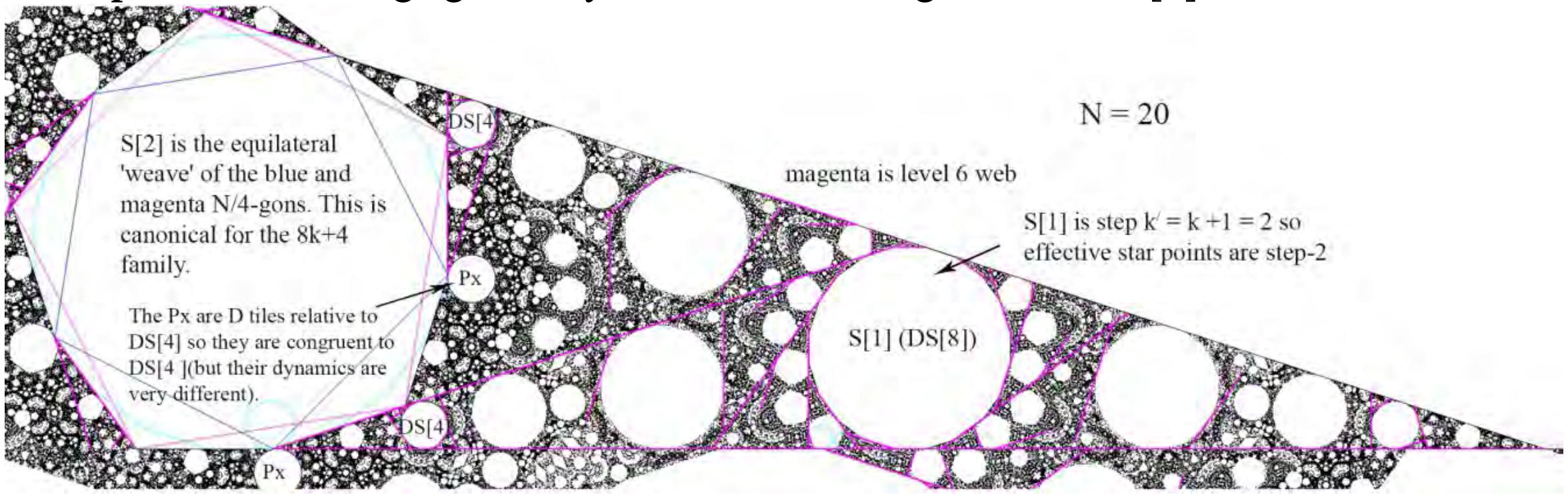

N = 20 is the second non-trivial member of the 8k+4 family. The Mutation Conjecture predicts that S[2], S[5] and S[6] will be mutated and the 8k+4 Conjecture above gives expected implications of the S[2] mutation. The 'odd' vertices of S[2] will always survive the mutation and the 'even' vertices will be extended outwards toward the predicted DS[4[ tiles. In the mod-16 subfamily anchored by N = 12 it appears that the DS[4]s will be in the First Family of Px 'parent' tiles centered at odd vertices of S[2], but for the N = 20 subfamily these Px tiles may be virtual beginning with N = 52 .Here with N = 20 the Px are still quite real and there is a surviving S[3] of Px with a promising web which could support future generations. But in general very little is known about the geometry at star[1] of S[2]. Even in the N = 12 subfamily where the parent Px tiles exist, a case like N = 44 has a Px where none of the First Family tiles survive the web, and there is no natural path towards recursion. A tile like S[2] with competing regular constituents is an example of a 'quasi-regular' polygon in Appendix G of [H3]. The special case of two congruent regular constituents is called 'semi-regular' in [S3].

## Appendix I   Exact calculations using symbolic dynamics

In [H2] (and in Example A3 below) we describe mappings from physics, astronomy, circuit theory and quantum mechanics which have singularity sets which may be locally conjugate to the outer-billiards map. Here we describe two such mappings and show how they can be used to perform exact calculations. The sample calculation will involve finding the parameters of the weakly conforming Mx tile of N = 11, as described in Example 5.4. We will first solve this problem with τ using 'corner sequences'. These sequences are examples of 'symbolic dynamics' as first formulated by George Birkhoff and Stephen Smale. See [De].

The other two maps in question are the 'digital-filter' map of Chau & Lin [CL] and a refinement of a 'dual-center' map of Arek Goetz [Go]. Like the outer-billiards map, these are piecewise isometries based on rational rotations - also known as affine piecewise rotations. These three maps appear to have conjugate webs and this may be due to the fact that each can be reduced to a form of shear and rotation. Even though the webs are congruent, the dynamics appear to be very different but we shall show that they have a consistent form of symbolic dynamics.

**Example A1** Use the evolution of the 'web' W to find the parameters of the Mx tile of N = 11 with (i) the outer-billiards map, (ii) the digital-filter map and (iii) a complex-valued Goetz map.

**Part (i): The outer-billiards map**. In Example 5.4, we noted that the Mx tile of N = 11 is a weakly conforming regular N-gon that occurs in the $2^{nd}$ generation for N = 11 on the edges of D. Therefore star[5] of Mx is star[1] of DS[1] as shown below. By the Two-Star Lemma, the parameters of Mx can be determined by finding another star point of Mx. In this example we will find the star[4] point of Mx by tracing the web evolution of the interval H1 - in the context of N with radius 1. (By reflection, this evolution can equally be studied on the edges of N = 22.)

The interval H1 lies on the horizontal base edge of N = 11 so there are 11 such intervals equivalent to H1 under rotation, and the local $\tau$- web determined by H1 includes the iteration of each of these rotated copies. Under $\tau$, these 11 regions map to each other and after 99 iterations of each interval, 8 segments land back at D[1] as shown here.

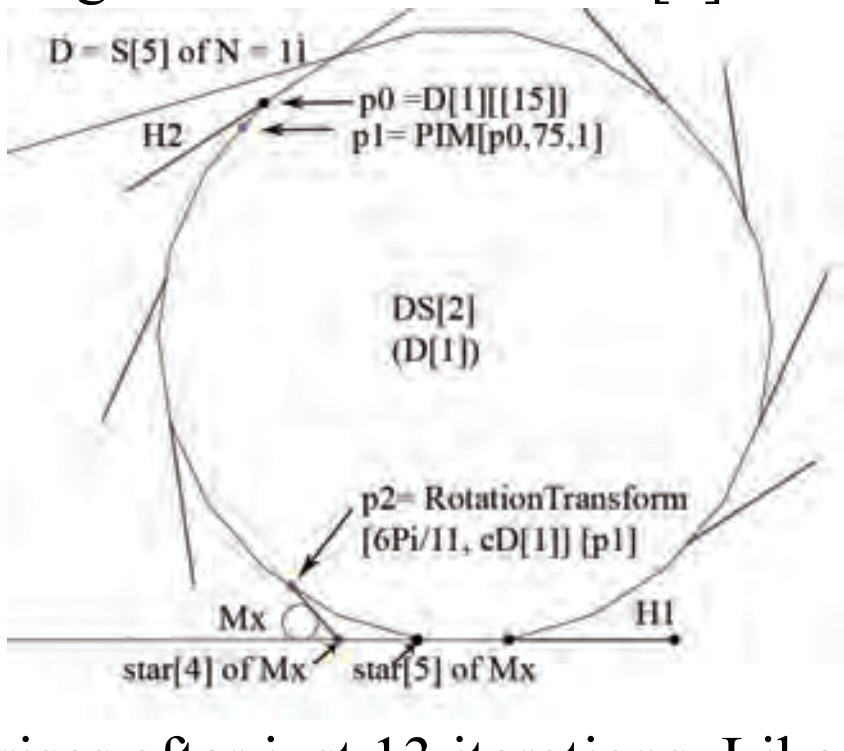

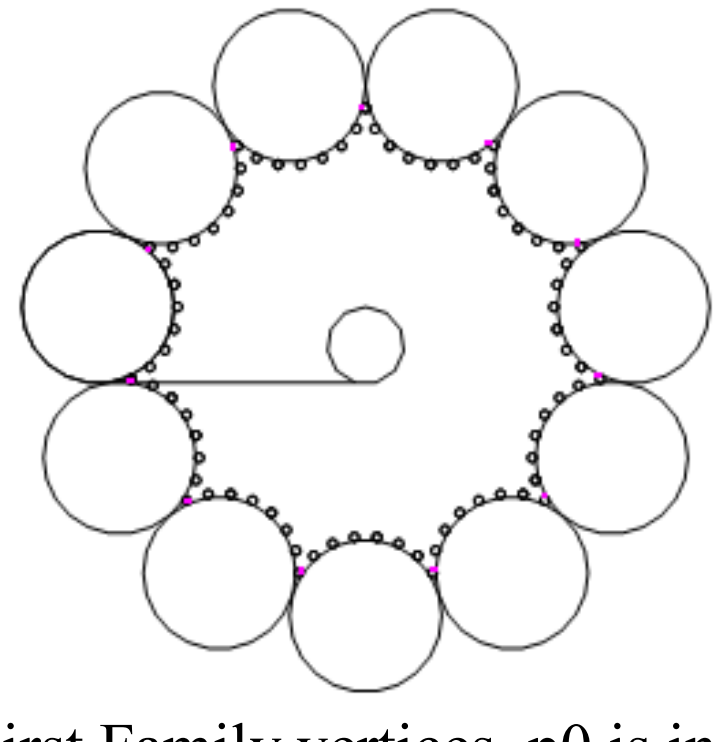

The interval H2 arises after just 13 iterations. Like all First Family vertices, p0 is in $\mathbb{Q}_{11}$ and we suspect that it maps to p1 (as a one-sided limit) - and hence determines the offset of p2. These vertices technically have no image under $\tau$, so we will find $p_1$ using the 'surrogate' orbit of a point that is close to p0 and on the interval H2. (Typically intervals like H2 will get truncated under iteration, but by inspection it is clear that the inner portion of H2 survives well beyond the few hundred iterations needed here.)

Here are the calculations using Mathematica:

(i) Since the image under $\tau$ of any edge is a parallel edge, the slope of H2 is known, so set p0N to be a point on H2 within 8 decimal places of p0. **Orbit = NestList[$\tau$, p0N, 200]** (an approximate orbit but initially reliable and any errors are easy to detect).

(ii) Since $\tau(p) = 2c_j - p$ for some vertex $c_j$ of N = 11, $\tau^k(p) = (-1)^k p + 2Q$ where Q is a +/– alternating sum of vertices. Every $\tau$-orbit determines a sequence $\{c_k\}$ of vertices and the matching indices are sufficient to find Q and determine the orbit. The study of these partition sequences is called 'symbolic dynamics' so we will call them S-sequences. The Mathematica module IND will use τ to find the S sequences to any depth (once again with possible error).

S[p0N,150] = **IND[p0N,150]** = {11,5,10,4,9,3,8,1,6,11,5,10,4,9,2,7,1,6,11,5,…}

Note that these indices initially advance by {5,5,5,5,5,5,4} (mod 11) because D is S[5] with step sequence {5} and D[1] has (periodic) step sequence {5,5,5,5,5,5,4}. Here this sequence will eventually break down. In general no web point can have a periodic orbit because these points have no inverse. We will use these indices in pairs, using the 'return' map $\tau^2(p) = p + 2(c_k - c_j)$.

(iii) P1 = **PIM[p0N,75,1]** will take IND and these 150 indices in pairs and reconstruct the orbit, while P3 = PIM[p0N,75,3] will construct a step-3 version of this orbit, which is called a 'projection' or algebraic graph as defined in [S2]. To get an exact orbit, simply use p0 instead of p0N.

P1 = PIM[p0,75,1] ; p1= P1[[75] = $\tau^{150}$(p0) ; p1[[1]] =

$$-6Cos[\frac{3\pi}{22}] - Cos[\frac{\pi}{11}]Cot[\frac{\pi}{22}] + 4Sin[\frac{\pi}{11}] + 8Sin[\frac{2\pi}{11}] + Cot[\frac{\pi}{22}]Sin[\frac{\pi}{11}]Tan[\frac{\pi}{11}] - Sec[\frac{\pi}{22}]Sin[\frac{\pi}{11}]Sin[\frac{5\pi}{22}]Tan[\frac{\pi}{11}]$$

(iv) Rotate by 6π/11 about the center of D[1]: p2= **RotationTransform[6·Pi/11, cD[1]][p1]**

(v) As indicated earlier, the slope of the web interval determined by p2 must match an edge of N. Here it has the same slope as the (right-side) star[4] edge of N, which we call slope4.
x1 = p2[[1]], y1 = p2[[2]]; b = y1– slope4·x1 so star[4][[1]] = (1– b)/slope4

(vi) By the Two-Star Lemma hMx = d/(Tan[5Pi/11]–Tan[4Pi/11]) where d is the horizontal displacement of star[4] and star[5].

(vii) Of course the displacement d depends of hN but that dependence vanishes when hMx is divided by hN, which here is Cos[π/11].

**AlgebraicNumberPolynomial[ToNumberField[hMx/hN,GenScale[11], x]** =

= $1 - 23x - \frac{27x^2}{2} + \frac{x^4}{2}$ where x = GenScale[11] = Tan[π/7]·Tan[π/14].

Any other hN and matching hMx must yield this same ratio so this is a fundamental polynomial for Mx. Any 'canonical' tile with scaling in $S_N$ will have such a polynomial. See Example A2.

**Part (ii) – The Digital Filter Map** The Df map ([CL] and [H2]) is only compatible with the outer-billiards map when N is even, so it is necessary to work inside N = 22 (with hN = 1). This is not a burden and actually simplifies the calculations. Except for scale, the First Family is unchanged from part (i) - but S[9] is now playing the part of N = 11, so we will show that hMx/hS[9] satisfies the polynomial above.

The Digital Filter map Df: $[-1,1)^2 \to [-1,1)^2$ is defined as Df[{x,y}]:={y, *f* (–x + *a*y)} where *f*(v) = Mod[v+1,2]-1 models a 2's complement sawtooth register. In matrix form (where *f*(y) ≡ y):

$$\begin{bmatrix} x_{k+1} \\ y_{k+1} \end{bmatrix} = f\left(\begin{bmatrix} 0 & 1 \\ -1 & a \end{bmatrix}\begin{bmatrix} x_k \\ y_k \end{bmatrix}\right) = \begin{bmatrix} y_k \\ f(-x_k + ay_k) \end{bmatrix}$$

Setting $a = 2\cos\theta$ the matching elliptical rotation

$$\begin{bmatrix} 0 & 1 \\ -1 & 2\cos\theta \end{bmatrix}$$

is conjugate to a true rotation

$$\begin{bmatrix} \cos\theta & \sin\theta \\ -\sin\theta & \cos\theta \end{bmatrix}$$

but this conversion is optional.

To generate the web of N = 22, set θ = 2π/22 and here we are interested in generating symbolic orbits. The S sequences will be much simpler than τ, because the function *f* distinguishes just 3 regions - so Df is a piecewise isometry with three primary regions (atoms) which can be labeled 1 (overflow), 0 (in bounds), or –1 (underflow).The equations for these atoms are given here.

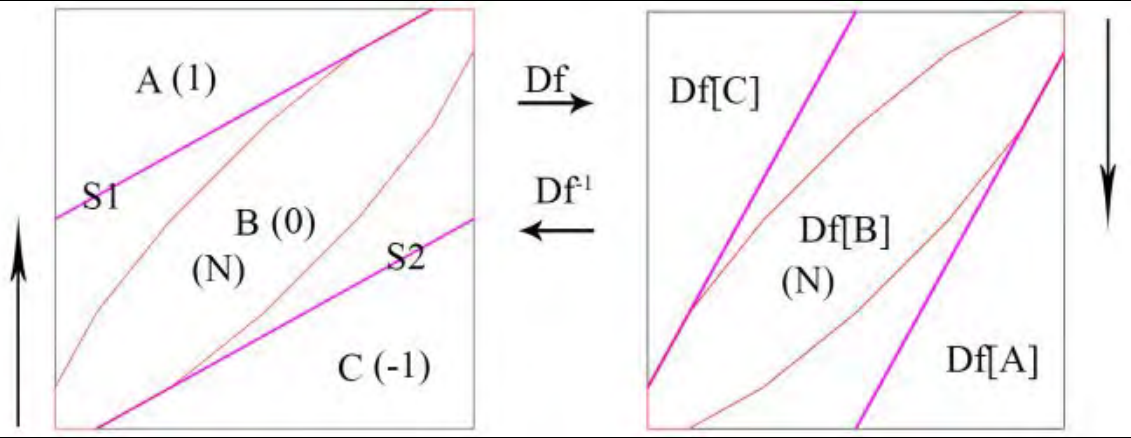

$$S[\{x,y\}] = \begin{cases} 1 & \text{if } -x + ay \ge 1 \\ 0 & \text{if } -1 \le -x + ay < 1 \\ -1 & \text{if } -x + ay < -1 \end{cases}$$

Under Df all points experience a rotation by θ and *f* determines the corresponding 'vertical' shear of –2, 0 and 2 respectively for A, B and C. The central B region is free of translation since it is 'in-bounds', so points will rotate by θ, and under iteration they will construct copies of N - one of which is shown here. The seperatices S1 and S2 define the maximal extent of this linear (elliptical) rotation, so they define the bounds of the three regions. Therefore if the S-sequence of a point p is known the Df map is simply DfS[{x,y},k_]:={y,–x+*a*y–2S[[k]]}

**Example**: The first 10 points in the Df orbit of cDS[1] for N = 22. Set p = TrToDf[cDS[1] ] where TrToDf is a change of coordinates between traditional Euclidean space and Df toral space.

Like all points adjacent to D, p ={x,y} will be in an 'overflow' so Df[p] = {y,–x+*a*y-2} as shown here. The point {y,–x+*a*y} is 'out of bounds' at top right and the {0,-2} shear is the *f* correction. Since Df[p] is 'in-bounds', it allows for 8 consecutive central rotations, but the last rotation will yield an 'underflow' which will generate a compensating displacement back to an overflow position to repeat the cycle. Therefore the S sequence will be {1,0,0,0,0,0,0,0,0, -1,1} and this periodic sequence determines the orbit of p.

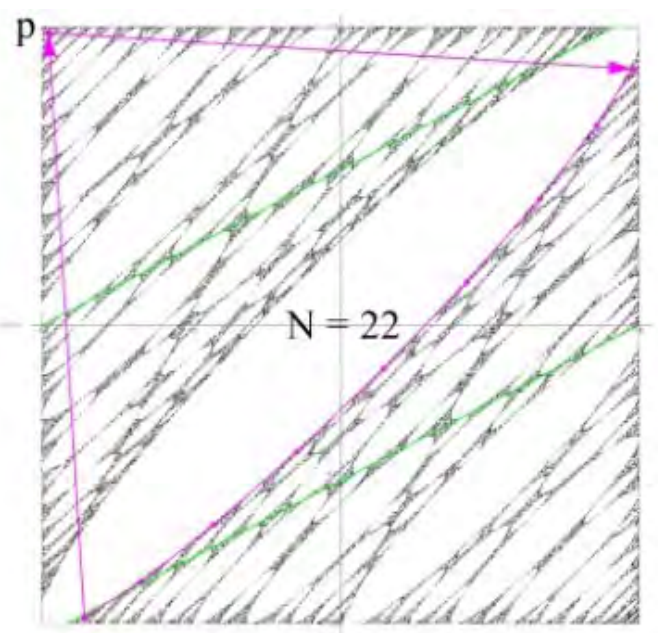

To get the level-k S sequence of any point p, first generate Orbit = NestList[Df,p,k] and then apply the function S to the elements of this orbit: Sequence = S/@Orbit.

As a guide we will use part (i) above to construct the expected Mx tile inside N = 22

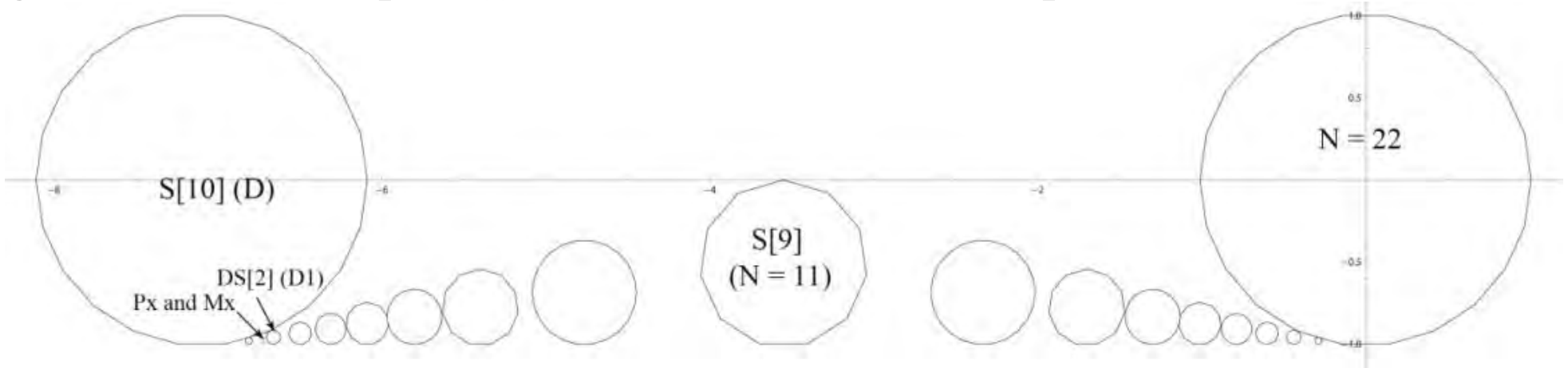

The closest connection between Mx and the First Family of D1 appears to be the (virtual) step-3 tile of D1 – which we call DS[3][2] or simply DS3. We will use DS3 to define an interval that maps to Mx. The green interval shown below is a rectified portion of an edge in W, and under Df this interval will map to the blue interval – when rectified. Therefore p0 maps to px (using Df) and p0 is exact. It only takes 136 iterations to accomplish this – but an exact calculation will be awkward with Df - so we will use surrogate orbits and DfS instead.

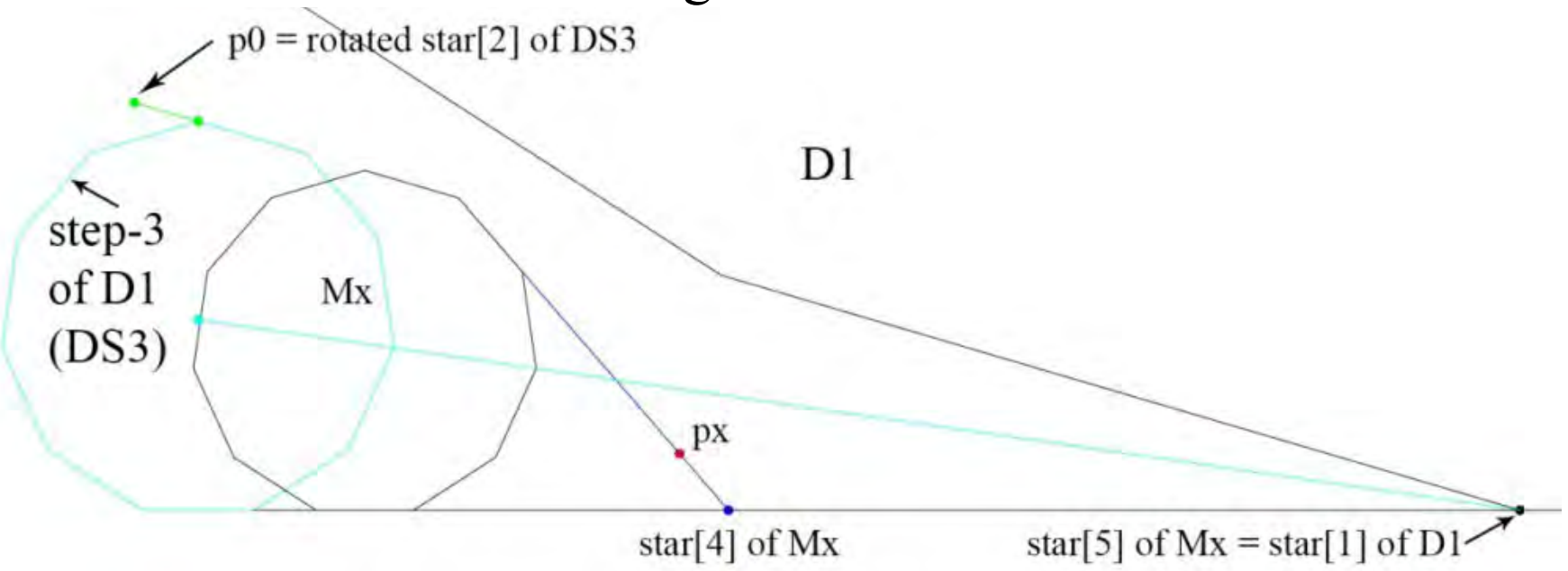

Here are the calculations :

(i) Use p1 = TrToDf[p0] to generate an approximate Df orbit of length 140 using w = 2Cos[2Pi/22] to 30 decimal places. Orbit=NestList[Df, p1,140]; S=S/@Orbit ={1,0,0,0,0,0,0,0,0,-1,1,0,0,0,…} (As indicated earlier the points at the foot of D are typically in overflow positions so sequences such as this are common. Underflow must map to overflow because –1 and 1 are only one bit apart in 2's complement.)

(ii) Generate an exact value of $Df^{137}$[p1] using DfS and exact w = 2Cos[2Pi/22]. Set q=Range[140]; q[[1]]= p1; **For[s=1, s<=140, s++,q[[s+1]]=Simplify[Re[DfS[q[[s]],s]]]];**

pxx (px unrectified) = q[[137]] ≈ {–6.679854552975, -0.9991754019926052}

px = Simplify[Re[DfToTr[pxx]]]; px[[1]] =

$$(-422+419(-1)^{1/11}-414(-1)^{2/11}+419(-1)^{3/11}-422(-1)^{4/11}+421(-1)^{5/11}-420(-1)^{6/11}+424(-1)^{7/11}-424(-1)^{8/11}$$
$$+420(-1)^{9/11}-421(-1)^{10/11})/((1+(-1)^{2/11})(-1+(-1)^{1/11}+(-1)^{6/11})^{2})$$

(iii) Set x1 = px[[1]]; y1 = px[[2]]; b = y1–slope4·x1 (where slope4 is the slope of right-side star[4] of N as in part (i)). Therefore star[4] [[1]] = (1–b)/slope4 =

$$\frac{i(-13-(-1)^{1/11}+20(-1)^{2/11}+14(-1)^{4/11}-14(-1)^{5/11}-20(-1)^{7/11}+(-1)^{8/11}+13(-1)^{9/11})}{(1+(-1)^{1/11})(1+(-1)^{2/11})^{4}}-2Cot[\frac{\pi}{11}]$$

(iv) As in part (i), hMx = d/(Tan[5Pi/11] –Tan[4Pi/11]) where d is the horizontal displacement of star[4] and star[5].

(v) This displacement d is relative to hN so the ratio hMx/hN is the scaling field $S_{11}$, but N is a 22-gon, so it makes more sense to use hMx/hS[9] where hS[9] = Tan[π/22]/Tan[π/11].

**AlgebraicNumberPolynomial[ToNumberField[hMx/hS[9], GenScale[11]],x]** gives

$1-23x-\frac{27x^2}{2}+\frac{x^4}{2}$ as in part (i) above.

**Part (iii) A complex-valued Goetz Map**

Perhaps the simplest map that reproduces the outer-billiards web for a regular N-gon is a y-axis version of a 'dual center' map of Arek Goetz. As described in [H2] this mapping has the form Dc[z] = Exp[– Iw]·(z – Sign[Im[z]) where w = 2π/N and for a scalar v , Sign[v] = 1, 0 or –1 iff v is positive, 0 or negative. Therefore Dc is a pure (clockwise) rotation for points on the x-axis but in general it is a 'shear and rotate' where the shear has magnitude -1 above the x-axis and +1 below (in a manner similar to W with sN = 1). So if Dc is applied to the interval [-1,1] it will form a perfect N-gon above the x-axis and its negative below. If this initial interval is expanded as shown below for N = 11 (with w = 2π/11), then the edges of –N will intersect the edges of N to simulate the interaction of the τ-domains. This juxtaposition actually occurs with W but not at the origin.

The web on the left below was generated by iterating the x-axis interval [–12,0], 200 times under Dc with w = 2π/11 (to 30 decimal places). In the limit, the region above (or below) the x-axis will be a perfect reproduction of the web for N = 11 with a side of 1 and star[1] at the origin. It is an easy matter to scale the First Family (and Mx) to use as guides to track the web development at D. On the right we show a portion of the ideal second generation in black and the hypothetical Mx in cyan. Unlike the Df map, the webs tend to evolve is a predictable fashion from intervals on the x-axis and this makes it easy to find intervals that will map to Mx. The chosen interval is shown in blue in the enlargement on the lower right, and the magenta interval is the negative of the image of this blue interval under $Dc^{564}$ .

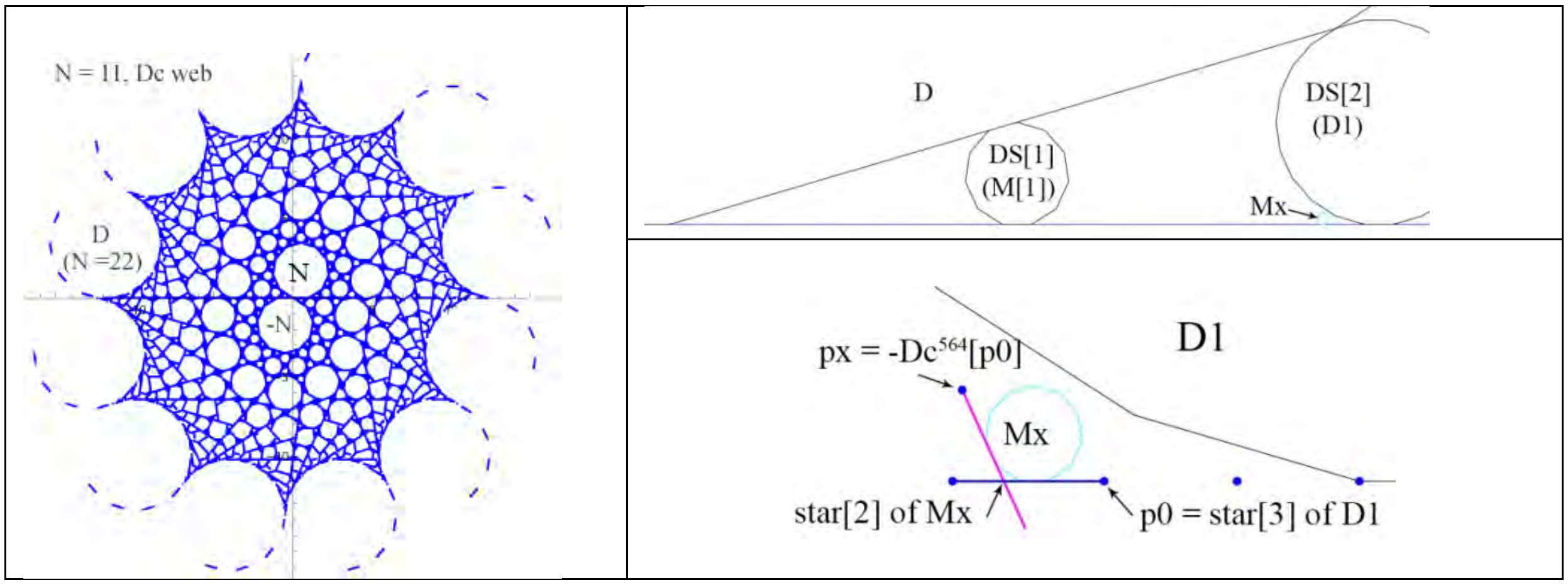

Since the slope of this magenta line is known, all that is needed to find star[2] of Mx is some point like px on this line. Here are the calculations.

(i) The initial point is p0=StarD1[[3]] =

$$\{\frac{5}{2}-Cot[\frac{\pi}{22}]Cot[\frac{\pi}{11}]+Tan[\frac{\pi}{22}]Tan[\frac{\pi}{11}]+\frac{1}{2}Cot[\frac{\pi}{22}]Cot[\frac{\pi}{11}](1-Tan[\frac{\pi}{22}]Tan[\frac{\pi}{11}](2+Tan[\frac{\pi}{22}]Tan[\frac{\pi}{11}]))-(\frac{1}{2}Cot[\frac{\pi}{11}]+\frac{1}{2}Cot[\frac{\pi}{11}](-1+Tan[\frac{\pi}{22}]Tan[\frac{\pi}{11}](2+Tan[\frac{\pi}{22}]Tan[\frac{\pi}{11}])))Tan[\frac{3\pi}{22}],0\}$$

(ii) Using p0N = N[p0,30] (30 decimal place approximation to p0), find the first 600 points in the (complex-valued) orbit: Orbit = NestList[Dc,p0N,600]

(iii **S[z_]:=Sign[Im[z]];** Sx=S/@Orbit = {0,1,1,1,1,1,1,–1,–1,–1,–1,–1,–1,1,1,1,1,1,…}

(iv) Define a 'literal' version of Dc based on Sx , namely **DcS[z_,k_]:=Exp[CI·w](z–Sx[[k]])** (with an exact w) to obtain the exact orbit of p0 using the Sx sequence from the surrogate orbit. Store the orbit in the sequence q:

**q=Range[600];q[[1]]=p2[[1]]; For[s=1,s<=600, s++, q[[s+1]]=Simplify[DcS[q[[s]],s]]];**

**Note:** By modifying F to allow approximate calculations of the Sign function, it is may be feasible to do calculations like this directly with NestList. One possible modified function is F[z] = Simplify[Exp[–I*w](z – I*IntegerPart[Sign[N[Re[z]]]])]. This tricks Mathematica into regarding the Sign output as exact. Normally Mathematica attempts to evaluate Sign in an exact fashion and these expressions get so complex that typically it fails after a few hundred iterations.

(v) q[[564]] = $\frac{6-8(-1)^{1/11}+5(-1)^{2/11}-2(-1)^{3/11}-(-1)^{4/11}+3(-1)^{5/11}-(-1)^{6/11}+5(-1)^{7/11}-3(-1)^{8/11}+(-1)^{9/11}+4(-1)^{10/11}}{-1+2(-1)^{1/11}-3(-1)^{2/11}+3(-1)^{3/11}-2(-1)^{4/11}+(-1)^{5/11}}$

The desired point is px = – q[[564]] ≈ –10.415046959467414 + 0.01643324914196498·I

(vi) x1 = Re[px]; y1 = Im[px]; The slope of star[2] of Mx is –slope2 as defined using star[2] of N, so b = y1+ slope2·x1 and star[2] [[1]] = b/slope2 ≈ –10.407542146047456045

(vii) Using the Two-Star Lemma with opposite sides hMx = d/(Tan[2Pi/11]+Tan[5Pi/11]) where d = star[5][[1]]–star[2][[1]].

(viii) Since N has side 1, hN = $\frac{1}{2}Cot[\frac{\pi}{11}]$. Use this to convert hMx to a scale in $S_{11}$.

**AlgebraicNumberPolynomial[ToNumberField[hMx/hN],GenScale[11]],x]**

= $1-23x-\frac{27x^2}{2}+\frac{x^4}{2}$ as in parts (i) and (ii)

Since the star[3]-star[4] interval of D1 generates Mx, it can be regarded as a 'mutated' DS3[2]. Likewise Px is generated by the star[6] to star[7] interval so it may be a modified DS6[2].

The historical connection between these maps is illustrated in Example A3 below - which shows that Df is equivalent to a sawtooth (tuned) version of the classical Standard Map while Dc is a sawtooth (tuned) version of a kicked harmonic oscillator. The connection is that the Standard Map is equivalent to a kicked 'free' rotor in zero gravity while a harmonic oscillator is assumed to be affected by gravity – and hence has a natural frequency of oscillation.

The classical Twist Map of J. Moser is also a 'free rotor' but it allows for the possibility of periodic kicks so the Standard Map illustrates the case of a Twist Map with these (ideal) periodic 'gravitational perturbations' turned on - and the issue is whether there will be a non-zero measure of initial conditions which yield 'integrable' solutions to the matching Hamiltonian.

It is still not clear what role gravity plays in quantum mechanics. In the late 50's P.B. Harper and others used a kicked harmonic oscillator with a natural frequency of oscillation to model both classical and quantum diffusion - based on a stationary Schrodinger equation. This is illustrated in Example A3 for $\omega = 2\pi/4$ (N = 4) which shows initial resonances similar to the Standard Map – and then diffusive breakdown and weak mixing as K increases.

The Df and Dc maps are very efficient ways to generate the local web for a regular N-gon but Dc has a number of advantages and we will be using it as an aid toward publishing a detailed '4K' catalog of edge geometry in [H6]. For this purpose the juxtaposition of N and –N is ideal and the natural +/– symmetry of W is augmented by reflective symmetry to yield very efficient maps.

Below is an example from N = 14 where we iterate 1,000 points in the interval H = {–2,–1} at a depth of 5,000. (The interval {–1,1} will generate N and –N in a period N orbit.) Here we crop these 5 million web points and their negatives and reflections to the desired region. (Less than 1 minute to generate and 1 minute to crop on a modest computer.)

**Example A2** (The edge geometry of N = 14) Dc[z_]:=Exp[–I*w]*(z–Sign[Im[z]]); w =N[ 2*Pi/14, 35];(35 decimal places); H=Table[x,{x, –2, –1, .001}]; Web = Table[NestList[Dc, H[[k]], 5000],{k, 1, Length[H]}]; RealWeb = {Re[#], Im[#]}&/@Web; WebPoints = Crop[Union[RealWeb, -RealWeb, Reflection[RealWeb]]]; (about 450,000 points) Graphics[{AbsolutePointSize[1.0], Point[WebPoints]}}]

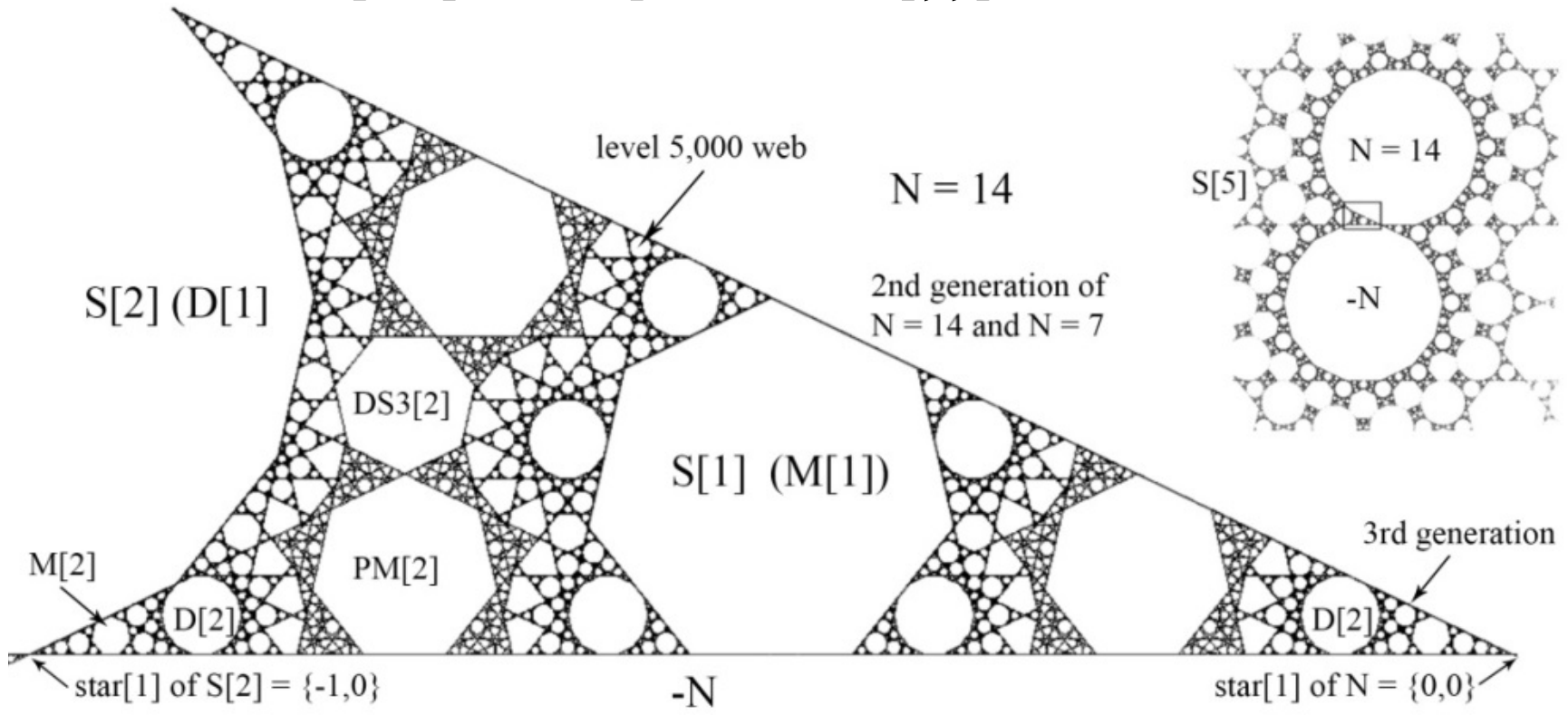

Example 5.3 discussed the N = 14 and N = 7 families, and here we examine the scaling. The scaling fields $S_7$ and $S_{14}$ are generated by x = GenScale[7] = Tan[π/7]·Tan[π/14]. Inside N = 14, S[5] is the surrogate N = 7, so by convention all heptagons are scaled relative to S[5]. M[1] (a.k.a.M1) is not a true 2nd generation N = 7 (except by scale) because the edge dynamics are very different with PMs in place of S[2]s. However M[2] is a valid 3rd generation 'matriarch' with the same edge geometry as N = 7. The even and odd generations appear to be self-similar.

| hM1/hS[5] | hD1/hN | hDS3[2]/hS[5] | hPM[2]/hS[5] | hS[5]/hN | hM2/hS[5] | hD2/hN |
|---|---|---|---|---|---|---|
| x | x | $\frac{1}{2}-4x-\frac{3x^2}{2}$ | $-\frac{3}{8}+\frac{17x}{4}+\frac{9x^2}{8}$ | $\frac{3}{7}+\frac{3x}{7}-\frac{x^2}{7}$ | $x^2$ | $x^2$ |

Example A3 - Some classical maps relevant to this paper

Twist Maps (zero gravity)

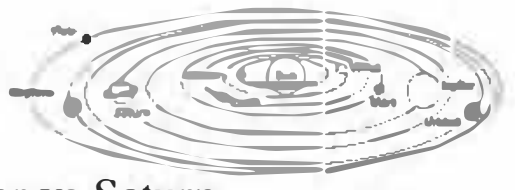

Jupiter vs Saturn

Poincare Section

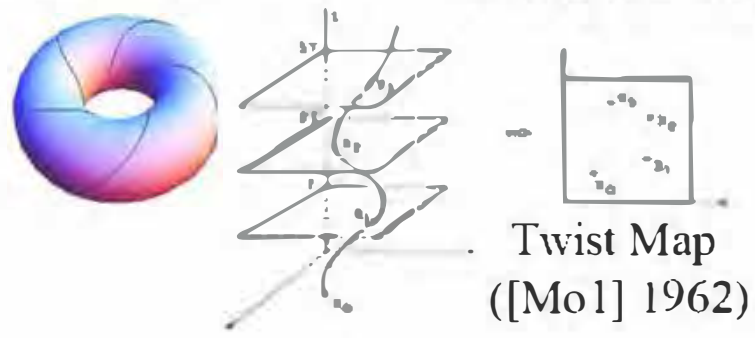

Twist Map ([Mo1] 1962)

$$\begin{bmatrix} x_{k+1} \\ y_{k+1} \end{bmatrix} = \begin{bmatrix} x_k \\ y_k + x_k \end{bmatrix}$$

(Here y is $\theta$ and x is momentum (constant))

Standard Map (Chirikov 1969) (K = .971635406)

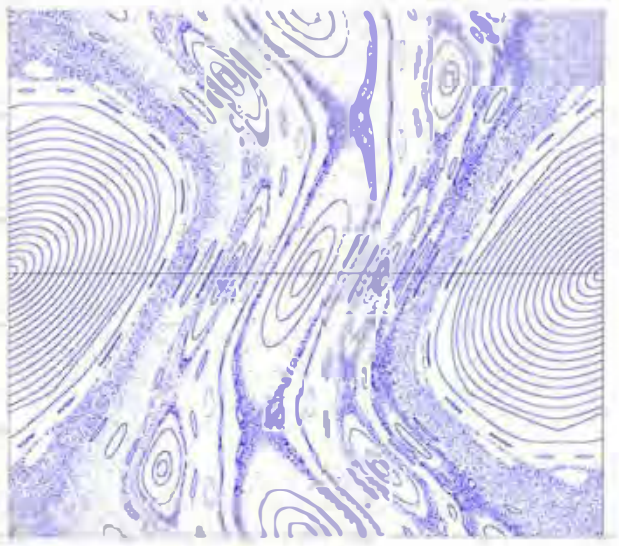

This is a 'Mod -1' periodically perturbed Twist Map showing resonant 'islands' of stability on all scales

$$\begin{bmatrix} x_{k+1} \\ y_{k+1} \end{bmatrix} = \begin{bmatrix} x_k + (K\sin 2\pi y_k)/2\pi \\ y_k + x_k + (K\sin 2\pi y_k)/2\pi \end{bmatrix}$$

(gravity turned on)

Sawtooth Standard Map ([A] 2001) ( k = 2Cos(2Pi/14)-2)

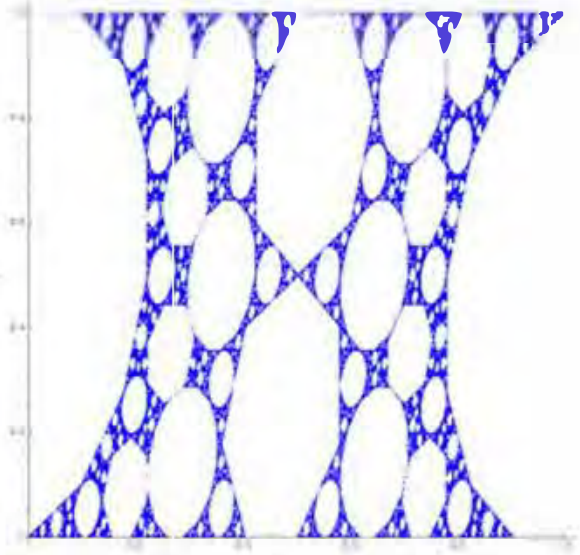

The S[k] of N = 14 are the primary resonances of this Standard Map tuned to 'Ch14'

$$\begin{bmatrix} x_{k+1} \\ y_{k+1} \end{bmatrix} = \begin{bmatrix} x_k + k\mathrm{Saw}(y_k) \\ y_k + x_k + k\mathrm{Saw}(y_k) \end{bmatrix}$$

(Saw(x) = x - 1/2 on [0,1] )

(The outer-billiards web W can be regarded as a sawtooth (tuned) version of the Standard Map)

These are equivalent under a change of variables

Harmonic Oscillator (non-zero gravity)

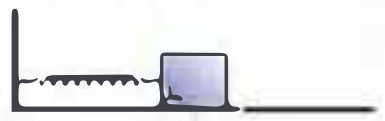

natural frequency of oscillation $\omega$

As in the Standard Map there may be periodic 'kicks' (and possible friction)

Kicked Harmonic Oscillator (KHO) Harper Map ($\omega$ = 2Pi/4) ([Ha]1955)

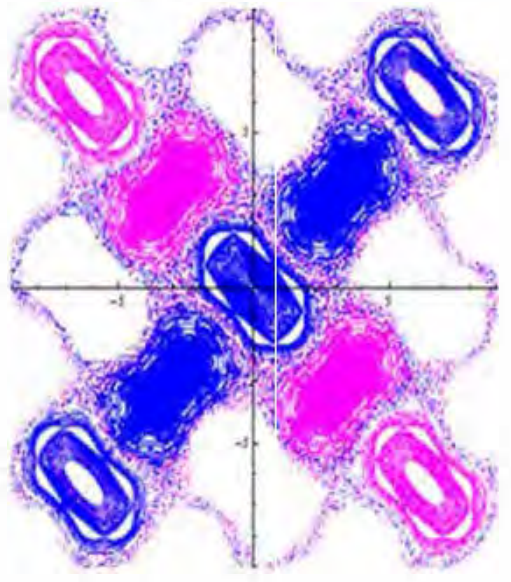

$$\begin{bmatrix} x_{k+1} \\ y_{k+1} \end{bmatrix} = \begin{bmatrix} \cos\omega & \sin\omega \\ -\sin\omega & \cos\omega \end{bmatrix} \begin{bmatrix} x_k + K\sin y_k \\ y_k \end{bmatrix}$$

(sometimes unequal kicks are applied to both x and y)

Dissipative KHO (DKHO) ([H2] 2012) ($\omega$ = 2Pi/7 (N = 7))

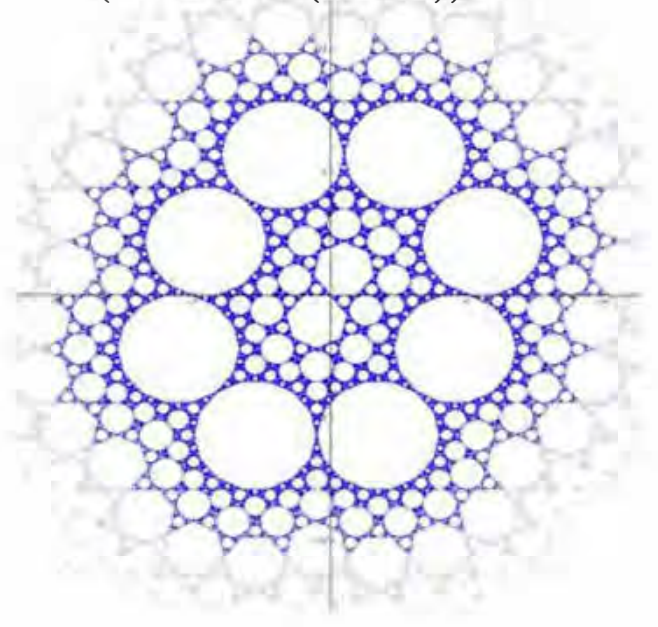

$$\begin{bmatrix} x_{k+1} \\ y_{k+1} \end{bmatrix} = \begin{bmatrix} \cos\omega & \sin\omega \\ -\sin\omega & \cos\omega \end{bmatrix} \begin{bmatrix} x_k - \mathrm{sgn}\, y_k \\ y_k \end{bmatrix}$$

Replacing sine with sgn (sign) yields a sawtooth (tunable) map that we call the dual-center map. In complex form: F[z] = Exp[-I$\omega$]*(z- Sign[Im[z]])

Digital Filter Map ([CL] 1988) ($a$ = 2Cos($\omega$) - here $\omega$ =2Pi/14

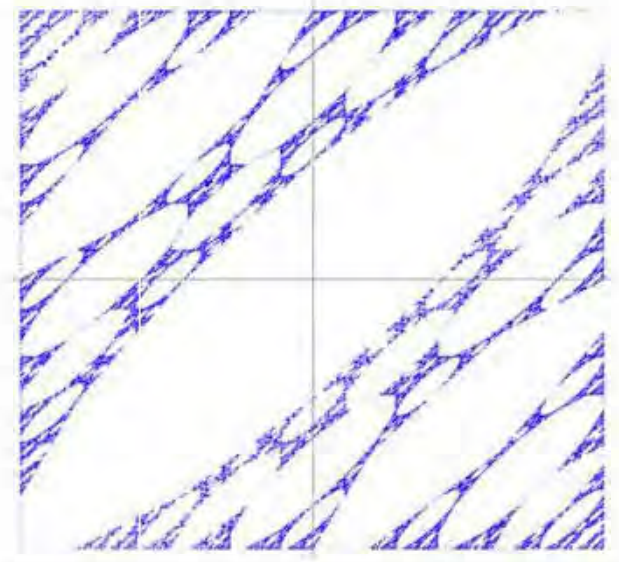

$$\begin{bmatrix} x_{k+1} \\ y_{k+1} \end{bmatrix} = f\left(\begin{bmatrix} 0 & 1 \\ -1 & a \end{bmatrix} \begin{bmatrix} x_k \\ y_k \end{bmatrix}\right) = \begin{bmatrix} y_k \\ f(-x_k + a y_k) \end{bmatrix}$$

Since $a$ = 2Cos($\omega$) this is conjugate to $\begin{bmatrix} \cos\omega & \sin\omega \\ -\sin\omega & \cos\omega \end{bmatrix}$

The (sawtooth) register overflow function is $f(x)$ = Mod[x+1,2] -1 and $f(y_k) = y_k$ because y is always in range

Overlay of Sawtooth Standard Map (magenta), Digital Filter Map (blue) and rectified N = 14 web in black

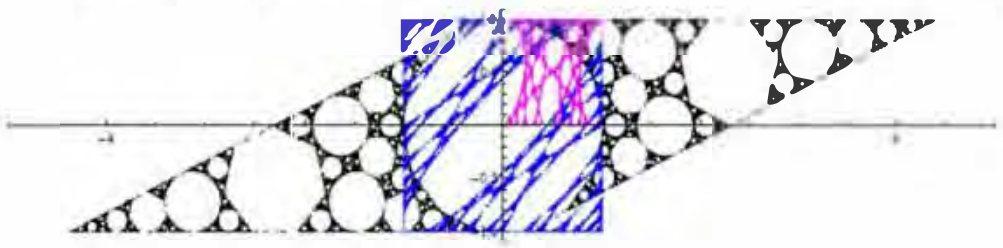

## Appendix II     Invariant Structures of Regular Polygons

| N = 20 (Twice-even ) | N = 14 (Twice-odd) | N = 11  (Odd) |
|---|---|---|
| 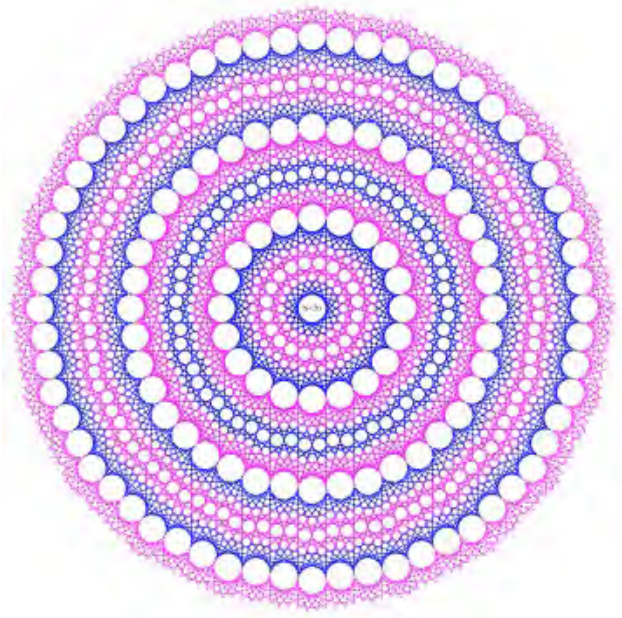 | 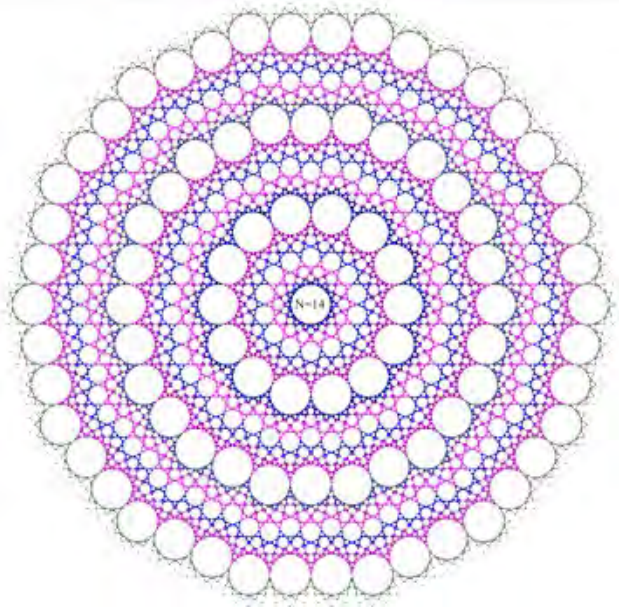 | 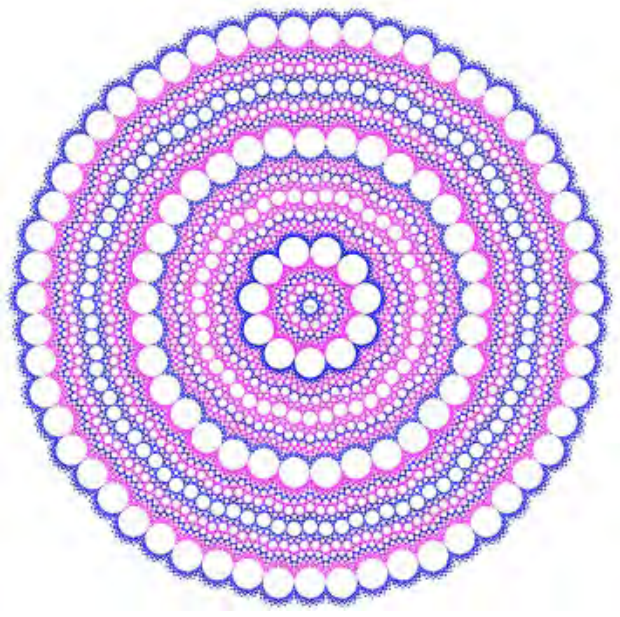 |

The periodic concentric rings of D tiles described in Shaidenko-Vivaldi [SV](1989) and later in [H3] *Outer billiards on Regular Polygons*, *arXiv:1311.6763* are a salient feature of the large scale evolution of the web W. In SV, the authors showed that if the rings are linked together in unbroken chains like 'necklaces' they will form dynamical barriers. Since $\tau$ is an area preserving mapping points on either side of a necklace cannot cross over. Therefore since the rings are periodic, all orbits are bounded. We believe that it is sufficient to study the dynamics inside the first ring of D tiles, which we call Ring 0. The chiral nature of N and the Reflection Principle still apply here so the regions between rings will support the $I_k$ structure. But the emphasis here is the global evolution.

We will first describe the evolution of +W and –W on the horizontal extended 'base' edge of a regular N-gon when the edges of N are not truncated. The Reflection Principle explains the nested rings above  since each D tile will generate two-sided families and guarantee periodic sequences of D tiles (and matching star points) on the horizontal edge. We will use these as 'seed' tiles to describe how the rings evolve.

-

In Figure 1 below we describe this evolution for the twice-odd case of N = 22 – which includes the N = 11 case. The twice-even case is very similar.

**Figure 1**  The first 3 iterations of the global web for N = 22 which includes N = 11 as the matching S[N/2-2] M tile. Here we combine cw and ccw evolution. On each iteration evenly spaced points on the extended edges of N are iterated under $\tau$.  The black forward edges show how the S[k] tiles evolve in the First Family Theorem (FFT) and the magenta trailing edges show the contribution of –W. Together they explain how new D tiles  are formed at periodic intervals .

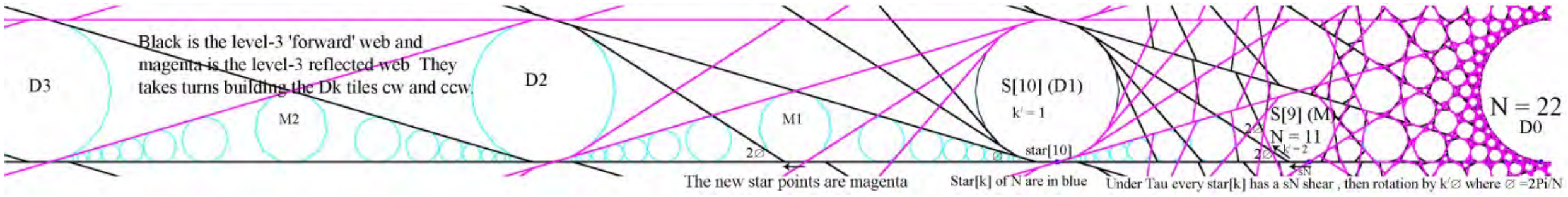

As in Example 4.3, for the black forward edges all the star[k] points on the horizontal axis L will evolve with an sN outward shear followed by a $k\phi'$ rotation  This horizontal strip is part of a domain of $\tau^{-1}$ based on star[1] of N, so for points on this strip $\tau^{-1}(p) = 2v_1 - p$ where $v_1$ = star[1] of N. Based on the magenta trailing edges of N, there will be N/2-1 star[k] hence N/2-1 sub-domains of $\tau$.

Now suppose a point p is circling N cw the right-side. As p approaches L it will 'see' $v_l$ so $\tau(p) = 2v_l - p$ will approach $2v_l$ which is sN. Therefore every point on L will experience a (one-sided) leftward shear of sN and in particular the blue star[9] of N will be extended outwards by sN, but now the discontinuity at star[9] implies a new reference base L′ which will be rotated by $2\phi$ from L, so this sN displacement forms part of next edge of S[9]. This will be repeated recursively on each iteration and the first 3 iterations are shown here. D at step-1 will replicate N so it will be a reflected version of N capable of generating its First Family on either side.

Meanwhile on the right side of N the trailing edges are evolving in a ccw fashion and we show in magenta what happens when this evolution is reflected to the left side of N As expected ,relative to these imported rays, D is evolving from the top down and this is perfect for generating the next D2. These trailing edges of D generate the new magenta star[k] points which will interact with L to generate the new M tile using the same sN offset and D2 will evolve in the same two-sided fashion as D to guarantee endless periodic seed D tiles on the horizontal axis.

**Figure 2** Below is the resulting geometry. Since $\Phi(22)/2 = \Phi(11)/2 = 5$, there will be 4 distinct regions on either side of the central M tiles. For the web W it does not matter whether N =22 or N = 11 is at the origin, but dynamically it does matter and here the rings are centered on N = 11.

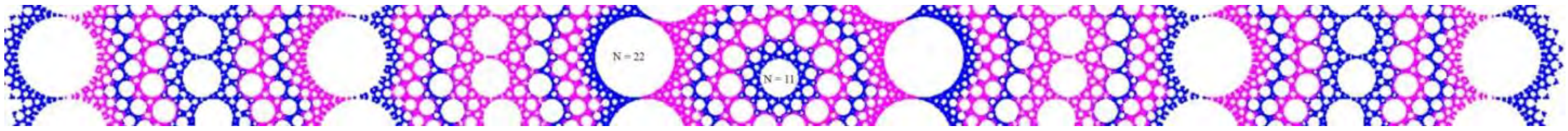

Based on the web evolution shown here, we will make the following assumptions about the 'seed' D[k] on the horizontal axis: *For N even the spacing of the $D_k$ tiles will be the same as N to D but for N odd it will be twice this*. Once the seed D tiles are known, iterating the centers under $\tau$ will determine their period and step-sequence. All the tiles in the orbit of the seed D will share the same step-sequence so the centers of the D tiles in the ring can be generated by iterating cD under $\tau^k$ where k is the sum of the step sequence.

**Case 1**; The Rings of N = 8 (with the first ring called Ring 0 with seed D=$D_0$)
Since $C_8 = 2$, there is just one invariant region right of M = S[2] and it is the black region shown here. $D_0$ is S[3] so it is step-3 and gcd(3,8) = 1 so the orbit of $cD_0$ will be period 8 as it circles N with winding number (mean steps on each iteration) 3/8 . It is easy to verify that $cD_1$ will have period 16 and step-sequence {4,3} as in the table, so Ring 1 must be constructed step 7 and the only way to do this is to first make a table of length 16*7 and then plot them mod 7. This is necessary because these 16 D tiles have to fit together in exactly the same way they are generated. For Ring 2 set Orbit = V[$cD_2$,24*11], then the Ring 2 tiles will have centers given byTable[Orbit[[k]],{k,1,254,11}.This guarantees that the 24 D tiles in this ring will share a vertex with their neighbors and also share a common S[1] dual tile because they are indeed 'step-11' neighbors of each other in the web.

| Ring | Step Sequence | Period of center | Winding |
|---|---|---|---|
| 0 | {3} | 8 | 3/8 |
| 1 | {4,3} | 16 | 7/16 |
| 2 | {4,4,3} | 24 | 11/24 |
| k | {3} + k*{4} | 8(k+1) | (3 + 4k)/8(k+1) |

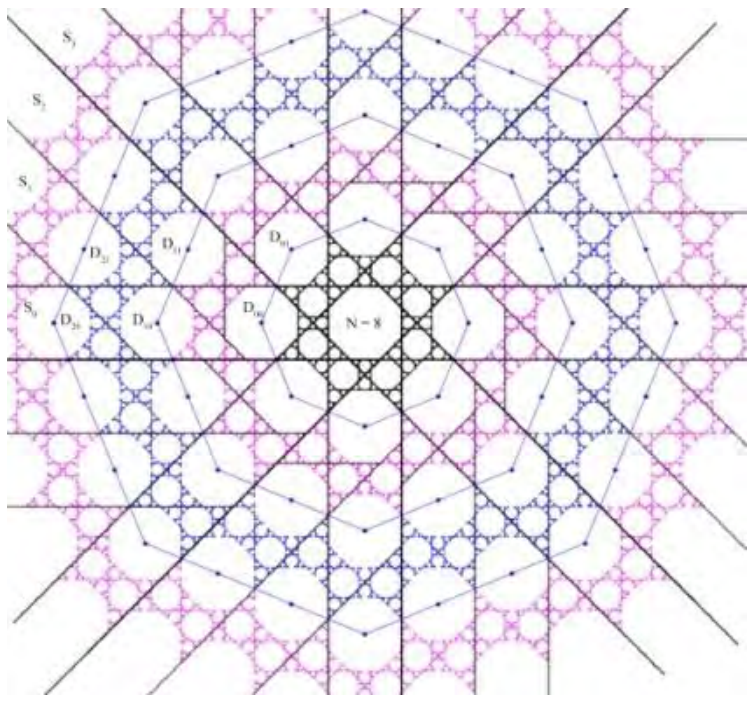

**Case 2** The rings of N = 14

There is always a strong connection between the dynamics of these rings and the overall dynamics of N. For N twice-odd the mutation condition of an S[k] is gcd(k′,N) > 1 so all the odd S]k] will experience a gender-change when they evolve in W. In the same fashion the orbits of the D tiles will decompose into 2 cycles when the step-sequence has a common factor with N. For example the seed $D_2$ in Ring 2 has 6+7+7 steps and gcd(20,14) = 2 so this ring will have 2 disjoint cycles of period 21 as shown by the blue and magenta 21-gons. The second cycle will be generated by a rotational neighbor of the seed D tile. This means that Ring 2 will be the 'Riffle' or weave of these two and hence a (regular) 42-gon. Ring 0 with step- sequence {6} is the Riffle of blue and magenta heptagons so it is a regular N-gon. The odd orbits do not decompose so they are easy to generate but it is necessary to know their step sequences to generate the matching 'invariant polygons' joining their centers. For example to generate Ring 3, over-scan with K = V[cD3, 56*27]; Then Ring 3 polygon = Table[K[[k]], {k, 1, 56*27, 27]. This guarantees that each D tile in Ring 3 will share a vertex with its step-27 neighbors.

| Step sequences and periods for D"s for N = 14 | | | |
|---|---|---|---|
| Ring | Step Sequence | Periods of centers | Winding Number - |
| 0 | {6} | 7 & 7 | 6/14 |
| 1 | {6,7} | 28 | 13/28 |
| 2 | {6,7,7} | 21 & 21 | 20/42 |
| 3 | {6,7,7,7} | 56 | 27/56 |
| k even | {6} + k*{7} | 14(k+1)/2 & 14(k+1)/2 | |
| k odd | {6} + k*{7} | 14(k+1) | |

Invariance can exist at all scales and in [VS] the authors discovered that some non-regular N-gons will have major rings or secondary rings. When N is twice-odd the S[N/2-2] „M" tiles have k′ = 2 and they will be surrogate N/2 tiles. For N = 14 this M is S[5] and we noted in the introduction that S[5] will combine together with its dual S[2], so form an unbroken „necklace" when rotated about N. S[4] and S[3] will also be dual and create another necklace, so there will be 2 invariant region right of M as shown in the middle below. Thd Edge Sharing Lemma to follow will show that any pair of adjacent dual tiles will share an edge.

| The first 3 rings of D tiles | Ring 0 has 2 secondary rings on each side of M. | In the Introduction we showed that for N = 14 there will be 3 pairs of dual tiles right of M, namely S[1]&N, S[3]&S[14] and S[5]&S[2]. Each of these has a matching star point which defines a separatrix $I_k$ between regions, so the 3 $I_k$ determine 2 invariant regions right of M and there are 3 symmetric $I_k$ left of M that determine equivalent regions. |
|---|---|---|
| 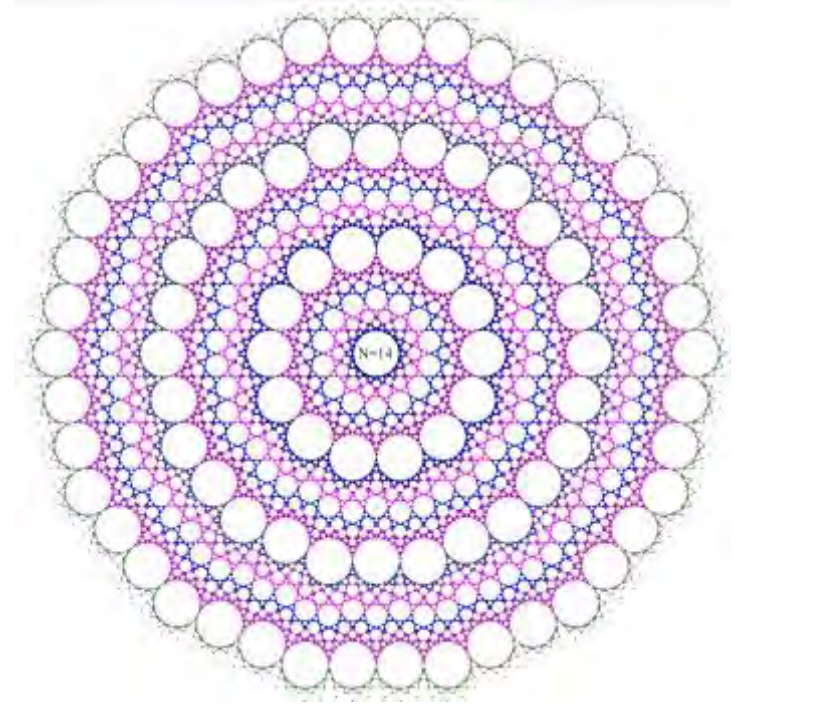 | 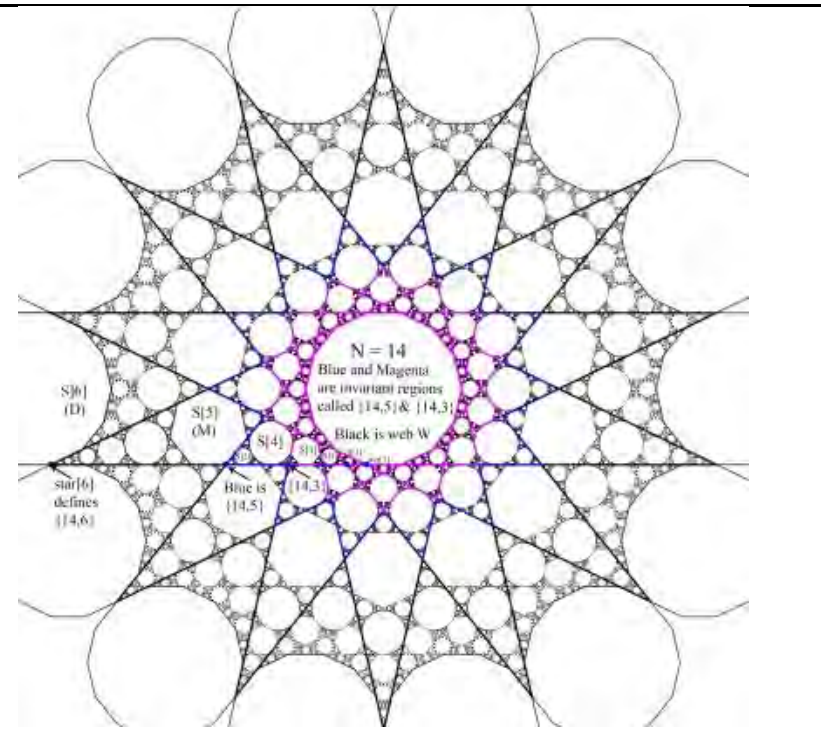 | |

**Case 3** The rings of N = 7

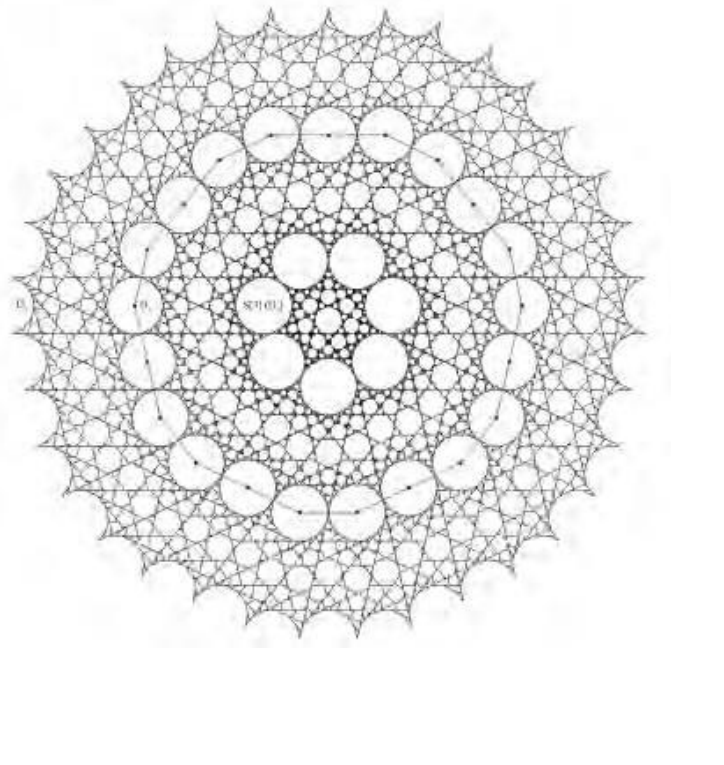

| Step sequences and periods for D's for N = 7 | | | |
|---|---|---|---|
| Ring | Step Sequence | Period of center | Winding Number |
| 0 | {3} | 7 | 3/7 |
| 1 | {3,3,4} | 21 | 10/21 |
| 2 | {3,3,4,3,4} | 35 | 17/35 |
| k | {3} + k*{3,4} | 7(2k+1) | 3(k+1 +4k) |

Here HalfN = 3 just like N = 8 so once again D has step sequence {3}, but now the ring growth is a little slower at {3,3,4} for Ring 1. There are no step-4 orbits inside Ring 0 and points that approach D from inside appear to form a Devils Staircase of winding numbers approaching {4}.

When N is odd the rings past Ring 0 are irregular because they are formed from the Riffle (weave) of odd and even cycles of the return map. For Ring 1 the cycle periods are 11 (odd) and 10 (even) so if K = V[cD1,21] is the $\tau$-orbit of cD1, then the odd and even cycles are M1 =Table [K[[k]], {k, 1, 21, 2}] and M2 = Table [K[[k]], {k, 2, 21, 2}]; Ring 1 polygon = Riffle[M1(blue), M2(magenta)] as shown on the right below. This explains why these resulting polygons are 2N-gons and why they have small winding numbers. Ring 0 is the regular Riffle of a non-regular triangle and an (odd-cycle) quadrilateral. This is reminiscent of the 'projections' from [H3].

Ring 1 can be used as a template for the remaining rings where each side increases by 1 diameter of D, so asymptotically they will be regular 2N-gons. Below are the first 4 rings where for display the polygons are rotated by one step relative to the web.

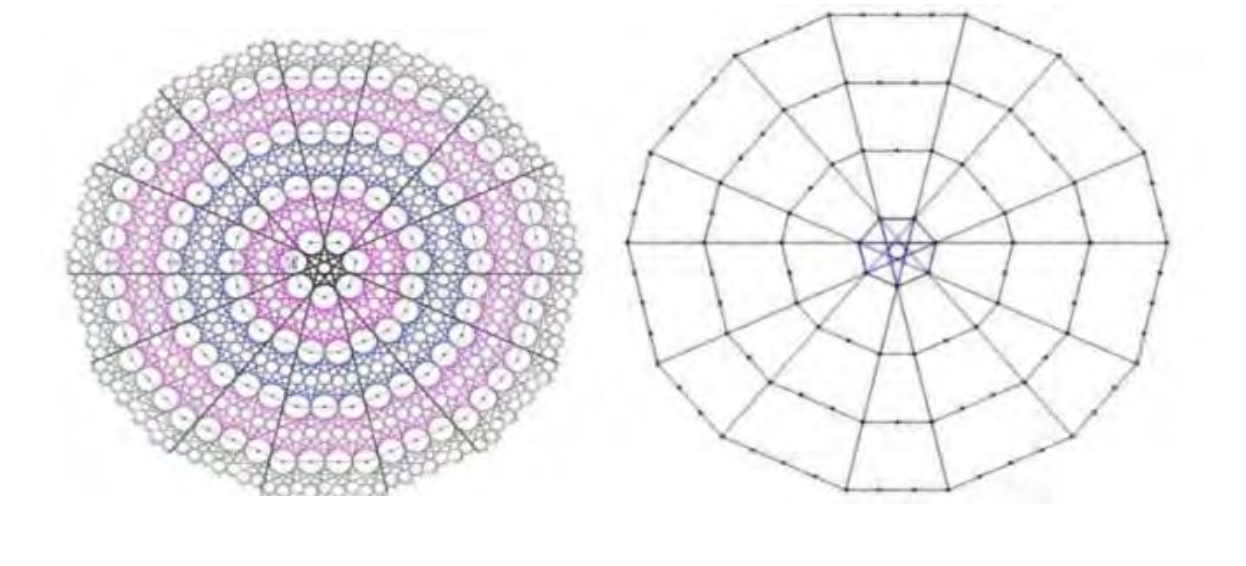

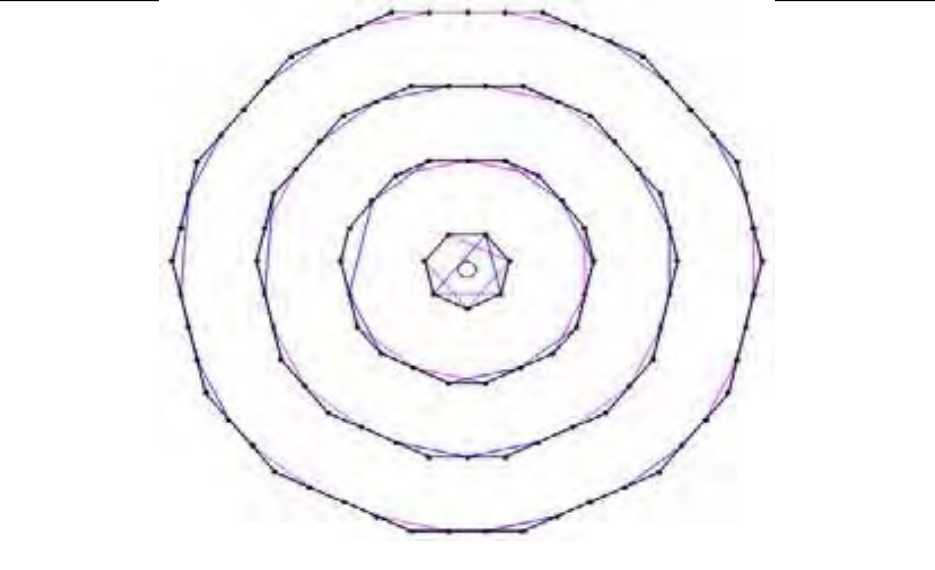

In [H4] we analyzed some of these 'semi-regular' N-gons to see how they would behave at the origin. As expected many were capable of generating rings of D tiles which showed they were bounded and 'stable', but their measure in the class of convex polygons was still zero so they did not disprove the R.S. conjecture that 'almost all' convex N-gons will be unstable. This remains an unresolved issue as there is certainly no KAM theory for regular N-gons.

**Secondary Invariant Regions of a Regular N-gon with Examples for N ≤ 50**

In the hierarchy of the singularity set W for a regular N-gon, the infinite chains of concentric rings of D tiles can be regarded as the primary resonances where there is no truncation of the extended edges of N. Since the regions between rings will be invariant each region would be expected to contain ‘secondary’ resonances in the form of invariant ‘tori’. These secondary regions should themselves decompose into invariant regions in an endless chain. We claim that it should be sufficient to restrict our investigations to the region inside the first ring that we call Ring 0, and all the dynamics outside of this ring will be congruent to dynamics inside.

In the Introduction we noted that the chiral nature of N implies that the default cw ‘forward’ web of N will be distinct from the ccw web so we call these W (or +W) and –W. We also claimed that the boundary between W and –W includes the $I_k$ closed curves which act as separatrices between invariant regions and we claimed that these $I_k$ curves are based on adjacent pairs of edge-sharing S[k]-S[k′] tiles where for N even k′ = N/2-k These are called ‘dual’ tiles and if at least one of the indices k or k′ is ‘primitive’ with gcd(k,N) = 1 then they will survive in W and form a barrier between invariant regions in the form of closed polygonal curves which we call $I_k$.

We illustrated this with N = 42 and that will be continued here. The Invariance Conjecture to follow will claim that for all regular N-gons, the number of distinct invariant regions will be $C_N$ - 1 where $C_N$ is the algebraic complexity of N, namely Φ(N)/2. For N = 42, $C_N$ = 6 so there should be 5 invariant regions to the right of M = S[N/2-2] acting as N = 21 as shown here.

**Figure 4** The 5 invariant regions right of M and matching $I_k$ curves for N = 42 with $C_N$ = 6

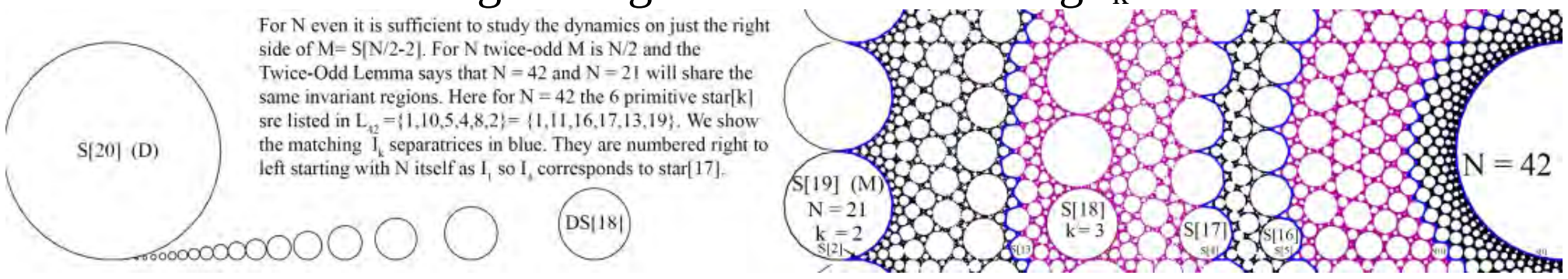

In the Introduction we claimed that these dual pairs will arise for two reasons:
(i) The chiral nature of N and Reflection principle implie that the combined web must evolve both cw and ccw as shown here.

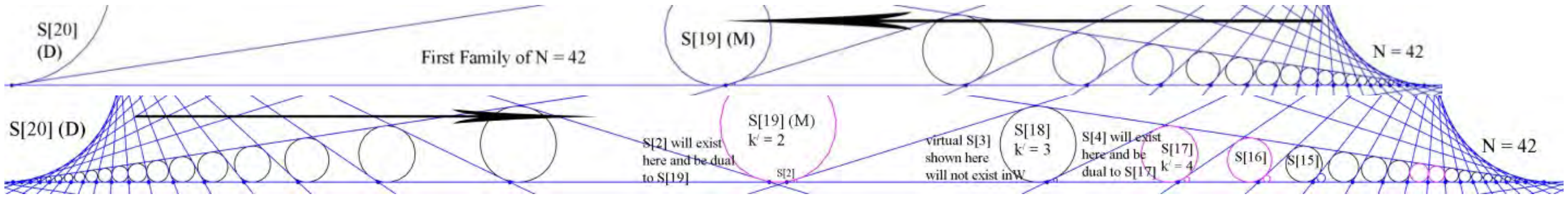

(ii). The Generalized First Family Theorem (GFFT) says that every First Family S[k] for k < N/2-1tile has the potential to construct S[k+1] on the left and the dual S[k′] on the right

So S[17] can generate S[18] on the left because every S[k] shares a star[k+1] of N on the left. On the right side every S[k] shares star[k] of N, but in order to guarantee that the S[k] are conforming to the star polygons, it was necessary that locally star[k′] of S[k] matches star[k] of N. So star[17] of N will form S[17] on the left with offset {sN/2,0} but (if space permits) the matching offset on the right will form an S[4] clone as part of the local family of S[17] . By the GFFT, MidpointS4 = star[4] + {sN/2,0} and cS4 = MidS4 + {0,hS[4]} where hS[4] is known from the FFT. Then Mathematica can make a copy of S[4] at the new location.

Theoretically this same method would generate an S[3] clone for N = 18 but as mentioned in the Introduction, both S[18] and S[3] are mutated and this implies that S[18] will in fact generate an S[2] that does not share an edge with S[18]. Here we look more carefully at S[17]

**Figure 5**

| S[17] and its dual S4. | A schematic showing synchronization of edges |
| --- | --- |
| 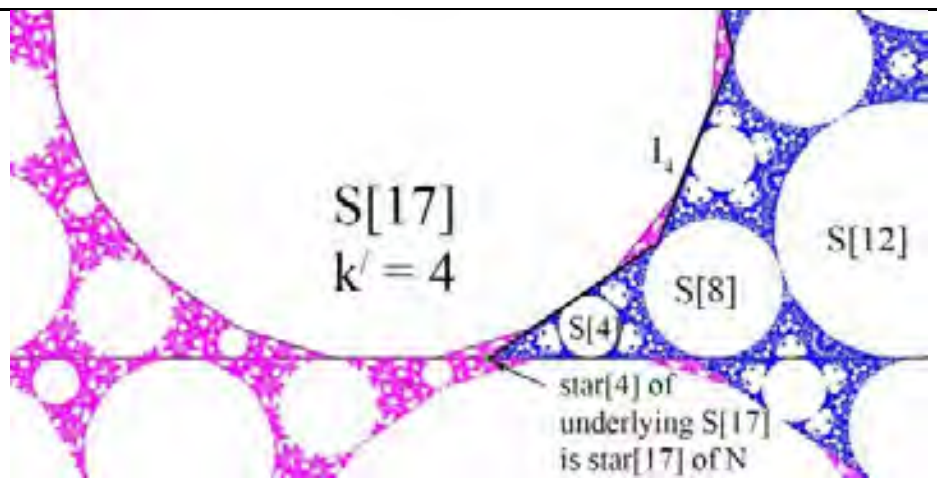 | 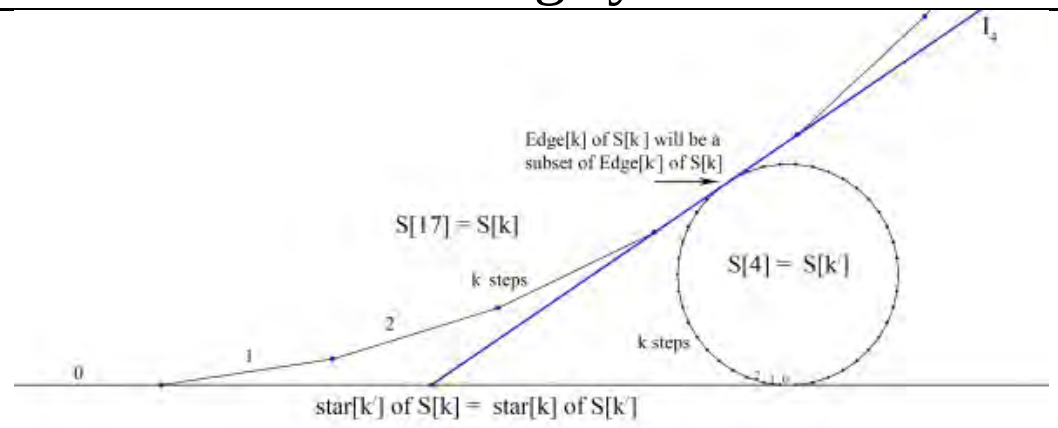 |

As noted above S[17] will have a dual relationship with a clone of S[4]. Therefore star[4] of S[17] is also star[17] of S[4]. This means they share the same blue tangent line and we can define the 'twist' of this line as 4 ccw steps relative to S[17] or 17 cw steps relative to S4 so edge E[17] of S[4] must be a subset of E[4] of S[17]. Clearly $4\pi/21$ and $-17\pi21$ and just two ways to define the same twist but our convention in Figure 22 will be to use the twist relative to the larger of S[k] and S[k′], so we will use $4\pi/21$. Note that S[17] will be mutated into an N/2-gon but that has no effect on the 'synchronization' with S[4] because the star points of the mutated S[17] are an even subset of the star points of the underlying US17 N-gon shown here. So the 4 steps of US[17] here are just 2 steps of the mutated S[k] – but we know this correspondence.

**Figure 6** For N even the S[k] tiles are born as N-gons, but when N is twice-odd, the S[k] (including N/2) must be compatible with N so the 'odd' S[k] will be mutated. But as noted above this 'gender change' is never an issue with analysis, because the center and height are unchanged and the star points of the N/2-case are an even subset of the original star points. Here we show S[1] in its mutated form and it will have no effect on the synchronization or edge sharing with N

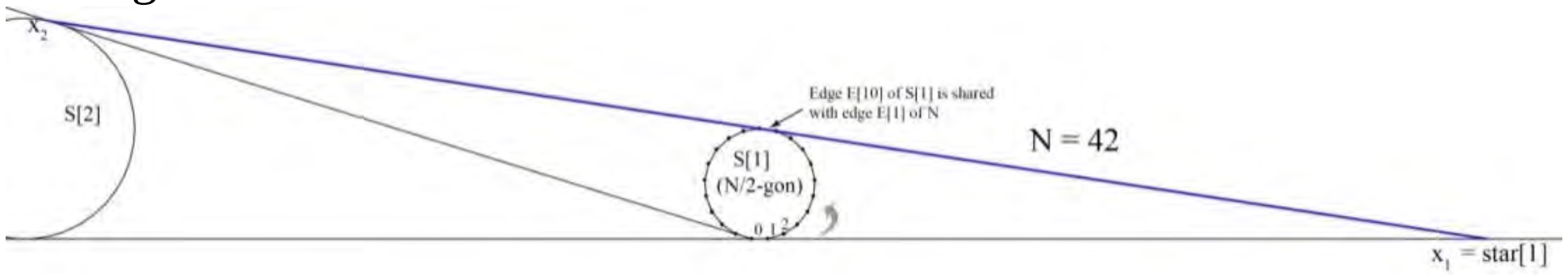

As with all dual pairs there are two opposing twists - one around N and the other around S[1]. From the perspective of N, the blue shared edge is just one cw step so it is $-\pi/21$. From the perspective of S[1] the blue edge is 10 ccw steps, but these are double steps because S[1] is an N/2-gon, so the twist is $20\pi/21$ and of course it is the 'same' twist as N. Once again we will defer to the larger tile and use $-\pi/21$. In this case it is easy to complete the $I_1$ seraratrix with one rotation around N. Since this rotation is negative, the product of twist and length is just $2\pi$ and this will be the $\omega$ value for N. In exactly solvable systems such as a Hamiltonian with a simple twist, it is common for 'solitons' or stationary waves to arise. They owe their stability to the integrals of motion which for N = 42 are the $I_k$. The 6 dual pairs act as stable solitons and indeed each pair defines a primitive star[k] with a well-defined 'twist' so there is a Hamiltonian form $H_k$ for each dual pair whose solution is the matching $I_k$. Each primitive star[k] acts a'stacked' boson and the height of the stacks are the eigenstates which for N = 42 are 1,2,4,5,8,10 . See Fig 22.

**The Synchronized Edge-Sharing Lemma**: For N even if dual S[k] and S[k′] are adjacent then we will show that they must have equal and opposite 'twists' and hence share an edge. There is no loss of generality in assuming that the larger tile S[k] is on the left as shown here - otherwise reflect about M for an equivalent pair. We illustrate here with the 'generic' case of N twice-even where both S[k] and S[k′] are regular N-gons and then show how to generalize this to N twice-odd where one of these two will be mutated into an N/2-gon.

(i) When neighboring S[k]-S[k′] are dual they must share the same blue 'tangent' line as shown here because star[k′] of S[k] = star[k] of S[k]. By definition relative to edge 0 on the horizontal axis, k′ cw steps of S[k] will align with the blue tangent line and on the right side, k ccw steps of S[k′] will also align with the tangent line so Edge[k] of S[k′] must be a subset of Edge[k′] of S[k].

**Figure 7** Showing dual S[k} and S[k′] as N-gons (here S[k] is S[7] of N = 20)

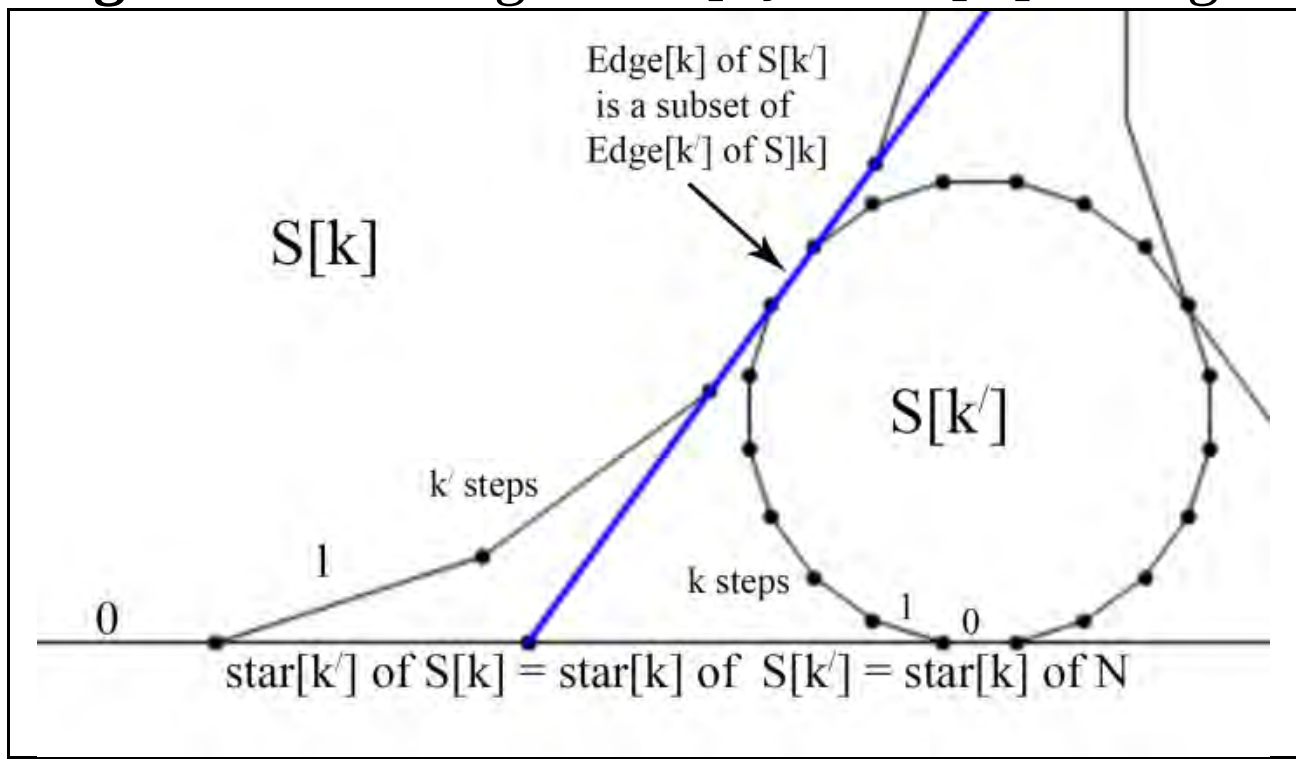

(ii) Now suppose N is twice-odd and the S[k] shown here is the underlying US[k] of the mutated S[k]. Then this diagram is still valid but the k′ steps of US[k] are now k′/2 of the doubled steps of S[k] so still synchronized.

(iii) If instead S[k′] is an N/2-gon and the S[k′] shown here is US[k′]. Then it will take 2k steps to reach the blue line but these are half-steps relative to S[k′] so still synchronized.

The Invariance Conjecture below will state that in W at least one of k or k′ will be primitive, but that assumption was not needed here. If N is twice-even the mutation condition for an S[k] is gcd(k′,N) >2 but when k is primitive k′ is also primitive and relatively prime with N so neither S[k] or S[k] can be mutated and for N twice-odd it is not possible for both to be mutated.

In all cases the geometry can be reduced to the S[1]-N case by regarding the smaller S[k′] as an S[1] relative to S[k]. By the FFT, hS[k]·hS[k′] = $(s_1/s_k)(s_1/s_{k'})$ but when both are N-gons. $s_{k'} = 1/s_k$ so hS[k]·hS[k′] = $s_1^2$ = hS[1], so if the heights of S[k] and S[k′] are scaled by hS[k]. then S[k] will have height 1 relative to S[k′] serving as S[1]. This will remain true for N twice-odd but now the scaling will have to include scale[2] of N which is the gender-change scaling.

The twice-odd family is the only one where neighboring S[k] can share an edge and we show here why the S[Floor[N/4]] tile will always be dual to its neighbor. Set k = ⌊N/4⌋ (Floor of N/4) and ask what tile is S[k] ? ⌊N/4⌋ = ⌊M/2⌋ where M = N/2 is odd so ⌊M/2⌋ = ⌊(2j+1)/2⌋ = j so k is simply the k in M = 2k +1. Then k′ = N/2-k = M-k = (2k+1) –k = k+1, so for N = 42, S[10] and S[11] are a dual and in general *for N twice-odd S[N/4⌋ and S[N/4⌋ +1* will be dual and hence share an edge. This S[N/4⌋ tile of N is what we call the 'shepherd' of the initial invariant region. It will be the S[k′] of its larger neighbor acting as S[k]. Since S[k] will be an N/2-gon the edge fit between them will not be as tight as it would be otherwise. See Fig. 8 below. (Note that S[N/4] for N twice-even is 'self-dual' and mutated but not a boundary tile.)

**Figure 8** The initial invariant region for N = 42 with 'shepherd' S[10] = S[Floor[N/4]]

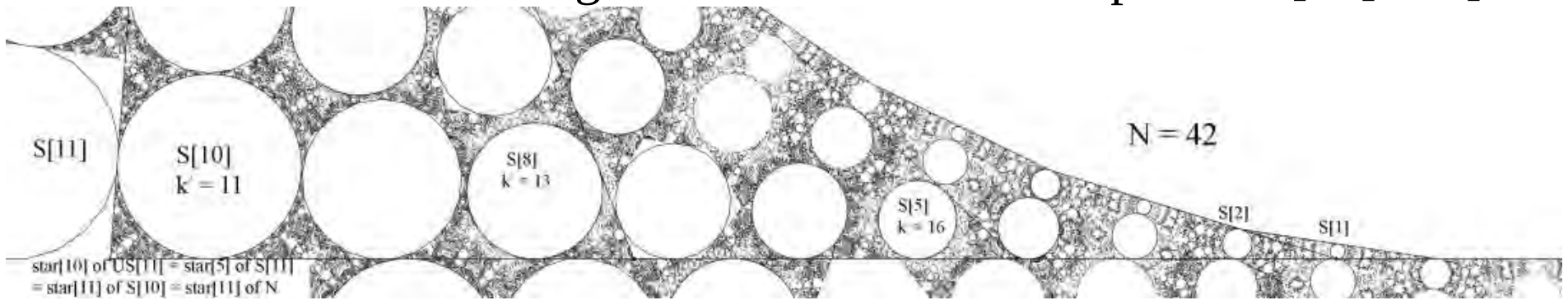

Our convention is for Mathematica to display S[11] as an N/2-gon with US[11] as the underlying N-gon, so duality implies star[10] of US[11] is star[11] of S[10] and therefore star[5] of S[11] shown here will match star[11] of S[10] and Edge[5] of S[11] will contain Edge[11] of S[10].

**Synchronization of the S[k] tiles with N**

Below are un-cropped versions of Figure 7 . On the left we show the $\tau$-orbit of S[k] and demonstrate how it can be defined by pure rotations without need for $\tau$. This is true because every S[k] First Family tile is 'synchronized' with N since they share star[k] of N. The tangent line linking S[k] and N implies that edge[k′] of S[k] will align with edge[k] of N – where we show the matching S[1] tile. Under rotation about N this edge will be step-k so it will have period N/(gcd(k,N) just like cS[k] and this provides an alternate proof of Lemma 4.1.

To track this orbit of S[k] we show the magenta star-polygon step-k orbit of star[k] which is defined to be Orbit = Table[RotationTransform[2kPi/N, {0,0}]][star[k]], {k,0, Nk,k}. The orbit of S[k] will be 'entrained' with this orbit as shown here by the first 3 iterations. By definition it is a valid $\tau$- orbit of cS[k] but the synchronization of N and S[k] means the orbit can be regarded as a simple rotation of edge[k] of N.

**Figure 7b**

| The orbit of S[k] as defined by N and not $\tau$ | When S[k] has a dual , edge[k] will align with edge[k] |
|---|---|
| 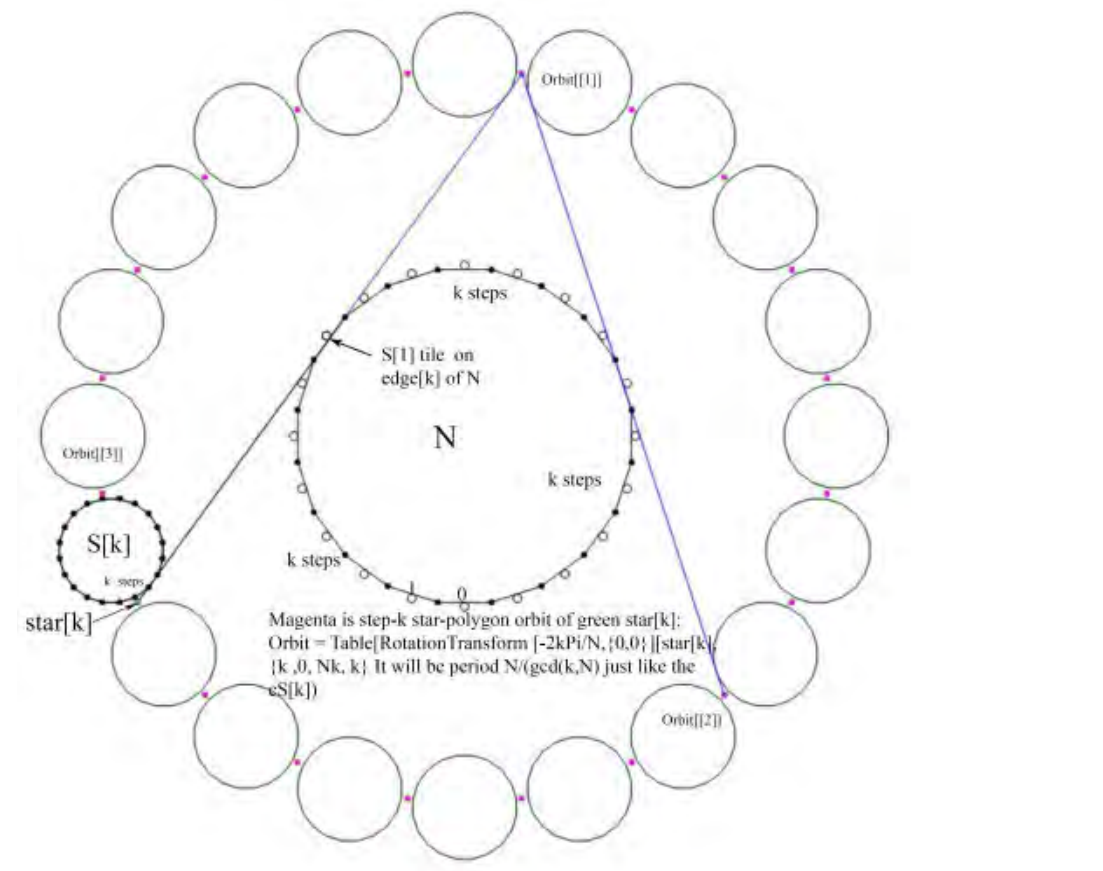 | 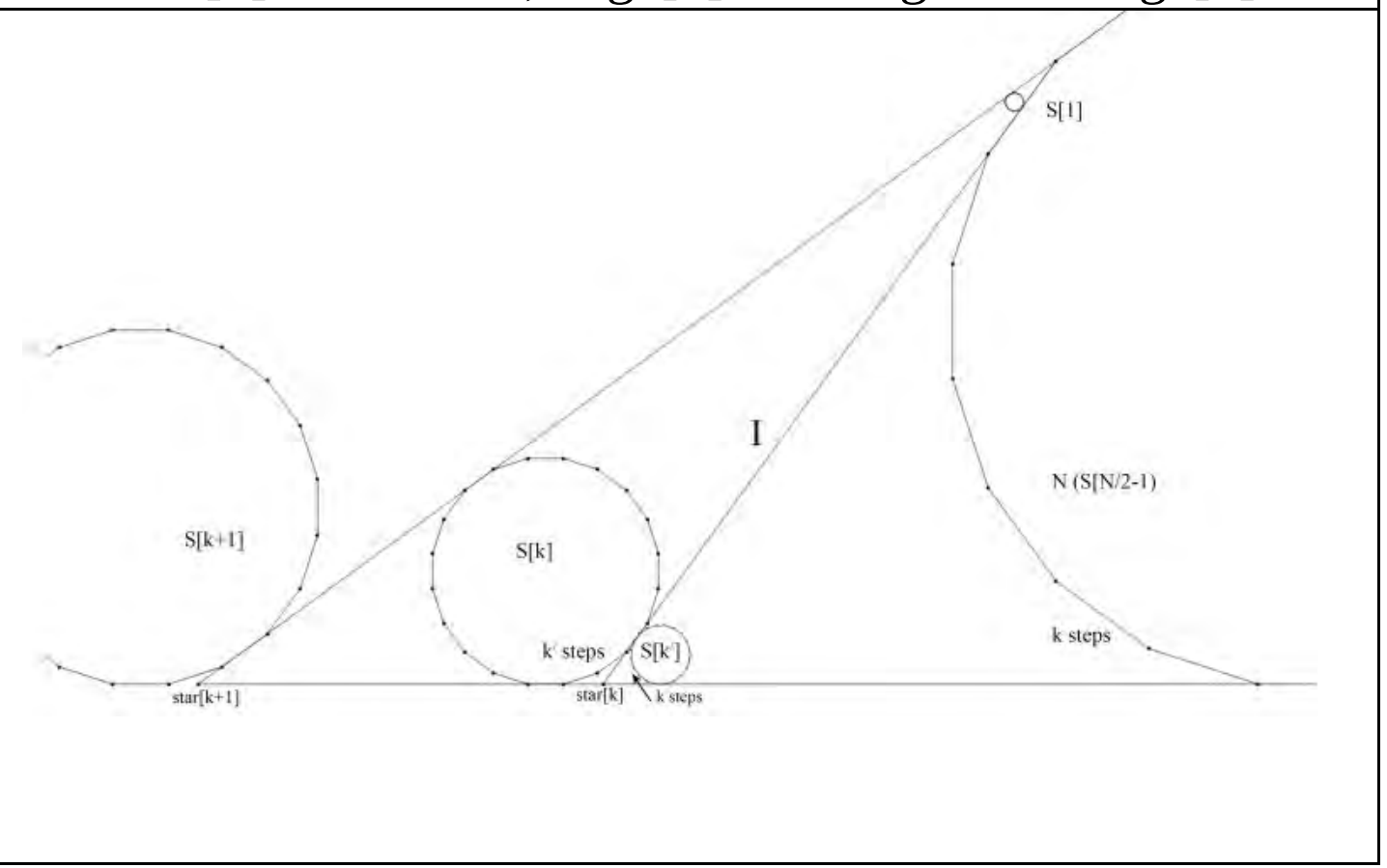 |

When S[k] has a dual S[k′] as on the right , the same linkage exists but now edge[k] of S[k′] is aligned with edge[k] of N. The the ratio of their heights is hS[k′] = 1/S[k] and if S[k] is scaled by 1/hS[k] to act as N the new hS[k′] will be 1/hS[k]$^2$ acting as a local S[1]. (It looks 'chubby' compared with S[1] but this is because the S[k] family is step[k′]. (Here N = 20 and k = 7.)

Since the S[k] dynamics are independent of $\tau$, the linkage between the S[k] and smaller S[k] clones is also independent, so S[k′] here is a satellite of S[k] not N. Since no S[k] are clones of N (except D), their families are only distantly related and each generation brings them further apart.

**The Invariance Conjecture**

Note: The examples of invariant regions to follow includes both even and odd N-gons, but the Twice-Odd Lemma implies that there is no loss of generality in regarding N/2 as embedded in N so we will concentrate here on the N-even case. Because N/2 and N share the same cyclotomic field they will share the 'same' invariant regions, but these regions will be reflections of each other as N/2 and N are swapped, so for example $L_{26}$ = {1,6,4,3,5,2} while $L_{13}$ = {2,5,3,4,6,1}.

For a regular N-gon with N even the central S[N/2-2] is known as M. This is the only tile that is in the First Family of both N and D and the Reflection Principle says that the geometry of N and D are equivalent under reflection about cM, so it is sufficient to study the geometry right of M.

For N even in this region right of M there will be N/2-1 star[k] points but just $C_N = \Phi(N)/2$ primitive star[k] points with gcd(k,N) = 1 and this is also the rank of the of the maximal real subfield $\mathbb{Q}_N^+$ which defines the algebra and geometry of N. We conjecture that there will be a 1-1 correspondence between the primitive star[k] of N and the 'fundamental resonances' of N because each primitive star[k] point can be regarded as an eigenvector of N with the 'twist' of star[k] as the matching eigenvalue (in radian measure).

This twist will match that of an S[k] - S[k′] dual pair whose shared edge defines star[k]. As shown in Figure 5 above every dual pair shares a common star point so the rotational steps of S[k] and S[k′] are synchronized and this defines the twist of star[k]. Every one of the $C_N$ surviving dual pairs is a 'soliton-like' resonance between the +W and –W webs where the S[k] and S[k′] lie on opposite sides and the shared edge defines the boundary. To escape mutations at least one of the indices k or k′ will have to be primitive.

Using the lists $L_N$ of primitive indices it is not difficult to construct the $I_k$ boundary curves as shown in Figure 4. In classical systems such curves were called 'integrals of motion' and they were signs of stability. Because of rotational symmetry this shared edge can first be rotated around local S[k] or S[k′] and then around N to form the $I_k$ dynamical barriers as in Fig 23-27.

Each $I_k$ will be a separatrix between invariant regions and $I_1$ is always based on the dual pair S[1] and N so it is the edges of N as the primary resonance. When N is twice-odd just one of k or k′ will be a primitive index, but when N is twice-even and k is primitive, k′ will also be primitive so the number of $I_k$ are reduced by half but now the right and left sides of M are algebraically distinct so the total number of distinct regions is still $C_N$ -1. (M here is no longer a proxy N/2.)

In the Introduction we note that this large-scale structure is inherit in N. Because of its chiral nature the phase space of every regular N-gon will consists of two overlapping but disjoint manifolds +W and –W bounded by the star polygons of N and in the interior of the phase space, they will meet at the $I_k$ formed from the extended shared edge of an S[k] tile of ±W and its dual S[k′] of the opposite polarity. The Jacobian of N can distinguish ±W because it is dependent on the direction of rotation chosen but of course J or –J will yield equivalent results.

In the spirit of Emmy Noether (1882-1935) *"every symmetry should have a matching invariance"*. It seems that along with the rotational symmetry of W, the cw-ccw lateral symmetry is behind these invariant regions. That symmetry is also inherent in the tan function that of C.S. Siegal found so interesting. Historian John Stillwell [St] said "*the tan function may be considered more fundamental than either cos or sin*". The critical Twice-Odd Lemma would not exist without the half-angle tan relationships and they also ensure that $\tan^2[\pi/N]$ or $\tan[\pi/N]\cdot\tan[2\pi/N]$ can be alternate generators of the scaling field $\mathbb{Q}_N^+$ in place of $2\cos[2\pi/N]$.

The twice-even case has no 'gender-change' mutations since the mutation condition is now gcd(k,N) >2. The extra symmetry of 4|N guarantees that $s_1^2 = hS[1]$ will be an algebraic unit and a good choice as generator of $\mathbb{Q}_N^+$ since the traditional generator $\lambda_N=2\text{Cos}(2\pi/N)$ may not even be an algebraic integer. But for N odd or twice-odd there is a gender-mismatch between S[1] and N so $s_1^2$ may fail to be an algebraic integer as witnessed by $\text{Tan}^2[\pi/6] = 1/3$. A simple solution is to use $hS[1]/hS[2] = s_1s_2 = \text{Tan}[\pi/N]\cdot\text{Tan}[2\pi/N]$ where $\text{Tan}[\pi/6]\cdot\text{Tan}[\pi/3] = 1$ as it should be.

When N is twice-even, the dual pairs are still dominant, but the ratio of primitive indices is doubled because all odd k are primitive and if k is odd, k′ is also odd. Therefore all primitive dual pairs have both indices primitive. This means that the regions right and left of M are now algebraically distinct so the total number of distinct invariant regions will again be $C_N$ -1, but now it applies to the combined right and left sides of M.

**Figure 9** For N = 20 with $C_N$= 4, the primitive indices in $Z_{10}$ are {1,3,7,9} but these are redundant so the short list is just $L_{20}$ = {1,3} which is equivalent to {9,7}. M is S[8] but it has no dual so there are just 3 regions and the four $I_k$ are $\pm I_1$ and $\pm I_2$ with sign change at cM.

| The first 3 rings | Ring 0 | The right-side of M | The S[1]-S[2] region [H6] |
|---|---|---|---|
| 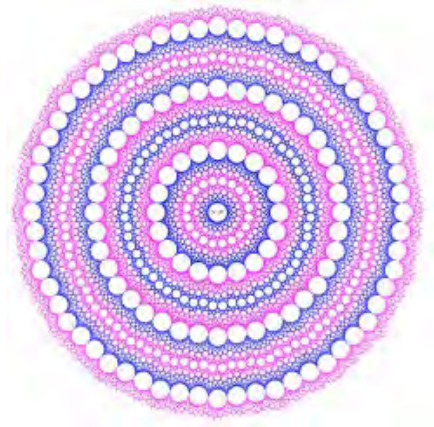 | 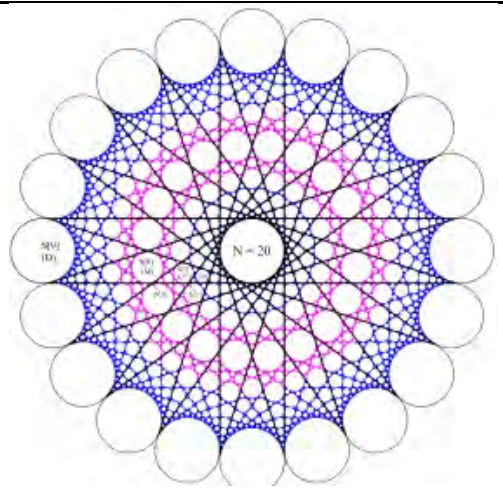 | 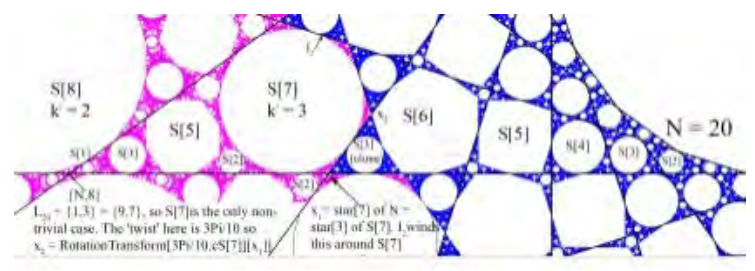 | 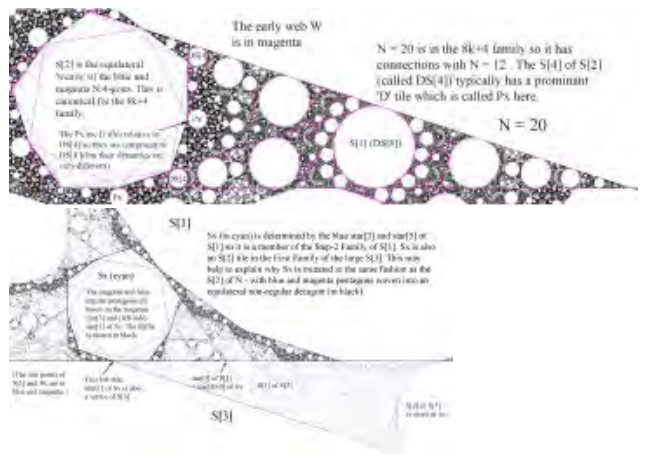 |

The extra symmetry of N|4, dictates that the mutation condition for S[k] is gcd(k′,N) > 2 but here S[6], S[5] and even S[2[] are mutated. The insert shows detail of S[2] and partially explains the connection with N = 12. Here the two distinct $I_k$ have 'stack' height 1 and 3 and these are the only eigen-values so a crystal for N = 20 will not be very friendly toward electrons.

To explain the algebraic difference between N twice-odd and N twice-even recall that our default convention for N was height (apothem) of 1 so the horizontal displacement of the star points are $s_k$ = Tan[kPi/N] = Tan[2kPi/2N]. Therefore the star[k] points divide the horizontal axis into 2N symmetric segments where the primitive $s_k$ are Galois conjugates of the minimal polynomial $T_N$ of Tan[Pi/N], but when 4|N there is extra ± symmetry in the Galois group of the extension and this reduces the algebraic complexity of Tan[2Pi/N] from $\Phi(N)$ to $\Phi(N)/2$ and now $T_N$ is odd so it is necessary to use both $T_N[x]$ and $T_N[-x]$ to obtain all the primitive $s_k$ and $-s_k$ can be regarded as algebraically distinct to $s_k$. See Figure 10 below.

**Figure 10** The primitive star points of S[8] are star[1[ and star[3] in black. Note that ±star[1] and ±star[3] are algebraically distinct but retain the geometric symmetry of $T_8[x]$\ and $T_8[-x]$

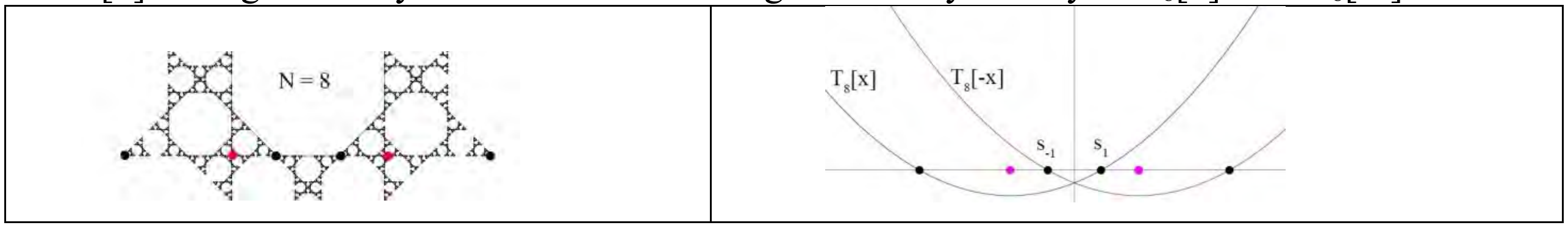

**Figure 11** N = 60 has $C_N = 8$ so there are 7 invariant regions on either side of N and $|L_{60}|$ will be just 4 so there will be 4 regions right (and left) of M as shown in the enlargement. These regions look identical but now each S[k]-S[k′] pair has both primitive so algebraically they are distinct. To find the primitive indices for N = 60, begin with $Z_{30}$ but this is still even so it is sufficient to use $Z_{15}$ which yields {1,7,11,13} so the {k, k′} dual pairs are {1 (29),7 (23),11 (19), 13 (17)} and by inspection the order is $L_{60}$ = {1,13,11,7}. Note that M has no dual S[2] so upon rotation about N it counts as just one region. This is why there are just 4 regions right of M (including M).

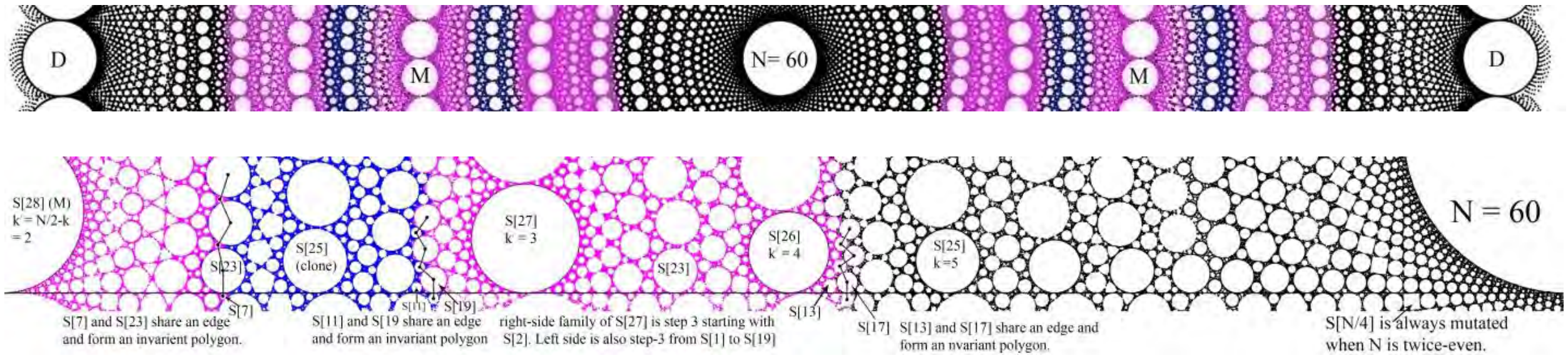

Section 2 describes the geometry and scaling of this region and we will repeat some of that analysis here. In Example 2.6 we applied this GFFT analysis to the S[13] and S[17] tiles which are a dual pair. Now the Edge Sharing Lemma shows that they share an edge since they are in the same step-4 GFFT family of S[26]. Note that they are really satellites of S[26] and it is clear that they can isolate S[26] locally so when S[16] is rotated around N, this isolation will be global.

**Figure 12** The S[25] region of N = 60. The only tiles here which are part of the (known First Family of N = 60 are the S[23]-S[26] tiles. The rest of these are displaced S[k] which are part of the GFFT local families of the S[k]. The blue lines of symmetry link all this to N and the local families of S[k] tiles always include N at S[N/2-1].

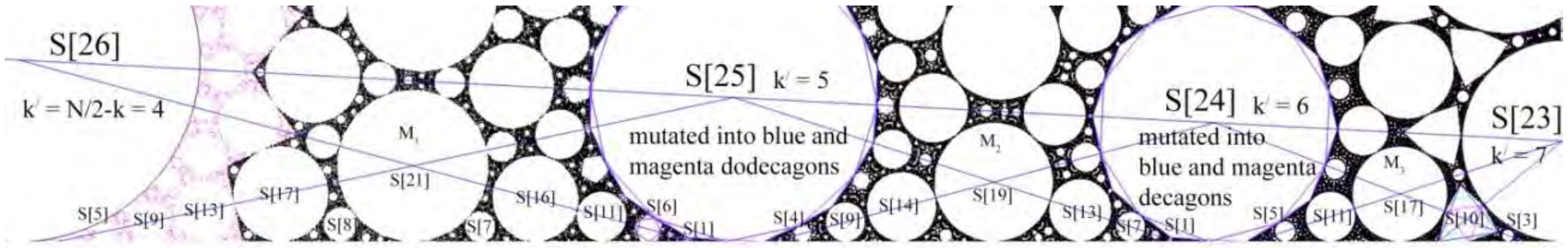

Here S[25] is the „shepherd" of the invariant black region and its step-5 right-side GFFT family includes S[24] and also N. The first member is S[4] and this is close to a dual match. The step-4 family of S[26] must also include N as well as S[25] so it has adjacent S[17] and S[13] tiles which share an edge. Note that as S[17] is rotated around S[26], S[13] is the perfect link between copies of S[17] and this remains true as S[26] is rotated around N. This local tiling is repeated at each primitive dual pair. When N was twice-odd these pairs always have different genders but there is no evidence showing that this sharing is more or less efficient than the twice-even case.

With careful analysis this dual pair could have been predicted because as the S[k] grow in size, there will be 'volunteers' forming between S[k] and they grow at a predictable rate.

The primary volunteers here are $M_1$, $M_2$ and $M_3$ and they have a special role as the only tile in the local families of neighboring S[k]. For N even the only tile in the First Family of N and D is M = S[N/2-2]. It will have 'effective' k′ that is the sum of N and D and that continues here so $M_1$ will be an S[21] with k′ = 5 + 4. This is what makes it possible to generate the local families and as N grows this can continue for multiple layers. A level-3 case for N = 300 is explained below.

The Mutation Conjecture says that for N twice-even S[k] will be mutated when N/gcd(k′ ,N) > 2 and here S[25] is the weave of 2 dodecagons because N/gcd(25,N) = 12 and likewise S[24] is the weave of 2 decagons. Even secondary tiles like S[6] of S[10] are mutated in the same fashion as their counterparts adjacent to N. We have no firm theory about how mutations evolve but in [H4] we generated W for some non-regular mutations and did not see much stability.

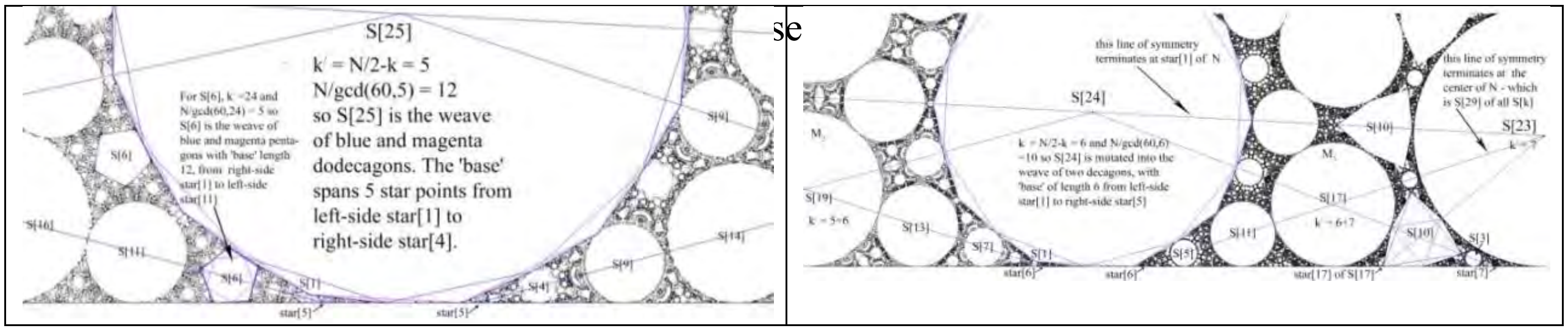

For N twice-even we have no general rule for the outermost 'shepherd' S[k]. Our best estimate in general is N/3 instead of N/4 as in the twice-odd family. For N = 60, S[25] is beyond the predicted S[20] and for N = 300 below the shepherd is S[128] well past the predicted S[100].

**Figure 14** The first few invariant regions for N = 300 where the 'shepherd' is S128].

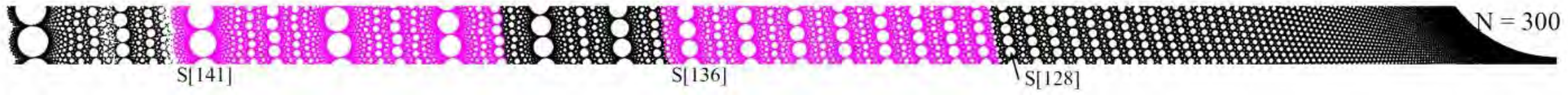

As N increases the S[k] are torn apart and most dual pairs occur in secondary families.

**Figure 15** The S[136] region for N = 300 is described here

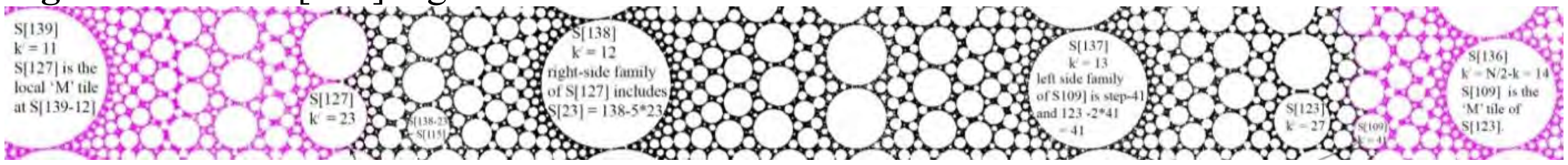

The S[136] transition is what we call level-3 because it entails two intermediate „M" tiles – namely S[123] and S[109]. The left-side magenta transition and the initial S[128] transition are both level 2 since they required just a single M tile. Here S[123] is the central M tile between the step-14 left-side family of S[136] and the step-13 right-side family of S[137] so it is the only tile in both families. S[123] has k′ = 27 so its right-side family has an S[136-27] = S[109] tile with k′ = 41. This matches the short left-side family of S109] which counts down twice from 123 mod 41. This small S[41] will share an edge with S[109].For black S[138] on the left the local M tile is S[127] with k′ =23 and it supports an S[23] because 138-5(23) = 23.

By way of full disclosure the next page will include a short practical guide toward generating these invariant regions and analyzing them using the GFFT. If the reader has a copy of Mathematica (any of the versions starting with 9.0 should work just fine) first download the FirstFamily.nb notebook from DynamicsOfPolygons.org and the only input needed is 'npoints'.

The default npoints is N = 14 but change it to 42 and from the Menu select Evaluate Notebook and Mathematica will quickly generate the symmetric First Families for N and D. To generate a rough web, Mathematica already knows the N/2-1 star points so set. **q1=star[[1]] ; q2=star[[20]]; WebPoints =W0[-q2[[1]]- s0/2, .10, -q1[[1]]- s0/2]]** This gives .1 spacing of points on each „forward‟ edge. Then **Web = V[WebPoints, 150]** will generate a level 150 web. Set bounds for the S[19] region of N with **left = c[19][[1]]; right = cS[18][[1]]; top = -.25; bottom =-1.1;** Then **Graphics[{AbsolutePointSize[1.0],Point[Web], poly/@FirstFamily, poly[S[19]},PlotRange[{{left,right},{bottom,top}}]**

If desired, open a new notebook with npoints = 21 and copy the First Family as a simple list. These tiles can then be imported into the N = 42 notebook with the same scale and location as the S[19] tile**. N21FF = TranslationTransform[cS[19]] /@(N21FirstFamily*hS[19])**

**Figure 16** The „M‟ region for N = 42 showing the imported N = 21 First Family

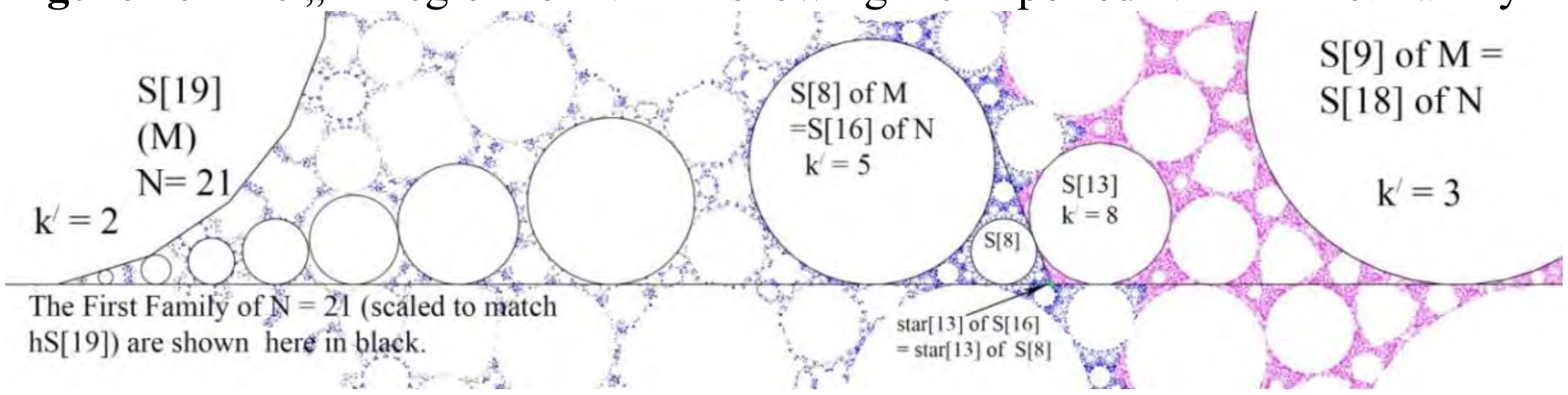

Since N is twice-odd, S[19] (and all „odd‟ S[k]) are mutated into N/2-gons. This matches the imported First Family of N = 21. These imported S[k] are S[2k] N-gons relative to N, so the local S[8] of M is an S[16] N-gon here and since it has k′ = 5, its local GFFT family counts down from S[18] to yield S[13] and S[8] which are dual. To construct them use star[8] and star[13] of S[16] wth offsets of sN/2 from the star points.

**MidS16 = {cS[16],-1}; StarS16 = Table[MidS16 + {hS[16]*Tan[k*Pi/42], 0}, {k,1,20}] ; MidS8 = {StarS16[[8]][[1]],-1} + {sN/2,0}, cS8 = MidS8 + {0,hS[8]}; S8 = TranslationTransform[cS8-cS[8]] [S[8]]** (then do the same for S[13] and S[3] if desired)

Like the star[k] of N, this star[13] of S8 has a star-polygon periodic orbit of the type shown below for star[11], but in general these orbits are just first approximations to the bounding orbits shown in Figure 21. The fit is pretty close for the canonical S[10]-S[11] pair shown below.

**Figure 17** The magenta period 42 star-polygon orbit of star[11[ of N

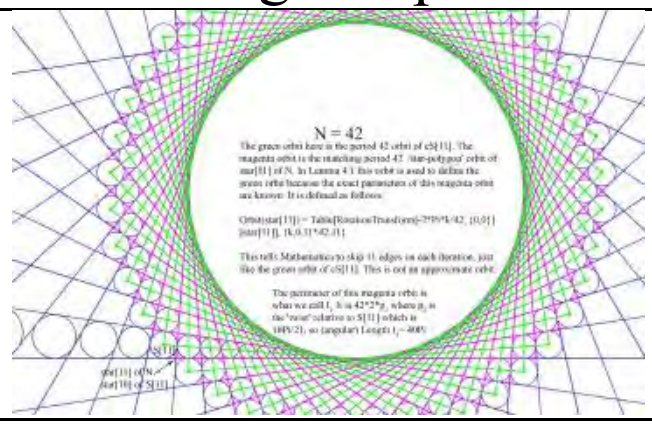

These star-polygon orbits were used in Lemma 4.1 to show that the matching green S[11] orbit has the same period. The star-polygon orbit here is exact. It is Orbit = Table[Rotation Transform [2kPi/N][star[11], {k,0,42*11,11}]. The boumdry is $I_2$ which is an 'integral curve' with integral (42)(2)(10$\pi$/21) = 40 $\pi$.

Below we will define the intrinsic angular 'twist' of star[11] to be the radian measure of the angle defined by S[11] and S[10]. Here that value is $10\pi/21$ and it would have units of 'velocity' as star[11] traces $I_2$ so integrating over the length 2*42 gives $E_2$ = energy = $40\pi$. See Fig.22.

The 'kite' that R.S. used in S[1] (2007) had a matching 'dart' and in 1974 Roger Penrose showed that together they could tile the plane in a non-periodic fashion. This was a startling discovery that no one had predicted. At about the same time the KAM theorem was being debated in the secrecy of the Iron Curtain, Penrose at Oxford had become interested in relativistic geometry and the stability problem of a collapsing star. In 1905 Karl Schwarzschild found a solution to Einstein's new equations that would allow a region so dense that space-time would fold in on itself, but no one believed this could actually happen. In the early 1960's a group of Russian scientists including Evgeny Lifschitz used the new perturbation theory to show that this kind of collapse of a star could never actually happen. Most scientists around the world agreed with them, but Roger had already shown how it was possible to trap light so in 1965 he published a short paper showing how a star might survive the collapse as a 'singularity'. Over the years these ideas gained support as 'cosmic rays' were detected with no known source. In 2020 Penrose was awarded a Nobel prize for this paper.

**Figure 20** The Schwarzschild 'point of no return' radius is about 1.5 km for the medium mass black hole off the screen on the left. The dynamics here are relativistic Hamiltonian which applies to much of the quantum world also. See ]No] and [H2].

<table>
<tr>
<td>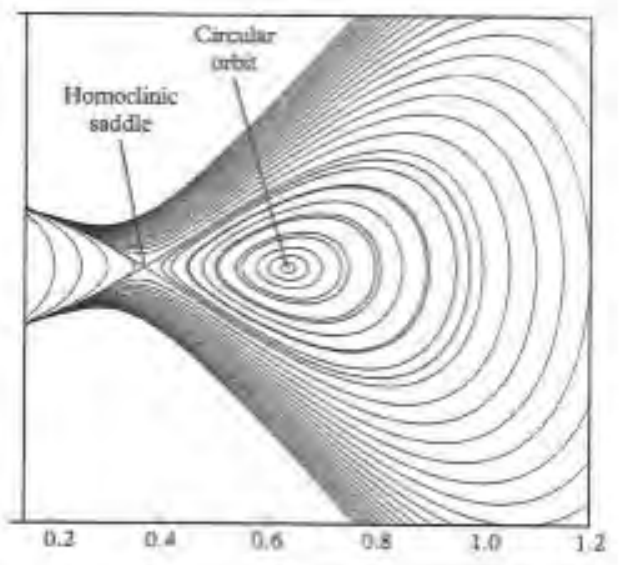
</td>
<td>In the vicinity of a black hole of mass M the Schwarzschild radius $r = 2GM/c^2$ where G is the gravitational constant and c is the speed of light. Here it is about .15 km as shown by the left side of this rectangle. The unstable homoclinic fixed point is where circular or elliptical orbits break down as they approach the black hole. These orbits can diverge at high speeds as observable 'radiation' which was detected in the 1960's by U.S. satellites which were looking for signs of Soviet nuclear tests in violation of treaties.</td>
</tr>
</table>

<table>
<tr>
<td>In this cosmic 'web' from Scientific American (2/26) each node is a collection of galaxies and on all scales there seem to be rotations.<br>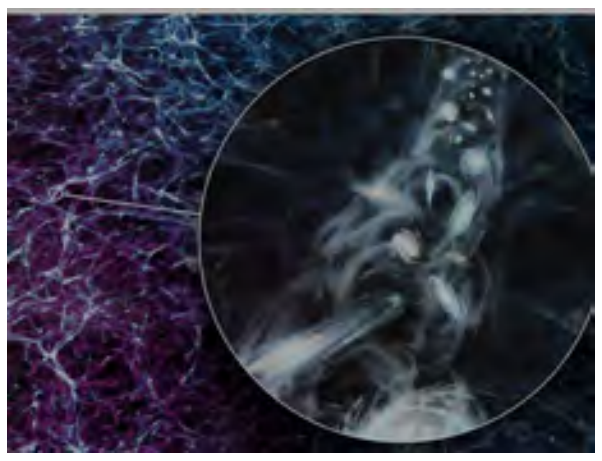</td>
<td>If a cloud of gas is irradiated by polarized light from early stars this could tend to polarize certain regions as cw and ccw rotations. On smaller scales molecular clouds would also be expected to react to radiation and if some of these molecules are 'prebiotic' with chiral structure they could be selected out by the same processes recently described in the Nobel Prize work of H. Kegan and K. Soai. See astrobiology and the last page below.</td>
</tr>
</table>

### The Invariant Regions and Indices $L_N$ for N-even ≤ 50 and N-odd ≤ 25

This list will cover the three cases of N twice-odd, N twice-even and N odd but from our remarks above the N-odd case can be regarded as a subset of the twice-odd case so we will only do the first 25 cases. In each example we will give the ordered minimal list $L_N$ which should be sufficient to reconstruct the invariant regions. This is an annotated list that will include comments.

**Case I - The twice-odd family**

This is the best-understood family which includes the matching N/2 case. The first non-trivial dual pair after S[1] and N, will be S[f4] and S[f4+1] where f4 = Floor[N/4] and this is the only known case where a dual pair is simply two neighboring S[k] tiles. In addition every S[k]-S[k′] dual pair has exactly one primitive index so the number of primitive pairs $|L_N|$, will always be equal to the algebraic complexity $C_N$ of N. The illustrations will only show the region right of the central S[N/2-2] 'M' tile because the left side is symmetric and the indices in $L_N$ will apply to this side with the understanding that the left-side indices are obtained by swapping S[k] of N with DS[k] of D. Therefore there will be $|L_N|$-1($C_N$-1) invariant regions right of M and the total number of invariant regions will be $2|L_N|$-2.

| | |
|---|---|
| N = 6 (3) $L_6$ = {1} with the only dual pair S[1] and N at DS[2], so no invariant regions right of S[1] = M. Complexity 1 | 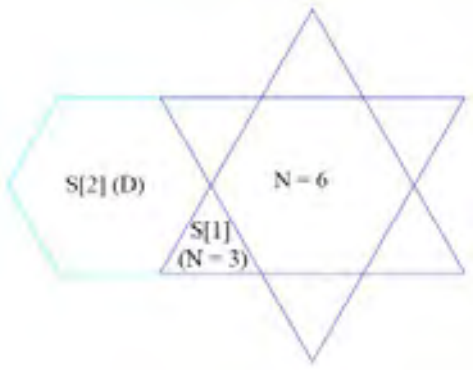  |
| N = 10 (5) $L_{10}$ = {1,2} with dual pairs S[1]-N and S[2]-S[3] These two transitions right of M = S[3] define just one region containing S[1] and S[2]. Since S[2] has gcd(k,N) = 2 its orbit will be period 5 instead of 10 and the blue-magenta S[2] decomposition will occur for N twice-odd. But these two regions are related by a simple rotation and can be regarded as a single region to match the black region left of M. Every periodic orbit defines a locally invariant region but this is the 'next-level' of invariance and can be very difficult to apply. In our detailed studies of the S[1]-S[2] region in [H6] and [H7] we make no attempt to separate their orbits but instead use the combined orbits. | 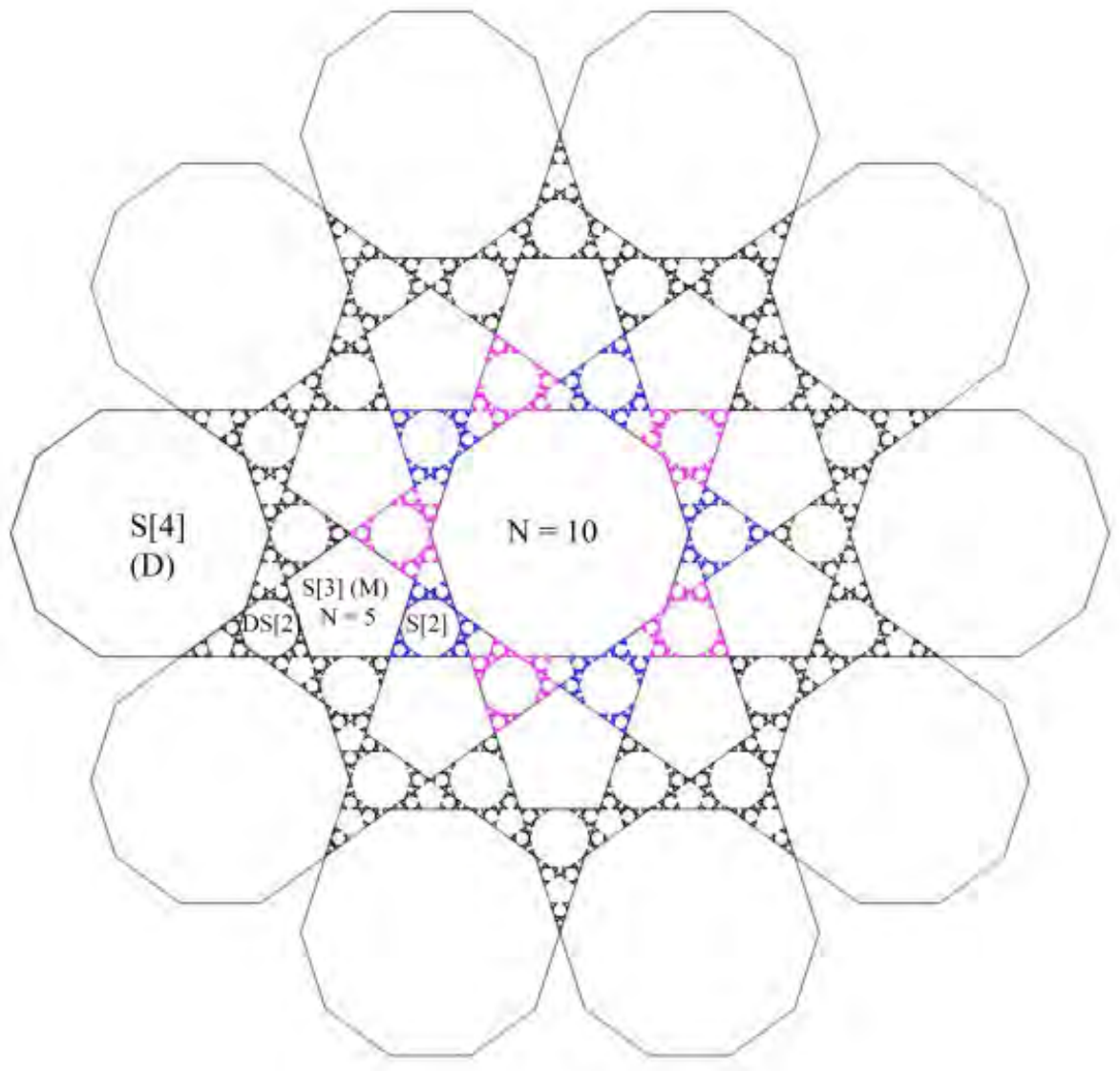  |

<table>
<tr>
<td>N = 14 (7) $L_{14}$ = {1,3,2} and these three transitions right of M are labeled $I_1$, $I_2$ and $I_3$. They correspond to the dual pairs S[1]-N, S[3]-S[4] and a S[2]-S[5] where S[2] is a clone of the original S[2]. Note that the primitive star[k] of N are indeed star[1],star[3] and star[5] but the Invariance Conjecture only implies that at least some clone of these will be a boundary. In the case of star[3] that clone is the star[3] of S[4] as noted in the Introduction. Note that for N twice-odd the inner region will always have bounding 'shepherd' tile S[Floor[N/4]].</td>
<td>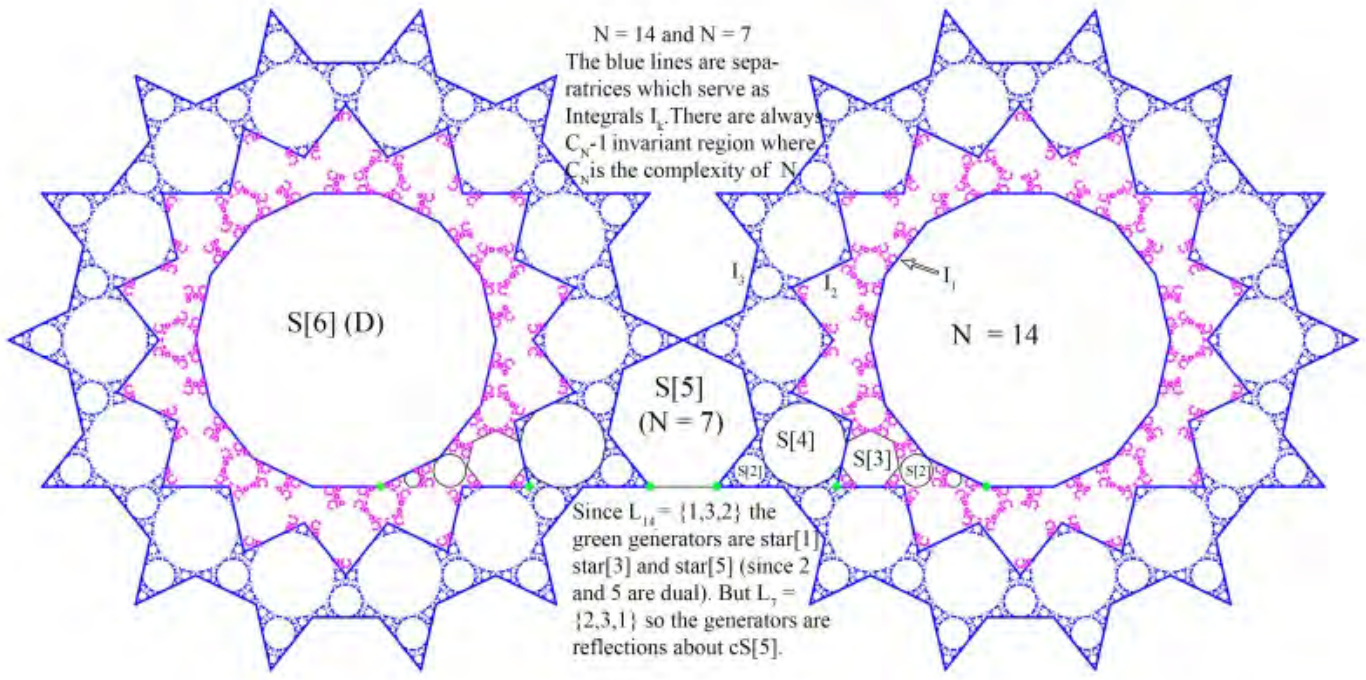
</td>
</tr>
<tr>
<td>N = 18 (9) $L_{18}$ = {1,4, 2}<br>Only 2 regions right of M but we show all 4 regions. The S[3]-S[6] pair are clearly not dual. Once again we show the various resonances inside the initial invariant region but make to attempt to isolate them. At D the distinction always vanishes and we will no longer show the individual regios in the plots below.</td>
<td>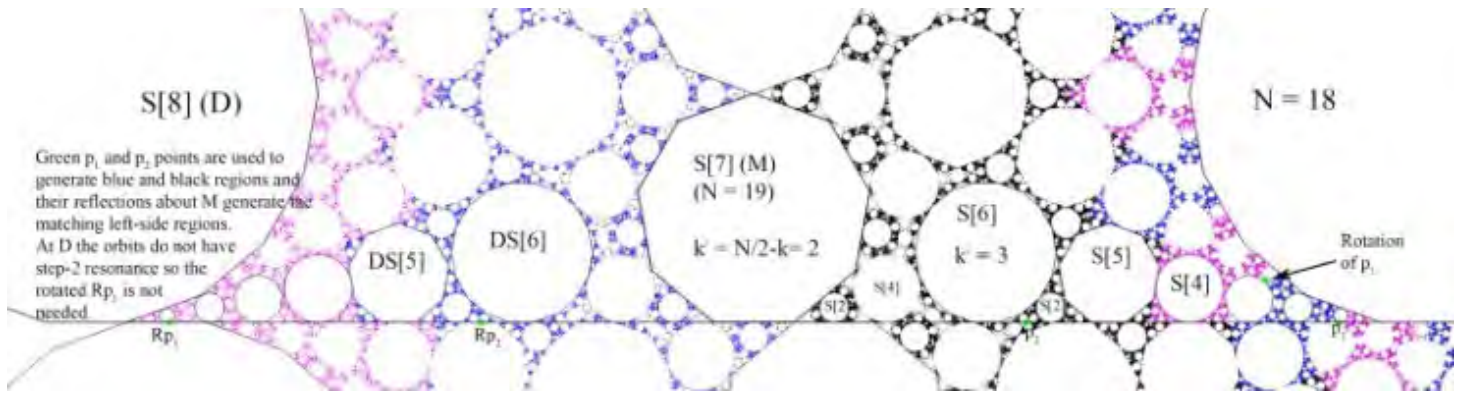
</td>
</tr>
<tr>
<td>N = 22 (11) $L_{22}$ ={1,5,4,3,2}<br>Maximal 4 regions right of M in a perfect 5,4,3 count-down. This is matched by N = 11 with $L_{11}$ = {2,3,4,5,1}. This gives a fresh perspective on this unique case of quintic algebraic complexity with 1-1 correspondence between invariant regions and primitive scales.</td>
<td>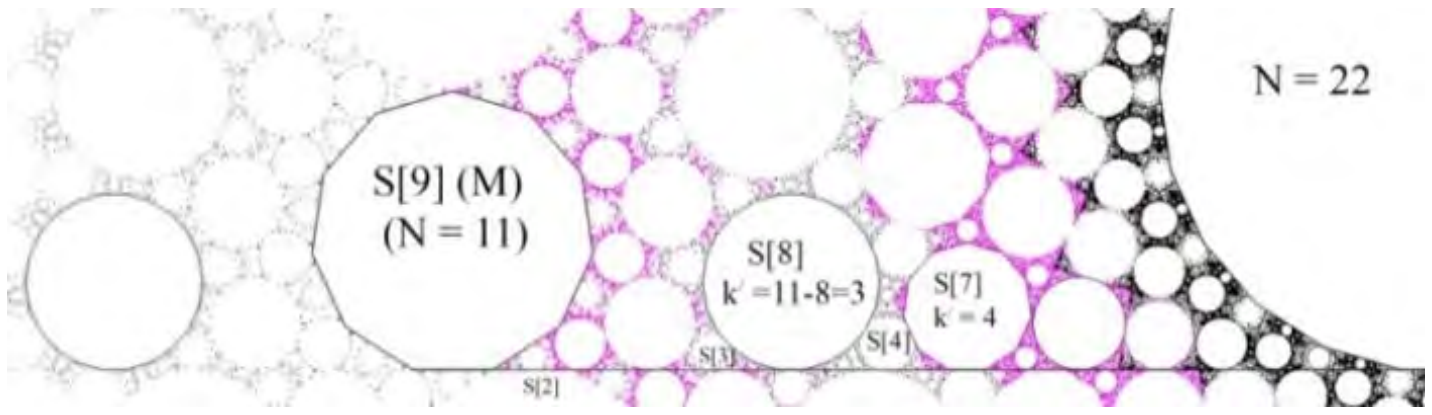
</td>
</tr>
<tr>
<td>N = 26 (13) $L_{26}$= {1,6,4,3,5,2}<br>Twice-prime so maximal 5 regions right of M. All possible pairs exist.</td>
<td>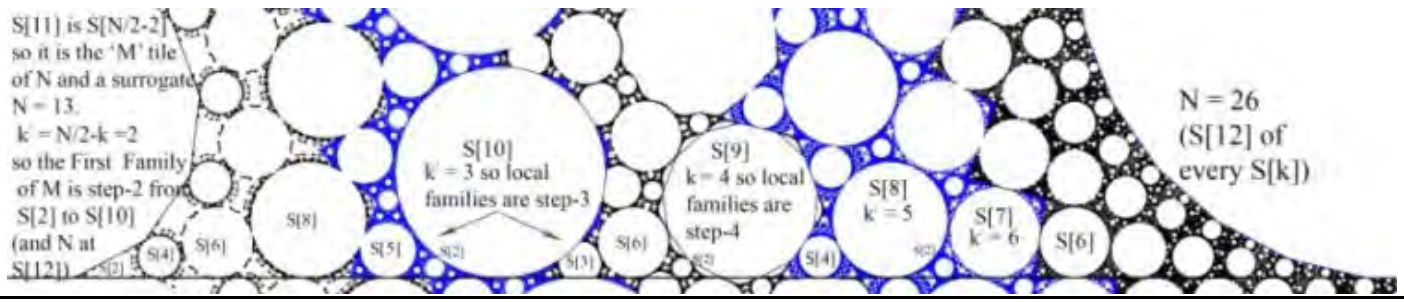
</td>
</tr>
<tr>
<td>N = 30 (15) $L_{30}$ = {1,7,4,2}<br>because these are the primitive indices in $Z_{15}$ , so just 3 regions right of M instead of the ideal 6</td>
<td>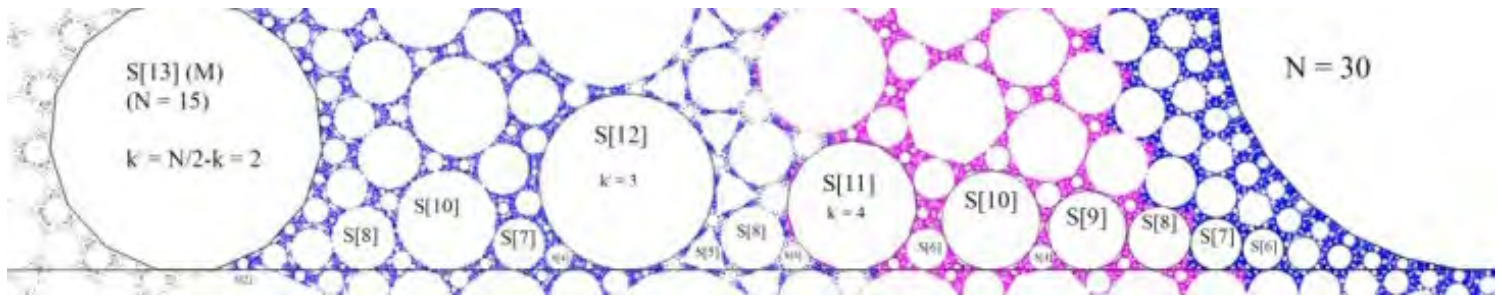
</td>
</tr>
</table>

| | |
|---|---|
| N = 34 (17) $L_{34}$ = {1,8,6,4,7,3,5,2} Twice-prime so maximal 7 regions right of M to match the complexity of N = 34. Upon reflection these regions match N = 17. | 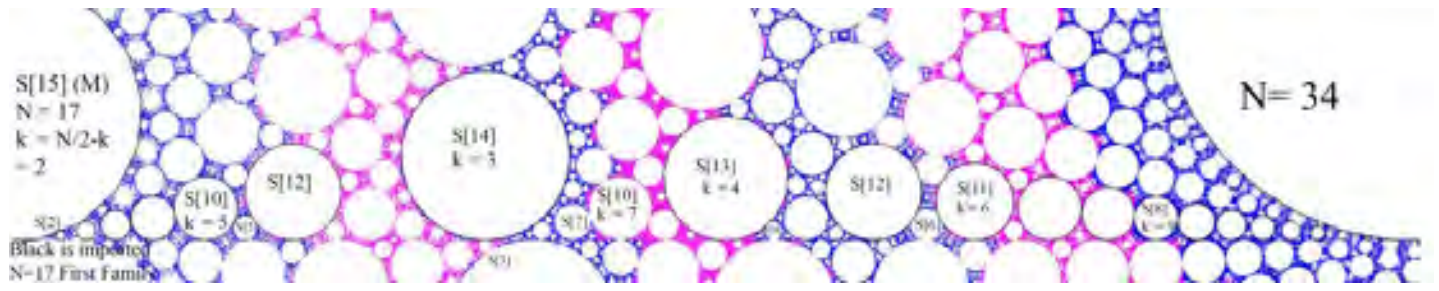 |
| N = 38 (19) $L_{38}$={1,9,6,5,4,3,8,7,2} Maximal 8 regions right of M. The arrows point to the S[k] involved. | 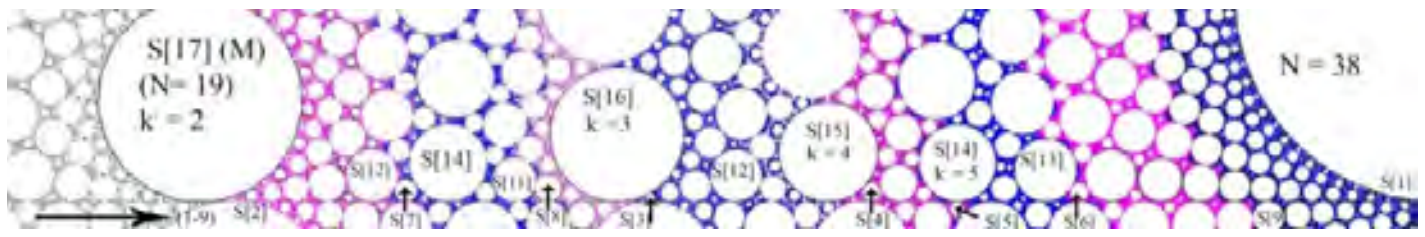 |
| N = 42 (21) $L_{42}$ ={1,10, 5, 4, 8, 2} devoid of multiples of 3 and 7 in $Z_{21}$ so just 5 regions right of M. | 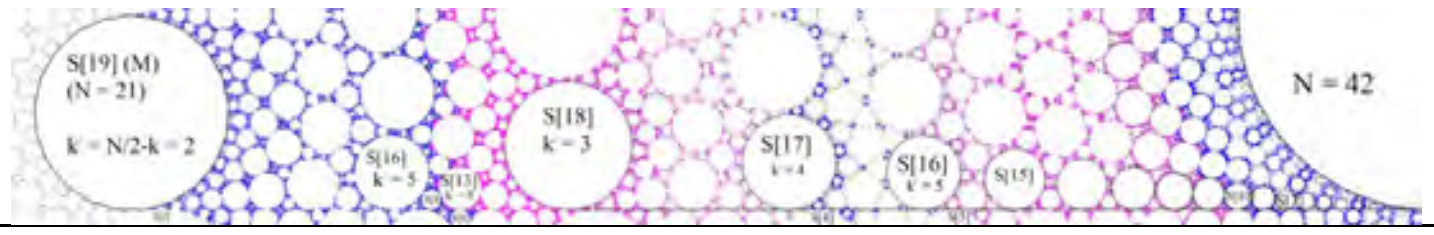 |
| N = 46 (23) $L_{46}$ = {1,11,8,6,9,4,10,3,5,7,2} Maximal 10 regions right of M. Secondary sequence from base of S[12] with an S[1] and the sequence grows step-2 to yield an S[8] duel to S[15] then it continues all to S]16] | 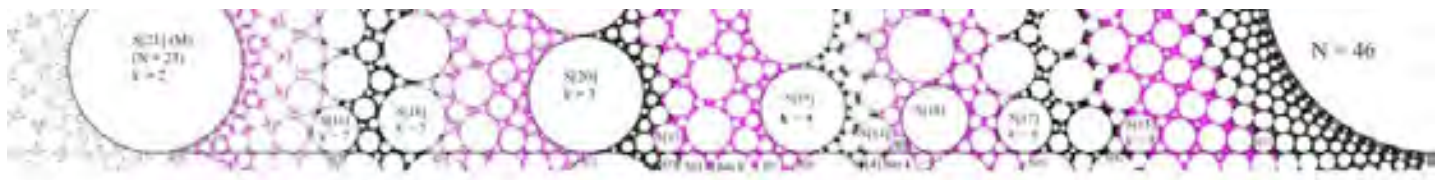 |
| N = 50 (25) {1,12,8,6,4,7,3,11,9,2} missing S[5] & S[10] as expected | 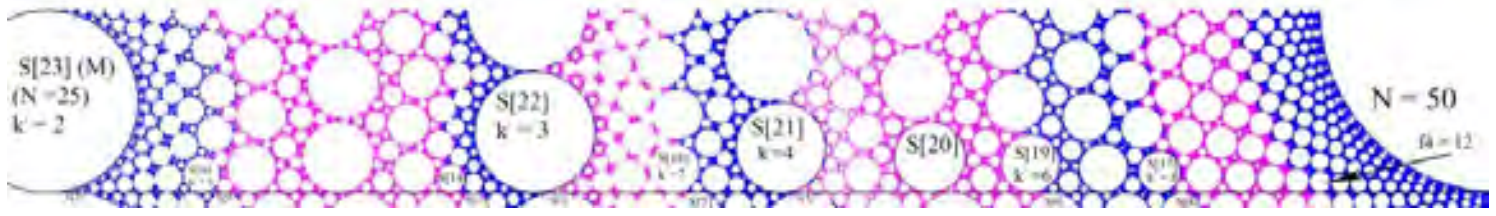 |
| N = 202 (101) $L_{202}$= {1,50,……,2} (length 50) This is the region between S[98] and S[99] which is M and a surrogate N = 101. There clearly are multiple invariant regions here but only a few are shown. The total number of invariant regions right of M is $\|L_{202}\|$-1 = 49 but the first region takes up 50 S[k] and the next region extends out to S[67]. This this leaves only 33 S[k] to host the remaining 48 regions, so there will be multiple regions between some S[k] as show here. | 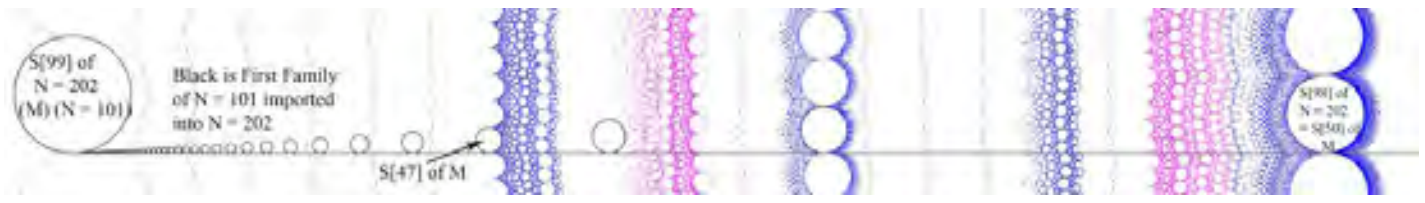 |

**Case 2 - N twice-even**. When N is N twice-even the canonical S[N/2-2] M tile has reflective symmetry just like the twice-odd case but it has no dual S[2] to create a dynamical barrier. Therefore M now defines just one local invariant region instead of two and the total number of invariant regions is $2|L_N|-1$ but the complexity $C_N$ is now $2|L_N|$, so the number of invariant regions is still $C_N$-1 just like the twice-odd case. This occurs because the left side of M is now algebraically independent of the right side even though it is symmetric with symmetric indices. Unlike the twice-odd case we no longer have a formula for the 'shepherd' of the first region, but the lines of symmetry shown here can be useful. As noted earlier the extra symmetry of the twice-even case means that when listing primitive pairs it is sufficient to use $Z_{N/4}$ instead of $Z_{N/2}$. So for N = 300 we use $Z_{75}$ and choose the dual from $Z_{150}$.

| | |
|---|---|
| 8 {1} $Z_{2\backslash}$ = {1,2} so the only dual pair is {1,N} = {1,DS[3]} = {1,3}. | 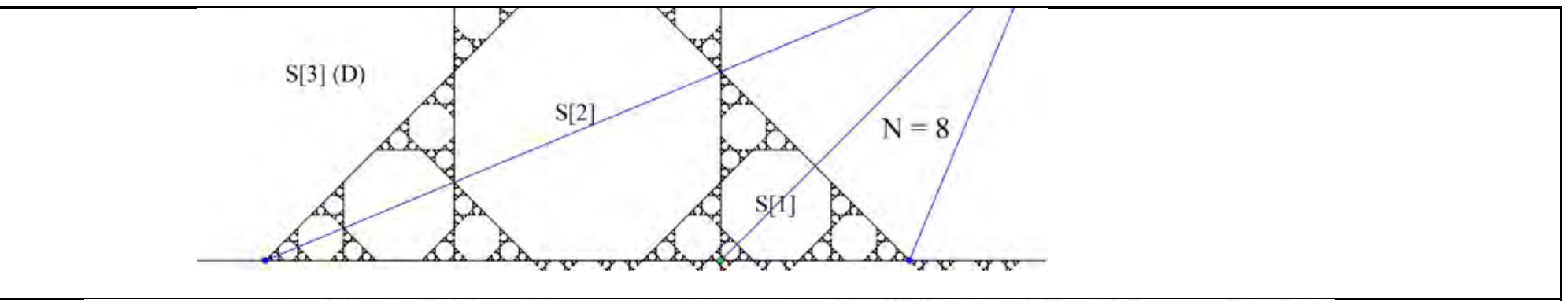  |
| 12 {1} $Z_3$ = {1,2,3} so the only dual pair is {1,5}. | 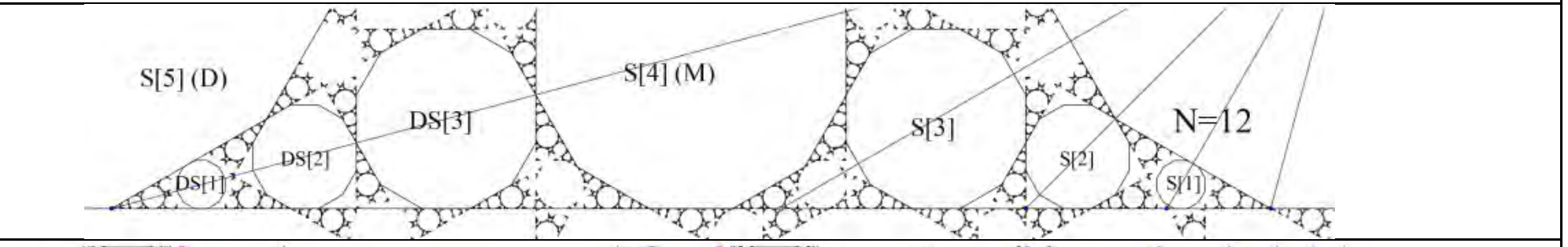  |
| 16 {1,3} $Z_4$ = {1,2,3,4} so two pairs {1,7} and {3,5} complexity 4 | 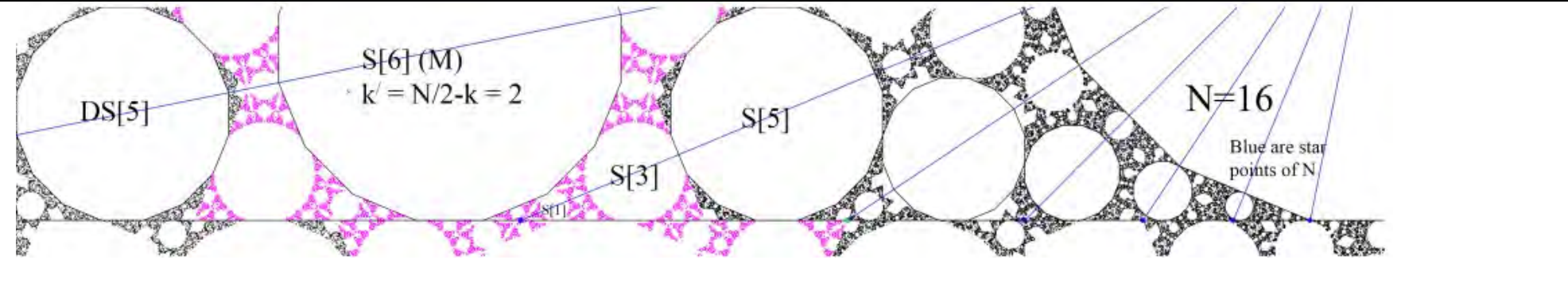  |
| 20 (1,3} The only non-trivial dual pair is S[7] and S[3] | 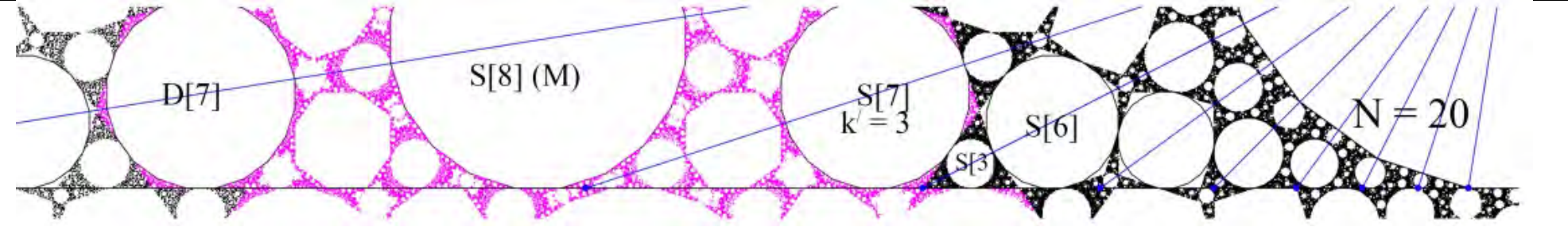  |
| 24 {1,5} Lines of symmetry pass through cS[9] the 'shepherd'. S[5] and S[7] are 'foundlings' at the foot of M but S[7] is rejected by M. | 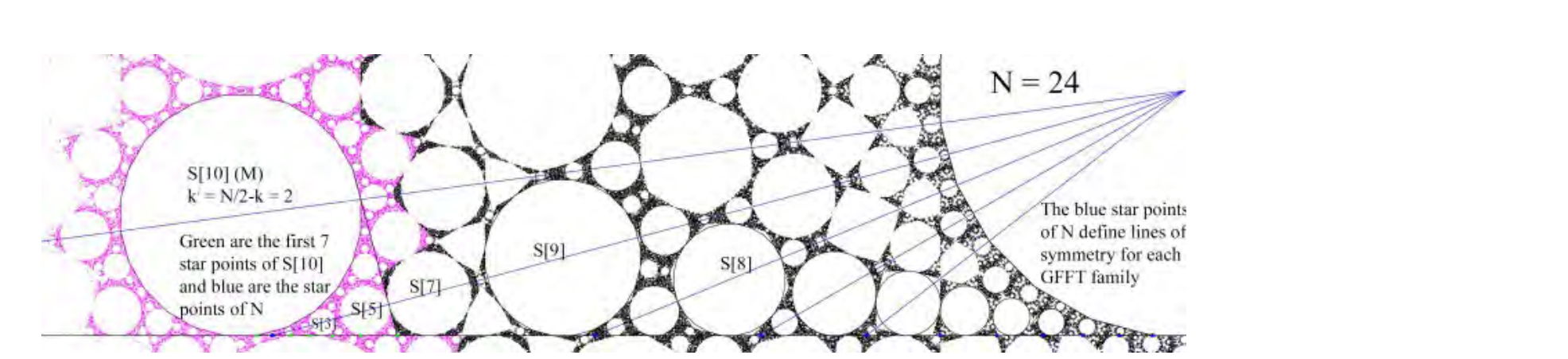  |
| 28 {1,5,3} Friendly world with S[9] at k′ = 5 next to S[10] | 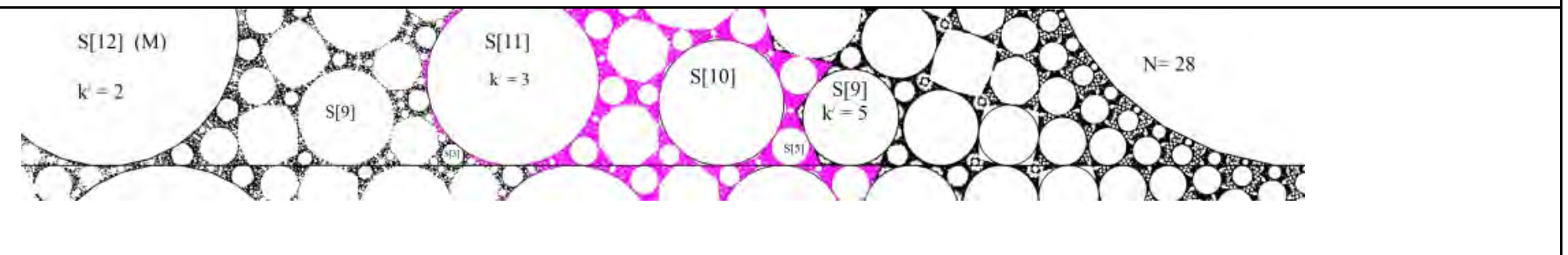  |
| 32 {1,5,3,7} Both sides of M shown here. | 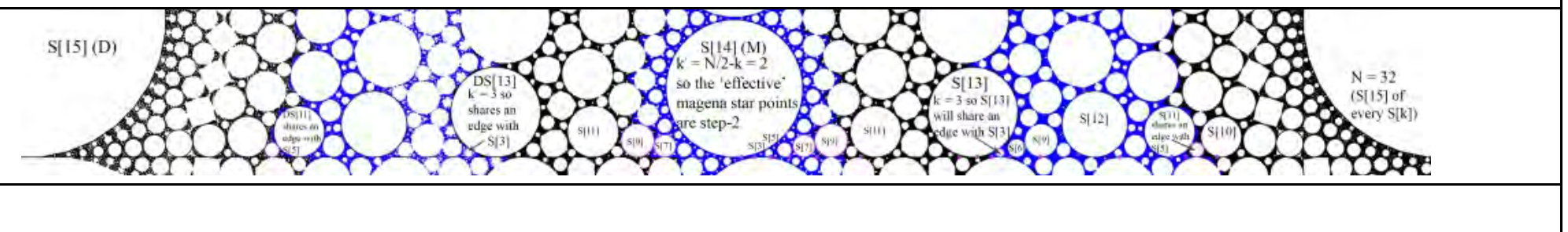  |
| | |
| 36 {1,7,5}<br>magenta S[13]<br>has $k' = 5$ | 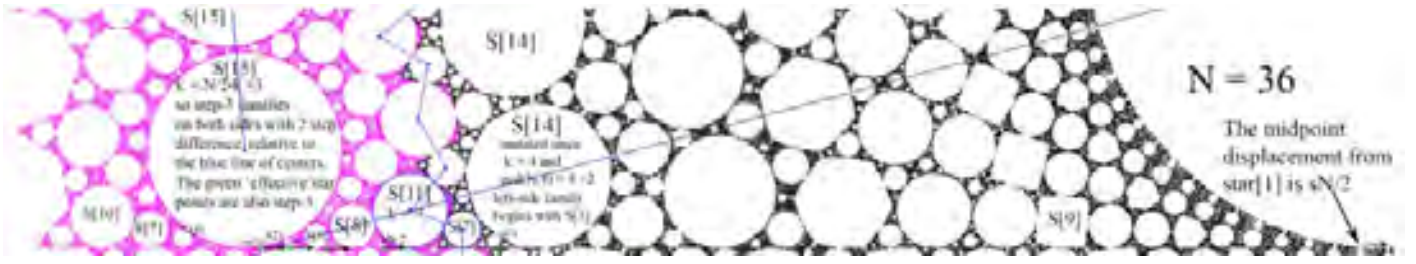 |
| 40 {1,7,3,9}<br>S[3] isolates S[17]<br>like a local M tile | 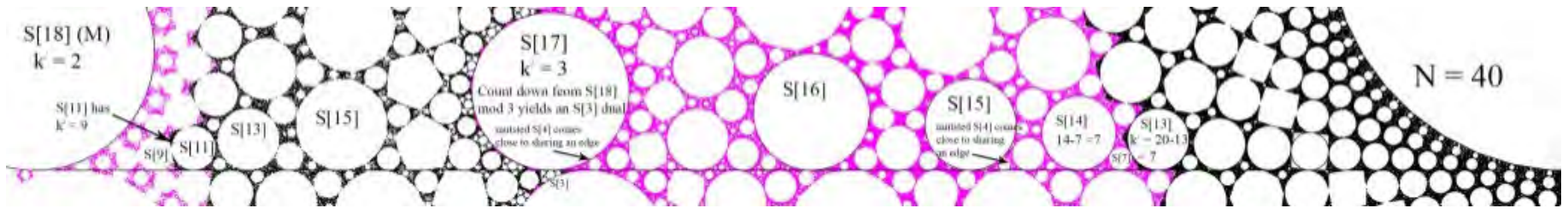 |
| 44 {1,7,9, 3, 5}<br>$Z_{11}$ gives<br>1,9,7,6,5,3} and<br>{7,15} is first | 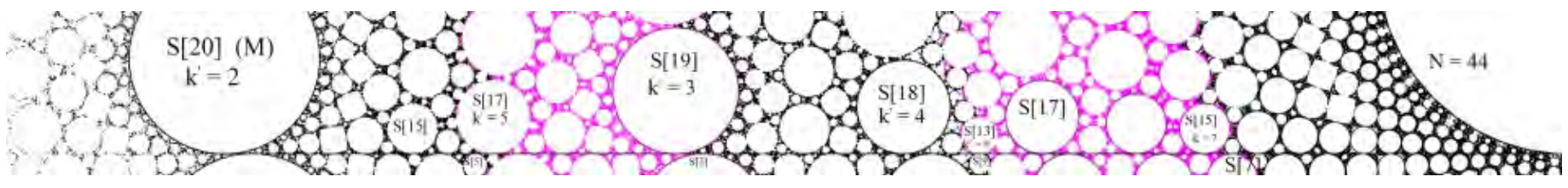 |
| 48 {1,5,7,11\|<br>$Z_{12}$ gives<br>{1,11,7,5} none at<br>level 1 so {5,19} is<br>first | 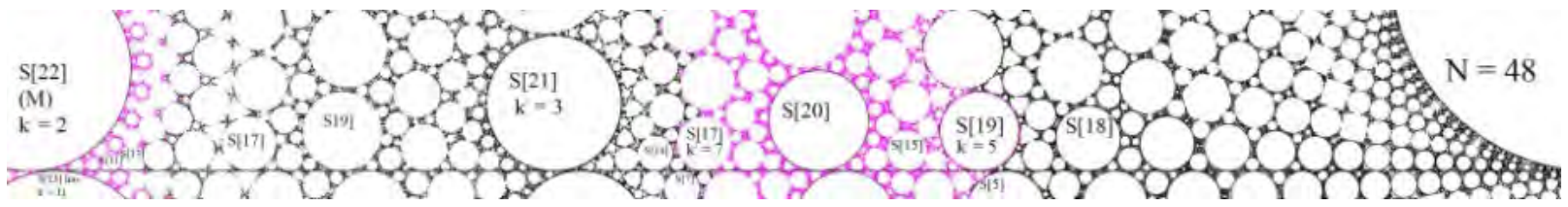 |
| 60 {1,13,11,7}<br>from $Z_{15}$ and now<br>S[13] and S[17] are<br>foundlings of S[26] | 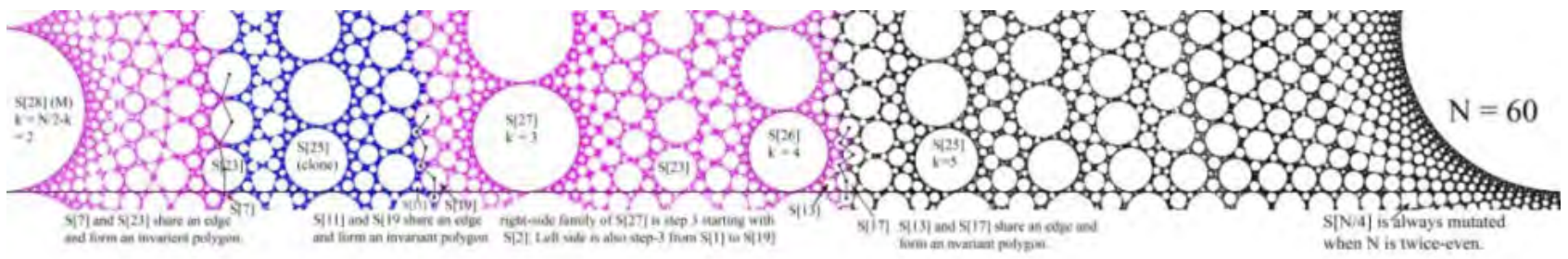 |
| 300 {1,43,41,23,..}<br>(from $Z_{75}$). $\|L_{300}\|$<br>should be 40 | 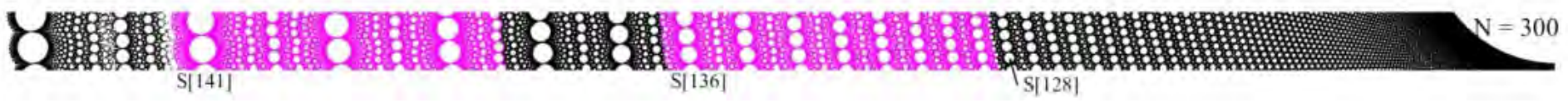 |

**Case 3** : **N-odd**

When an odd N-gon is at the origin all the S[k] are 2N-gons which are congruent to DS[2k[ of D. In particular S[1] of N is congruent to a DS[2] of D and all analysis of the S[k] of N can be carried out in the context of D with N as the DS[N-2] 'M' tile of D. The $L_N$ sequence will show the number of transitions from N to D, so it will be identical to $L_{2N}$ but the order will be reversed. The total number of invariant regions will be $|L_N| -1 = C_N -1$.

| | |
|---|---|
| 5 {2,1} 'On' then 'off'<br>This is the inverse of N<br>= 10 which is {1,2}. | 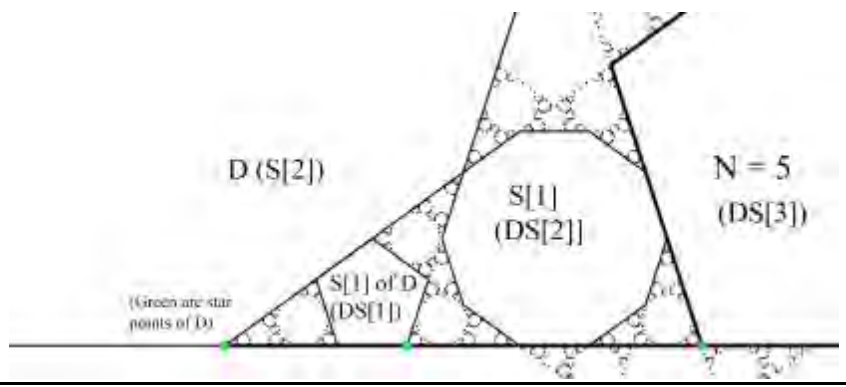 |
| 7 (2,3,1} This DS[3]-<br>DS[4] edge sharing is<br>the twice-odd signature.<br>One transition for every<br>DS[k] from 1 to 3.<br>Complexity 3 | 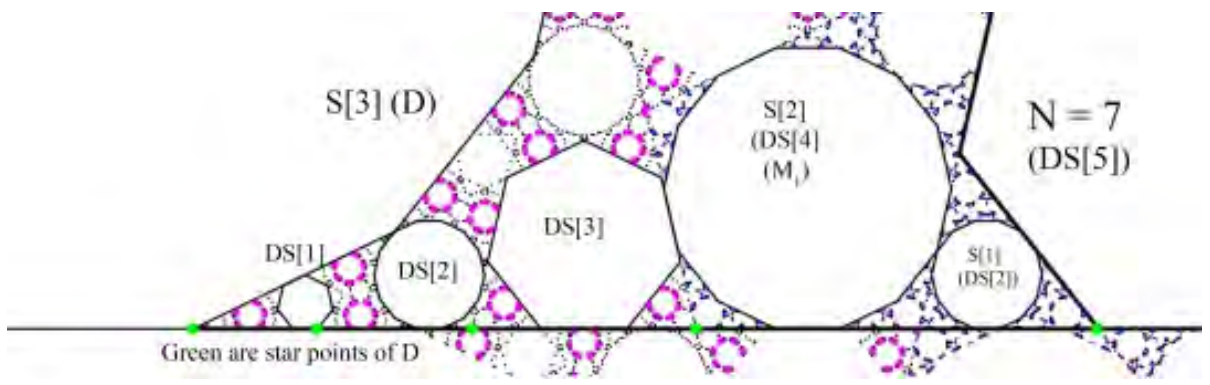 |
| 9 {2,4,1} One transition for each DS[k] from 1 to 4 – excluding 3 | 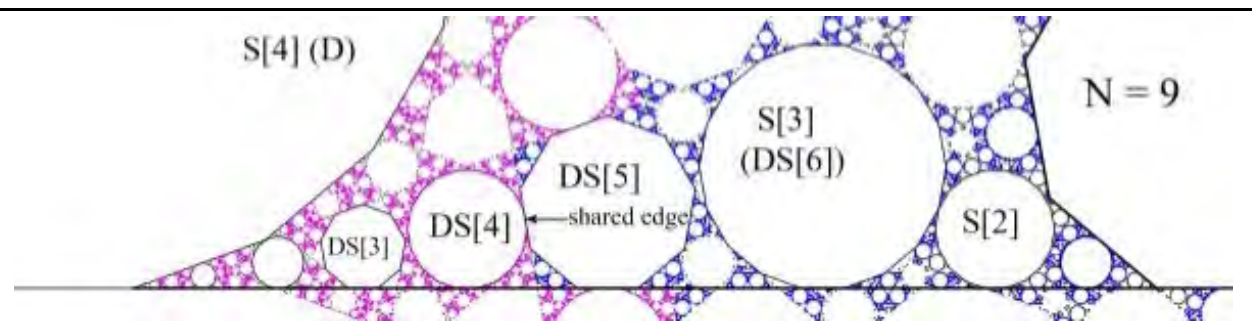  |
| 11 {2,3,4,5,1} See [H3] for a 3D map of this region. | 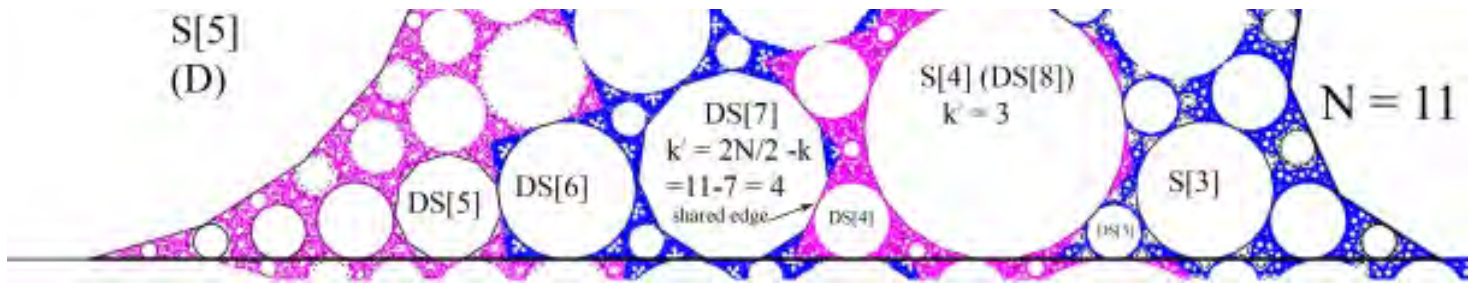  |
| 13 {2,5,3,4,6,1} the inverse of N = 26. One transition for each DS[k] from 1 to 6 | 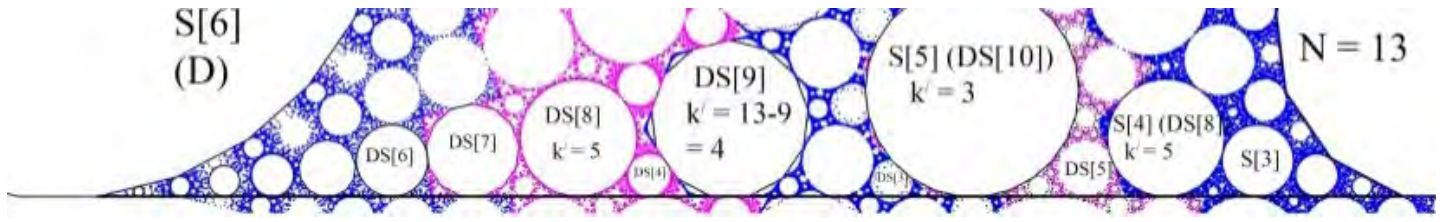  |
| 15 {2,4,7,1} Missing 3 and 5 | 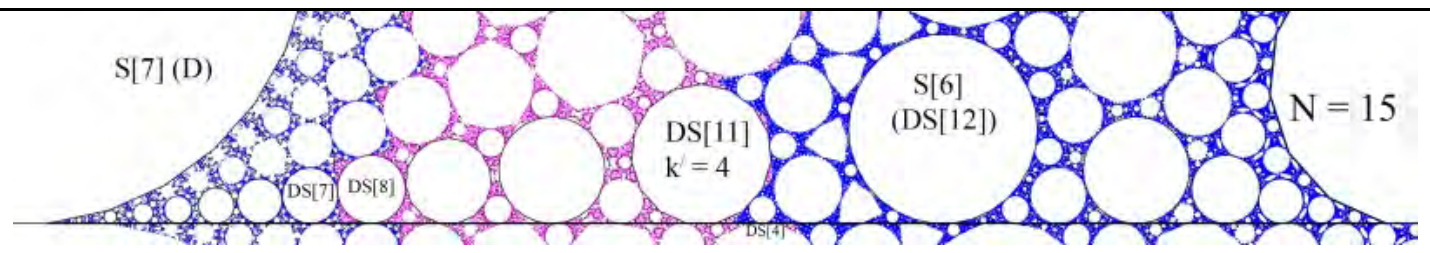  |
| 17 {2,5,3,7,4,6,8,1} 7 regions right of D. Here we track the 'invariant curves' as they circle N'. | 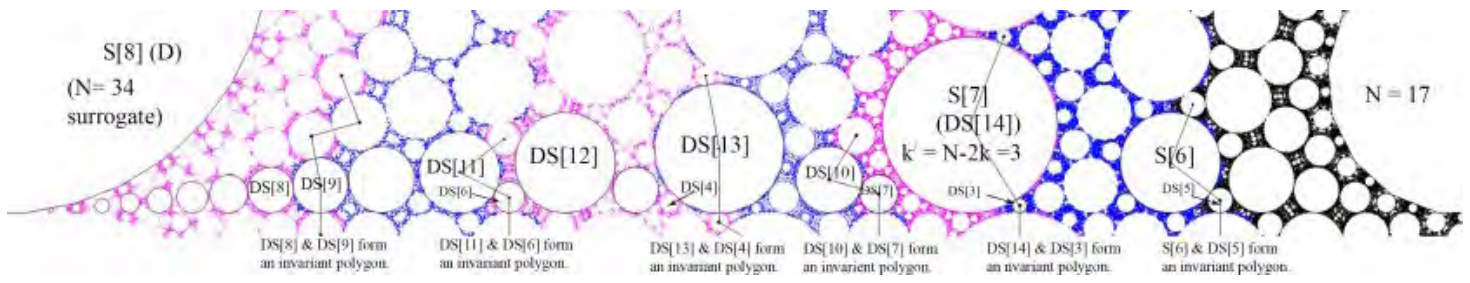  |
| 19 {2,7,8,3,4,5,6,9,1} One transition for every DS[k] from 1 to 9 | 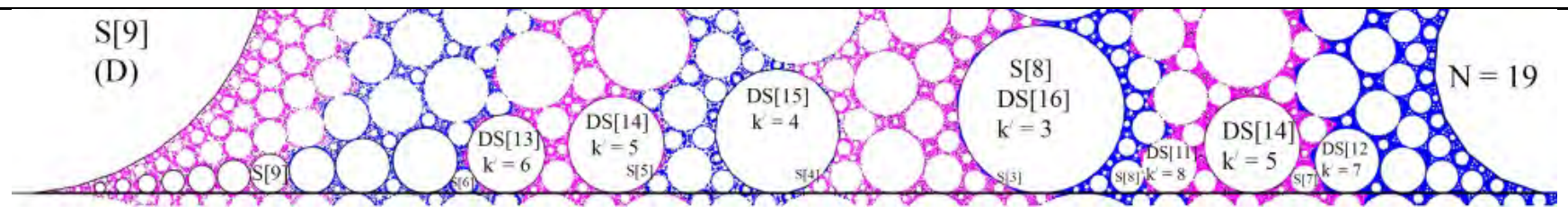  |
| 21 {2,8,4,5,10,1} One transition for each DS[k] from 1 to 10 excluding 3,6,9 and 7 | 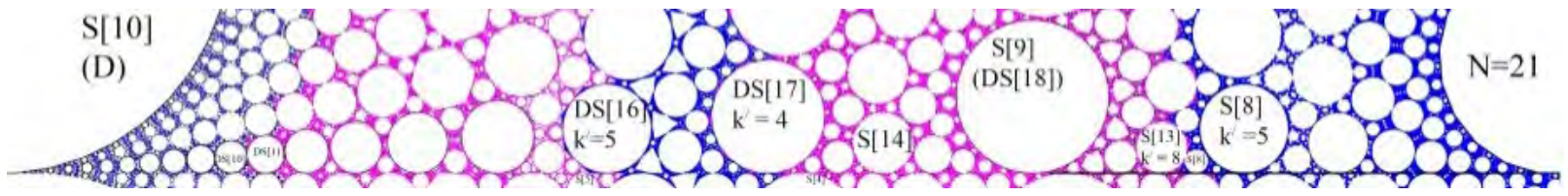  |
| 23 {2,7,5,3,10,4,9,6,8,11,1} One transition for each DS[k] from 1 to 11 | 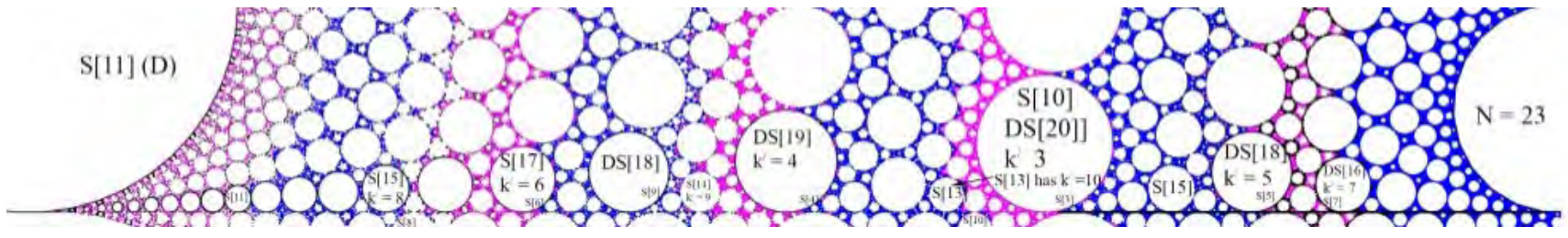  |
| 25 {2,9,11,3,7,4,6,8,12,1} (5 and10 missing) $\|L_{25}\|$ = 10 so 9 regions. The local $M_1$ is DS[22] with k′ = 3 and it defines two regions. | 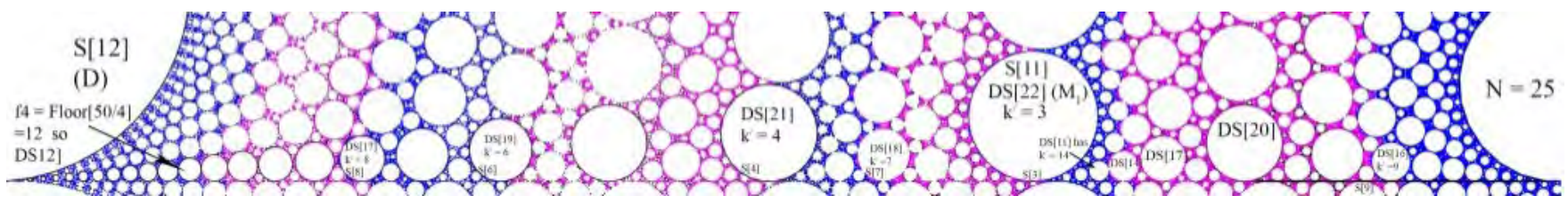  |

## How to define a discrete Hamiltonian for a regular N-gon

It can be a challenge to track the boundaries of the invariant regions but it can always be done with the help of the $L_N$ lists. Solving for these curves using Det(J) = 0 could be very difficult. Below are the 6 bounding curves for N = 42 . We call these $I_k$ because they correspond to 'integrals of motion' for classical Hamiltonian dynamics. By the Invariance Conjecture each $I_k$ curve is based on an S[k]-S[k′] primitive dual pair whose matching star[k] point is the base for $I_k$. These dual pairs are similar to the 'soliton' standing waves of a Hamiltonian which also have their origin in integral curves, so the Hamiltonian for N = 42 will have 6 solitons based on the primitive star[k] - which will be the only 'elements' of H.

The Arnold-Liouville Theorem of classical physics shows that knowing n integrals of a Hamiltonian with n degrees of freedom would be sufficient to determine the dynamics because analytical methods like Green's Theorem can be used to define the motion between integrals. For regular N-gons the Hamiltonian H will generate the $I_k$ integrals, but the invariant regions between the $I_k$ will not have a simple solution because every region represents a distinct 'eigen-state' of H (and N). The Hamiltonian we will use is based on a distinct twist map for each region but putting them all together will be a major challenge. The eigen-states of a regular N-gon are determine by the divisors of N, so to simulate the structure of N = 42 in material science would require solid state crystals that are opaque to frequencies based on 2, 3 and 7 . Such crystals can be formed from layered 'superlattices' and it is not unusual for solitons to exist in such crystals.

**Figure 21** The 6 $I_k$ seraratrices for N = 42 (and N = 21)

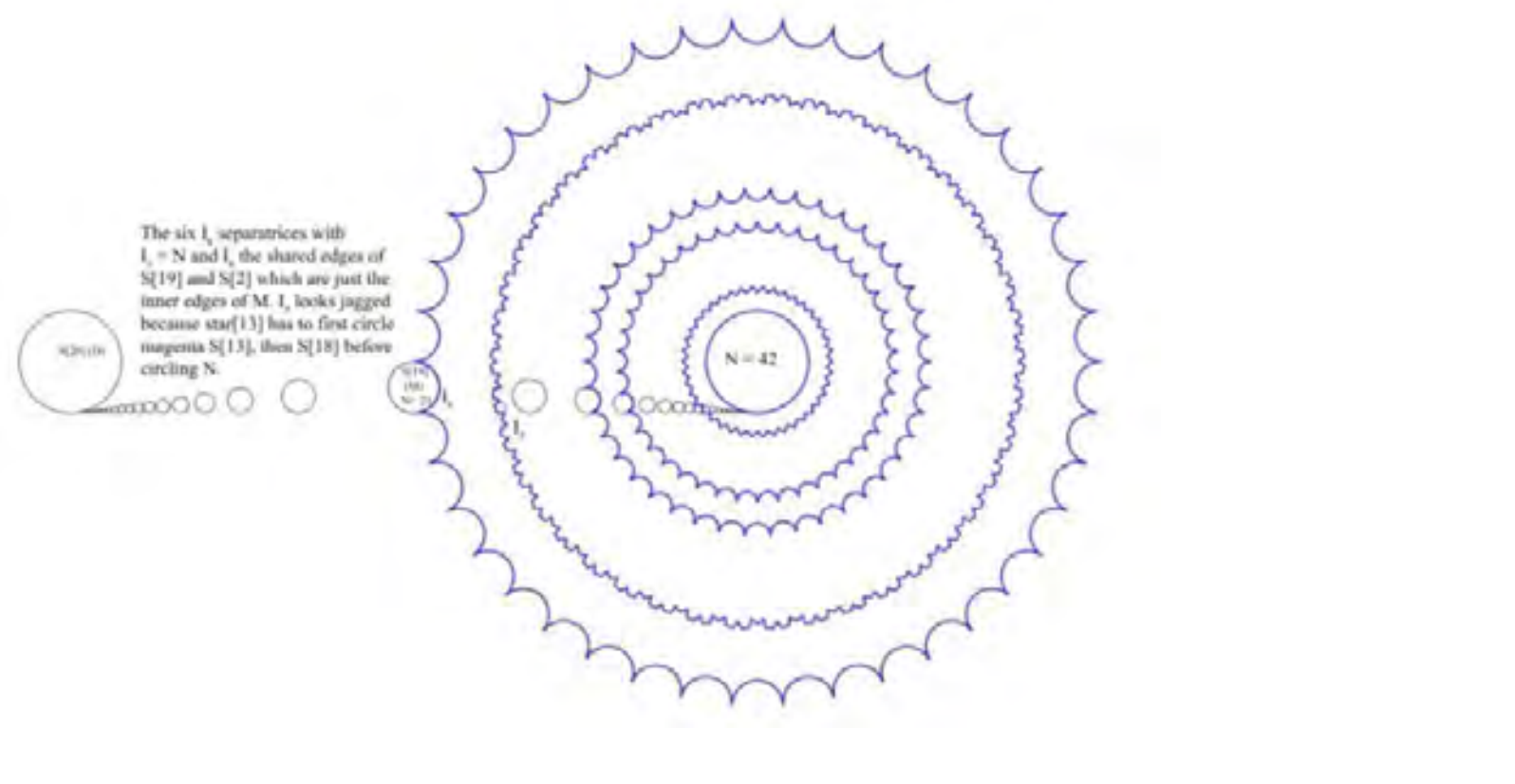

We will describe a discrete Hamiltonian H that will mimic N = 42 at least as far as generating the $I_k$. For a well-behaved integrable Hamiltonian like the one shown here, these $I_k$ are called integrals of motion and they may be sufficient to determine the dynamics.

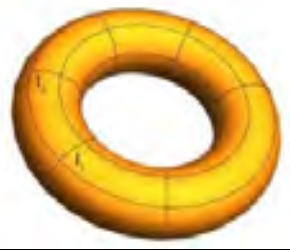

In the early 1970's physicists such as Donald Sutherland [Su] at Berkeley and Francisco Calogero [Ca] in Rome realized that the Schrodenger wave function might have more fine structure than expected as they found many-body dim-1 Hamiltonian systems that were exactly solvable in both the classical and quantum mechanical sense. Calogero [Ca] (1969) describes such a system $H=\frac{1}{2}\sum_{1}^{n} p_i^2 + V$ where $V = \sum_{k<j}(x_k - x_j)^{-2}$ for k = 1,,n where the $x_j$ are real and the matching $p_j$ are momentum. In Mo2 (1975) Moser supported the work of Calogero and Sutherland by showing that the Lax pairs of the classical Lagrangian can (in theory) be used to find the necessary eigenvalues in the form of polygonal curves on the $p_k$ or $x_k$.. It is no coincidence that the primitive star[k] of N = 42 can also be regarded as igenvalues of N because they are roots of the minimal polynomial of $\mathrm{Tan}^2[\pi/42]$ which is a generator of $\mathbb{Q}_N^+$.

Hamiltonians of this type are now known as Calogero-Moser systems. Calegero's example was based on the classical n-body problem but it was no surprise that there are similar system with potential of the form $V = 1/\sin^2(x_1 - x_j)$ which could be used to model the Pauli exclusion of the quantum mechanical case. This was also described by Moser and once again the integrals of motion were polynomials in $p_k$ derived from eigenvalues  This may be the best match for a regular N-gon because the primitive star[k] have a natural 'twist' which can be used to model the intrinsic angular momentum $p_k$. The general Hamiltonian will have the disjoint form

$H = \sum_{1}^{n} H_i(x_i, p_i) + \sum_{1}^{n} V_i$ where $V_i = 1/\mathrm{Sin}^2(|x - x_i|)$ Initially we will define H without the need for V, but this potential will be consistent with the formulation and it is no surprise that any regular N-gon can be regarded as a quantum oscillator with the obvious natural frequencies.

The n particles in Hamiltonian phase space will be the 6 star points determined by the primitive S[k] pairs, so each star[k] is of the form $[x_k, p_k]$ where $x_k$ is the horizontal coordinate of the shared star[k] and $p_k$ is the momentum in the form of the 'twist' of $x_k$ which has direction equal to the angle star[k] makes with the x-axis (in radians) so the minimal twist is $p_1 = \omega = 2\pi/N$ at star[1] and the upper bound is $\pi$. The star[k] have displacement $s_k = \mathrm{Tan}[k\pi/N]$ so the $x_k$ are known and the $p_k$ twists are also known since they have the form $k\pi/N$.

As noted earlier there is a natural measure of the stack height for each star[k] as shown in the table below. It is the integer energy level between two extended edges of N measured as $np_k$ namely 1,2,4,5,8,10. These are the permitted eigen-states of N (and H), so an electron with these energies would find a natural resonance in a solid state crystal based on N but the combined nature of such a crystal might be very hard to predict. It would be easier to predict the behavior of a crystal based on a prime (or twice prime) N-gon because the consecutive eigen-states would reduce to the single twist as in the KAM Theorem where the $I_k$ are nested tori.

**Figure 22** The parameters of the $I_k$ for N = 42

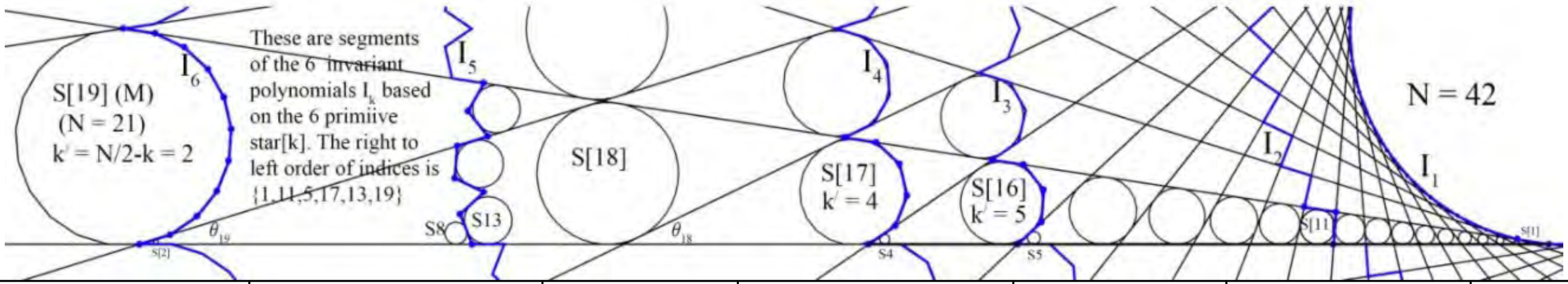

| $I_k$ | $I_1$ | $I_2$ | $I_3$ | $I_4$ | $I_5$ | $I_6$ |
|---|---|---|---|---|---|---|
| Position ($x_k$) | Star[1] | Star[11] | Star[5] of S[16] | Star[17] | Star[13] of S8 | Star[19] |
| Momentum $p_k$ (radians) | $p_1 = -(-\pi/21))$ | $p_{11} = 10\pi/21$ | $p_5 = 5\pi/21$ | $p_{17} = 4\pi/21$ | $p_{13} = -(-8\pi/21))$ | $p_{19} = 2\pi/21$ |
| length of $I_k$ | N | 2N | 4N | 5N | 8N | 10N |
| $E_k$ (angular) | $p_1(N) = 2\pi = \omega$ | $p_{11}(2N)$ | $p_5(4N)$ | $p_{17}(5N)$ | $p_{13}(8N)$ | $10p_{19}$ |
| Total $E_k$ | $2\pi$ | $40\pi$ | $40\pi$ | $40\pi$ | $128\pi$ | $40\pi$ |
| Height of stack | 1 | 2 | 4 | 5 | 8 | 10 |

The examples to follow will include N = 14, 22 and 34 in the twice-prime family and there the $p_k$ will always be consecutive $k\omega$ where $\omega = p_1 = 2\pi/N$, so they will just be 'overtones' of $\omega$ and this implies that prime N-gons are ‘harmonic’ oscillators with a single fundamental twist - like KAM with the twist being the ratio of the orbital periods of Jupiter and Saturn.

For N = 42 the ‘outlier’ is star[13] of the S[8] clone. There is a gap here in the progression of dual pairs because S[18] has a non-trivial mutation with 3 a divisor of N. This is similar to what happened earlier with S[14] and S[7] which fails to exist because N is divisible by 7. For S[18] much of the underlying step-3 left-side family survives and contributes to S8 and S13. Their shared edge defines star[13] of S8 and this shows that a ‘primitive’ star[k] may be primitive only with respect to the S[k]. But star[13] of N is dual to star[8] with equivalent twist.

The best way to describe the motion here is to convert from [x, p] Hamiltonian phase-space co-ordinates to [$\theta$,I] ‘action-angle’ coordinates . To prove that this change is ‘canonical’ it is necessary to show that area elements are preserved so $d\theta \wedge dI = dx \wedge dp$ . This will follow if the Poisson bracket {$\theta$,I} = 1. By definition {$\theta$,I}(with respect to x,p) $= \dfrac{\partial\theta}{\partial x}\dfrac{\partial I}{\partial p} - \dfrac{\partial\theta}{\partial p}\dfrac{\partial I}{\partial x}$

Now setting I = p and $\theta$ = x

$\dfrac{\partial\theta}{\partial x} = \dfrac{\partial x}{\partial x} = 1 \quad \dfrac{\partial\theta}{\partial p} = \dfrac{\partial x}{\partial p} = 0 \quad \dfrac{\partial I}{\partial x} = \dfrac{\partial p}{\partial x} = 0 \quad \dfrac{\partial I}{\partial p} = \dfrac{\partial p}{\partial p} = 1$ so {$\theta$,I} = 1·1 - 0·0 = 1 Therefore in action-angle coordinates H has a very simple form $H(\theta, I) = \sum_{k=1}^{n} H_k(\theta_k, I_k)$ and the flow on each $H_k$ will be the ‘twist’ $p_k$ of $x_k$ given by $\theta_{n+1} = \theta_n + p_k$ with $\theta_0 = 0$ *and* $I_k$ *fixed.*

Such Hamiltonians will be individually ‘integrable’ since each $H_k$ must preserve $I_k$ , but there is no known theory about what the combined solution will look like .The dimension must be $C_N$ to describe the invariant regions but the direct sum of independent twists will not be a simple issue. The N prime case with simple overtone twists should mimic the issues with KAM.

In material science the missing energy levels in H will correspond to ‘quantum band gaps’ where a solid-state silicon crystal will not pass electrons with specific frequencies. Such crystals are of interest in the construction of semiconductors , superconductors or quantum computers (where the synchronization of the S[k] and S[k′] bosons creates a stable entanglement). As noted above, these crystals can be formed using ‘superlattices’ where two or more physical crystal grids are nested inside each other, layer by layer like nested cyclotomic fields for the divisors of N. If N ‘friendly’ bosons are injected into a very cold Boise-Einstein gas they will form ‘bright’ solitons with height depending on N so bosons can ‘hug’ just like S[k] and S[k′].

Examples of superlattices

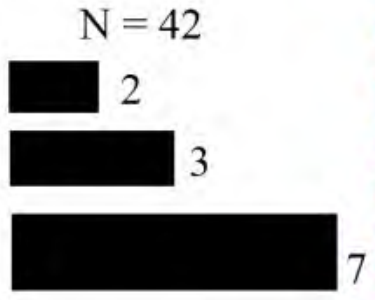

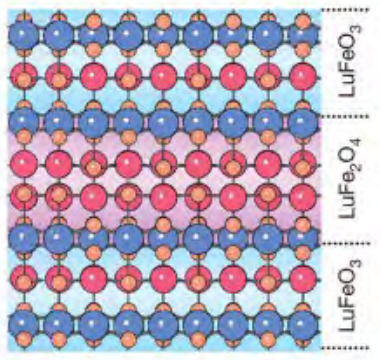

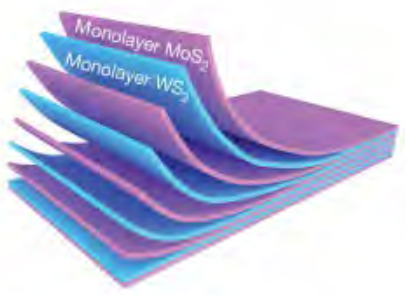

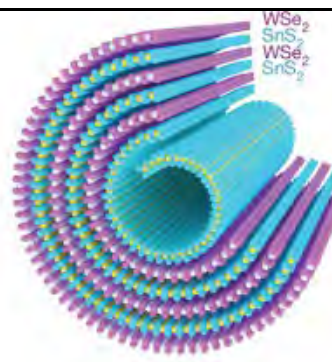

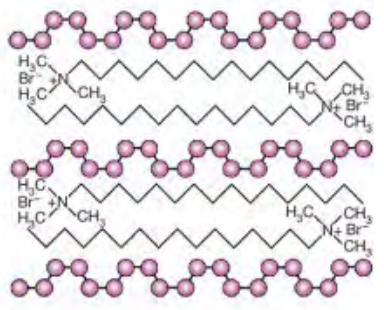

## Examples of $I_k$ and matching $E_k$ for N = 14, 20, 22, 34 and 60

**Figure 23** The $I_k$ for N = 14 and N = 7

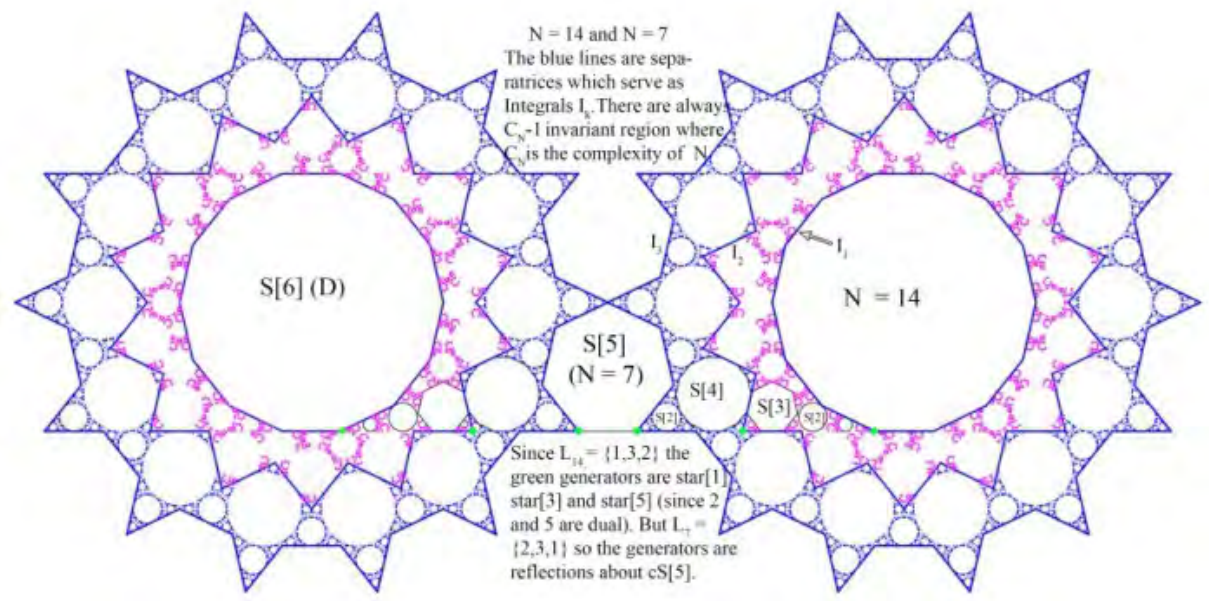

Heights of the $I_k$ are 1,2,3

| | $I_1$ | $I_2$ | $I_3$ |
|---|---|---|---|
| position | star[1] | star[4] | star[5] |
| momentum | $-\pi/7$ | $3\pi/7$ | $2\pi/7$ |
| Length | -N | 2N | 3N |
| $E_k$ | $2\pi$ | $12\pi$ | $12\pi$ |

**Figure 24** N = 20 with $L_{20}$ ={1,3} = {9,7}

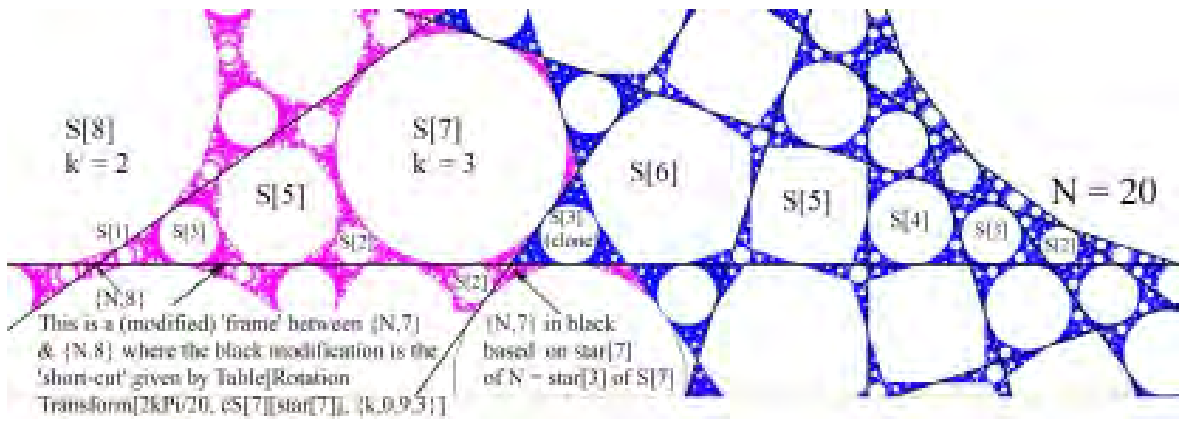

Heights of $I_k$ (and $-I_k$) are 1,3

| | $I_1$ | $I_2$ |
|---|---|---|
| position | star[1] | star[7] |
| momentum | $-\pi/N$ | $3\pi/10$ |
| Length | -N | 3N |
| $E_k$ | $2\pi$ | $18\pi$ |

For N twice-even the number of distinct $I_k$ is $C_N/2$ but technically the reflections are distinct.

**Figure 25** N = 22 with $L_{22}$ = {1,5,4,3,2}

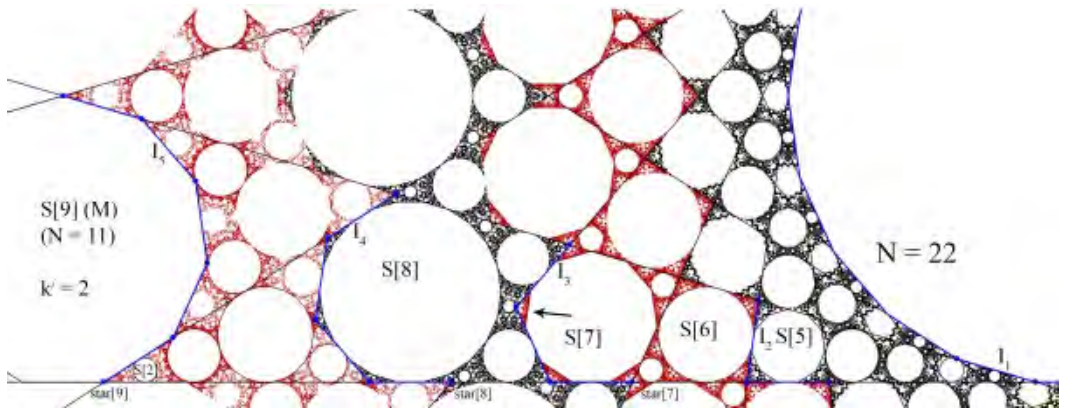

Heights of $I_k$ are 1,2,3,4,5

| | $I_1$ | $I_2$ | $I_3$ | $I_4$ | $I_5$ |
|---|---|---|---|---|---|
| position | Star[1] | Star[5] | Star[7] | star[8] | Star]9] |
| momentum | $-\pi/11$ | $5\pi/11$ | $4\pi/11$ | $3\pi/11$ | $2\pi/11$ |
| Length | -N | 2N | 3N | 4N | 5N |
| $E_k$ | $2\pi$ | $20\pi$ | $24\pi$ | $24\pi$ | $20\pi$ |

The arrow points to the magenta 'island' on the extended edges of S[7]. These islands typically have heights larger than the 'host' tile and that is true for the large S[8] background island also. The S[k] are well behaved step-k so their height as 'average' winding number is just k/N. Therefore in the density plot S[7] in the foreground is lowest at 7/22 and its island reaches up to S[8] at 8/22. It is no coincidence that this is also the intrinsic twist 8Pi/22 of star[7] so it is the fixed height of $I_3$ and an eigenstate of N. As $I_k$ approaches S[4] and the shared edge, the $I_3$ height seems to diminish but actually the neighbors get higher and the center of the edge could be what is called a Dirac point with equal energy on both sides of the border.

| Detail | The 'density' plot | |
|---|---|---|
| 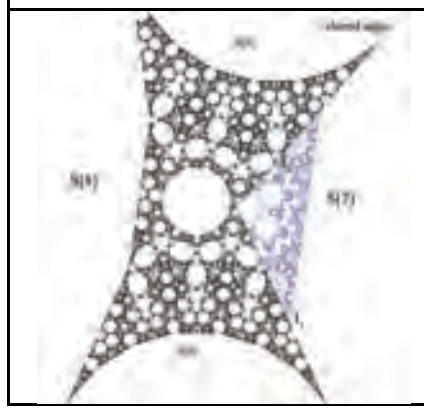 | 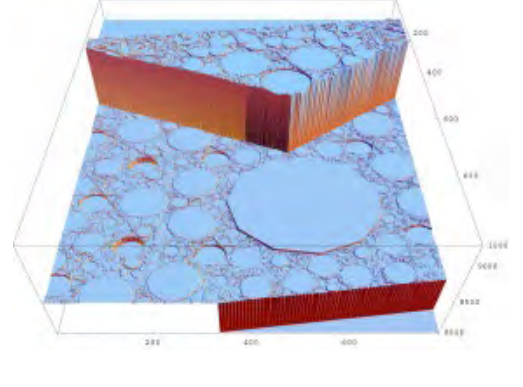 | When $I_3$ meets the extended vertex of S[7], theory says that diffractive tunneling could allow communication across the boundary in a manner similar to a Dirac point. |

**Figure 26** N = 34 and N = 17

$L_{34}$ = [1,8,6,4,7,3,5,2]

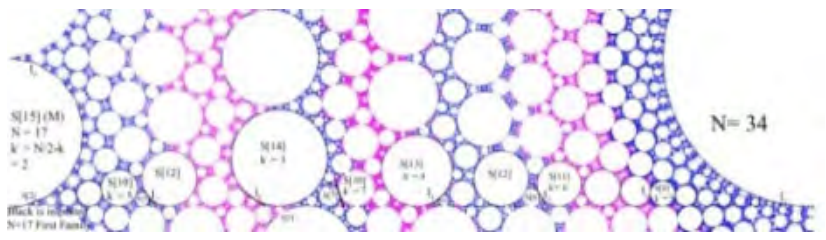

Heights of the $I_k$ are 1,2,3,…,8

| | $I_1$ | $I_2$ | $I_3$ | $I_4$ |
|---|---|---|---|---|
| position | star[1] | star[9] | star[11] | star[13] |
| momentum | $-\pi/17$ | $8\pi/17$ | $6\pi/17$ | $4\pi/17$ |
| Length | --N | 2N | 3N | 4N |
| $E_k$ | $2\pi$ | $32\pi$ | $36\pi$ | $32\pi$ |

| $I_5$ | $I_6$ | $I_7$ | $I_8$ |
|---|---|---|---|
| star[7] | star[3] | star[5] | star[15] |
| $7\pi/17$ | $3\pi/17$ | $5\pi/17$ | $2\pi/17$ |
| 5N | 6N | 7N | 8N |
| $70\pi$ | $36\pi$ | $70\pi$ | $32\pi$ |

Even though the primitive star[k] are not consecutive, the boson stack heights are consecutive. They represent energy gaps of the $I_k$ which are measured by the number of steps between two extended edges of N. These are what we call the eigenstates of H and N and they are of the form $k\omega$ with $\omega = p_1 = 2\pi/N$. so N resembles a harmonic oscillator with 'base' frequency $\omega$.

**Figure 27** N= 60 with $L_{60}$ = {1,13,11,7}

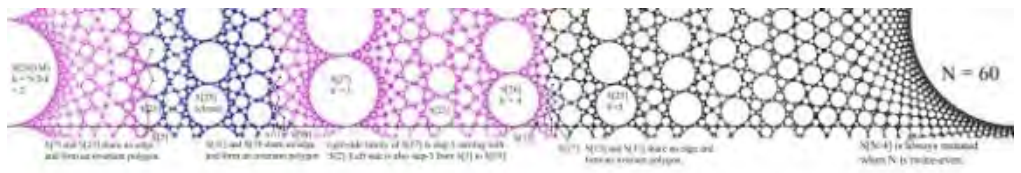

Heights of $I_k$ are 1,7,11,13

| | $I_1$ | $I_2$ | $I_3$ | $I_4$ |
|---|---|---|---|---|
| position | Star[1] | Star[13] | Star[11] | Star[7] |
| momentum | $-\pi/30$ | $13\pi/30$ | $11\pi/30$ | $7\pi/30$ |
| Length | -N | 7N | 11N | 13N |
| $E_k$ | $2\pi$ | $182\pi$ | $242\pi$ | $182\pi$ |

These star points are all 'proxies' but their energy levels still progress naturally.

As shown below N =60 has so many divisors that the primitive count is just 1,13,11 and 7 (and their duals). The plot given here has enough detail that it is easy to take each of the 4 primitive dual pairs and count the steps needed to bridge one extended edge of N. These are what we call the eigenstates of N because they represent the energy needed to transit an $I_k$ and in this sense they represent the natural energy states of N. Since the primitive star[k] with Det(J) = 0 are the Galois conjugate families determined by $\mathbb{Q}_N^+$ they represent the 'resonant' modes of N = 60 as a harmonic oscillator. These modes are all based on the nested cyclotomic fields that make up $\mathbb{Q}_N$ and these 'factor' polygons are shown below.

**Figure 28** The 'factor polygons' for N = 60 showing shared star[k] (larger star[k] not shown)

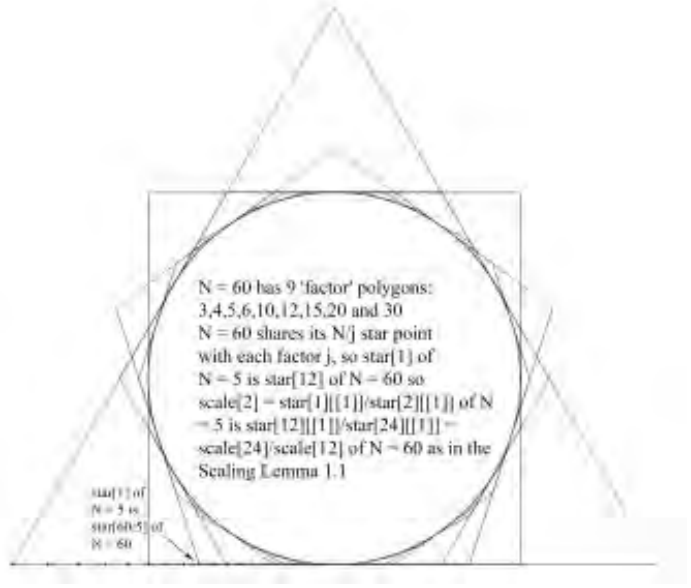

As a solid-state crystal N = 60 should be pretty 'opaque' with only a short list of eigenstates and just a few lonely $I_k$ to represent them. Since N is twice-even each $-I_k$ will also be an eigenstate where signs change left of M so technically 8 eigenvalues but looks like 4. This is in keeping with the 'chiral' nature of N.

We continue to probe the differences and similarities between the classical KAM results and regular N-gons. The classical case is summarized in the illustration below from Moser based on a periodically perturbed twist map that generates a mixture of rational (dotted) and irrational tori.

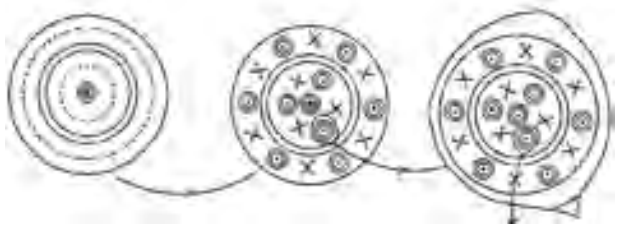

To generate these tori Moser relied on a 'monotone twist condition' which says that as the initial point moves outward the angular twist p varies smoothly so there is a continuous gradient of twists which can be used to find action-minimizing KAM tori. In the 1980's when John Mather and Serge Aurby generalized this in A-M Theory they had to retained the monotonic twist condition to generate a spectrum of what they called Lagrangian curves. Audry showed that it would be sufficient to have continuity and just minimize the action so Sorrentino [So] (2015) and others were able to weaken the continuity conditions and devise a Weak KAM Theory.

Bur even Weak KAM Theory required a monotone twist and hence does not apply to regular polygons or to any C-M discrete Hamiltonian which is based on polynomial Lagrangian curves. When Moser says they are 'integrable' it is because these curves exist but like N = 42 no one knows how to solve the combined Hamiltonian so the invariant regions are still a mystery

Taking a hint from Aubry who was a solid-state physist we are worked backwards from the known S[k] pairs and would like to know their true source. In the Introduction we showed how a dual version of the First Family Theorem can explain the emergence of the dual pairs assuming the Reflexive Principle which says that the web W must have cw-ccw Reflexive symmetry and we noted in Definition 4.2 that both $\tau$ and N do indeed have reflective symmetry.

Below are the $D_5$ and, $D_6$ Dihedral groups each with 2N elements equally divided between rotations $r$ and reflections $f$. They are 'achiral' because a simple rotation would overlay one on its reflected image, but this would fail if the edges had orientation as shown in $D_3$. Now the dihedral symmetry group is non-Abelian with $f \circ r = -r \circ f$ and the final $D_3$ is chiral with ± 'spin'

<table>
<tr><th>Unlabeled $D_5$ and $D_6$ are achiral</th><th>$D_3$ showing 3 reflection axes reduced to 0 by labeling neighbors</th><th>N = 20 showing how $\pi$ gives S[7] and S[3] opposing spins</th></tr>
<tr><td>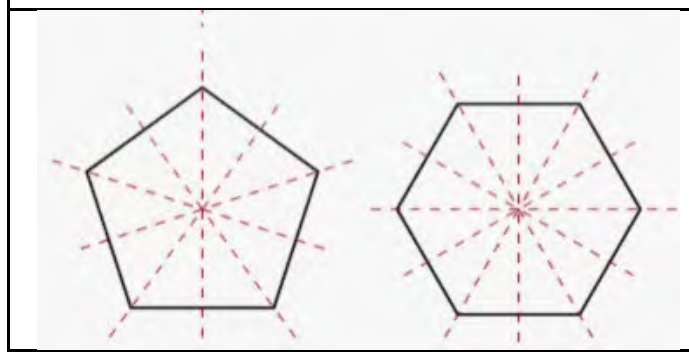</td><td>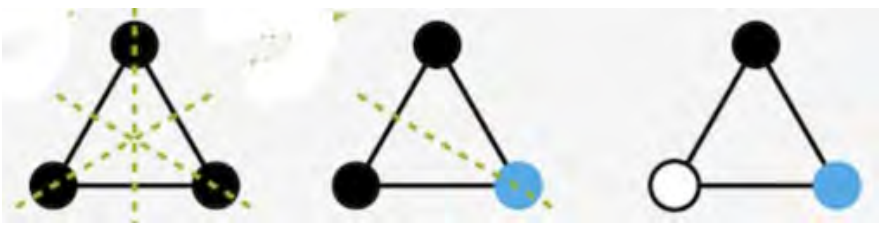</td><td>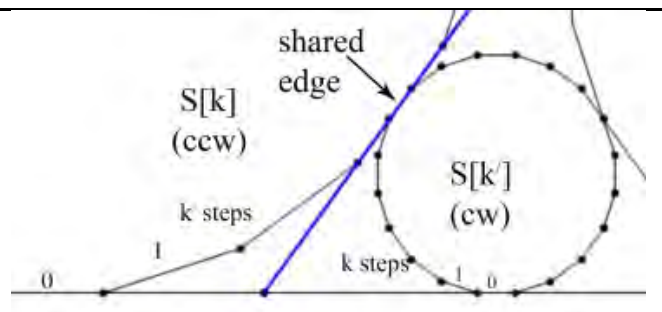
</td></tr>
</table>

So regular N-gons are 'chiral' but over the years of evolution of this paper we had 'directional bias' where +W was emphasized over –W. But about 3 years ago the curious dual S[k] tiles arose which had spin in opposing directions and we realized –W had an equal say in W.

**Theorem**: Every regular N-gon has $C_N$ distinct dual pairs which survive in W and each pair determines an $I_k$ separatrix between invariant regions so there will be $C_N$ -1 distinct invariant regions and the $I_k$ polygonal curves represent the eigenstates of N as a harmonic oscillator.

These many examples give credence to the notion that a regular N-gon can be regarded as a harmonic oscillator with $C_N$ natural resonant modes. In the quantum world this would correspond to $C_N$ particles in a Hamiltonian H with these resonant eigenvalues as opposing twists with known angular momentum. The Pauli exclusion principle would allow bosons with opposing twists to form stacks of arbitrary height and mimic molecules. Indeed the macro form of these stacks is known to produce solitons as ‘bright’ sources in Boise-Einstein gas.

It seems that there are multiple ways to view a regular N-gon: (i) As a geometric object with Dihedral symmetry group $D_N$ (ii) As a classical or quantum mechanical Hamiltonian H with $C_N$ elements which represent resonant energy gaps (iii) As a solid state crystal, super-lattice or molecular compound with $C_N$ ‘quantum band gaps’.

Therefore it might be possible to learn something about the real world using 2D N-gon models which are easy to construct and manipulate – yet they hide a complexity which is probably multi-fractal. The S[k] dual pairs which arise are very real and can be witnessed by a student on their cell-phone or graphing calculator. Geometrically they have easily discernable opposing ‘twists’ which are a basic property of many quantum particles.

If Molecules like theses did exist in nature as a mixed compound it is easy to separate this prebiotic soup into their homochiral modes +W and –W using $\tau$ as a polarizing catalyst because $\tau$ is a local isometry preserving area and polarity. (Globally it is biased by N but not locally.)

<table>
<tr><th>N =42 as a collection of prebiotic tiles</th><th>After applying $\tau$, blue = cw and magenta = ccw</th></tr>
<tr><td>
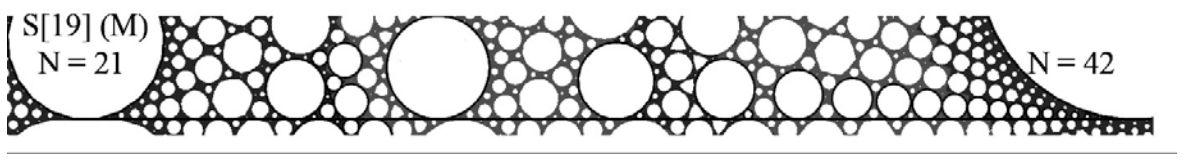

</td><td>
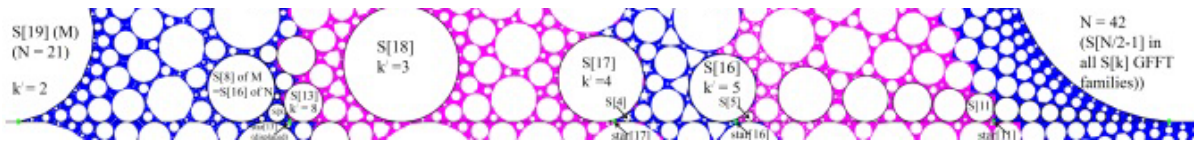

</td></tr>
</table>

In [KS] H.Kagan and K.Soai won the 2026 Nobel prize in Chemistry by showing how to do this in the real world. The exciting and innovative aspect of their work is to show that homochiralty can occur through natural processes. The Soai reaction used to form pure homochiral molecules is an iterative process and here in N = 42 Euclidean ‘phase-space’ $\tau$ polarizes a point p by iterating k times and noting that points in the ‘orbit’ $\{\tau^k(p)\}$ must have the same polarity ±W. But in the spirit of Kagan we believe that N can be ‘auto-catalytic’ and independent of $\tau$.

It is not uncommon for young stars to output circularly polarized light and in photochemistry UV light has been used to create asymmetric photolysis. Many such experiments were described by T. Fukue in 2010 followed by a more in-depth 2022 study by Kamada et al. It is a challenge to detect potential enantiomers in inter-stellar space but a natural place to look is the central region of our galaxy where the radiation is intense. Recently NASA astrobiology released the image below from the Green Bank Radio Telescope as part of the Prebiotic Interstellar Molecular Survey. The compound shown here is propylene oxide $CH_3CHOCH_2$ the most complex enantiomer candidate yet recognized.

<table>
<tr><th>Circularly polarized light can be polarized in space</th><th>Propylene oxide in a molucular cloud near Sgr B2</th></tr>
<tr><td>
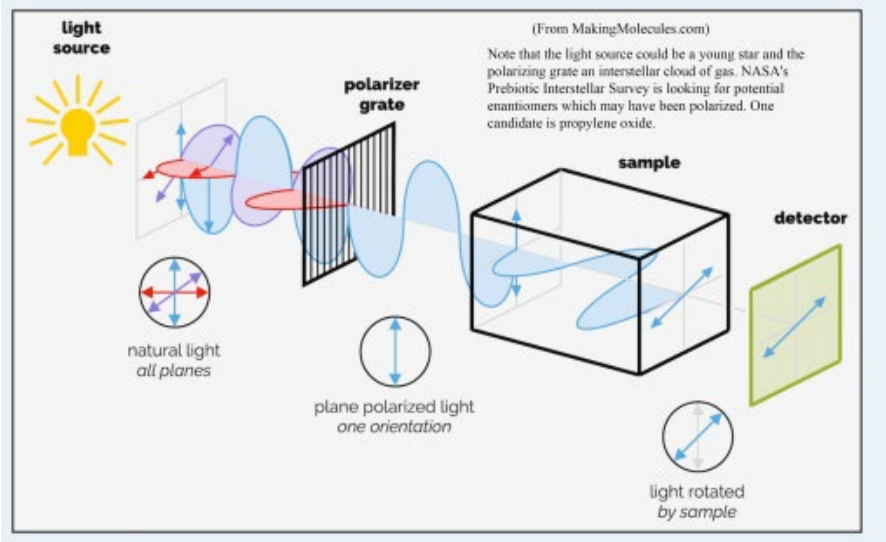

</td><td>
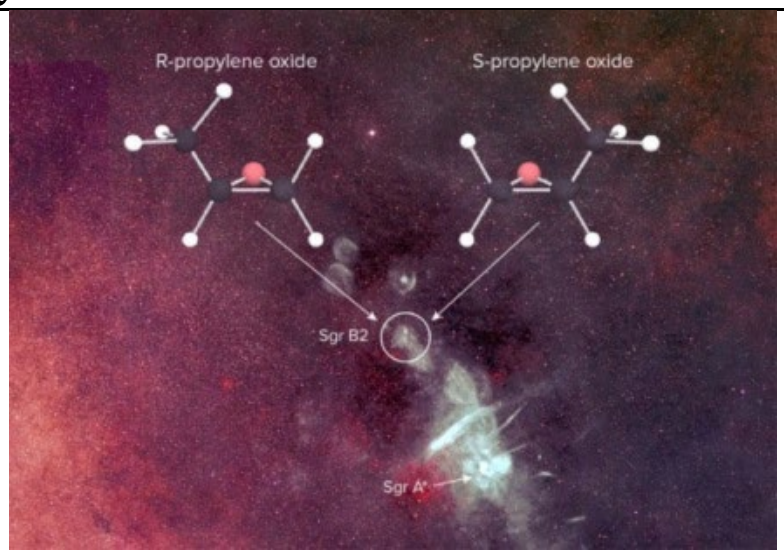

</td></tr>
</table>

**Links**

▪The author's web site at DynamicsOfPolygons.org is devoted to the outer billiards map and related maps from the perspective of a non-professional. The items below are available there.

Under Software a Mathematica notebook called FirstFamily.nb will generate the First Family and related star polygons for any regular polygon. It is also a full-fledged outer billiards notebook which works for all regular polygons. This notebook includes the Digital Filter map and the Dual Center map. This notebook is not necessary to implement the Digital Filter or Dual Center maps, but it may be useful to have a copy of the matching First Family as reference.

▪ Under Software there are also notebooks which include the Non-regular case.

▪ For someone willing to download the free Mathematica CDF Reader there are many 'manipulates' available at the Wolfram Demonstrations site - including an outer billiards manipulate of the author and two other manipulates based on the author's results in [H2]. It may be necessary to explicitly enter 'outer-billiards' in the search box.

▪ Some mp4 animations using Projections as proposed in [S2] and [H6].

| | | | | | | |
|---|---|---|---|---|---|---|
| N30P5 | N32p1P8 | N36pxP9A | N36s2P9 | N60s2P10 | N100s2P25 | N202P2 |
| N32p10P12 | N32p10P8 | N36P10-P13 | N45s7P6 | N60s2P15 | N202P17 | |

▪ An annotated Chronology of the KAM Theorem and the outer-billiards

map ▪ A graphical summary of maps relevant to this paper